# HF=HM III: Holomorphic curves and the differential for the ech/Heegaard Floer correspondence

Cagatay Kutluhan, Yi-Jen Lee and Clifford Henry Taubes

ABSTRACT: This is the third of five papers that construct an isomorphism between the Seiberg-Witten Floer homology and the Heegaard Floer homology of a given compact, oriented 3-manifold. The isomorphism is given as a composition of three isomorphisms; the first of these relates a version of embedded contact homology on an an auxillary manifold to the Heegaard Floer homology on the original. This paper describes the relationship between the differential on the embedded contact homology chain complex and the differential on the Heegaard Floer chain complex. The paper also describes the relationship between the various canonical endomorphisms that act on the homology groups of these two complexes.

The Heegaard Floer homology of a given, compact and oriented three manifold as defined by Peter Ozsváth and Zoltan Szabó [OS1] is computed using a suitably chosen Morse function and associated pseudogradient vector field. This given manifold is denoted by M. The second paper in this series [KLTII] used this Heegaard Floer data for M to construct geometric data on the connect sum of M with a certain number of copies of $S^1 \times S^2$. What follows uses Y to denote this connect sum but with orientation opposite from that on M. The geometric data on Y was used in [KLTII] to define a stable Hamiltonian version of Michael Hutching's embedded contact homology [H1]. The generators of the latter chain complex were described in [KLTII]. In particular, the $\mathbb{Z}$-module that serves as the chain complex for the embedded contact homology was written in [KLTII] as the tensor product of the $\mathbb{Z}$-module that serves as the Heegaard Floer chain complex on M and a second, canonical factor. The pseudoholomorphic curves that are used to define the embedded contact homology differential were also described in [KLTII]. This paper uses as input the latter's description of these curves to relate the differential on Y's embedded contact homology chain complex to the differential on M's Heegaard Floer chain complex. The relationship between the two differentials leads to an isomorphism between the embedded contact homology chain complex on Y and a tensor product of two factors, one being the Heegaard Floer homology on M and the other a certain canonical $\mathbb{Z}$ module.

The first paper in this series [KLTI] uses the isomorphism provided by this paper as one of a triad of isomorphisms which compose to define an isomorphism between the Heegaard Floer homology on M and the Seiberg-Witten Floer homology on M. The input from this paper to [KLTI] is summarized by the latter's Theorem 2.3. Theorem 2.3 in [KLTI] is restated here as part of Theorem 1.1. The proof of Theorem 1.1 and thus of Theorem 2.3 in [KLTI] constitutes almost all of this article.

What follows is a table of contents for this article.

Notation, definitions, constructions and results from [KLTII] are used here freely, and so the reader should be familiar with [KLTII]. Note in particular the following convention: Section numbers, equation numbers, and other references from [KLTII] are distinguished from those in this paper by the use of the Roman numeral II as a prefix. For example, 'Section II.1' refers to Section 1 in [KLTII]. Note also that the convention here as in [KLTII] is to use $c_0$ to denote a constant in $(1, \infty)$ whose value is independent of all relevant parameters. The value of $c_0$ can increase between subsequent appearances. A second convention used here and in [KLTII] concerns a function that is denoted by $\chi$. It is a fixed, non-increasing function on $\mathbb{R}$ that equals 1 on $(-\infty, 0]$ and equals 0 on $[1, \infty)$.

**Acknowledgements**

The authors owe much to Michael Hutchings for sharing his thoughts about embedded contact homology. Another debt is owed to Robert Lipshitz for sharing thoughts about his version of Heegaard Floer homology. The authors also acknowledge the generous support of the Mathematical Sciences Research Institute where the three authors did most of the work reported here. The second and third authors were also supported in part by the National Science Foundation. The third author was also supported by a David Eisenbud Fellowship.

**1. Embedded contact homology on Y and Heegaard Floer homology on M**

The first four subsections that follow provides a very brief summary of what is done in [KLTII]. Section 1a describes Y and its relevant geometry and Section 1b describes the $\mathbb{Z}$-module that serves for the embedded contact homology chain complex on Y. The Section 1c supplies a short primer on the almost complex geometry of $\mathbb{R} \times Y$. Section 1d briefly describes how the corresponding pseudoholomorphic curves are used to define the differential for the embedded contact geometry chain complex. This subsection also describes how these curves are used to define certain important endomorphisms of the of this chain complex.

Section 1e states the central result of this article, this being Theorem 1.1. It characterizes the embedded contact homology differential in terms of the differential that defines the Heegaard Floer homology of M. This theorem likewise characterizes the afore-mentioned endomorphisms of the embedded contact homology chain complex in terms of their analogs for the Heegaard Floer chain complex on M. As noted in the introductory remarks, Theorem 1.1 restates Theorem 2.3 in [KLTI]. Section 1f gives a brief look ahead at the proof of Theorem 1.1.

The last subsections, Sections 1g and 1h, supply some additional background from [KLTII] for use in the proof of Theorem 1.1

### a) The manifold Y and its geometry

The manifold Y is diffeomorphic to the connect sum of M and $G+1$ copies of $S^1 \times S^2$. A useful handle decomposition of Y is contructed from data on M that is used to define the M's Heegaard Floer homology. The first item from this data set is a self-indexing Morse function, this denoted by $f$. The map $f$ has image [0, 3], it has one index 0 critical point, one index 3 critical point, some $G \geq 1$ index 1 critical points and the same number of index 2 critical points. The latter are on the respective $f = 1$ and $f = 3$ level sets. The level set $f = \frac{3}{2}$ is denoted by $\Sigma$; this being the Heegaard surface, a surface of genus $G$. The second item from the data set is the choice of a class in $H^2(M; \mathbb{Z})$ which defines a homomorphism from $H_2(M; \mathbb{Z})$ to $2\mathbb{Z}$. This class is denoted in what follows by $c_{1M}$. A $\text{Spin}^{\mathbb{C}}$ structure will be chosen momentarily, and its first Chern class will play the role of $c_{1M}$. The third item from this data set is a fiducial point in $\Sigma$, this denoted here by $z_0$. The final item consists of an appropriate pseudogradient vector field for $f$. This vector field is denoted by $\mathfrak{v}$; it is defined on the complement of $f$'s critical points and it is such that $\mathfrak{v}(f) = 1$. This vector field is chosen to obey various constraints; these are described in Section II.1b, II.1c and II.1d. Note in particular that $\mathfrak{v}$ is constrained so as to give what Ozsváth and Szabó in [OS1] deem to be a strongly admissable Heegaard diagram for the chosen class $c_{1M}$ and the point $z_0$. The data consisting of $(f, c_{1M}, z_0, \mathfrak{v})$ is said in what follows to be the Heegaard Floer data. Constants that depend on just this data are said to *depend solely on the Heegaard Floer data*.

The construction of Y requires the choice of a pairing between the set of index 1 critical points of $f$ and the set of index 2 critical points of $f$. The resulting set of $G$ pairs is denoted by $\Lambda$. An element $\mathfrak{p} \in \Lambda$ is written as an ordered pair of points with it understood that the first entry is the index 1 critical point of $f$ and the second entry is the index 2 critical point of $f$.

The definition of Y also requires the choice of 2 additional positive numbers, these denoted by $\delta_*$ and R. The constant $\delta_*$ is from $(0, 1)$ and it is determined solely by the Heegaard Floer data. The constant R has the lower bound $-100 \ln \delta_*$. This constant R has no apriori upper bound, and the freedom to take R as large as needed is exploited in [KLTII] and in the constructions to come in this article.

The construction of the geometry needed for the $\mathbb{Z}$-module that serves as the embedded contact geometry chain complex requires the choice of two additional positive numbers, these denoted by $\delta$ and $x_0$. The latter with R are not determined by the Heegaard Floer data. The trio $(\delta, x_0, R)$ are constrained by the requirements that $\delta < \delta_*$, $x_0 < \delta^3$ and $R \geq -c_0 \ln x_0$. Note in particular that the choice of $\delta$ determines an upper bound for $x_0$, and that the choice of $x_0$ subject to this upper bound then determines a lower bound for R. Constants $\delta$, $x_0$, and R that satisfy these bounds are said to be *appropriate*.

The remaining parts of this subsection describe first Y and then the geometry that is needed to define the embedded contact homology $\mathbb{Z}$ module.

*Part 1*: As noted above, the manifold Y is diffeomorphic to the connect sum of M with $G+1$ copies of $S^1 \times S^2$. Section II.1a supplies a useful handle decomposition of Y as $M_\delta \cup \mathcal{H}_0 \cup_{\mathfrak{p} \in \Lambda} \mathcal{H}_\mathfrak{p}$ where $M_\delta$ is the complement in M of a certain set of $2(G+1)$ disjoint balls whose centers are the critical points of $f$; and where $\mathcal{H}_0$ and each $\mathfrak{p} \in \Lambda$ version of $\mathcal{H}_\mathfrak{p}$ is a copy of $[-1,1] \times S^2$. What follows summarizes from Section II.1a how these 1-handles are attached.

THE HANDLES $\{\mathcal{H}_\mathfrak{p}\}_{\mathfrak{p} \in \Lambda}$: Fix $\mathfrak{p} = (p_+, p_-) \in \Lambda$. The constant $\delta_*$ is chosen so that there are respective coordinate charts centered on the index 1 critical point $p_+$ and index 2 critical points $p_-$ with coordinates (x, y, z) defined where $|x|^2 + |y|^2 + |z|^2 \leq (10\delta_*)^2$ and such that $f$ appears as

$$f = 1 + x^2 + y^2 - 2z^2 \quad and \quad f = 2 - (x^2 + y^2 - 2z^2)$$
(1.1)

Use $(r_+, (\theta_+, \varphi_+))$ to denote the standard spherical coordinates for the Euclidean coordinate chart centered on $p_+$, and likewise use $(r_-, (\theta_-, \varphi_-))$ to denote the spherical coordinates for the coordinate chart centered on $p_-$. When $d \in (0, 10\delta_*)$, the ball in M given by $r_+ < d$ is said to be the *radius d coordinate ball* centered on $p_+$, and the corresponding $r_- < d$ ball is said to be the radius d coordinate ball centered on $p_-$.

The handle $\mathcal{H}_\mathfrak{p}$ is given coordinates $(u, (\theta, \phi))$ where $(\theta, \phi)$ are the standard spherical coordinates on the $S^2$ factor, and where $u \in [-R-\ln(7\delta_*), R+\ln(7\delta_*)]$ is the coordinate for the interval factor. The handle $\mathcal{H}_\mathfrak{p}$ is attached to the complement in M of the radius $e^{-2R}(7\delta_*)^{-1}$ coordinate balls centered on $p_+$ and $p_-$ via the identifications given by

$$(r_+ = e^{-(R-u)}, (\theta_+ = \theta, \varphi_+ = \phi)) \quad and \quad (r_- = e^{-(R+u)}, (\theta_- = \pi - \theta, \varphi_- = \phi)) \ .$$
(1.2)

The part of $\mathcal{H}_\mathfrak{p}$ where $(1 - 3\cos^2\theta) > 0$ is denoted by $\mathcal{H}^+_\mathfrak{p}$. Any given constant u slice of $\mathcal{H}^+_\mathfrak{p}$ is an annular neighborhood of the equator in $S^2$.

THE HANDLE $\mathcal{H}_0$: The constant $\delta_*$ is chosen so that respective coordinate charts centered on the index 0 and index 1 critical points of $f$ have coordinates (x, y, z) that are defined where the coordinate functions obey $|x|^2 + |y|^2 + |z|^2 \leq (10\delta_*)^2$ and are such that $f$ appears as

$$f = x^2 + y^2 + z^2 \quad and \quad f = 3 - (x^2 + y^2 + z^2)$$
(1.3)

Use $(r_+, (\theta_+, \varphi_+))$ to denote the standard spherical coordinates for the Euclidean coordinate chart centered on the index 0 critical point of $f$, and use $(r_-, (\theta_-, \varphi_-))$ to denote the spherical coordinates for the coordinate chart centered on the index 0 critical point of $f$. When $d \in (0, 10\delta_*)$, the ball in M given by $r_+ < d$ is said to be the radius d coordinate ball about the index 0 critical point of $f$. The corresponding $r_- < d$ ball is said to be the radius d coordinate ball about the index 3 critical point of $f$.

The handle $\mathcal{H}_0$ is given coordinates $(u, (\theta, \phi))$ where $(\theta, \phi)$ are the standard spherical coordinates on the $S^2$ factor, and where $u \in [-R-\ln(7\delta_*), R+\ln(7\delta_*)]$ is the coordinate for the interval factor. The handle $\mathcal{H}_0$ is attached to the complement in M of the radius $e^{-2R}(7\delta_*)^{-1}$ coordinate balls centered on the index 0 and index 3 critical points of $f$ by rule in (1.2).

The constant $\delta_*$ is chosen so that the respective radius $10\delta_*$ coordinate balls about any two distinct critical points of $f$ are disjoint. Given $r \in (e^{-2R}(7\delta_*)^{-1}, 10\delta_*)$, the complement in M of the union of the radius r coordinate balls centered on M is denoted by $M_r$. The description of Y just given identifies $M_r$ with a subset in Y. The latter is denoted also by $M_r$. The just described identification is used implicitly in what follows to view these two incarnations of $M_r$ as one and the same 3-manifold with boundary. In particular, this identification defines $f$ as a function on the $M_r$ part of Y, the latter also denoted by $f$.

*Part 2*: A stable Hamiltonian structure on Y consists of a pair $(a, w)$ where $a$ is a 1-form, $w$ is a 2-form, and these are such that $dw = 0$ and $da \in \text{Span}(w)$. Moreover, $a \wedge w$ is nowhere zero; and this 3-form defines the orientation for Y that is opposite to that defined for Y by M's orientation of $M_\delta$. Being nowhere zero, the 2-form $w$ defines a homomorphism from TY to T*Y whose kernel is a real line bundle over Y. The 1-form $a$ is non-zero on this line subbundle. This understood, let $v$ henceforth denote the vector field that spans the kernel of $w$ and has pairing 1 with $a$.

Sections II.1b-e describe a stable Hamiltonian structure for Y that is defined using the Heegaard Floer data $(f, c_{1M}, \mathfrak{v})$ and appropriate constants $\delta$, $x_0$ and R. The salient features of $a$, $w$ and $v$ are summarized momentarily. This summary restates what is said in Section II.1.e.

By way of notation, the upcoming formulae use functions $x, \chi_+$ and $\chi_-$ of $u \in \mathbb{R}$ given by $x = x_0\chi(|u|-R-\ln\delta-12)$ and $\chi_+ = \chi(-u+\frac{1}{4}R)$ and $\chi_- = \chi(u-\frac{1}{4}R)$. The formulae also imploy functions f and g of the variable u given by

$$f = x + 2(\chi_+ e^{2(u-R)} + \chi_- e^{-2(u+R)}) \quad and \quad g = (\chi_+ e^{2(u-R)} - \chi_- e^{-2(u+R)}).$$

(1.4)

Their respective derivatives are denoted by f´ and g´.

The bullets that follow supply the promised description of $a$, $w$ and $v$.

- ON $M_\delta$: *The 2-form w on $M_\delta$ is nowhere zero on the kernel of the 1-form df and $v$ here is the pseudogradient vector field $\mathfrak{v}$.*

- IN THE HANDLE $\mathcal{H}_0$: *The 2-form w and the vector field v on $\mathcal{H}_0$ are*

$$w = \sin\theta\, d\theta \wedge d\phi \quad \text{and} \quad v = \frac{1}{2(\chi_+ e^{2(u-R)} + \chi_- e^{-2(u+R)})} \frac{\partial}{\partial u}.$$

(1.5)

- IN THE HANDLES $\{\mathcal{H}_\mathfrak{p}\}_{\mathfrak{p}\in\Lambda}$: *Fix $\mathfrak{p}\in\Lambda$. The trio $a$, $w$ and $v$ on $\mathcal{H}_\mathfrak{p}$ are*

$$a = (x + g´)(1 - 3\cos^2\theta)\,du - \sqrt{6}\,f\cos\theta\sin^2\theta\, d\phi + 6g\cos\theta\sin\theta\, d\theta,$$
$$w = 6x\cos\theta\sin\theta\, d\theta \wedge du - \sqrt{6}\,d\{f\cos\theta\sin^2\theta\, d\phi\},$$
$$v = \alpha^{-1}\{f(1-3\cos^2\theta)\partial_u - \sqrt{6}\,x\cos\theta\,\partial_\phi + f´\cos\theta\sin\theta\,\partial_\theta\}.$$

(1.6)

Here, $\alpha$ is a certain positive function of the pair $(u,\theta)$.

The next bullet concerns the cohomology class of the form $w$. This bullet refers to the direct sum decomposition

$$H_2(Y;\mathbb{Z}) = H_2(M;\mathbb{Z}) \oplus H_2(\mathcal{H}_0;\mathbb{Z}) \oplus (\oplus_{\mathfrak{p}\in\Lambda} H_2(\mathcal{H}_\mathfrak{p};\mathbb{Z}))$$

(1.7)

that comes via Mayer-Vietoris by writing $Y = \mathcal{M}_\delta \cup \mathcal{H}_0 \cup (\cup_{\mathfrak{p}\in\Lambda}\mathcal{H}_\mathfrak{p})$. The summands in (1.7) that correspond to the various 1-handles are isomorphic to $\mathbb{Z}$; and any oriented, cross-sectional sphere is a generator. The convention in what follows is to orient these spheres with the 2-form $\sin\theta\, d\theta\, d\phi$.

- THE COHOMOLOGY CLASS OF $w$: *Integration of the 2-form w defines the linear map from $H_2(Y;\mathbb{Z})$ to $\mathbb{Z}$ that has value 2 on the generator of $H_2(\mathcal{H}_0;\mathbb{Z})$; it has value zero on each $\mathfrak{p}\in\Lambda$ version of $H_2(\mathcal{H}_\mathfrak{p};\mathbb{Z})$; and it acts on the $H_2(M;\mathbb{Z})$ summand in (1.7) as the pairing with the chosen class $c_{1M}$.*

A particular integral curve of the vector field $v$ plays a distinguished role in the embedded contact homology story. This curve is described next.

- THE CURVE THROUGH $z_0$: *There is a closed integral curve of $v$ in $M_\delta \cup \mathcal{H}_0$ whose intersection with $\Sigma$ is the chosen fiducial point $z_0$. This is curve is denoted by $\gamma^{(z_0)}$. It also has a single intersection with each cross-sectional sphere in $\mathcal{H}_0$.*

The final bullets introduce a pair of auxilliary 1-forms on Y that play central roles. The definition of the first of these 1-forms refers to the function $f_*$ that is defined on any given $\mathfrak{p} \in \Lambda$ version of $\mathcal{H}_\mathfrak{p}$ by the rule

$$f_* = (\chi_+ e^{2(u-R)} - \chi_- e^{-2(u+R)})(1 - 3\cos^2\theta).$$

(1.8)

The definition of the second of these 1-forms refers to the function $\chi_\delta$ that is defined on any given $\mathfrak{p} \in \Lambda$ version of $\mathcal{H}_\mathfrak{p}$ by the rule $\chi_\delta = \chi(|u| - R - \ln\delta - 10)$.

- THE 1-FORM $\upsilon_\diamond$: *The 1-form $\upsilon_\diamond$ is closed and is such that $\upsilon_\diamond \wedge w \geq 0$. Furthermore, $\upsilon_\diamond \wedge w = 0$ only where both $u = 0$ and $1 - 3\cos^2\theta = 0$ on each $\mathfrak{p} \in \Lambda$ version of $\mathcal{H}_\mathfrak{p}$. This 1-form equals $df$ on $M_\delta$, it is given by $\upsilon_\diamond = 2(\chi_+ e^{2(|u|-R)} + \chi_- e^{-2(|u|+R)}) du$ on $\mathcal{H}_0$, and it is given by $df_*$ on any given $\mathfrak{p} \in \Lambda$ version of $\mathcal{H}_\mathfrak{p}$.*

- THE 1-FORM $\hat{a}$: *The 1-form $\hat{a}$ has pairing 1 with $v$ and is such that $\hat{a} \wedge w > 0$. This 1-form is equal to $\upsilon_\diamond$ on $M_\delta \cup \mathcal{H}_0$ and it is equal to $\chi_\delta a + (1 - \chi_\delta)\upsilon_\diamond$ on any given $\mathfrak{p} \in \Lambda$ version of $\mathcal{H}_\mathfrak{p}$.*

### b) The embedded contact homology $\mathbb{Z}$-module

This subsection describes the $\mathbb{Z}$-module that serves as the chain complex for embedded contact homology. The subsection has four parts that briefly summarize material from Section II.1f and Section II.2

*Part 1*: Fix a Spin$^{\mathbb{C}}$ structure on M and use $c_{1M}$ now to denote the associated first Chern class in $H^2(M;\mathbb{Z})$. This class is used to construct the strongly admissible Heegaard diagram that is used to define the Heegaard Floer chain complex on M.

The $\mathbb{Z}$-module that serves for the Heegaard Floer chain complex on M for the chosen Spin$^{\mathbb{C}}$ structure can be defined with the help of a finite set that is denoted by $\mathcal{Z}_{HF}$. Any given element in $\mathcal{Z}_{HF}$ is viewed here and in [KLTII] as a suitably constrained, unordered G-tuple of integral curves of $\mathfrak{v}$. Let $\hat{\upsilon}$ denote an element from $\mathcal{Z}_{HF}$. There are three constraints on $\hat{\upsilon}$: First, each constituent integral curve from $\hat{\upsilon}$ runs from an index 1 critical point of $f$ to an index 2 critical point of $f$. Second, no two distinct constituents share the same index 1 critical point or the same index 2 critical point. This being the case, $\hat{\mathfrak{v}}$ defines a pairing between the set of index 1 critical points of $f$ and the set of index 2 critical points of $f$. The third constraint demands that the G points that comprise $\cup_{\mathfrak{v} \in \hat{\mathfrak{v}}} (\mathfrak{v} \cap \Sigma)$ with the point $z_0$ define the chosen $\text{Spin}^\mathbb{C}$ structure in the manner that is described in [OS1].

The Heegaard Floer chain complex is the free $\mathbb{Z}$-module generated by the elements of the set $\mathcal{Z}_{HF} \times \mathbb{Z}$. This module is denoted by $\mathbb{Z}(\mathcal{Z}_{HF} \times \mathbb{Z})$. This interpretation of the Heegaard Floer chain complex is used by Robert Lipshitz in [L] to reformulate Heegaard Floer homology.

*Part 2*: The class $c_{1M}$ is used to define Y and its stable Hamiltonian data $(a, w, v)$. The $\mathbb{Z}$-module for the relevant version of embedded contact homology on Y is defined with the help of a set that is denoted by $\mathcal{Z}_{ech,M}$ and whose elements are suitably constrained, finite sets of closed integral curve of $v$ that lie entirely in the union of the various $\mathfrak{p} \in \Lambda$ versions of $\mathcal{H}_\mathfrak{p}$ and the $f \in (1,2)$ part of $M_\delta$. The set $\mathcal{Z}_{ech}$ is described in the upcoming Part 3 of this subsection. What follows directly summarizes some of what is said in Section II.2 about the closed integral curves of $v$ that lie entirely in the subset of Y just described.

CLOSED CURVES IN $\cup_{\mathfrak{p} \in \Lambda} \mathcal{H}_\mathfrak{p}$: Fix $\mathfrak{p} \in \Lambda$. There are precisely two integral curves of $v$ that lie entirely in $\mathcal{H}_\mathfrak{p}$. These constitute the two components of the locus where both $u = 0$ and $1 - 3\cos^2\theta = 0$. The curve with $\cos\theta = \frac{1}{\sqrt{3}}$ is denoted by $\hat{\gamma}_\mathfrak{p}^+$ and that where $\cos\theta = -\frac{1}{\sqrt{3}}$ is denoted by $\hat{\gamma}_\mathfrak{p}^-$.

INTERSECTIONS WITH $M_\delta$: There exists a purely Heegaard Floer dependent constant $\kappa \geq 1$ whose significance is described in what follows. Construct Y with $\delta < \kappa^{-1}\delta_*^3$. Let $\gamma$ denote a closed integral curve of $v$ in $M_\delta \cup (\cup_\mathfrak{p} \mathcal{H}_\mathfrak{p})$ that intersects $M_\delta$. Then $\gamma \cap M_\delta$ consists of a finite set of segments of integral curves of $\mathfrak{v}$ in the $f^{-1}((1,2))$ part of $M_\delta$. Each such segment lies in the radius $\kappa\delta$ tubular neighborhood of an integral curve of $\mathfrak{v}$ that runs from an index 1 critical point of $f$ to an index 2 critical point of $f$.

INTERSECTIONS WITH $\cup_{\mathfrak{p} \in \Lambda} \mathcal{H}_\mathfrak{p}$: Let $\gamma$ denote a closed integral curve of $v$ in $M_\delta \cup (\cup_\mathfrak{p} \mathcal{H}_\mathfrak{p})$ that intersects $M_\delta$. Fix $\mathfrak{p} \in \Lambda$. The intersection of $\gamma$ with $\mathcal{H}_\mathfrak{p}$ consists of a finite set of segments. Let $\gamma_*$ denote any one such segment. The following is true:

- $\gamma_*$ *sits where* $1 - 3\cos^2\theta > 0$
- $\gamma_*$ *runs from the* $u = -R - \ln(7\delta_*)$ *end of* $\mathcal{H}_\mathfrak{p}$ *to the* $u = R + \ln(7\delta_*)$ *end*.

- *The function $ĥ = f(u)\cos\theta \sin^2\theta$ is constant on $\gamma_*$.*
- *The coordinate $u$ restricts as an affine coordinate to $\gamma$.*
- *The angle $\phi$ on $\gamma$ changes according to the rule:* $\frac{d\phi}{du} = -\sqrt{6}\,\frac{x(u)}{f(u)}\,\frac{\cos\theta(u)}{1-3\cos^2\theta(u)}$.

(1.9)

The assertion about the closed integral curves in $\cup_{\mathfrak{p}\in\Lambda}\mathcal{H}_\mathfrak{p}$ summarizes Lemma II.2.1; the assertion about $\gamma \cap M_\delta$ summarizes Lemma II.2.4; and the assertion about $\gamma \cap (\cup_{\mathfrak{p}\in\Lambda}\mathcal{H}_\mathfrak{p})$ summarizes Lemma II.2.2 and (II.2.5).

*Part 3*: Fix O denote the four element set $\{0, 1, -1, \{1, -1\}\}$. The elements in $\mathcal{Z}_{\text{ech},M}$ enjoy a 1-1 correspondence

$$\mathcal{Z}_{\text{ech},M} = \mathcal{Z}_{\text{HF}} \times (\times_{\mathfrak{p}\in\Lambda}(\mathbb{Z}\times\text{O})).$$

(1.10)

This correspondence is canonical when the various $\mathfrak{p} \in \Lambda$ factors of $\mathbb{Z}$ are viewed as affine spaces modeled in $\mathbb{Z}$. A base point in a given $\mathfrak{p} \in \Lambda$ version is determined by the choice of a lift to $\mathbb{R}$ of the $\mathbb{R}/(2\pi\mathbb{Z})$ value $\phi$ coordinate in $\mathcal{H}_\mathfrak{p}$.

What follows describes the geometric meaning of the correspondence in (1.10). Write a given element in $\mathcal{Z}_{\text{HF}} \times (\times_{\mathfrak{p}\in\Lambda}(\mathbb{Z}\times\text{O}))$ as $(\hat{\upsilon}, (\mathfrak{k}_\mathfrak{p}, \text{o}_\mathfrak{p})_{\mathfrak{p}\in\Lambda})$ with $\hat{\upsilon}$ from $\mathcal{Z}_{\text{HF}}$ and with any given $\mathfrak{p}\in\Lambda$ version of $(\mathfrak{k}_\mathfrak{p}, \text{o}_\mathfrak{p})$ in $\mathbb{Z}\times\text{O}$. Let $\Theta$ denote the corresponding element in $\mathcal{Z}_{\text{ech},M}$. As noted in Part 2, each element in $\mathcal{Z}_{\text{ech},M}$ is a finite set of closed integral curves of $\mathfrak{v}$ that lie in $M_\delta \cup (\cup_{\mathfrak{p}\in\Lambda}\mathcal{H}_\mathfrak{p})$. With this in mind, consider first the significance of the entry $\hat{\upsilon}$. The intersection of $\cup_{\gamma\in\Theta}\gamma$ with $M_\delta$ has G components, each being a segment of an integral curve of $\mathfrak{v}$ that runs from the boundary of the radius $\delta$ coordinate ball about an index 1 critical point of $f$ to the boundary of the radius $\delta$ coordinate ball centered on an index 2 critical point of $f$. The components $M_\delta \cap (\cup_{\gamma\in\Theta}\gamma)$ enjoy a 1-1 correspondence with the integral curves from $\hat{\upsilon}$ with the correspondence such that a given segment from $M_\delta \cap (\cup_{\gamma\in\Theta}\gamma)$ lies in the radius $c_0\delta$ tubular neighborhood of composed of its partner from $\hat{\upsilon}$. This version of $c_0$ depends only on the Heegaard Floer data.

To say more about the curves in $\Theta$, fix $\mathfrak{p}\in\Lambda$. The intersection of $\cup_{\gamma\in\Theta}\gamma$ with $\mathcal{H}_\mathfrak{p}$ has precisely one component that crosses $\mathcal{H}_\mathfrak{p}$ from the $u = -R - \ln(7\delta_*)$ end to the end where $u = R + \ln(7\delta_*)$. More is said about this component momentarily. The remaining components (if any) are determined by $\text{o}_\mathfrak{p}$ using the following rule:

- *If $\text{o}_\mathfrak{p} = 0$, then $\Theta$ contains neither $\hat{\gamma}_\mathfrak{p}^+$ nor $\hat{\gamma}_\mathfrak{p}^-$.*
- *If $\text{o}_\mathfrak{p} = 1$ or $\text{o}_\mathfrak{p} = -1$, then $\Theta$ contains $\hat{\gamma}_\mathfrak{p}^+$ or $\hat{\gamma}_\mathfrak{p}^-$ respectively, but not both of them.*
- *If $\text{o}_\mathfrak{p} = \{1, -1\}$, then $\Theta$ contains both $\hat{\gamma}_\mathfrak{p}^+$ and $\hat{\gamma}_\mathfrak{p}^-$.*

(1.11)

The end points of the segment of $(\cup_{\gamma \in \Theta} \gamma) \cap \mathcal{H}_\mathfrak{p}$ that crosses $\mathcal{H}_\mathfrak{p}$ intersects the respective $u = R + \ln\delta$ and $u = -R - \ln\delta$ spheres in $\mathcal{H}_\mathfrak{p}$ at a point whose spherical coordinates $(\theta, \phi)$ differ by at most $c_0\delta$ from the coordinates of the sphere's intersection with $\cup_{\upsilon \in \hat{\upsilon}} \upsilon$. Note that this version of $c_0$ depends solely on the Heegaard Floer data. The entry $\mathfrak{k}_\mathfrak{p}$ in $\Theta$'s label gives an indication of the total change in the coordinate $\phi$ along this segment. To elaborate, view $\mathfrak{p}$'s entry of $\mathbb{Z}$ in (1.10) as an affine space modeled on $\mathbb{Z}$. Suppose that $\Theta$ and $\Theta'$ are any two elements from $\mathcal{Z}_{ech,M}$. Label these two elements as in (1.10) and assume that they have identical $\mathcal{Z}_{HF}$ factor. Write their respective entries in $\mathfrak{p}$'s factor of the affine copy of $\mathbb{Z}$ as $\mathfrak{k}_\mathfrak{p}$ and $\mathfrak{k}_\mathfrak{p} + k$ with $k \in \mathbb{Z}$. Let $\Delta$ denote the total change in $\phi$ along part of the curve from $\Theta$ that crosses $\mathcal{H}_\mathfrak{p}$, and let $\Delta'$ denote the analogous $\Theta'$ angle change. Then

$$\Delta' - \Delta = k + \mathfrak{e}$$

(1.12)

where $|\mathfrak{e}| \leq c_0\delta$. As before, this version of $c_0$ depends only on the Heegaard Floer data.

See Section II.2 and in particular Corollary II.2.7 and Proposition II.2.8 for an expanded version of what was just said in this Part 3.

*Part 4*: The $\mathbb{Z}$ module for the relevant version of embedded contact homology is freely generated by a certain principle $\mathbb{Z}$ bundle over $\mathcal{Z}_{ech,M}$. This bundle is denoted in what follows by $\hat{\mathcal{Z}}_{ech,M}$ and the free $\mathbb{Z}$-module generated by $\hat{\mathcal{Z}}_{ech,M}$ is denoted by $\mathbb{Z}(\hat{\mathcal{Z}}_{ech,M})$. Elements of this $\mathbb{Z}$-module are finite, integer weighted formal sums of elements in $\hat{\mathcal{Z}}_{ech,M}$.

The definition of $\hat{\mathcal{Z}}_{ech,M}$ is given momentarily. By way of preliminaries, note that any given integral curve of $v$ is oriented by $v$. This being the case, each closed integral curve of $v$ defines a closed 1-cycle in Y. When $\gamma$ denotes the closed integral curve, then $[\gamma]$ is used to denote the corresponding cycle. The set of the closed integral curves from any given $\Theta \in \mathcal{Z}_{ech,M}$ likewise defines a closed 1-cycle in Y, this being the cycle $\sum_{\gamma \in \Theta} [\gamma]$. The latter cycle is denoted by $[\Theta]$. The Poincaré dual of $2[\Theta]$ defines via the tautological pairing the homomorphism from $H_2(Y;\mathbb{Z})$ to $2\mathbb{Z}$ that acts as follows: It sends the $H_2(\mathcal{H}_0;\mathbb{Z})$ summand in (1.7) to 0, it sends the generator of each $\mathfrak{p} \in \Lambda$ labeled summand in (1.7) to 2, and it acts on the $H_2(M;\mathbb{Z})$ summand as the pairing with the given class $c_{1M}$.

A somewhat non-canonical description of $\hat{\mathcal{Z}}_{ech,M}$ requires the choice of a fiducial element $\Theta_0 \in \mathcal{Z}_{ech,M}$. This done, any given element element $\hat{\Theta}$ can be viewed as an equivalence class of pairs $(\Theta, Z)$ with $\Theta \in \mathcal{Z}_{ech,M}$ and with Z an element in the $\mathbb{Z}$-module $H_2(Y;[\Theta]-[\Theta_0])$. The equivalence relation is defined with the help of the closed integral curve $\gamma^{(z_0)}$ that is described in the fifth bullet of Part 2 in Section 1b. Pairing with the Poincare dual of $\gamma^{(z_0)}$ defines a homomorphism from the $\mathbb{Z}$-module of closed 2-cycles to $\mathbb{Z}$. This pairing is denoted by $[\gamma^{(z_0)}]^{Pd}(\cdot)$. The equivalence relation has $(\Theta, Z) \sim (\Theta´, Z´)$ if and only if $\Theta = \Theta´$ and also $[\gamma^{(z_0)}]^{Pd}(Z - Z´) = 0$. The principal bundle projection map sends an equivalence class $(\Theta, Z)$ to $\Theta$. The element $1 \in \mathbb{Z}$ acts to send $(\Theta, Z)$ to $(\Theta, Z + [S_0])$ where $[S_0]$ is the $u = 0$ sphere in $\mathcal{H}_0$.

A different choice for $\Theta_0$ produces a different principal $\mathbb{Z}$ bundle over $\mathcal{Z}_{ech,M}$, but the new one and the original are canonically isomorphic. To elaborate, suppose that $\Theta_0´$ is a second choice. Any cycle in $H_2(Y; [\Theta_0] - [\Theta_0´])$ has a well defined intersection pairing with the curve $\gamma^{(z_0)}$. Adding a suitable multiple of the $f = \frac{3}{2}$ level set in $M_\delta$ will give a cycle in $H_2(Y; [\Theta_0] - [\Theta_0´])$ with zero intersection pairing against $\gamma^{(z_0)}$. Let $Z_0$ denote such a cycle. The isomorphism in question sends an the equivalence class of $(\Theta, Z)$ in the $\Theta_0$ version of $\hat{\mathcal{Z}}_{ech,M}$ to that of $(\Theta, Z + Z_0)$ in the $\Theta_0´$ version. A different intersection pairing zero choice gives the same equivalence class and thus the same isomorphism

As explained next, the existence of cycles with zero intersection pairing against $\gamma^{(z_0)}$ can be exploited to construct a canonical principal $\mathbb{Z}$-bundle isomorphism

$$\mathcal{Z}_{ech,M} \times \mathbb{Z} \to \hat{\mathcal{Z}}_{ech,M} .$$

(1.13)

The isomorphism depicted here is defined by a certain section of $\hat{\mathcal{Z}}_{ech,M}$ whose image corresponds via (1.13) to $\mathcal{Z}_{ech,M} \times \{0\}$. This section sends any given $\Theta \in \mathcal{Z}_{ech,M}$ to the equivalence class of a pair $(\Theta, Z)$ where $Z$ can be any 2-cycle in $H_2(Y; [\Theta] - [\Theta_0])$ that has pairing 0 against $\gamma^{(z_0)}$. A different choice of $Z$ with intersecting pairing zero against $\gamma^{(z_0)}$ defines the same equivalence class and so the same section. This is why (1.13) is canonical. In fact, (1.13) is canonical in the following stronger sense: The canonial isomorphism between any two $\Theta_0$ and $\Theta_0´$ versions of $\hat{\mathcal{Z}}_{ech,M}$ intertwines their respective versions of (1.13).

The image via (1.13) of the set $\mathcal{Z}_{ech,M} \times \{-\infty, \ldots, -1\}$ defines a sub-fiber bundle in $\hat{\mathcal{Z}}_{ech,M}$. The latter is denoted by $\hat{\mathcal{Z}}^0_{ech,M}$. The free $\mathbb{Z}$-module generated by the elements in $\hat{\mathcal{Z}}^0_{ech,M}$ plays a central role in Theorem 2.3 of [KLTI]. This submodule is denoted here by $\mathbb{Z}(\hat{\mathcal{Z}}^0_{ech,M})$.

Pairing with the class $c_{1M}$ defines a linear functional from $H_2(M; \mathbb{Z})$ to $2\mathbb{Z}$. Let $p_M \in 2\mathbb{Z}$ denote the divisibility of the subgroup defined by the image. Rules laid out by Hutchings (see [H2]) can be used here to give each generator of $\hat{\mathcal{Z}}_{ech,M}$ a relative $\mathbb{Z}/(p_M\mathbb{Z})$ degree and so give $\mathbb{Z}(\hat{\mathcal{Z}}_{ech,M})$ a relative $\mathbb{Z}/(p_M\mathbb{Z})$ grading.

### c) The almost complex geometry of $\mathbb{R} \times Y$

An endomorphism of $\mathbb{Z}(\hat{\mathcal{Z}}_{ech,M})$ that serves as the embedded contact homology differential is defined using certain sorts of submanifolds in $\mathbb{R} \times Y$. The latter are pseudoholomorphic for a chosen almost complex structure. Part 1 of this subsection describes the allowed almost complex structures. Part 2 of the subsection summarizes some standard definitions.

*Part 1*: Section II.3a and (II.6.1) describe the constraints that delineate the set of allowed almost complex structures on $\mathbb{R} \times Y$. The first two constraints are the $\mathbb{R} \times Y$ versions of standard constraints that are used in all contact and symplectic versions of Floer homology. The remaining constraints are special to the situation at hand. By way of notation, the Euclidean coordinate on the $\mathbb{R}$ factor of $\mathbb{R} \times Y$ is denoted by $s$.

Let J denote a given almost complex structure on $\mathbb{R} \times Y$. This almost complex structure is allowed if it has the properties listed in the seven bullets that follow.

- J *maps the Euclidean tangent vector $\partial_s$ to the $\mathbb{R}$ factor of $\mathbb{R} \times Y$ to* $v$.
- J *is not changed by constant translations along the $\mathbb{R}$ factor of $\mathbb{R} \times Y$.*
- J *preserves the kernel of the 1-form $\hat{a}$; and its restriction to this 2-plane field defines the orientation given by w.*

The next two bullets concern the restriction of J to any given $\mathfrak{p} \in \Lambda$ version of $\mathbb{R} \times \mathcal{H}_\mathfrak{p}$. The statement of the second refers to the vector fields

$$e_1 = -6g\cos\theta\sin\theta\, \partial_u + (x + g')(1 - 3\cos^2\theta)\, \partial_\theta \quad and \quad e_2 = \partial_\phi + \sqrt{6}\chi_\delta f \cos\theta \sin^2\theta\, v \,.$$
(1.14)

These span the kernel of $\hat{a}$ where both u and $(1 - 3\cos^2\theta)\sin\theta$ are not zero.

- J *is unchanged by constant, $\mathbb{R}/2\pi\mathbb{Z}$ translations of the coordinate $\phi$.*
- $Je_1 = \sigma^{-1}e_2$ *where $\sigma$ is a positive function of* u *and* $\theta$.

The final two bullets concerns the restriction of J to $\mathbb{R} \times M_\delta$. The first of these refers to two sets of pairwise disjoint annuli in the Heegaard surface $\Sigma$. The annuli in the first set are labeled by the index 1 critical points of $f$, and those in the second are labeled by the index 2 critical points of $f$. Let p denote a given index 1 or index 2 critical point. The corresponding annulus is denoted respectively by $T_{p+}$ or $T_{p-}$.

To say more about these annuli, let p denote an index 1 critical point of $f$. The annulus $T_{p+}$ is the image via Lie transport along the integral curves of $\mathfrak{v}$ of the annulus in the radius $\delta_*$ coordinate ball centered on p where $1 - 3\cos^2\theta_+ > 0$. By way of a review from [KLTII], the image of the central, $\theta_+ = \frac{\pi}{2}$, circle in this annulus is denoted by $C_{p+}$; this is the intersection between the Heegaard surface $\Sigma$ and the ascending disk from the critical point p. Let $\hat{h}_+$ denote the function $2e^{2t_+}\cos\theta_+ \sin^2\theta_+$ on the radius $8\delta_*$ coordinate ball centered on p. Lie transport by $\mathfrak{v}$ gives $T_{p+}$ coordinates $(\varphi_+, \hat{h}_+)$ with the former being $\mathbb{R}/2\pi\mathbb{Z}$ valued. The restriction of $w$ to $T_{p+} \subset \Sigma$ is given using these coordinates by $\sqrt{6}\,d\varphi_+ \wedge d\hat{h}_+$.

Let p now denote an index 2 critical point of $f$. The annulus $T_{p-}$ is the image via Lie transport via $\mathfrak{v}$ of the annular region where $1 - 3\cos^2\theta_- > 0$ in the boundary of the radius $\delta_*$ coordinate ball centered on p. The image in $T_{p-}$ of the central $\theta_- = \frac{\pi}{2}$ circle is denoted by $C_{p-}$; it is the intersection between $\Sigma$ and the descending disk from p. Set $\hat{h}_- = 2e^{2t_-}\cos\theta_- \sin^2\theta_-$. Lie transport by $\mathfrak{v}$ identifies gives $T_{p-}$ the coordinates $(\varphi_-, \hat{h}_-)$. The 2-form $w$ on $T_{p-}$ is $-\sqrt{6}\,d\varphi_- \wedge d\hat{h}_-$.

There is one more point to note regarding an intersection $T_{p+} \cap T_{p'-}$. The respective coordinates $(\varphi_+, \hat{h}_+)$ for $T_{p+}$ and $(\varphi_-, \hat{h}_-)$ for $T_{p-}$ are related on this intersection by the rule $(d\varphi_+, d\hat{h}_+) = \pm(d\hat{h}_-, d\varphi_-)$ with the + sign taken when the pair of vectors $(\frac{\partial}{\partial\varphi_+}, \frac{\partial}{\partial\varphi_-})$ define an oriented basis for $T\Sigma$ at the corresponding point in $C_{p+} \cap C_{p'-}$.

The union of the annuli in the set labeled by the index 1 critical points of $f$ is denoted by $T_+$, and of the union of the annuli from the set labeled by the index 2 critical points of $f$ is denoted by $T_-$. The union of the index 1 critical point versions of $C_{p+}$ is denoted by $C_+$ and the corresponding union of the index 2 critical point versions of $C_{p-}$ is denoted by $C_-$.

The bullet that follows identifies the $f \in (1, 2)$ part of $M_\delta$ with $(1, 2) \times \Sigma$ in the manner just described.

- $J\frac{\partial}{\partial\varphi_+} = \frac{\partial}{\partial\hat{h}_+}$ *on* $\mathbb{R} \times (1, 2) \times T_+$ *and* $J\frac{\partial}{\partial\varphi_-} = -\frac{\partial}{\partial\hat{h}_-}$ *on* $\mathbb{R} \times (1, 2) \times T_-$.

These two conditions are compatible on $T_+ \cap T_-$ because $(d\varphi_+, d\hat{h}_+) = \pm(d\hat{h}_-, d\varphi_-)$ on any given component of $T_+ \cap T_-$.

The final bullet refers to a certain residual subset in the $C^\infty$ Fréchét space of almost complex structures that obey the preceding bullets. This is the subset $\mathcal{J}_{ech}$ from Theorem II.A.1.

- *J comes from the residual set* $\mathcal{J}_{ech}$.

The membership in the residual set $\mathcal{J}_{ech}$ guarantees the vanishing of the cokernel of a Fredholm operator that is associated to certain sorts of pseudoholomorphic submanifolds.

An almost complex structure that obeys the first three bullets is compatible with the 2-form $\hat{\omega} = ds \wedge \hat{a} + w$. Said differently: The bilinear form $\hat{\omega}(\cdot, J(\cdot))$ defines a Riemannian metric on $\mathbb{R} \times Y$ when J obeys the first three bullets. Note that this metric gives both $\frac{\partial}{\partial s}$ and $v$ norm 1, it makes them mutually orthogonal, and it makes both orthogonal to the kernel of $\hat{a}$. This metric also makes J into an orthogonal endomorphism of $T(\mathbb{R} \times Y)$. With J given, the metric $\hat{\omega}(\cdot, J(\cdot))$ is used implicitly in what follows to define norms and covariant derivatives on the various tensor bundles over $\mathbb{R} \times Y$.

An almost complex structure that obeys all seven of these bullets will be said to be a member of $\mathcal{J}_{ech}$.

*Part 2*: Let J denote an almost complex structure on $\mathbb{R} \times Y$. Assume for the moment that J obeys only the constraints from the first three bullets of Part 1. A proper subset $C \subset \mathbb{R} \times Y$ is said in what follows to be a J-holomorphic subvariety if it has the following properties:

- *C has no isolated points and the complement of a finite set in C is a submanifold with J-invariant tangent space.*
- *The integral of w over C is finite.*

(1.15)

A J-holomorphic subvariety is said to be irreducible if the complement of any given finite set is connected.

A J-holomorphic subvariety may or may not be compact. If not, these conditions have the various standard implications ([HWZ], [S], [HT]) about the large $|s|$ part of the subvariety. To say more, let C denote a given, non-compact pseudoholomorphic subvariety. There exists $s_0 > 1$ such that the $|s| \geq s_0$ part of C is a disjoint union of embedded cylinders. The 1-form $ds$ is non-zero on the tangent space of each such cylinder. A component cylinder of the $|s| \geq s_0$ part of C is said to be an *end* of C. An end of C where $s \geq s_0$ is said to be *positive* and an end where $s \leq -s_0$ is said to be *negative*. A constant $|s| \geq s_0$ slice of any given end is an embedded circle in Y. This circle appears as a braid in a small radius tubular neighborhood of a closed integral curve of $v$. As $|s|$ increases, the circle in question moves via an ambient isotopy so as to converge pointwise as $|s| \to \infty$ as a multiple cover of the central integral curve of $v$. The closed integral curve in question is said to be *associated* to the given end.

The set of J-holomorphic subvarieties is given the topology that associates to any given J-holomorphic subvariety a basis of open neighborhoods of the following sort: Let

C denote the given subvariety. The sets of C's neighborhood basis are labeled by (0, 1). A J-holomorphic subvariety C´ is in a given $\varepsilon \in (0, 1)$ member of this basis when

- $\sup_{z \in C} \text{dist}(z, C´) + \sup_{z \in C´} \text{dist}(C, z) < \varepsilon$.
- *Let $\mu$ denote a smooth 2-form on $\mathbb{R} \times Y$ with $|\mu| \leq 1$, with $|\nabla \mu| < \varepsilon^{-1}$ and with compact support where $|s| < \varepsilon^{-1}$. Then $|\int_C \mu - \int_{C´} \mu| < \varepsilon$.*

(1.16)

The resulting topological space is called the *moduli space* of J-holomorphic subvarieties. The group $\mathbb{R}$ has a continuous action on the moduli space, this given by the constant translations along the $\mathbb{R}$ factor of $\mathbb{R} \times Y$. An irreducible, $\mathbb{R}$-invariant J-holomorphic subvariety is the product of $\mathbb{R}$ with a closed integral curve of $v$.

Of particular interest in much of what follows are the moduli space components that contain elements that are characterized as follows: Let C denote a member.

- *C is embedded.*
- *Distinct ends of C have distinct associated closed integral curves of $v$. This is also the case for distinct negative ends of C.*
- *The constant $|s|$ slice of any given end are isotopic in the tubular neighborhood of the associated integral curve of $v$ to this central integral curve.*
- *The set of integral curves of $v$ that are associated to the postive ends of C defines an element of $\mathcal{Z}_{\text{ech},M}$. This is also true for the negative ends.*

(1.17)

A J-holomorphic subvariety that is described by the second, third and fourth bullets of (1.17) is said to be an *ech-subvariety*. An ech-subvariety is said to be an *ech-HF subvariety* if it lacks irreducible components that intersect $\mathbb{R} \times M_\delta$ in an $f$ = constant level set, or that intersect $\mathbb{R} \times \mathcal{H}_0$ in a u = constant level set, or that intersect some $\mathfrak{p} \in \Lambda$ version of $\mathbb{R} \times \mathcal{H}_\mathfrak{p}$ in the u = 0 level set. Any such forbidden irreducible component is described completely by one of Propositions II.3.1-II.3.4. An ech-HF subvariety is said here to be an *ech-HF* submanifold if it obeys all four of the bullets in (1.17)

Let C denote an ech-subvariety. The element in $\mathcal{Z}_{\text{ech},M}$ that comes from the positive ends of C via the fourth bullet of (1.17) is denoted by $\Theta_{C+}$; the analogous negative end element in $\mathcal{Z}_{\text{ech},M}$ is denoted by $\Theta_{C-}$. The smooth part of C is oriented by J, and so C's image in Y via the projection defines a 2-cycle with boundary $[\Theta_{C+}] - [\Theta_{C-}]$. This 2-cycle is denoted by $[C]_Y$. Meanwhile, Hutchings (see [Hu2]) gives rules for assigning an integer to C, its *ech index*. This ech index is denoted here by $I_{\text{ech}}(C)$.

Let $\hat{\Theta}´$ and $\hat{\Theta}$ denote a given pair from $\hat{\mathcal{Z}}_{\text{ech},M}$ and let k denote a given integer. Use $\mathcal{M}_k(\hat{\Theta}´, \hat{\Theta})$ to denote the set of ech subvarieties with membership characterized as follows: Write $\hat{\Theta}$ as a pair $(\Theta,Z)$ with $\Theta \in \mathcal{Z}_{ech,M}$ and $Z \in H_2(Y;[\Theta] - [\Theta_0])$. Write in the analogous fashion as $\hat{\Theta}' = (\Theta',Z')$. The subvariety C is a member of $\mathcal{M}_k(\hat{\Theta}',\hat{\Theta})$ when

$$\Theta_{C+} = \Theta, \quad \Theta_{C-} = \Theta', \quad Z' = Z - [C]_Y \quad and \quad I_{ech}(C) = k. \tag{1.18}$$

The set $\mathcal{M}_k(\hat{\Theta}',\hat{\Theta})$ is a union of components of the moduli space of J-holomorphic subvarieties. By way of a parenthetical remark, the set $\mathcal{M}_k(\hat{\Theta}',\hat{\Theta})$ is empty unless the sum of k $mod(p_M)$ and the $\mathbb{Z}/(p_M\mathbb{Z})$ grading of $\hat{\Theta}'$ equals the $\mathbb{Z}/(p_M\mathbb{Z})$ grading of $\hat{\Theta}$.

### d) The differential and the geometric endomorphisms of $\mathbb{Z}(\hat{\mathcal{Z}}_{ech,M})$

Part 1 of this subsection supplies a brief description of the differential on $\mathbb{Z}(\hat{\mathcal{Z}}_{ech,M})$ that defines the relevant version of embedded contact homology. As noted in Appendix II.Aa, rules laid out by Hutchings [HS] can be used to define an action of the algebra $\mathbb{Z}[\mathbb{U}] \otimes (\wedge^*(H_1(Y;\mathbb{Z})/torsion))$ on the embedded contact homology $\mathbb{Z}$-module. Part 2 of this subsection describes the generators of this $\mathbb{Z}[\mathbb{U}] \otimes (\wedge^*(H_1(Y;\mathbb{Z})/torsion))$ action. Part 3 talks about the grading of this $\mathbb{Z}$-module.

All that is said in what follows assumes that J comes from $\mathcal{J}_{ech}$.

*Part 1*: An endomorphism of $\mathbb{Z}(\hat{\mathcal{Z}}_{ech,M})$ is given by its action on the generators and the action on any given generator $\hat{\Theta} \in \hat{\mathcal{Z}}_{ech,M}$ results in a formal sum of the form

$$\hat{\Theta} \to \sum_{\hat{\Theta}' \in \hat{\mathcal{Z}}_{ech,M}} N_{\hat{\Theta}',\hat{\Theta}} \hat{\Theta}', \tag{1.19}$$

with each coefficient an integer, and where only finitely many coefficients are non-zero. The collection of integers $\{N_{\hat{\Theta}',\hat{\Theta}}\}_{\hat{\Theta}',\hat{\Theta} \in \hat{\mathcal{Z}}_{ech,M}}$ defines the endomorphism.

Theorem II.A.1 asserts that the endomorphism of $\mathbb{Z}(\hat{\mathcal{Z}}_{ech,M})$ that serves as the differential for embedded contact homology can be defined according to the rules laid out by Hutchings (see [HS], [Hu1], [HT]). These rules are summarized in Part 1 of Section 9b to come. Suffice it to say here that each $\{\hat{\Theta}',\hat{\Theta}\}$ version of the relevant version of $N_{\hat{\Theta}',\hat{\Theta}}$ is computed using the components of $\mathcal{M}_1(\hat{\Theta}',\hat{\Theta})$.

Theorem II.A.1 implies that $\mathcal{M}_1(\hat{\Theta}',\hat{\Theta})$ is a smooth manifold with a finite set of components, each $\mathbb{R}$-equivariantly diffeomorphic to $\mathbb{R}$, and that each component of this space contributes either +1 or -1 to a sum that gives $N_{\hat{\Theta}',\hat{\Theta}}$.

The endomorphism that defines the differential for embedded contact homology is denoted in what follows by $\partial_{ech}$.

*Part 2*: As noted above, the homology of $\partial_{ech}$ is a $\mathbb{Z}$-module with certain canonical endomorphisms that generate an action of $\mathbb{Z}(\mathbb{U}) \otimes (\wedge^*(H_1(Y;\mathbb{Z})/\text{torsion}))$. The generators of this action are defined by endomorphisms of $\mathbb{Z}(\hat{\mathcal{Z}}_{ech,M})$. As explained in Part 3 of Appendix II.Aa, the endomorphism of $\mathbb{Z}(\hat{\mathcal{Z}}_{ech,M})$ that supplies the action of $\mathbb{U}$ on the homology is defined with the help of a chosen point in either $\mathcal{H}_0$ or the part of $M_\delta$ where $f \in (0,1) \cup (2,3)$. Let y denote such a point. It follows from (II.A.6) and Theorem II.A.1 that a given $\hat{\Theta}', \hat{\Theta} \in \hat{\mathcal{Z}}_{ech,M}$ version of the coefficient $N_{\hat{\Theta}', \hat{\Theta}}$ is zero unless $\hat{\Theta}'$ and $\hat{\Theta}$ are related as follows: Write $\hat{\Theta}$ as $(\Theta, Z)$ with $\Theta \in \mathcal{Z}_{ech,M}$ and with $Z \in H_2(Y;[\Theta] - [\Theta_0])$. Then $\hat{\Theta}' = (\Theta, Z - [S])$ where $[S]$ here denotes the u = 0 sphere in $\mathcal{H}_0$. The coefficient $N_{\hat{\Theta}', \hat{\Theta}}$ in this case is 1.

What follows is also a consquence of Theorem II.A.1. The endomorphisms of $\mathbb{Z}(\hat{\mathcal{Z}}_{ech,M})$ that generate the action of $\wedge^*(H_1(Y;\mathbb{Z})/\text{torsion})$ on the homology are defined with the help of a chosen, suitably generic basis of cycles that generate $H_1(Y;\mathbb{Z})/\text{torsion}$. Fix such a basis and let $\hat{\imath} \subset Y$ denote a chosen basis element. Any given $\hat{\Theta}', \hat{\Theta} \in \hat{\mathcal{Z}}_{ech,M}$ coefficient $N_{\hat{\Theta}', \hat{\Theta}}$ in the corresponding version of (1.19) is computed using the submanifolds from $\mathcal{M}_1(\hat{\Theta}', \hat{\Theta})$ that intersect $\{0\} \times \hat{\imath}$ and the corresponding intersection points. If γ is suitably generic, then the set of pairs consisting of a submanifold in $\mathcal{M}_1(\hat{\Theta}', \hat{\Theta})$ and an intersection point with $\{0\} \times \hat{\imath}$ is a finite set. Moreover, each such intersection point contributes either +1 or -1 to a sum that gives $N_{\hat{\Theta}', \hat{\Theta}}$. The upcoming Section 9c explains how these ±1 contributions are determined.

The upcoming Theorem 1.1 refers to an *M-adapted* 1-cycle basis for $H_1(Y;\mathbb{Z})/\text{torsion}$  The definition of term M-adapted requires the introduction from Part 7 in Section II.1c of a certain finite set in the interior of $\Sigma - (T_- \cup T_+)$. The set contains the fiducial point $z_0$ and $\dim(H^1(M;\mathbb{Z}))$ additional points. This set is denoted by ¥. Each $z \in$ ¥ is the intersection point of $\Sigma$ with a closed integral curve of v. The latter curve is denoted by $\gamma^{(z)}$. Pairing with the Poincaré duals of the homology classes of the cycles in the set $\{[\gamma^{(z)}] - [\gamma^{(z_0)}]\}_{z \in ¥ - z_0}$ generates the dual in $\text{Hom}(H_2(Y;\mathbb{Z});\mathbb{Z})$ of the $H_2(M;\mathbb{Z})$ summand in (1.7).

An M-adapted basis is characterized as follows: The basis contains the cycle $[\gamma^{(z_0)}]$, it contains a set of cycles that can be labeled $\{\hat{\imath}^{(z)}\}_{z \in ¥ - z_0}$, and it is rounded out by a set of G cycles that can be labled $\{\hat{\imath}_{\mathfrak{p}}\}_{\mathfrak{p} \in \Lambda}$. A given $z \in ¥ - z_0$ version of $\hat{\imath}^{(z)}$ lies entirely in the $M_{7\delta_*}$ part of Y. It is homologous to $[\gamma^{(z)}] - [\gamma^{(z_0)}]$ and it is obtained from the latter by first truncating the $\mathcal{H}_0$ portions of the curves $\gamma^{(z)}$ and $\gamma^{(z_0)}$ and then reconnecting the respective endpoints by arcs on the boundary of the radius $7\delta_*$ coordinate balls about the index 0 and index 3 critical points of $f$. A given $\mathfrak{p} \in \Lambda$ version of $\hat{\imath}_{\mathfrak{p}}$ is disjoint from the $f \in [1, 2]$ part of $M_{7\delta_*}$, and it intersects the rest of $M_{7\delta_*}$ and $\mathcal{H}_0$ as a smooth curve that is transverse to the level sets of $f$ in $M_\delta$ and the constant u spheres in $\mathcal{H}_0$; the orientation is such that it has intersection number 1 with the u = 0 sphere in $\mathcal{H}_0$. Meanwhile, $\hat{\imath}_\mathfrak{p}$ intersects $\cup_{\mathfrak{p}' \in \Lambda} \mathcal{H}_{\mathfrak{p}'}$ as the $\theta = 0$ arc in $\mathcal{H}_\mathfrak{p}$; and the orientation is such that it has intersection number -1 with each u = 0 sphere in $\mathcal{H}_\mathfrak{p}$.

**e) The Heegaard Floer equivalence**

A three part digression follows directly to set the notation for Theorem 1.1.

*Part 1*: Use (1.10) and (1.13) to write $\hat{\mathcal{Z}}_{\text{ech},M}$ as

$$\hat{\mathcal{Z}}_{\text{ech},M} = \mathcal{Z}_{\text{ech},M} \times \mathbb{Z} = (\mathcal{Z}_{\text{HF}} \times (\times_{\mathfrak{p} \in \Lambda}(\mathbb{Z} \times \text{o}))) \times \mathbb{Z} \ .$$

(1.20)

The principle $\mathbb{Z}$ bundle action of $\mathbb{Z}$ acts on the right most $\mathbb{Z}$. This factor is now moved next to the $\mathcal{Z}_{\text{HF}}$ factor to write (1.20) as

$$\hat{\mathcal{Z}}_{\text{ech},M} = (\mathcal{Z}_{\text{HF}} \times \mathbb{Z}) \times (\times_{\mathfrak{p} \in \Lambda}(\mathbb{Z} \times \text{o})) \ .$$

(1.21)

The identification in (1.21) induces the tensor product decomposition

$$\mathbb{Z}(\hat{\mathcal{Z}}_{\text{ech},M}) = \mathbb{Z}(\mathcal{Z}_{\text{HF}} \times \mathbb{Z}) \otimes (\otimes_{\mathfrak{p} \in \Lambda} \mathbb{Z}(\mathbb{Z} \times \text{o})) \ .$$

(1.22)

This representation of $\mathbb{Z}(\hat{\mathcal{Z}}_{\text{ech},M})$ is used implicitly by Theorem 1.1

*Part 2*: The factor $\mathbb{Z}(\mathcal{Z}_{\text{ech},M} \times \mathbb{Z})$ is the $\mathbb{Z}$-module for the Heegaard Floer homology on M. The endomorphism that supplies the differential for this homology is denoted by $\partial_{\text{HF}}$. Theorem 1.1 describes the differential $\partial_{\text{ech}}$ on $\mathbb{Z}(\hat{\mathcal{Z}}_{\text{ech},M})$ in terms of $\partial_{\text{HF}}$ and certain endomorphisms that are induced on (1.22) by corresponding endomorphisms of the various $\mathfrak{p} \in \Lambda$ factors of $\mathbb{Z}(\mathbb{Z} \times \text{o})$. A given $\mathfrak{p} \in \Lambda$ version is denoted by $\partial_\mathfrak{p}$. All are the same endomorphism of $\mathbb{Z}(\mathbb{Z} \times \text{o})$, this being the endomorphism $\partial_*$ that acts on the generating set as follows:

- $\partial_*(k,0) = 0$ *for each* $k \in \mathbb{Z}$.
- $\partial_*(k,1) = (k,0) + (k+1,0)$ *for each* $k \in \mathbb{Z}$.
- $\partial_*(k,-1) = (k,0) + (k-1,0)$ *for each* $k \in \mathbb{Z}$.
- $\partial_*(k,\{1,-1\}) = (k,-1) + (k+1,-1) - (k,1) - (k-1,1)$ *for each* $k \in \mathbb{Z}$.

(1.23)

As noted by Lemma 2.5 in [KLT1], the homology of the chain complex $(\mathbb{Z}(\mathbb{Z}\times O), \partial_*)$ is $\mathbb{Z} \oplus \mathbb{Z}$, and generators are the closed elements (0,0) and (0,1)-(1,-1).

The various versions of Heegard Floer homology enjoy an action of the algebra $\mathbb{Z}[\mathbb{U}] \otimes \wedge^*(H_1(M;\mathbb{Z})/\text{torsion})$ whose generators can be defined by endomorphisms of $\mathbb{Z}(\mathcal{Z}_{ech,M} \times \mathbb{Z})$. Note in this regard that the generator of the action of $\mathbb{Z}[\mathbb{U}]$ sends any given pair $(\hat{\upsilon}, k) \in \mathcal{Z}_{ech,M} \times \mathbb{Z}$ to $(\hat{\upsilon}, k-1)$.

There is one other Heegaard Floer endomorphism that plays a role in what follows. The latter is defined by its action on the generators, and in doing so, it acts solely on the $\mathcal{Z}_{HF}$ and ignores the $\mathbb{Z}$ factor. This is the endomorphism that appears in Theorem 4.1 of [OS1] and Definition 8.1 of [OS2]. The latter is denoted here by $\partial_{HF0}$.

The upcoming formula for $\partial_{ech}$ and the other endomorphisms of (1.22) use the following convention: Suppose that E and E´ are graded chain complexes and that $\Delta$ and $\Delta´$ are respective graded endomorphisms of E and E´. The latter induce on $E \otimes E´$ an endomorphism, $\Delta + \Delta´$, that is defined by the following action on the reducible elements: Let $\mathfrak{e}$ and $\mathfrak{e}´$ denote respective elements of E and E´. Then $(\Delta + \Delta´)(\mathfrak{e} \otimes \mathfrak{e}´)$ is defined to be $\Delta\mathfrak{e} \otimes \mathfrak{e}´ + (-1)^{\deg ree(\Delta´)degree(\mathfrak{e})} \mathfrak{e} \otimes \Delta´\mathfrak{e}´$.

*Part 3*: Let $p_M \in 2\mathbb{Z}$ again denote the greatest divisor of the image of $H_2(M;\mathbb{Z})$ in $\mathbb{Z}$ via the pairing homomorphism with $c_{1M}$. As noted previously, the $\mathbb{Z}$-module $\mathbb{Z}(\hat{\mathcal{Z}}_{ech,M})$ has a relative $\mathbb{Z}/(p_M\mathbb{Z})$ grading that is induced by a relative grading of its generators. The grading difference between given generators $\hat{\Theta}´, \hat{\Theta}$ is denoted in what follows by $gr_{ech}(\hat{\Theta}´) - gr_{ech}(\hat{\Theta})$.

As explained in [OS], the $\mathbb{Z}$-module $\mathbb{Z}(\mathcal{Z}_{HF} \times \mathbb{Z})$ has a relative $\mathbb{Z}/(p_M\mathbb{Z})$ grading that is induced by a relative $\mathbb{Z}/(p_M\mathbb{Z})$ grading of the set $\mathcal{Z}_{HF}$. The difference between the respective gradings of given elements $\hat{\upsilon}´, \hat{\upsilon} \in \mathcal{Z}_{HF}$ is denoted by $\deg_{HF}(\hat{\upsilon}´) - \deg_{HF}(\hat{\upsilon})$. Granted this notation, the difference between the gradings of corresponding elements $(\hat{\upsilon}´, k´)$ and $(\hat{\upsilon}, k)$ from $\mathcal{Z}_{HF} \times \mathbb{Z}$ is $\deg_{HF}(\hat{\upsilon}´) - \deg_{HF}(\hat{\upsilon}) + 2(k´ - k)$.

The module $\mathbb{Z}(\mathbb{Z} \times O)$ has an absolute $\mathbb{Z}$ grading with values in the 3-element set $\{0, 1, 2\}$. The latter grading is induced by a grading of the generators that depends only on the factor $O = \{0, 1, -1, \{1, -1\}\}$: The element 0 has grading zero, the elements -1 and 1 have grading 1, and the element $\{-1, 1\}$ has grading 2. The resulting grading map from $\mathbb{Z} \times O$ to $\{0, 1, 2\}$ is denoted by $gr_o()$.

With the preceding as background, what follows is this paper's central result.

**Theorem 1.1**: *Identify $\mathbb{Z}(\hat{\mathcal{Z}}_{ech,M})$ with $\mathbb{Z}(\mathcal{Z}_{HF} \times \mathbb{Z}) \otimes (\otimes_{\mathfrak{p} \in \Lambda} \mathbb{Z}(\mathbb{Z} \times O))$ as in (1.22).*
- *The differential $\partial_{ech}$ appears as $\partial_{ech} = \partial_{HF} + \sum_{\mathfrak{p} \in \Lambda} \partial_\mathfrak{p}$.*

- *The $\mathbb{U}$-map acts as the map $((\hat{\upsilon},k),(\mathfrak{k}_p,o_p)_{p\in\Lambda}) \to ((\hat{\upsilon},k-1),(\mathfrak{k}_p,o_p)_{p\in\Lambda})$.*
- *Use the M-adapted 1-cycle basis $\{[\gamma^{(z_0)}],\{\hat{\imath}^{(z)}\}_{z\in¥-z_0},\{\hat{\imath}_p\}_{p\in\Lambda}\}$ to define endomorphisms of $\mathbb{Z}(\hat{\mathcal{Z}}_{ech,M})$.*
    a) *The endomorphism defined by $[\gamma^{z_0}]$ acts as $\partial_{HF}-\partial_{HF0}$.*
    b) *The endomorphisms defined by cycles from $\{\hat{\imath}^{(z)}\}_{z\in¥-z_0}$ act only on the $\mathbb{Z}(\mathcal{Z}_{HF}\times\mathbb{Z})$ factor. In doing so, they induce a set of generators of the $\wedge^*(H_1(M;\mathbb{Z})/\text{torsion})$ action on the Heegaard Floer homology.*
    c) *The endomorphism that is defined by any given $\mathfrak{p}\in\Lambda$ version of $\hat{\imath}_p$ acts as $\mathbb{I}_p+\hat{\imath}_p$ where $\mathbb{I}_p$ acts as the identity on the factors $\mathbb{Z}(\mathcal{Z}_{HF}\times\mathbb{Z})\otimes(\otimes_{p'\in\Lambda-p}\mathbb{Z}(\mathbb{Z}\times o))$ and $\hat{\imath}_p$ is the degree -1 endomorphism that acts only on $\mathfrak{p}$'s factor of $\mathbb{Z}(\mathbb{Z}\times o)$. It acts on this factor as the endomorphism that sends $(\mathfrak{k}_p,o_p)$ to $(\mathfrak{k}_p,o_p')$ with coefficient either 1 or 0. The coefficient 1 appears if and only if both $o_p=1$ and $o_p'=0$, or both $o_p=\{1,-1\}$ and $o_p'=-1$.*
- *Let $\hat{\Theta}'=((\hat{\upsilon}',k'),(\mathfrak{k}_p',o_p')_{p\in\Lambda})$ and $\hat{\Theta}=((\hat{\upsilon},k),(\mathfrak{k}_p,o_p)_{p\in\Lambda})$ denote any two elements. Then*

$$gr_{ech}(\hat{\Theta}')-gr_{ech}(\hat{\Theta})=gr_{HF}(\hat{\upsilon}')-gr_{HF}(\hat{\upsilon})+2(k'-k)+\sum_{p\in\Lambda}(gr(o_p')-gr(o_p))$$

The subsequent sections in this article contain the proof of Theorem 1.1. The next subsection gives an indication of what the proof involves.

**f) A look ahead at the proof**

Three fundamental observations serve as the foundation for the proof of Theorem 1.1. The first is provided by Robert Lipshitz [L] and his theorem to the effect that the differential for Heegaard Floer homology can be defined using certain sorts of pseudoholomorphic subvarieties that reside in the $f^{-1}(1,2)$ part of $\mathbb{R}\times M_\delta$. These are described in Section II.6 and their properties are summarized in the next subsection. The second observation is supplied by Propositions II.7.2 and II.7.3. The latter assert that the $\mathbb{R}\times M_\delta$ part of any of the relevant ech-HF subvariety looks very much like a subvariety of the sort considered by Lipshitz. The third observation is jointly supplied by Propositions II.4.5 and II.5.8. These two propositions jointly hint at a canonical form for the $\mathbb{R}\times(\cup_{p\in\Lambda}\mathcal{H}_p)$ part of any given ech-HF subvariety. The subsequent proof of Theorem 1.1 uses this view of an ech-HF subvariety as the union of a Heegaard Floer looking $\mathbb{R}\times M_\delta$ part and a roughly canonical $\mathbb{R}\times(\cup_{p\in\Lambda}\mathcal{H}_p)$ part to derive the decomposition given by the first bullet of Theorem 1.1, and likewise to prove the assertions of the remaining bullets.

A proof of Theorem 1.1 along the preceding lines must address the following fundamental question:

> *Fix a subvariety from* [L]. *As noted above, there is some set of ech-HF subvarieties that look much like it on* $\mathbb{R} \times M_\delta$. *What can be said about this set; in particular, can enough be said to justify the claims of Theorem 1.1?*

As it turns out, only submanifolds need be considered, and the upcoming sections study the question just posed with regards to submanifolds. This is done by constructing the appropriate set of ech-HF submanifolds from a given submanifold from [L]. The construction has two parts. The following two parts of this subsection says a few things about the two parts of the construction and how they lead to Theorem 1.1.

*Part 1*: The first part of the construction starts with a submanifold from [L] and a suitably compatible pair of elements $\hat{\Theta}'$, $\hat{\Theta} \in \hat{\mathcal{Z}}_{\text{ech},M}$; it then uses this data to build a canonical approximation to what would be an ech-HF submanifold from the moduli space $\mathcal{M}_1(\hat{\Theta}', \hat{\Theta})$. This approximation exploits the dichotomy between what is said in Propositions II.4.5 and II.5.8 and what is said in Proposition II.7.2 about the $\mathbb{R} \times (\cup_\mathfrak{p} \mathcal{H}_\mathfrak{p})$ and $\mathbb{R} \times M_\delta$ parts of an ech-HF submanifold. In particular, the approximation consists of a set $\mathcal{C}_0 = \{C_{S0}, \{C_{p0}\}_{\mathfrak{p} \in \Lambda}\}$ where $C_{S0}$ denotes a J-holomorphic submanifold with boundary in $\mathbb{R} \times M_\delta$, and where any given $\mathfrak{p} \in \Lambda$ version of $C_{p0}$ denotes a J-holomorphic submanifold with boundary in $\mathbb{R} \times \mathcal{H}_\mathfrak{p}$. The submanifold $C_{S0}$ looks very much like one of the subvarieties from [L]; and each $\mathfrak{p} \in \Lambda$ version of $C_{p0}$ is described by Propositions II.4.5 and II.5.8. The submanifold $C_{S0}$ has 2G boundary components, one on a certain $f \in (1+\delta^2, 1+\delta_*^2)$ level set in each $\mathfrak{p} \in \Lambda$ version of $\mathbb{R} \times \mathcal{H}_\mathfrak{p}$ and the other on a certain $f \in (2-\delta_*^2, 2-\delta^2)$ level set in each $\mathfrak{p} \in \Lambda$ version of $\mathbb{R} \times \mathcal{H}_\mathfrak{p}$. Meanwhile, any given $\mathfrak{p} \in \Lambda$ version of $C_{p0}$ has two boundary components, one on each of these same level sets of $f$ in $\mathbb{R} \times \mathcal{H}_\mathfrak{p}$. However, the boundary components of $C_{p0}$ need not agree with the corresponding $C_{S0}$ boundary components on the relevant level sets of $f$.

The set of such approximations to would-be ech-HF submanifolds can be used to define an ersatz version of $\mathcal{M}_1(\hat{\Theta}', \hat{\Theta})$. This ersatz version can then be used to define coefficients of endomorphisms of $\mathbb{Z}(\hat{\mathcal{Z}}_{\text{ech},M})$ using Hutching's rules. To say a bit more, remark that the definition of the endomorphism coefficients using honest ech-submanifolds is along standard symplectic field theory lines in the sense that a family of Fredholm operators and a certain tautological $\mathbb{R}$ action play the central roles. The space $\mathcal{M}_1(\hat{\Theta}', \hat{\Theta})$ parametrizes the right sort of family; and the tautological $\mathbb{R}$ action is induced by the constant translations along the $\mathbb{R}$ factor of $\mathbb{R} \times Y$. Granted this remark about the definitions, what follows is a key point: The corresponding ersatz moduli space that is constructed from the canonical approximations to ech-HF submanifolds has an analogous family of Fredholm operators and an analogous action of $\mathbb{R}$. This being the case, Hutching's rules can also be used with the ersatz moduli spaces to define endomorphisms $\mathbb{Z}(\hat{\mathcal{Z}}_{ech,M})$. Meanwhile, the canonical nature of the construction guarantees that the resulting versions of the endomorphisms relevant to Theorem 1.1 satisfy the conclusions of Theorem 1.1.

Section 2 describes in detail the canonical approximations to ech-HF submanifolds; Sections 2-6 construct them.

*Part 2*: Part 2 of the construction builds a cobordism between the ersatz version of a given $\mathcal{M}_1(\hat{\Theta}', \hat{\Theta})$ and the version with honest ech-HF submanifolds. The cobordism maps to the interval [0, 1] with the inverse image of 0 giving the ersatz moduli space and that of 1 giving the version with honest ech-submanifolds. The cobordism defines a smooth manifold with boundary such that the map to [0, 1] is proper and smooth. A key point here is that the relevant family of Fredholm operators extends across the cobordism, as does the relevant $\mathbb{R}$ action. Given that the approximation versions of Theorem 1.1's endomorphisms obey Theorem 1.1's conclusions, these last facts are seen to imply that Theorem 1.1's conclusions must hold with its endomorphisms defined using honest ech-submanifolds.

To say a bit more about this cobordism, consider for a moment one of the approximation sets, $\mathcal{C}_0 = \{C_{S0}, \{C_{\mathfrak{p}0}\}_{\mathfrak{p} \in \Lambda}\}$. As noted in Part 1, its elements are J-holomorphic manifolds with boundary with the boundaries lying on certain level sets of $f$. The boundary of $C_{S0}$ is determined solely by the given subvariety from [L]. The boundary of any given $\mathfrak{p} \in \Lambda$ version of $C_{\mathfrak{p}0}$ is constrained in part by that of $C_{S0}$. Keeping this in mind, let $\tau$ denote the parameter in [0, 1]. The inverse image of $\tau$ in the cobordism consists of a set of the form $\mathcal{C} = \{C_S, \{C_\mathfrak{p}\}_{\mathfrak{p} \in \Lambda}\}$ where $C_S$ is a J-holomorphic submanifold with boundary in $\mathbb{R} \times M_\delta$, and where each $\mathfrak{p} \in \Lambda$ version of $C_\mathfrak{p}$ is a J-holomorphic submanifold with boundary in $\mathbb{R} \times \mathcal{H}_\mathfrak{p}$. The submanifold $C_S$ has 2G boundary components, these on the afore-mentioned level sets of $f$ in $\cup_{\mathfrak{p} \in \Lambda}(\mathbb{R} \times \mathcal{H}_\mathfrak{p})$. Meanwhile, each $\mathfrak{p} \in \Lambda$ version of $C_\mathfrak{p}$ has two, one each on the $\mathbb{R} \times \mathcal{H}_\mathfrak{p}$ parts of these level sets. The parameter $\tau$ indicates the extent to which the two boundary components of any given $\mathfrak{p} \in \Lambda$ version of $C_\mathfrak{p}$ agree with the relevant pair of boundary components of $C_S$. In the case when $\tau = 1$, they match up and so define an honest ech-subvariety. This is not necessarily true for $\tau < 1$.

By way of a hint as to the nature of the family of Fredholm operators, the operator for a given $\tau \in [0, 1]$ version of $\mathcal{C} = \{C_S, \{C_\mathfrak{p}\}_{\mathfrak{p} \in \Lambda}\}$ is viewed as a set of G+1 operators, with one defined by $C_S$ and one by each $\mathfrak{p} \in \Lambda$ version of $C_\mathfrak{p}$. Thus, each is defined on a manifold with boundary and as such, its definition requires the specification of some boundary conditions. These are local (as opposed to spectral) boundary conditions that couple the $C_S$ operator to those defined by the various $\mathfrak{p} \in \Lambda$ versions of $C_\mathfrak{p}$ so as to associate to $\mathcal{C}$ a single Fredholm operator. The parameter $\tau$ determines the degree of coupling.

Section 7 constructs the cobordism space that interpolates between the ersatz moduli space and the space of ech-HF submanifolds.

Section 8 supplies the background needed to use the cobordism to compute the differential and other endomorphisms that appear in Theorem 1.1. Section 9 uses the properties of the cobordism space to complete the proof of Theorem 1.1

**g) The subvarieties used by Lipshitz**

This subsection summarizes some of what is said in Section II.6 about the subvarieties that are used by Lipshitz in [L]. These are subvarieties in the $f \in (1, 2)$ part of $\mathbb{R} \times M$ that are best described by viewing this part of $M$ as $(1,2) \times \Sigma$ via the identification given by Lie transport with the pseudogradient vector field $\mathfrak{v}$.

The relevant subvarieties in $\mathbb{R} \times (1,2) \times \Sigma$ are pseudoholomorphic for an almost complex structure with certain special properties. These are described in Section 1 of [L]. Let J denote an almost complex structure on $\mathbb{R} \times Y$ that obeys the constraints in Part 1 of Section 1c. The restriction of J to the $f \in (1,2)$ part of $\mathbb{R} \times M_\delta$ can be extended to the whole of $\mathbb{R} \times (1, 2) \times \Sigma$ so as to give an almost complex structure of the sort considered by Lipshitz, and in particular, of the sort that is described by (II.6.1). Conversely, a suitably generic almost complex structure on $\mathbb{R} \times (1,2) \times \Sigma$ that obeys (II.6.1) will serve for Lipshitz. Moreover, such an almost complex structure will restrict to the $\mathbb{R} \times M_\delta$ part of $\mathbb{R} \times (1,2) \times \Sigma$ as the restriction of an almost complex structures on $\mathbb{R} \times Y$ that obeys the constraints in Section 1c. This understood, let J denote an almost complex structure on $\mathbb{R} \times Y$ that obeys the constraints in Part 1 of Section 1c and let $J_{HF}$ denote an almost complex structure on $\mathbb{R} \times (1,2) \times \Sigma$ that obeys (II.6.1). To say more about $J_{HF}$, note that the stable Hamiltonian 2-form $w$ appears on the $M_\delta$ part of $(1,2) \times \Sigma$ as the pull-back via the projection of an area form on $\Sigma$. Denote the latter by $w_\Sigma$. This form extends in the obvious way to the whole of $(1,2) \times \Sigma$. Let $t$ denote the Euclidean coordinate on $(1,2)$. The almost complex structure $J_{HF}$ maps $\partial_s$ to $\partial_t$, it preserves the level sets of $t$ and it is compatible with the symplectic from $ds \wedge dt + w_\Sigma$. It also commutes with the $\partial_s$ Lie derivative and it obeys the constraint given by the sixth bullet in Part 1 of Section 1c.

Lipshitz considers $J_{HF}$-holomorphic subvarieties in $\mathbb{R} \times (1,2) \times \Sigma$ with eight special properties that are listed momentarily. For the purposes at hand, it is sufficient to consider the case where the subvariety in question is a smooth submanifold. The closure in $\mathbb{R} \times [1,2] \times \Sigma$ of a submanifold with these properties is said here to be a *Lipshitz submanifold*. Let $S_0$ denote interior of a Lipshitz submanifold

PROPERTY 1: *The integral over $S_0$ of $w_\Sigma$ is finite. This is also the case for the integral of $ds \wedge dt$ over any subset of $S_0$ with bounded image in the $\mathbb{R}$ factor of $\mathbb{R} \times (1,2) \times \Sigma$.*

The second property refers to the G circles in $\Sigma$ that comprise the latter's intersection with the ascending disks from the index 1 critical points of $f$, and the corresponding set of G circles that comprise $\Sigma$'s intersection with the descending disks from the index 2 critical points of $f$. If p is a given index 1 or index 2 critical point of $f$, the corresponding circle is denoted by $C_{p+}$ or $C_{p-}$ as the case may be. The union of the index 1 critical point versions of $C_{p+}$ is denoted by $C_+$ and the union of the index 2 critical point versions is denoted by $C_-$.

PROPERTY 2: *The $J_{HF}$-holomorphic submanifold $S_0$ is the interior of a properly embedded surface $\mathbb{R} \times [1,2] \times \Sigma$ with 2G boundary components. Half of the boundary components are in $\mathbb{R} \times \{1\} \times C_+$ and no two of these lie in the same component. The other half are in $\mathbb{R} \times \{2\} \times C_-$ and likewise, no two are in the same component.*

The surface with boundary in PROPERTY 2 is denoted by S. If p is an index 1 or index 2 critical point of $f$, then the corresponding boundary component of S is denoted by $\partial_p S$. It is a properly embedded copy of $\mathbb{R}$ in $\mathbb{R} \times \{1\} \times C_{p+}$ or $\mathbb{R} \times \{2\} \times C_{p-}$ as the case may be.

With regard to notation, Lipshitz and also Section II.6 view what is denoted here by S as the image of a complex surface via a $J_{HF}$-holomorphic map, $u$. What is denoted by S here is denoted in [L] and in Section II.6 by $u(S)$.

The third property refers to elements from the set $\mathcal{Z}_{HF}$. By way of a reminder, an element of $\mathcal{Z}_{HF}$ consists of a set of G integral curves of $\mathfrak{v}$ with each starting at an index 1 critical point of $f$ and ending at an index 2 critical point of $f$. Moreover, distinct curves from such a set have distinct starting points and distinct ending points.

PROPERTY 3: *The surface S is the complement of 2G points in a compact surface with boundary. The function s on S increases with no finite limit on sequences that limit to G of these points, and it decreases with no finite limit on sequences that limit to the remaining G points.*

This surface with boundary in question is denoted by $\underline{S}$. The G points of $\underline{S}-S$ with neighborhoods where s is unbounded from above are said to be *positive points*, and the remaining points are said to be *negative* points.

The remaining properties of S are all consequences of the first three. The next two properties restate Lemmas II.6.2 and II.6.3. They refers to the coordinates $(\varphi_+, h_+)$ for any given component of $T_+$ and to the coordinates $(\varphi_-, h_-)$ for any given component of $T_-$. The $h_+ = 0$ or $h_- = 0$ locus in the relevant component of $T_+$ or $T_-$ is the corresponding component of $C_+$ or $C_-$.

PROPERTY 4: *Let p denote either an index 1 or index 2 critical point of f. The closure in $\underline{S}$ of the corresponding boundary component $\partial_p S$ adds one positive point and one negative point from $\underline{S}-S$. Meanwhile, $\partial_p S$ appears in $\mathbb{R} \times C_{p+}$ or $\mathbb{R} \times C_{p-}$ as a graph over $\mathbb{R}$ of the form $x \to (s = x, \varphi_+ = \varphi^{S,p}(x))$ or $x \to (s = x, \varphi_- = \varphi^{S,p}(x))$ as the case may be. In either case, $\varphi^{S,p} \colon \mathbb{R} \to \mathbb{R}$ is a smooth map with bounded derivatives to any given order. Moreover, the $x \to \pm\infty$ limits of $\varphi^{S,p}$ exist and both are in $C_+ \cap C_-$.*

The next property describes the behavior of S near any given boundary component.

PROPERTY 5: *There exists $z_S > 0$ and $\kappa_S > 1$ with the following significance: Let p denote either an index 1 or index 2 critical point of f. Then a neighborhood of $\partial_p S$ in S appears as the image of a map from $\mathbb{R} \times (0, z_S)$ to $\mathbb{R} \times (1, 1+z_S) \times T_{p+}$ or $\mathbb{R} \times (2 - z_S, 2) \times T_{p-}$ as the case may be. This map has the form*

- $(x, z) \to (s = x, t = 1+z, \varphi_+ = \varphi(x, z), h_+ = \varsigma(x, z))$ *when p has index 1.*
- $(x, z) \to (s = x, t = 2 - z, \varphi_- = \varphi(s, z), h_- = \varsigma(s, z))$ *when p has index 2,*

*where $\varphi(\cdot)$ and $\varsigma(\cdot)$ are maps from $\mathbb{R} \times [0, z_p)$ to $\mathbb{R}$ that obey*

- $|\varsigma(x, z)| + z^{-1} |\varphi(x, z) - \varphi^{S,p}(x)| < \kappa_S z$,
- $|\partial_x \varsigma(s, z)| + z^{-1} |(\partial_x \varphi)(s, z) - (\partial_x \varphi^{S,p})(x)| < \kappa_S z$.

*In addition, the pair $\varphi$ and $\varsigma$ have bounded derivatives to any given order on $\mathbb{R} \times [0, z_S)$.*

The next property says more about the large $|s|$ part of S.

PROPERTY 6: *There exists $\kappa_S > 1$ such that the $s \leq -\kappa_S$ and $s \geq \kappa_S$ portions of S are disjoint unions of G half open rectangles. Those where $s \leq -\kappa_S$ are properly embedded submanifolds in $(-\infty, \kappa_S] \times [1, 2] \times (T_- \cap T_+)$ that appear as a graph over $(-\infty, \kappa_S] \times [1, 2]$ of a map to a component of $T_- \cap T_+$ with the following properties: Let q denote the point in $C_- \cap C_+$ that lies in the given component and let $\psi$ denote the map. Then $\mathrm{dist}(\psi, q) \leq e^{-|s|/\kappa_S}$. The derivatives of $\psi$ to any given order are also bounded by a constant times this* *same exponential factor. Meanwhile, the components of the $s \geq \kappa_S$ have an analogous description as a graph over $[\kappa_S, \infty) \times [1, 2]$.*

This property and PROPERTY 2 lead directly to the next property.

PROPERTY 7: *The set of constant $s \in \mathbb{R}$ slices of S when viewed in M converge pointwise as $s \to \infty$ to define an element in $\mathcal{Z}_{HF}$. This is also true of the $s \to -\infty$ limit of the constant s slices.*

These two elements in $\mathcal{Z}_{HF}$ are denoted respectively by $\hat{\upsilon}_+$ and $\hat{\upsilon}_-$.

The set of Liphsitz submanifolds is given the topology whereby an open neighborhood of a given subvariety S has a basis of open sets labeled by the positive numbers. The open set labeled by $\varepsilon \in (0, \infty)$ is characterized as follows: A submanifold S´ is a member when

- $\sup_{z \in S} \text{dist}(z, S´) + \sup_{z \in S´} \text{dist}(S, z) < \varepsilon$
- *Let $\mu$ denote a smooth 2-form on $[-\frac{1}{\varepsilon}, \frac{1}{\varepsilon}] \times [1,2] \times \Sigma$ with compact support, with supremum norm 1 and with $|\nabla \mu| \leq \frac{1}{\varepsilon}$. Then $|\int_{S´} \mu - \int_S \mu| \leq \varepsilon$.*

(1.24)

The resulting topological space is denoted by $\mathcal{A}_{HF}$.

The group $\mathbb{R}$ acts continuously on $\mathcal{A}_{HF}$ via its action on $\mathbb{R} \times (1,2) \times \Sigma$ as the group of constant translations along the $\mathbb{R}$ factor. This $\mathbb{R}$ action is free on the complement of the set $\{\mathbb{R} \times (\cup_{\upsilon \in \hat{\upsilon}} \upsilon)\}_{\hat{\upsilon} \in \mathcal{Z}_{HF}}$ of 1-point components $\mathcal{A}_{HF}$. Lemmas II.6.6 and II.6.7 say more about the structure of $\mathcal{A}_{HF}$. These lemmae refer to a certain $\mathbb{R}$-linear, Fredholm incarnation of the $\bar{\partial}$-operator that is canonically associated to any given Lipshitz subvariety. The relevant operator is described in II.6e for the case when the variety in question is a submanifold. The operator for a Lipshitz submanifold S is denoted by $D_S$. Let $N_S \to S$ denote the complex normal bundle of S, with the complex structure defined by $J_{HF}$ and with the Hermitian structure and thus holomorphic structure defined by the metric $w_\Sigma(\cdot, J_{HF}(\cdot))$. Let $T^{0,1}S$ denote the (0, 1) part of $T^*S \otimes \mathbb{C}$. The operator $D_S$ maps sections of $N_S$ to sections of $N_S \otimes T^{0,1}S$ by the rule

$$\eta \to D_S \eta = \bar{\partial}\eta + \nu\eta + \mu\bar{\eta}$$

(1.25)

where $\nu$ denotes a certain section of $T^{0,1}S$ and $\mu$ denotes a section of $N_S^2 \otimes T^{0,1}S$.

To say something about the Fredholm domain, remark that PROPERTY 4 can be used as in Section II.6e to identify the bundle $N_S$ along the boundary of S with $T\Sigma$ along

C$_+$ and C$_-$. This understood, the Fredholm domain is the L$^2_1$ completion of the space of compactly supported sections of N$_S$ that obey

- η ∈ TC$_+$ *on the t = 1 part of the boundary of* S.
- η ∈ TC$_-$ *on the t = 2 part of the boundary of* S.

(1.26)

Note in this regard that the definition of the L$^2_1$ norm uses the metric $w_\Sigma(\cdot, J_{HF}(\cdot))$ to define the integration measure, inner products and covariant derivative on all tensor bundles constructed from N$_S$ and TS. Meanwhile the range space for this Fredholm incarnation of D$_S$ is the space of square integrable sections of N$_S \otimes T^{0,1}S$.

The final property speaks of this operator D$_S$.

PROPERTY 8: *The operator* D$_S$ *has trivial cokernel.*

Given this last property, it follows from Lemma II.6.7 that the subspace of Lipshitz submanifolds in $\mathcal{A}_{HF}$ has the structure of a smooth manifold whose dimension near any given submanifold S is the Fredholm index of D$_S$.

### h) Coordinates for the 1 - 3cos$^2$θ > 0 part of $\mathbb{R} \times \mathcal{H}_\mathfrak{p}$

The upcoming construction of ech-HF submanifolds exploits the parametrization of the 1 - 3cos$^2$θ > 0 part each $\mathfrak{p} \in \Lambda$ version of $\mathbb{R} \times \mathcal{H}_\mathfrak{p}$ from Part 1 of Section II.4c. The parametrization is denoted by $\Psi_\mathfrak{p}$. The three parts of this subsection that follow define $\Psi_\mathfrak{p}$ and list some of its important features.

*Part 1*: Fix $\mathfrak{p} \in \Lambda$. The upcoming description of $\Psi_\mathfrak{p}$ requires introducing the coordinates (u, θ, φ) for $\mathcal{H}_\mathfrak{p}$ and the function $\hat{h}$ of the variables u and θ given by

$$\hat{h} = f(u)\cos\theta \sin^2\theta$$

(1.27)

with f as defined in (1.4). The 1-form $d\hat{h}$ is nowhere zero where 1 - 3cos$^2$θ > 0. This function is also annihilated by the vector field $v$ and so it has constant value along $v$'s integral curves in $\mathcal{H}_\mathfrak{p}$.

The definition of $\Psi_\mathfrak{p}$ also involves the J-holomorphic submanifolds from Proposition II.3.2's space $\mathcal{M}_\Sigma$ and Proposition II.3.4's space $\mathcal{M}_{\mathfrak{p}0}$. By way of a reminder, the space $\mathcal{M}_\Sigma$ is $\mathbb{R}$-equivariantly diffeomorphic to $\mathbb{R} \times (1, 2)$. Each element is a compact submanifold that is diffeomorphic to Σ. A given $(s, t) \in \mathbb{R} \times [1+\delta^2, 2-\delta^2]$ element is the $(s = s, f = t)$ slice of $\mathbb{R} \times M_\delta$. An element parameterized by $\mathbb{R} \times (1, 1+7\delta_*^2)$ intersects the u > 0 portion of each $\mathfrak{p} \in \Lambda$ version of $\mathbb{R} \times \mathcal{H}_\mathfrak{p}$ where $1 - 3\cos^2\theta > 0$; and an element in the $\mathbb{R} \times (2 - 7\delta_*^2, 2)$ part of $\mathcal{M}_\Sigma$ intersects the u < 0 and $1 - 3\cos^2\theta > 0$ part of each $\mathfrak{p} \in \Lambda$ version of $\mathbb{R} \times \mathcal{H}_\mathfrak{p}$. In each case, the intersection is a properly embedded annulus that can be parametrized by the functions $\phi$ and $\hat{h}$. This parametrization is such that the range of $\hat{h}$ is symmetric with respect to multiplication by -1; and such that functions s and u restrict as ±1 symmetric functions of $\hat{h}$. Thus, each constant $\hat{h}$ slice of the annulus is a circle in some constant (s, u) sphere in $\mathbb{R} \times \mathcal{H}_\mathfrak{p}$. Taken together, these annuli from $\mathcal{M}_\Sigma$ foliate the part of $\mathbb{R} \times \mathcal{H}_\mathfrak{p}$ where u ≠ 0 and $1 - 3\cos^2\theta > 0$.

The space $\mathcal{M}_{\mathfrak{p}0}$ is $\mathbb{R}$-equivariantly diffeomorphic to $\mathbb{R}$. An element in the space $\mathcal{M}_{\mathfrak{p}0}$ is a properly embedded annulus in the part of $\mathbb{R} \times \mathcal{H}_\mathfrak{p}$ where $1 - 3\cos^2\theta > 0$ and u = 0. The pair $(\phi, \hat{h})$ restrict as coordinates to this annulus such that $\hat{h}$ defines a proper map to $(-(x_0 + 2e^{-2R})\frac{2}{3\sqrt{3}}, (x_0 + 2e^{-2R})\frac{2}{3\sqrt{3}})$. The coordinate s on the annulus is a symmetric function of $\hat{h}$ that is unbounded from above on both ends of its domain. The member parametrized by $0 \in \mathbb{R}$ intersects the s = 0, u = 0 slice of $\mathbb{R} \times \mathcal{H}^+_\mathfrak{p}$ as the $\theta = \frac{\pi}{2}$ circle in $S^2$. The annuli from $\mathcal{M}_{\mathfrak{p}0}$ foliate the $1 - 3\cos^2\theta > 0$ part of the u = 0 slice of $\mathbb{R} \times \mathcal{H}_\mathfrak{p}$.

*Part 2*: Introduce $\mathcal{H}^+_\mathfrak{p}$ to denote the $1 - 3\cos^2\theta > 0$ and $|u| < R + \ln(\delta_*)$ part of $\mathcal{H}_\mathfrak{p}$. The inverse of $\Psi_\mathfrak{p}$ is an $\mathbb{R}$-equivariant embedding of $\mathbb{R} \times \mathcal{H}^+_\mathfrak{p}$ into

$$\mathbb{R} \times (-R - \ln(\delta_*), R + \ln(\delta_*)) \times (\mathbb{R}/2\pi\mathbb{Z}) \times (-\tfrac{4}{3\sqrt{3}}\delta_*^2, \tfrac{4}{3\sqrt{3}}\delta_*^2) \ .$$

(1.28)

The image is denoted by $\mathbb{R} \times \mathcal{X}$. The coordinate functions on the space depicted in (1.28) and thus on $\mathbb{R} \times \mathcal{X}$ are denoted by $(x, \hat{u}, \hat{\phi}, h)$. The rules that follow define $\Psi_\mathfrak{p}$.

- $\Psi_\mathfrak{p}$ *sends a given* $(x, \hat{u} \neq 0, \hat{\phi}, h)$ *point in* $\mathbb{R} \times \mathcal{X}$ *to the* $\phi = \hat{\phi}$, $\hat{h} = h$ *point on the subvariety from* $\mathcal{M}_\Sigma$ *that intersects the* $\theta = \frac{\pi}{2}$ *slice of* $\mathbb{R} \times \mathcal{H}^+_\mathfrak{p}$ *where* $s = x, u = \hat{u}$.
- $\Psi_\mathfrak{p}$ *sends a given* $(x, \hat{u} = 0, \hat{\phi}, h)$ *point in* $\mathbb{R} \times \mathcal{X}$ *to the* $\phi = \hat{\phi}$, $\hat{h} = h$ *point on the subvariety from* $\mathcal{M}_{\mathfrak{p}0}$ *that intersects the* $\theta = \frac{\pi}{2}$ *slice of* $\mathbb{R} \times \mathcal{H}_\mathfrak{p}^+$ *where* $s = x, u = 0$.

(1.29)

Formulas for the $\Psi_\mathfrak{p}$-pushforwards of the coordinate vector fields $\partial_x, \partial_{\hat{u}}, \partial_{\hat{\phi}}$ and $\partial_h$ as given in (II.4.4) can be written as

- $\Psi_{\mathfrak{p}*}\partial_x = \partial_s$ ,
- $\Psi_{\mathfrak{p}*}\partial_{\hat{u}} = \upsilon(\nu + \alpha^{-1}\sqrt{6}x\cos\theta\,\partial_\phi + \varpi\,\partial_s)$ ,

- $\Psi_{p*}\partial_{\hat{\phi}} = \partial_{\phi}$ ,
- $\Psi_{p*}\partial_{h} = -\beta^{-1}(e_1 - \sigma^{-1}\sqrt{6}\chi_\delta f \cos\theta \sin^2\theta \, \partial_s)$ .

(1.30)

Here, $\alpha$ is from (1.6), both $\upsilon$ and $\beta$ are certain positive functions of the pair $(u, \theta)$, the function $\varpi$ depends only on $(u, \theta)$, the function $\sigma$ is from the fifth bullet in Part 1 of Section 1c, and the vector field $e_1$ is defined by (1.14).

*Part 3*: The definition just given endows $\Psi_p$ with the properties listed momentarily. Let p denote the index 1 critical point of $f$ from $\mathfrak{p}$ and let p´ the corresponding index 2 critical point. The list refers to the annuli $T_{p+}$ and $T_{p'-}$ in $\Sigma$ and their respective coordinates $(\varphi_+, h_+)$ and $(\varphi_-, h_-)$. These are introduced in Part 1 of Section 1c. The list also writes the $f \in (1,2)$ part of $\mathbb{R} \times M_\delta$ as a subset of $\mathbb{R} \times (1,2) \times \Sigma$; and it uses $(s, t)$ to denote the Euclidean coordinates on $\mathbb{R} \times (1,2)$. What follows next is the promised list.

- *The constant $(x, \hat{u})$ surfaces in $\mathbb{R} \times X$ are mapped by $\Psi_p$ to J-holomorphic submanifolds.*
- *The map $\Psi_p$ is equivariant with respect to the $\mathbb{R}$ actions on $\mathbb{R} \times X$ and $\mathbb{R} \times \mathcal{H}^+_p$ along their $\mathbb{R}$ factors.*
- *The map $\Psi_p$ is equivariant with respect to the $\mathbb{R}/2\pi\mathbb{Z}$ action that translates the coordinate $\hat{\phi}$ on $\mathbb{R} \times X$ and translates the coordinate $\phi$ on $\mathbb{R} \times \mathcal{H}^+_p$.*
- *The $\hat{u} \geq R + \ln\delta$ part of $X$ is $(R+\ln\delta, R+\ln\delta_*) \times \mathbb{R}/2\pi\mathbb{Z} \times (-\frac{4}{3\sqrt{3}}\delta_*^2, \frac{4}{3\sqrt{3}}\delta_*^2)$; and $\Psi_p$ maps this part of $\mathbb{R} \times X$ diffeomorphically onto $\mathbb{R} \times [1+\delta^2, 1+\delta_*^2) \times T_{p+}$ by the rule*

$$(x, \hat{u}, \hat{\phi}, h) \to (s = x, t = e^{2(\hat{u}-R)}, \varphi_+ = \hat{\phi}, h_+ = h).$$

- *The $\hat{u} \leq -R - \ln\delta$ part of $X$ is $(-R-\ln\delta, -R-\ln\delta_*) \times \mathbb{R}/2\pi\mathbb{Z} \times (-\frac{4}{3\sqrt{3}}\delta_*^2, \frac{4}{3\sqrt{3}}\delta_*^2)$; and $\Psi_p$ maps this part of $\mathbb{R} \times X$ diffeomorphically onto $\mathbb{R} \times (2-\delta_*^2, 2-\delta^2) \times T_{p+}$ by the rule*

$$(x, \hat{u}, \hat{\phi}, h) \to (s = x, t = e^{-2(\hat{u}+R)}, \varphi_+ = \hat{\phi}, h_+ = -h).$$

(1.31)

The first three bullets of this list reproduce (II.4.3).

## 2. The approximations

As noted in Part 1 of Section 1f, each Lipshitz submanifold can be used to construct a corresponding set of approximations to ersatz ech-HF submanifolds. The set in question is parametrized by a subset $\hat{\mathcal{Z}}^S \subset \hat{\mathcal{Z}}_{\text{ech,M}} \times \hat{\mathcal{Z}}_{\text{ech,M}}$ which is invariant with respect to the diagonal action of $\mathbb{Z}$. Each element in $\hat{\mathcal{Z}}^S$ determines a corresponding ersatz ech-HF subvariety, this being a collection of $G + 1$ submanifolds with boundary that is denoted by $\mathcal{C}_0 = \{C_{S0}, \{C_{\mathfrak{p}0}\}_{\mathfrak{p} \in \Lambda}\}$. By way of a reminder, $C_{S0}$ is a submanifold with boundary in the $f \in (1, 2)$ part of $\mathbb{R} \times M_\delta$ and each $\mathfrak{p} \in \Lambda$ version of $C_{\mathfrak{p}0}$ is a submanifold with boundary in $\mathbb{R} \times \mathcal{H}^+_\mathfrak{p}$. The upcoming Section 2b defines $\hat{\mathcal{Z}}^S$; and the remaining subsections describe the ersatz ech-HF submanifold that is associated to any given element in $\hat{\mathcal{Z}}^S$. Section 2a describes the data needed to construct this association.

**a) The parameters $(\delta, x_0, R)$ and a new parameter, $z_*$**

Section 1f does not mention one important point: The desired set of ersatz ech-HF submanifolds can be constructed from a given Lipshitz submanifold only if the parameter $\delta$ from the data set $(\delta, x_0, R)$ that defines Y and its stable Hamiltonian geometry is sufficiently small. In particular, the chosen almost complex structure $J_{HF}$ and the orbit in $\mathcal{A}_{HF}/\mathbb{R}$ of the chosen Lipshitz submanifold jointly determine an upper bound on $\delta$. As noted in Section 1a, the latter determines an upper bound for $x_0$, and then $x_0$ determines an upper bound for R.

The construction of the ersatz set of ech-HF submanifolds requires the specification of an additional parameter, this denoted by $z_*$. This $z_*$ is a positive number whose maximum allowed value is determined by the orbit in $\mathcal{A}_{HF}/\mathbb{R}$ of the chosen Lipshitz submanifold. In any event, $z_*$ is less than $e^{-32}\delta_*^2$. The choice of $z_*$ must be made prior to choosing $\delta$ since the constructions that follow require the maximum allowed value for $\delta$ be less than $e^{-16} z_*^{1/2}$.

Additional refinements for allowed maximum of $z_*$ and $\delta$ are stated as they are needed in the constructions to come. In any event, both are determined solely by the $\mathbb{R}$-orbit of the chosen Lipshitz surface. What follows are further comments on this issue that are of specific concern with regards to the proof of Theorem 1.1.

The first point is perhaps self evident: A given finite set in $\mathcal{A}_{HF}/\mathbb{R}$ determine maximum values for $z_*$ and $\delta$ such that the constructions to follow can be made using any Lipshitz surface from this chosen set of $\mathbb{R}$-orbits if $z_*$ and $\delta$ are less than their allowed maxima.

As it turns out, the set of $\mathbb{R}$-orbits need not be finite to obtain this same conclusion. Saying more requires a digression to introduce the notion of a *weakly compact set*. Let $\mathcal{K} \subset \mathcal{A}_{HF}$ denote an $\mathbb{R}$-invariant set of submanifolds. This set is said to be weakly compact when the following two requirements are met: First, integration of the 2-form $w_\Sigma$ over the Lipshitz surfaces maps $\mathcal{K}$ to a bounded subset in $\mathbb{R}$. Here is the second requirement: Fix any sequence from $\mathcal{K}$ and there exists a surface $S \in \mathcal{K}$ and subsequence of the given sequence which, after renumbering as $\{S_n\}_{n=1,2,\ldots}$, obeys

- $\sup_{z \in S \cap ([-1/n, 1/n] \times [1,2] \times \Sigma)} \text{dist}(z, S_n) + \sup_{z \in S_n \cap ([-1/n, 1/n] \times [1,2] \times \Sigma)} \text{dist}(S, z) < \frac{1}{n}$.
- *Let $\mu$ denote a smooth 2-form on $[-\frac{1}{n}, \frac{1}{n}] \times [1,2] \times \Sigma$ with compact support, with supremum norm 1 and with $|\nabla \mu| \leq \frac{1}{n}$. Then $|\int_S \mu - \int_{S_n} \mu| \leq \frac{1}{n}$.*

(2.1)

Suppose $\mathcal{K}$ is a given, weakly compact subset of Lipshitz submanifolds. As it turns out, maximum values of $z_*$ and $\delta$ can be chosen so that the constructions to come can be done using any surface from $\mathcal{K}$ and values for $z_*$ and $\delta$ that are less than these $\mathcal{K}$-dependent maxima. A data set $D = (z_*, \delta, x_0, R)$ that can be used for all Lipshitz submanifolds in $\mathcal{K}$ is said in what follows to be *$\mathcal{K}$-compatible*. Note that the applications to the proof of Theorem 1.1 require only finite sets of $\mathbb{R}$-orbits of Lipshitz submanifolds.

The proof that $\mathcal{K}$-compatible data sets exist is a straightforward affair given how the maxima for $z_*$ and $\delta$ are subsequently determined from any given Lipshitz submanifold. The proof is left to the reader save for what is said in the two parts that follow.

*Part 1*: Fix a Lipshitz submanifold S. The upper bounds for $z_*$ and $\delta$ are determined by certain data that can be associated to S. This first two element of this data set come from PROPERTY 5 of Section 1e. These are the constants $z_S$ and $\kappa_S$. In particular $z_*$ is constrained to be less than $e^{-32} z_S$. The third element also comes via PROPERTY 5 of Section 1e. This is a bound for the $C^6$ norm over any length 1 interval in $\mathbb{R}$ of any index 1 and index 2 critical point version of the map $\varphi^{S,p}$. The data set also contains the $C^6$ norms of the intrinic and extrinsic curvatures of S, and a maximum for the allowed diameter of a tubular neighborhood of S in $\mathbb{R} \times [1,2] \times \Sigma$.

The final element in the data set is a norm for a certain inverse of the operator $D_S$. To say more about what this means, recall from Section II.6e that $D_S$ defines an $\mathbb{R}$-linear, Fredholm map from a certain Hilbert subspace of $L^2_1$ sections of the complex normal bundle of S to the $L^2$ Hilbert space of sections of the tensor product of this normal bundle with the (0,1) cotangent bundle of S. PROPERTY 5 of Section 1e implies that this map is surjective. As a consequence, the operator $D_S$ has an inverse that maps the range Hilbert space to the $L^2$-orthogonal complement in the domain Hilbert space of the kernel of $D_S$. The latter map is continuous and so bounded; it is the desired inverse.

*Part 2*: The apriori bound on the integral of $w_\Sigma$ over the submanifolds that comprise $\mathcal{K}$ implies that all of the data listed in Part 1 lie in a compact set. The point being that the convergence criteria for membership in $\mathcal{K}$ and the bound on the integral of $w_\Sigma$ implies that the subsequence $\{S_n\}_{n=1,2,...}$ converges to the submanifold S in the $C^\infty$ topology on compact subsets of $\mathbb{R} \times [1,2] \times \Sigma$. Note in this regard that the constraint on $J_{HF}$ given by the sixth bullet in Part 1 of Section 1c has the following consequence: When written as in the proof of Lemma II.6.3, the equations that define a Lipshitz subvariety in $\mathbb{R} \times (1,2) \times \Sigma$ are $\mathbb{C}$-linear equations on neighborhoods of $\mathbb{R} \times \{1\} \times T_+$ and $\mathbb{R} \times \{2\} \times T_-$. This linearity is exploited in Part 3 of Section II.6b. In particular, only slight modifications to the arguments used in the proof of Lemma II.6.3 establish $C^\infty$ convergence for the parts of $\{S_n\}$ near the boundary of $\mathbb{R} \times [1, 2] \times \Sigma$. Meanwhile, the $C^\infty$ convergence in the interior of $\mathbb{R} \times [1, 2] \times \Sigma$ is proved using standard arguments about sequences of pseudoholomorphic curves. See, for example [MS].

The $C^\infty$ convergence on compact subsets of $\mathbb{R} \times (1,2) \times \Sigma$ to a surface in $\mathcal{K}$ implies the desired apriori bound on all but one element of any $S \in \mathcal{K}$ version of the data set given in Part 1. The one element missing is the norm of the inverse of $D_S$. The needed bound on this norm can be derived using (II.6.15)–(II.6.17) to deal with the case when $\{S_n\}_{n=1,2,...}$ does not converge pointwise on the whole of $\mathbb{R} \times [1, 2] \times \Sigma$ to its limit. In the latter case the large n versions of $S_n$ will have long, nearly $\mathbb{R}$-invariant cylinders (a consequence of Lemma II.5.6). Even so, (II.6.15)–(II.6.17) supply an $S \in \mathcal{K}$ independent constant $c_0$ such that $\|D_S\eta\|_{L^2} \geq c_0^{-1}\|\eta\|_{L^2_1}$ if $\eta$ is in the domain of $D_S$ and has compact support on such a cylinder. This sort of bound plus the $C^\infty$ convergence on compact subsets of $\mathbb{R} \times [1, 2] \times \Sigma$ implies the desired $S \in \mathcal{K}$ independent bound on the norm of the inverse of $D_S$.

**b) The set $\hat{\mathcal{Z}}^S$**

Fix a Lipshitz submanifold, S. As noted at the outset, the set of ersatz ech-HF submanifolds that are constructed from S is indexed by a subset $\hat{\mathcal{Z}}^S \subset \hat{\mathcal{Z}}_{ech,M} \times \hat{\mathcal{Z}}_{ech,M}$. The two parts of this subsection describe the set $\hat{\mathcal{Z}}^S$.

*Part 1*: The diagonal action of $\mathbb{Z}$ on $\hat{\mathcal{Z}}_{ech,M} \times \hat{\mathcal{Z}}_{ech,M}$ preserves $\hat{\mathcal{Z}}^S$ and gives it the structure of a principle $\mathbb{Z}$ bundle over its image in $\mathcal{Z}_{HF} \times \mathcal{Z}_{HF}$. The image is described in the upcoming Part 2 of this subsection. To say more about the fiber over this image, first write $\hat{\mathcal{Z}}_{ech,M}$ using (1.13) as $\mathcal{Z}_{ech,M} \times \mathbb{Z}$. Introduce next $\mathfrak{n}_S$ to denote the intersection number between S and the $J_{HF}$-holomorphic subvariety $\mathbb{R} \times (1,2) \times z_0$. This is a non-negative integer. A given element $((\Theta_-, k_-), (\Theta_+, k_+))$ from $\hat{\mathcal{Z}}_{\text{ech},M} \times \hat{\mathcal{Z}}_{\text{ech},M}$ sits in $\hat{\mathcal{Z}}^S$ only if $k_+ = k_- + n_S$.

*Part 2*: Let $\hat{\upsilon}_-$ and $\hat{\upsilon}_+ \in \mathcal{Z}_{HF}$ denote the elements that are defined by S as described in PROPERTY 6 of Section 1g. These sets define via (1.10) a corresponding set of elements in $\mathcal{Z}_{\text{ech},M}$, thus a subset in $\mathcal{Z}_{HF} \times \mathcal{Z}_{HF}$ of the form

$$(\hat{\upsilon}_- \times (\times_{\mathfrak{p} \in \Lambda} (\mathbb{Z} \times 0))) \times (\hat{\upsilon}_+ \times (\times_{\mathfrak{p} \in \Lambda} (\mathbb{Z} \times 0))) .$$
(2.1)

The set $\hat{\mathcal{Z}}^S$ will sit over a subset in (2.1). The latter is denoted in what follows by $\mathcal{Z}^S$.

The elements in (2.1) that lie in $\mathcal{Z}^S$ are characterized by G conditions, one for each $\mathfrak{p} \in \Lambda$. As explained momentarily, a given pair $(\hat{\upsilon}_-, (\mathfrak{k}_{\mathfrak{p}-}, o_{\mathfrak{p}-})_{\mathfrak{p} \in \Lambda})$ and $(\hat{\upsilon}_+, (\mathfrak{k}_{\mathfrak{p}+}, o_{\mathfrak{p}+})_{\mathfrak{p} \in \Lambda})$ from (2.1) defines a corresponding set of integers, this denoted by $\{\mathfrak{m}_{\mathfrak{p}}\}_{\mathfrak{p} \in \Lambda}$. The given pair defines an element in $\mathcal{Z}^S$ if and only if each $\mathfrak{p} \in \Lambda$ version of $\mathfrak{m}_{\mathfrak{p}}$, $o_{\mathfrak{p}-}$ and $o_{\mathfrak{p}+}$ obey

- $\mathfrak{m}_{\mathfrak{p}} = 0$ *and one of the following:*
    a) $o_{\mathfrak{p}-} = o_{\mathfrak{p}+} = \{0\}$
    b) $o_{\mathfrak{p}-} = \{0\}$ *and* $o_{\mathfrak{p}+} = \{-1, 1\}$.
- $\mathfrak{m}_{\mathfrak{p}} = -1$ *and* $o_{\mathfrak{p}-} = \{0\}$ *and* $o_{\mathfrak{p}+} = \{1\}$,
- $\mathfrak{m}_{\mathfrak{p}} = 1$ *and* $o_{\mathfrak{p}-} = \{0\}$ *and* $o_{\mathfrak{p}+} = \{-1\}$.

(2.2)

Fix $\mathfrak{p} = (p, p') \in \Lambda$. To say more about $\mathfrak{m}_{\mathfrak{p}}$, let $\gamma_{\mathfrak{p}+}$ and $\gamma_{\mathfrak{p}-}$ denote the respective segments of the integral curves of $v$ in $\mathcal{H}^+_{\mathfrak{p}}$ that are defined by the data $(\hat{\upsilon}_-, (\mathfrak{k}_{\mathfrak{p}-}, o_{\mathfrak{p}-})_{\mathfrak{p} \in \Lambda})$ and $(\hat{\upsilon}_+, (\mathfrak{k}_{\mathfrak{p}+}, o_{\mathfrak{p}+})_{\mathfrak{p} \in \Lambda})$. Fix $r \in [\frac{1}{2}\delta_*, \delta_*]$ so that projection map from $\mathbb{R} \times (1,2) \times \Sigma$ to $(1,2) \times \Sigma$ restricts to S so as to define a map that is transverse to surfaces in $(1,2) \times \Sigma$ that correspond to the $f \in (1,2)$ part of the boundary of the radius r coordinate balls centered at p and p'. These respective intersections define a pair of embedded arcs, one in the boundary of the radius r coordinate ball centered at p and the other in the boundary of the radius r coordinate ball centered at p'. The former starts at the point where the relevant integral curve of $\mathfrak{v}$ from $\hat{\upsilon}_+$ intersects the radius r coordinate ball centered at p and ends at the point where the relevant integral curve of $\mathfrak{v}$ from $\hat{\upsilon}_-$ intersects this radius r coordinate ball. Denote this arc by $\upsilon_{\mathfrak{p}1}$. The second arc starts from the point where the relevant integral curve of $\mathfrak{v}$ from $\hat{\upsilon}_-$ intersects the radius r coordinate ball centered on p' and it ends where the relevant integral curve of $\mathfrak{v}$ from $\hat{\upsilon}_+$ intersects the radius r coordinate ball centered on $p_2$. Denote this second arc by $\upsilon_{\mathfrak{p}2}$.

As noted by Corollary II.2.6, the starting point of $\upsilon_{\mathfrak{p}1}$ has distance no greater than $c_0 \delta$ from the point where $\gamma_{\mathfrak{p}+}$ intersects the boundary of the radius r coordinate ball centered at p and its ending point has distance no greater than $c_0 \delta$ from the point where $\gamma_{\mathfrak{p}-}$ intersect the boundary of the radius r coordinate ball centered at p. There is an analogous observation about the starting and ending points of $\upsilon_{p2}$.

Granted these observations, what follows defines a closed 1-cycle in $\mathcal{H}^+_p$. Start where $\gamma_{p+}$ intersects the boundary of the radius r coordinate ball centered on p´ and proceed along $\gamma_{p+}$ until it intersects the boundary of the radius r coordinate ball centered on p. Then proceed along a geodesic arc in this sphere of length $c_0\delta$ or less to the starting point of $\upsilon_{p1}$. Proceed along $\upsilon_{p1}$ to its endpoint and then along the geodesic arc in the sphere to the its intersection point with $\gamma_{p-}$. Return to the boundary of the radius r coordinate ball centered on p´ by traversing backwards along $\gamma_{p-}$. Then proceed along the short geodesic in this sphere to the starting point of $\upsilon_{p2}$, follow $\upsilon_{p2}$ to its end, and then follow the short geodesic in this sphere to the nearby intersection point $\gamma_{p+}$.

This closed 1-cycle defines a class in $H_1(\mathcal{H}^+_p; \mathbb{Z})$. The latter group is isomorphic to $\mathbb{Z}$ with generator the equatorial circle in the u = 0 slice with the orientation given by $\frac{\partial}{\partial \phi}$. This understood, the closed 1-cycle defines an integer. This integer is $\mathfrak{m}_p$.

### c) The submanifold $C_{S0}$

Fix a Lipshitz submanifold S; or if needed, take S from some chosen, weakly compact subset $\mathcal{K} \subset \mathcal{A}_{HF}$ of Lipshitz submanifolds. In any event, with S chosen, fix a pair $(\hat{\Theta}_-, \hat{\Theta}_+) \in \hat{\mathcal{Z}}^S$. This data labels an ersatz ech-HF submanifold, $\mathcal{C}_0 = \{C_{S0}, \{C_{p0}\}_{p\in\Lambda}\}$. This subsection describes $C_{S0}$. The description is given in the first four parts of this subsection. An existence/uniqueness assertion is stated in Part 6 by Proposition 2.1. This proposition gives an indication of the role played by S. Part 5 of the subsection sets some background for Proposition 2.1.

*Part 1*: Use Lemmas II.6.3 and II.6.4 to find a constant, $z_S < e^{-32}\delta_*^2$ so that the conclusions of Lemma II.6.3 holds and so that Lemma II.6.4 holds when $z \leq z_S$. With regards to Lemma II.6.4, choose $z_S$ so as to guarantee the following: The composition of first projection from $\mathbb{R} \times (1,2) \times \Sigma$ to $(1,2) \times \Sigma$ and then the identification of the latter with the $f \in (1,2)$ part of M sends the $f \leq 1 + z_S$ portion of S into the union of the radius $e^{-16}\delta_*$ coordinate balls centered on the index 1 critical points of $f$; and it sends the portion where $f \geq 2 - z_S$ into the union of the radius $e^{-16}\delta_*$ coordinate balls centered on the index 2 critical points of $f$. Note that $z_S$ can be taken to be $\mathcal{K}$-compatible when $\mathcal{K}$ is specified.

Fix $z_* \in (0, e^{-64}z_S)$ and then $\delta < e^{-16}z_*^{1/2}$. Some additional purely S-dependent ($\mathcal{K}$-compatible) constraints on the upper bounds for $z_*$ and $\delta$ are given subsequently.

*Part 2*: The element $C_{S0}$ from $\mathcal{C}_0$ is a properly embedded, J-holomorphic submanifold with boundary in $\mathbb{R} \times [1+z_*, 2-z_*] \times \Sigma$. This surface has 2G boundary components, with G on $\mathbb{R} \times \{1\} \times \Sigma$ and G on $\mathbb{R} \times \{2\} \times \Sigma$. These former set are mapped via the projection to $\Sigma$ into pairwise distinct components of $T_+$, and the latter are mapped via this projection into pairwise distinct components of $T_-$. In any event, each boundary component is an embedded copy of $\mathbb{R}$.

*Part 3*: This part describes the large $|s|$ part of $C_{S0}$. To this end, let $\Theta_-$ and $\Theta_+$ denote the respective base points in $\mathcal{Z}_{ech,M}$ for the chosen elements $\hat{\Theta}_-, \hat{\Theta}_+ \in \hat{\mathcal{Z}}^S$. What follows first describes the $s \ll -1$ behavior.

There exists $s_1 \geq 1$ which is such that the $s \leq -s_1$ part of $C_{S0}$ is a disjoint union of G graphs over $(-\infty, -s_1] \times [1+z_*, 2-z_*]$. Each such graph has the form

$$(s, t) \to (s, t, \psi(s, t))$$
(2.3)

where $\psi$ is a map from $(-\infty, -s_1] \times [1+z_*, 2-z_*]$ to $T_+ \cap T_-$. The images of these G maps are in distinct components; and each such component contains the intersection with $\Sigma$ of an integral curve of $v$ from $\Theta_-$. Let $q_*$ now denote such an intersection point and let $\psi_*$ denote the map from (2.3) with image in the $q_*$ component of $T_+ \cap T_-$. Then

$$\text{dist}(\psi_*(s, \cdot), q_*) \leq c\, e^{-|s|/c}$$
(2.4)

where $c \geq 1$ is a purely S-dependent (or $\mathcal{K}$-compatible) constant. Finally, if the given component of $T_+ \cap T_-$ is parametrized using the relevant index 1 critical point versions of the functions $(\varphi_+, \hat{h}_+)$, and if $\psi$ is written with respect to these coordinates in terms of functions $(\varphi_+ = \varphi_*(s,t), \hat{h}_+ = \varsigma_*(s,t))$, then the pair $(\varphi_*, \varsigma_*)$ obey the Cauchy-Riemann equations $\partial_s \varphi_* - \partial_t \varsigma_* = 0$ and $\partial_s \varsigma_* + \partial_t \varphi_* = 0$.

The $s \geq s_1$ part of $C_{S0}$ has the analogous description with $\Theta_+$ replacing $\Theta_-$.

*Part 4*: This part describes the boundary behavior of $C_{S0}$ near any given boundary component. This involves a purely S-dependent (or $\mathcal{K}$-compatible) constant, $\kappa_{S*}$, which is greater than 100. If $z_S$ is chosen less than $\kappa_{S*}^{-2}$ then what follows holds is true. Let p denote either an index 1 or index 2 critical point of $f$. When p has index 1, use $(\varphi_+, \hat{h}_+)$ to parametrize $T_{p+}$, and when p has index 2, use $(\varphi_-, \hat{h}_-)$ to parametrize $T_{p-}$. Then the part of $C_{S0}$ in $\mathbb{R} \times [1+z_*, 1+z_S] \times T_{p+}$ or in $\mathbb{R} \times [2-z_S, 2-z_*] \times T_{p-}$ is diffeomorphic to $\mathbb{R} \times [z_*, z_S]$ and parametrized via a map of the form

- $(x, z) \to (s = x, t = 1+z, \varphi_+ = \varphi^{S0}(x, z), \hat{h}_+ = \varsigma^{S0}(x,z))$ *when p has index 1*.

- $(x, z) \to (s = x, t = 2 - z, \varphi_- = \varphi^{S0}(s, z), \hat{h}_- = \varsigma^{S0}(s, z))$ when p *has index* 2.

(2.5)

The functions $\varphi^{S0}$ and $\varsigma^{S0}$ that appear here are $\mathbb{R}$-valued functions that obey the Cauchy-Riemann equations $\partial_x \varphi^{S0} - \partial_z \varsigma^{S0} = 0$ and $\partial_x \varsigma^{S0} + \partial_z \varphi^{S0} = 0$. In addition, the first and higher derivatives of these functions to any given order are bounded uniformly on $\mathbb{R} \times [z_*, z_S]$ by a constant that depends only on the given order and S (it is $\mathcal{K}$-compatible when $\mathcal{K}$ is given). Finally $|\varsigma^{S0}(\cdot, z_*)| \leq \kappa_S z_*$.

*Part 5*: The upcoming Proposition 2.1 states an existence/uniqueness assertion for $C_{S0}$. This part of the subsection supplies some background for this proposition.

Proposition 2.1 views the $f \in (1, 2)$ part of M as $(1, 2) \times \Sigma$ so as to describe $C_{S0}$ as a submanifold with boundary in $\mathbb{R} \times [1 + z_*, 2 - z_*] \times \Sigma$. The submanifold $C_{S0}$ in this guise is the $t \in [1 + z_*, 2 - z_*]$ portion of a properly embedded, $J_{HF}$-holomorphic submanifold in $\mathbb{R} \times [1, 2] \times \Sigma$, this denoted by $S_*$. The submanifold $S_*$ has 2G boundary components, one in each index 1 critical point version of $\mathbb{R} \times \{1\} \times T_{p+}$; and likewise, one in each index 2 critical point version of $\mathbb{R} \times \{2\} \times T_{p-}$.

The submanifold $S_*$ is isotopic to S in a small radius tubular neighborhood of S. The description in Proposition 2.1 identifies this tubular neighborhood with a disk bundle in the normal bundle of S using an exponential map of the sort that is described in Section II.6e. What follows reviews some aspects of this sort of exponential map.

To start, recall from Section II.6e that S has a complex normal bundle, $N_S \to S$ and an exponential map the embeds a disk subbundle as a tubular neighborhood of S. The exponential map is denoted by $\mathfrak{e}_S$ and the disk subbundle by $N_0$. The latter has radius $\rho_S$ and its image in $\mathbb{R} \times [1, 2] \times \Sigma$ is a tubular neighborhood that contains the set of points with distance $c^{-1} \rho_S$ from S. The map $\mathfrak{e}_S$ embeds each fiber disk as a $J_{HF}$-holomorphic disk. If $\mathcal{K}$ is a previously specified, weakly compact set of Lipshitz submanifolds, then $\rho_S$ and $c$ can be taken to be $\mathcal{K}$-compatible; and as can the derivatives to any given order of the exponential map $\mathfrak{e}_S$.

Although not stated as such in Section II.6e, the map $\mathfrak{e}_S$ can be chosen so as to respect the graph structure described in Section II.6c and PROPERTY 5 of Section 1g near the boundaries of S. In particular, $\mathfrak{e}_S$ can be chosen so that it maps any given fiber of $N_0$ over the $t \in (1, 1 + z_S)$ and $t \in (2 - z_S, z_S)$ portions of S as follows: The graph structure indicated by PROPERTY 5 of Section 1g identifies the bundle $N_S$ over this part of S with the restriction to S of T$\Sigma$. In particular, the 1-forms $(d\varphi_+, d\hat{h}_+)$ when $t \in (1, 1 + z_S)$ and with $(d\varphi_-, -d\hat{h}_-)$ when $t \in (2 - z_S, 2)$ with the underlying real bundle defines an orientation preserving isomorphism to the product $\mathbb{R}^2$ bundle. Given this isomorphism, $\mathfrak{e}_S$ on these parts of S can and should be chosen so as to send any point in S parameterized by (x, z) and a pair (a,b) in the $\mathbb{R}^2$ factor of the product bundle to the point given by one of

- $(s = x, t = 1+z, \varphi_+ = \varphi(x,z) + a, \hat{h}_+ = \varsigma(x, z) + b)$.
- $(s = x, t = 2 - z, \varphi_- = \varphi(x,z) + a, \hat{h}_- = \varsigma(x, z) - b)$.

(2.6)

The pair $(d\varphi_+, d\hat{h}_+)$ as defined by any given index 1 critical point of $f$ are also used to write $N_S$ as the product bundle over the corresponding $s \ll -1$ and $s \gg 1$ parts of S. By way of a reminder, the $s \ll -1$ part of S has G components, each projecting to a distinct component of $[1, 2] \times (T_+ \cap T_-)$, and these G components can be labeled by the index 1 critical points of $f$. The same can be said for the $s \gg 1$ part of S. In any event, the exponential map $\mathfrak{e}_S$ can and should be chosen so that it acts as in the top bullet of (2.6) on these large $|s|$ parts of S.

The construction of a map $\mathfrak{e}_S$ of this sort can be made using the techniques that are used to prove Lemma 5.4 in [T1].

Proposition 2.1 refers to the Fredholm operator, $D_S$, that is described in Section II.6e; it is depicted in (1.25). By way of a reminder, this operator maps a certain Hilbert space of sections of $N_S$ to the space of square integrable sections of $N_S \otimes T^{0,1}S$. The Hilbert space for the domain is the Sobolev $L^2_1$ norm completion of the subspace of sections that obeys the constraints in (II.6.12). The *kernel* of $D_S$ refers to the sections of $N_S$ in the domain Hilbert space that are annihilated by $D_S$. The $L^2$ inner product on sections of $N_S$ is defined using the fiber metric on $N_S$ and the integration measure on S that comes from the metric induced by its embedding in $\mathbb{R} \times [1, 2] \times \Sigma$.

Write the pair $\hat{\Theta}_-$ and $\hat{\Theta}_+$ from the chosen element in $\hat{\mathcal{Z}}^S$ as $(\Theta_-, k_-)$ and $(\Theta_+, k_+)$. Proposition 2.1 refers to a number that is associated to each index 1 and each index 2 critical point of $f$ by $\Theta_+$ and another that is determined by $\Theta_-$. When p is used to denote the critical point in question, the corresponding two numbers are denoted respectively by $h_{p+}$ and $h_{p-}$. When p is an index 1 critical point of $f$, the numbers $h_{p+}$ and $h_{p-}$ denote the respective $\hat{h}_+$ coordinates of the $T_{p+}$ intersection point of an integral curve of $v$ from $\Theta_+$ and $\Theta_-$ with the $t = 1+z_*$ slice of $(1, 2) \times \Sigma$. When p is an index 2 critical point of $f$, the numbers $h_{p+}$ and $h_{p-}$ denote the the respective $\hat{h}_-$ coordinates of the $T_{p-}$ intersection point of an integral curve of $v$ from $\Theta_+$ and $\Theta_-$ with the $t = 2-z_*$ slice of $(1,2) \times \Sigma$. By way of a parenthetical remark, it follows from what is said in Section II.2 that $|h_{p+}|$ and $|h_{p-}|$ are both bounded by $c_0\delta^2$.

With regards to $\Theta_-$ and $\Theta_+$, Proposition 2.1 uses $\hat{\upsilon}_-$ and $\hat{\upsilon}_+$ to denote the respective HF-cycles that are used in (2.1) for their definition. Let $q \in \Sigma$ denote a given intersection point with an integral curve of $\upsilon$ from either $\hat{\upsilon}_-$ or $\hat{\upsilon}_+$. Let $q_*$ denote the corresponding, nearby intersection point of the corresponding segment of an integral curve of $v$ from $\Theta_-$ or $\Theta_+$ as the case may be. Note that $q_*$ has distance at most $c_0\delta$ from q. The point q is an element in $C_+ \cap C_-$ and therefore in some index 1 critical point component of $C_+$. This critical point labels a corresponding $s \ll -1$ or $s \gg 1$ component of S. The latter component is denoted by $\mathcal{E}_{Sq}$. Proposition 2.1 writes the exponential map $\mathfrak{e}_S$ over $\mathcal{E}_{Sq}$ as in the top line of (2.6) so as to view $q_*$ as a section of $N_S$ over $\mathcal{E}_{Sq}$.

Here is one final item of notation: The function $x \to \chi(x)$ maps $\mathbb{R}$ to [0, 1]; it is nonincreasing, equal to 1 where $x \leq 0$ and equal to zero where $x \geq 1$.

*Part 6*: This final part of the subsection first states and then proves the existence/uniqueness proposition about $C_{S0}$.

**Proposition 2.1**: *Fix a Lipshitz submanifold S or one from some specified weakly compact subset $\mathcal{K} \subset \mathcal{A}_{HF}$ of Lipshitz submanifolds. There exist $\kappa_{S*} \geq 100$ and $z_S \in (0, \kappa_{S*}^{-2})$ that depend only on S (and are $\mathcal{K}$-compatible if relevant) such that what follows is true. Fix $z_* \in (0, e^{-32} z_{S*})$ and then $\delta \in (0, e^{-16} z_*^{1/2})$ and $x_0$ and R. There exists a unique section $\eta_*$ of $N_0$ that is characterized by:*
- *The $C^4$-norm of $\eta_*$ is bounded by $\kappa_{S*}\delta$.*
- *The restriction of $\eta_*$ to the $t \in [1+z_*, 2-z_*]$ part of S is $L^2$-orthogonal to the corresponding restriction of the elements in the kernel of the operator $D_S$.*
- *Let p denote either an index 1 or index 2 critical point of f. The pairing of the section $\eta_*$ along the corresponding boundary component of S with the relevant 1-form $d\hat{h}_+$ or $d\hat{h}_-$ is the function on $\mathbb{R}$ given by $\hat{h}_{p+}(1-\chi) + \hat{h}_{p-}\chi$.*
- *Let $q \in \Sigma$ denote an intersection point with an integral curve of $v$ from either $\hat{v}_-$ or $\hat{v}_+$ and let $q_*$ denote the corresponding section of $N_S$ over $\mathcal{E}_{Sq}$. The pointwise norm of $\eta_* - q_*$ converges to zero as $|s| \to \infty$ on $\mathcal{E}_{Sq}$. Moreover, given $k \geq 0$, there exists a purely S-dependent (or $\mathcal{K}$-compatible) constant c such that the derivatives to order k on $\mathcal{E}_{Sq}$ are bounded by $c\,e^{-|s|/c}$.*
- *With $(1,2) \times \Sigma$ viewed now as the $f \in (1,2)$ part of M, use $C_{S0} \subset \mathbb{R} \times M_\delta$ to denote the $t \in [1+z_*, 2-z_*]$ part of $S_* = \mathfrak{e}_S \circ \eta_*(S)$. This version of $C_{S0}$ obeys the properties listed in the preceding Parts 1-4 of this subsection.*

*Proof of Proposition 2.1*: The proof that follows has three steps.

<u>Step 1</u>: Suppose that $\eta_1$ is a smooth section of $N_0$ that has the same large $|s|$ asymptotics and boundary behavior as the desired $\eta_*$. Assume that the pointwise norm of $\eta_1$ and those of its derivatives to sixth order are bounded by $c_0\delta$. In addition, require that the surface $\mathfrak{e}_S \circ \eta_1$ be $J_{HF}$ holomorphic where $|s| \geq c$ where $c > 1$ is a purely S dependent (or $\mathcal{K}$-compatible) constant. A section $\eta_1$ satisfying these requirements is given in Step 2.

What follows constructs a section, $\eta_2$, of $N_0$ in the Fredholm domain of the operator $D_S$ such that $\eta_* = \eta_1 + \eta_2$ obeys all of the bullets of the proposition. To this end, keep in mind what is said by Part 2 of Section II.6e and by (II.6.10): If $\eta = \eta_1 + \eta_2$ is a section of $N_S$, then $\mathfrak{e}_S \circ \eta$ is $J_{HF}$ holomorphic if and only if $\eta_2$ obeys an equation that has the schematic form

$$\bar{\partial}\eta_2 + \mathfrak{r}_{1*}(\eta_2) \cdot \partial\eta_2 + \mathfrak{r}_{0*}(\eta_2) = \mathfrak{r}_*$$

(2.7)

where the notation is as follows: First, $\mathfrak{r}_*$ is a smooth section of $N_S \otimes T^{0,1}S$ with compact support where $|s| \leq c$ and $C^5$ norm and $L^2$ norm bounded by $c\delta$. Here, $c > 1$ again denotes a purely S dependent (or $\mathcal{K}$-compatible) constant. This term $\mathfrak{r}_*$ is determined by $\eta_1$. Second, $\mathfrak{r}_{1*}$ and $\mathfrak{r}_{0*}$ are analogous to their counterparts (II.6.10). They differ from the latter by virtue of a dependence on $\eta_1$, but even so, this difference has $C^5$ norm bounded by $c\delta$ with $c$ as just described. In particular, they obey $|\mathfrak{r}_{1*}(b)| \leq c|b|$ and $|\mathfrak{r}_{2*}(b) - \upsilon b - \mu \bar{b}| \leq c|b|^2$ where $\upsilon$ and $\mu$ are from (1.25) and $c$ is as described above. In addition, their derivatives to any given order are bounded by purely S (or $\mathcal{K}$-compatible) constants.

Granted these last remarks, the equation for $\eta_2$ can be written as

$$D_S\eta_2 + \mathfrak{z}(\eta_2) = \mathfrak{r}_*$$

(2.8)

where $|\mathfrak{z}(\eta_2)| \leq c(|\eta_2|^2 + |\eta_2||\nabla\eta_2|)$. Here again, $c \geq 1$ is a purely $s$ dependent (or $\mathcal{K}$-compatible) constant. With this last fact understood, and given the afore-mentioned bounds on the higher derivatives of $\mathfrak{r}_{0*}$, $\mathfrak{r}_{1*}$ and $\mathfrak{r}_*$, the existence and uniqueness of the desired solution to (2.8) follows via a standard application of the implicit function theorem.

Step 2: Consider now $\eta_1$. What follows here describes $\eta_1$ on a component of the $s \ll -1$ part of S. To do this, return to the notation used in Part 3. Let $q \in C_+ \cap C_-$ denote the relevant point. As noted in PROPERTY 6 of Section 1g, the end $\mathcal{E}_{Sq}$ can be viewed as a graph of a map from the $s \leq -c$ part of $\mathbb{R} \times [1,2]$ into q's component of $T_+ \cap T_-$. Use the coordinates $(\varphi_+, h_+)$ to write q as the origin in $\mathbb{R}^2$ and the corresponding map $\psi$ to $T_+ \cap T_-$ as a map to $\mathbb{R}^2$. With the normal bundle $N_S$ identified with $\mathbb{R}^2$ as in the top line of (2.6), the section $\eta_1$ where $s \leq -c(s_1 + |\ln\delta|)$ is $q_* - \psi$. Note that this formula is such that $\mathfrak{e}_S \circ \eta_1$ is the $J_{HF}$-holomorphic surface $\mathbb{R} \times (1,2) \times q_*$ on this part of S. There is an analogous formula for $\eta_1$ where $s \geq c(s_1 + |\ln\delta|)$. The desired behaviour of $\eta_1$ near the boundary of S can be obtained in a straightforward fashion using the description of S given in PROPERTY 5 of Section 1g. A similarly straightforward use of 'cutoff' functions will extend the section of $N_S$ that is defined by the resulting formula for $\eta_1$ near the boundary of S and that given above on the $|s| \gg c$ on S so as to define the desired version of $\eta_1$ over the whole of S.

Step 3: The demands of the fifth bullet follow directly given the third and fourth bullets and given that $\eta_2$ obeys (2.8). As $\eta_2$ is in the Fredholm domain of $D_S$, the third bullet follows if $\eta_2$ is smooth up to the boundary. This can be proved using slightly modified versions of the arguments that are used to prove Lemma II.6.3. Meanwhile, the assertion made by the fourth bullet is proves using arguments that are little different from those used in Section II.6c.

### d) The submanifolds $\{C_{\mathfrak{p}0}\}_{\mathfrak{p} \in \Lambda}$

The subsequent three parts of this subsection describe the salient features of the submanifolds that comprise the subset $\{C_{\mathfrak{p}0}\}_{\mathfrak{p} \in \Lambda}$ from $\mathcal{C}_0$. To this end, fix $\mathfrak{p} \in \Lambda$ so as to focus on the corresponding element $C_{\mathfrak{p}0}$. The subsequent description uses $\mathcal{H}^+_{\mathfrak{p}*}$ to denote the $e^{-2(R-|u|)}(1 - 3\cos^2\theta) \leq z_*$ part of $\mathcal{H}^+_\mathfrak{p}$. Part 4 of the subsection states and existence/uniqueness assertion about $C_{\mathfrak{p}0}$.

*Part 1*: What is denoted by $C_{\mathfrak{p}0}$ is a properly embedded submanifold with boundary in $\mathbb{R} \times \mathcal{H}^+_{\mathfrak{p}*}$ with J-holomorphic interior. There are two boundary components, one on the $u > 0$ component of the boundary of $\mathbb{R} \times \mathcal{H}^+_{\mathfrak{p}*}$ and the other on the $u < 0$ component. Define $\Delta_\mathfrak{p} \in \{0, 1, 2\}$ as follows: If $\mathfrak{m}_\mathfrak{p} = 0$ and Item a) of the first bullet of (2.2) is relevant, than $\Delta_\mathfrak{p} = 0$. If Item b) is relevant, than $\Delta_\mathfrak{p} = 2$. If $\mathfrak{m}_\mathfrak{p} = 1$ or $\mathfrak{m}_\mathfrak{p} = -1$, then $\Delta_\mathfrak{p} = 1$. In the case $\Delta_\mathfrak{p} = 0$, the submanifold $C_{\mathfrak{p}0}$ is diffeomorphic to the product of $\mathbb{R}$ with a closed interval. When $\Delta_\mathfrak{p} = 1$, the submanifold $C_{\mathfrak{p}0}$ is diffeomorphic to the complement of a single interior point in the product of $\mathbb{R}$ with a closed interval. When $\Delta_\mathfrak{p} = 2$, it is diffeomorphic to the complement of two interior points in such a product.

To describe the large $|s|$ behavior of $C_{\mathfrak{p}0}$, introduce $\gamma_{\mathfrak{p}-}$ and $\gamma_{\mathfrak{p}+}$ to denote the respective integral curve segments in $\mathcal{H}^+_\mathfrak{p}$ that come from $\Theta_-$ and $\Theta_+$, and introduce $\hat{\gamma}^+_\mathfrak{p}$ and $\hat{\gamma}^-_\mathfrak{p}$ to denote the respective closed integral curves of $v$ in the $u = 0$ slice of $\mathcal{H}_\mathfrak{p}$ that comprise the loci where $\cos\theta = \frac{1}{\sqrt{3}}$ and $\cos\theta = -\frac{1}{\sqrt{3}}$.

- *Each constant $s \ll -1$ slice of $C_{\mathfrak{p}0}$ is a properly embedded arc in $\mathcal{H}^+_\mathfrak{p}$; and these arcs converge pointwise as $s \to -\infty$ to $\gamma_{\mathfrak{p}-}$.*

- *If $\Delta_p = 0$, then each constant $s \gg 1$ slice of $C_{p0}$ is a properly embedded arc in $\mathcal{H}^+_p$; and these arcs converge in an isotopic fashion in $\mathcal{H}_p$ as $s \to \infty$ to $\gamma_{p+}$.*
- *If $\Delta_p = 1$, then each constant $s \gg 1$ slice of $C_{p0}$ has two components.*
  a) *One component is a properly embedded arc in $\mathcal{H}^+_p$; and these arcs converge in an isotopic fashion in $\mathcal{H}_p$ as $s \to \infty$ to $\gamma_{p+}$.*
  b) *The other component is an embedded circle; and these circles converge pointwise in $\mathcal{H}_p$ as $s \to \infty$ to $\hat{\gamma}^+_p$ when $\mathfrak{m}_p = -1$, and they converge in an isotopic fashion in $\mathcal{H}_{p*}$ as $s \to \infty$ to $\hat{\gamma}^-_p$ when $\mathfrak{m}_p = 1$.*
- *If $\Delta_p = 2$, then each constant $s \gg 1$ slice of $C_{p0}$ has three components. One component is a properly embedded arc in $\mathcal{H}^+_p$; and these arcs converge in an isotopic fashion in $\mathcal{H}_p$ as $s \to \infty$ to $\gamma_{p+}$. The other two components are embedded circles. One s-parametrized set of these circles converges pointwise in an isotopic fashion in $\mathcal{H}_p$ as $s \to \infty$ to $\hat{\gamma}^+_p$; the other set converges in an isotopic fashion in $\mathcal{H}_{p*}$ as $s \to \infty$ to $\hat{\gamma}^-_p$.*

(2.9)

*Part 2*: This part says more about about how $C_{p0}$ sits in $\mathbb{R} \times \mathcal{H}^+_{p*}$. To this end, reintroduce from Section 1h the parametrization of $\mathcal{H}^+_p$ by the map $\Psi_p$. The domain of this map $\Psi_p$ is an open subset in $\mathbb{R} \times (-R - \ln\delta_*, R + \ln\delta_*) \times (\mathbb{R}/2\pi\mathbb{Z}) \times (-\frac{4}{3\sqrt{3}}\delta_*^2, \frac{4}{3\sqrt{3}}\delta_*^2)$. As in Section 1h, the coordinates for the latter are written as $(x, \hat{u}, \hat{\phi}, h)$. By way of a reminder, the domain of $\Psi_p$ is an open subset of the form $\mathbb{R} \times X$. The $\Psi_p$-inverse image of $\mathbb{R} \times \mathcal{H}^+_{p*}$ is the subset of the domain $\mathbb{R} \times X$ where $|\hat{u}| \leq R + \frac{1}{2}\ln z_*$. This being the case, it proves useful to introduce $I_*$ to denote the interval $[-R - \frac{1}{2}\ln z_*, R + \frac{1}{2}\ln z_*]$ and restrict $\Psi_p$ to the $\hat{u} \in I_*$ part of its domain.

The $\Psi_{p*}$-inverse image of $C_{p*}$ is given as the image of a proper map from a certain domain in $\mathbb{R} \times I_*$ to $\mathbb{R} \times X$. This map has the form

$$(x, \hat{u}) \to (x, \hat{u}, \hat{\phi} = \varphi^{p0}(x, \hat{u}), h = \varsigma^{p0}(x, \hat{u}))$$

(2.10)

where $(\varphi^{p0}, \varsigma^{p0})$ is a map from a domain in $\mathbb{R} \times I_*$ to $\mathbb{R}/(2\pi\mathbb{Z}) \times (-\frac{4}{3\sqrt{3}}\delta_*^2, \frac{4}{3\sqrt{3}}\delta_*^2)$. The domain of this map is $\mathbb{R} \times I_*$ if $\Delta_p = 0$, it is the complement of a single point in the $\hat{u} = 0$ slice of $\mathbb{R} \times I_*$ if $\Delta_p = 1$, and it is the complement of two points in the $\hat{u} = 0$ slice if $\Delta_p = 2$.

*Part 3*: This part describes the behavior of $C_{p0}$ near its boundaries. Consider first the boundary on the $u > 0$ component of the $e^{-2(R-|u|)}(1 - 3\cos^2\theta) = z_*$ locus. To set the notation, note that the $f \geq 1 + \delta^2$ part of $M_\delta \cap \mathcal{H}^+_p$ can be parametrized using the coordinates $(t, \varphi_+, h_+)$ from the relevant component of the $(1+\delta^2, 1+\delta_*^2) \times T_+$ portion of $(1,2) \times \Sigma$.

The inverse image via $\Psi_p$ of the $f \in (\delta^2, z_*]$ part of $\mathbb{R} \times \mathcal{H}_p^+$ corresponds to the part of $\mathbb{R} \times \mathcal{X}$ where $\hat{u} \in [R+\ln\delta, R+\frac{1}{2}\ln z_*]$. The coordinate functions $(x, \hat{u}, \hat{\phi}, h)$ on the $\hat{u} \in (R+\ln\delta, R+\frac{1}{2}\ln z_*]$ part of $\mathbb{R} \times \mathcal{X}$ are related to the coordinate functions $(s, t, \varphi_+, h_+)$ using the rule $(s = x, t = 1 + e^{-2(R-|\hat{u}|)}, \varphi_+ = \hat{\phi}, h_+ = h)$.

Granted the preceding, it follows that the functions $(\varphi^{p0}, \varsigma^{p0})$ that appear in (2.10) can be viewed as functions of the coordinates $(s, z)$. Doing so writes the part of $C_{p0}$ in the $f \in (1+\delta^2, 1+z_*]$ part of $\mathbb{R} \times (M_\delta \cap \mathcal{H}_p^+)$ as the graph over $\mathbb{R} \times (\delta^2, z_*]$ given by the rule

$$(x, z) \to (s = x, t = 1+z, \varphi_+ = \varphi^{p0}, h_+ = \varsigma^{p0}).$$
(2.11)

The fact that $C_{p0}$ is J-holomorphic implies that the pair $(\varphi^{p0}, \varsigma^{p0})$ obey the Cauchy-Riemann equations: $\partial_x \varphi^{p0} - \partial_z \varsigma^{p0} = 0$ and $\partial_x \varsigma^{p0} + \partial_z \varphi^{p0} = 0$.

There is a corresponding picture of $C_{p0}$ on the $t \in [2-z_*, 2-\delta_2)$ portion of $\mathbb{R} \times (M_\delta \cap \mathcal{H}_{p*}^+)$. This part of $\mathbb{R} \times \mathcal{H}_{p*}^+$ corresponds via $\Psi_p$ to the $\hat{u} \in [-R-\frac{1}{2}\ln z_*, -R-\ln\delta)$ part of $\mathbb{R} \times \mathcal{X}$. It is parametrized by coordinates $(x, z, \varphi_-, h_-)$ with $z \in (\delta^2, z_*]$ related to the coordinate $t$ by the rule $t = 2-z$ and $z$ related to $\hat{u}$ by the rule $z = e^{-2(R+|\hat{u}|)}$. The part of $C_{p0}$ here is parametrized by viewing $(\varphi^{p0}, \varsigma^{p0})$ as functions of $(x, z)$ and writing $\varphi_- = \varphi^{p0}$ and $h_- = -\varsigma^{p0}$. The pair $(\varphi^{p0}, \varsigma^{p0})$ obey here the Cauchy-Riemann equations when written as functions of $(x, -z)$.

What follows is now a crucial point: The functions $\varphi^{p0}$ is constrained on both boundary components of $C_{p0}$ as follows: Write $\varphi^{p0}$ on either boundary as a function of the coordinate $x \in \mathbb{R}$. Meanwhile, write the function $\varphi^S$ on the relevant $\mathbb{R} \times \{1+z_*\} \times T_+$ or $\mathbb{R} \times \{2-z_*\} \times T_-$ part of the boundary of $C_{S0}$ as a function of x also. Then

$$\varphi^{p0}(x, z_*) = \varphi^{S0}(x, z_*).$$
(2.12)

There is no apriori constraint on the value of $\varsigma^{p0}$ on the boundaries of $\partial C_{p0}$ but for what is implied by (2.9).

*Part 4*: The proposition given below states the fundamental existence/uniqueness theorem for $C_{p0}$. The proposition refers to the preceding Parts 1-3.

**Proposition 2.2**: *There exists a purely $S$-dependent (or $\mathcal{K}$-compatible) constant $\kappa > 1$ with the following significance: Define the geometry of $Y$ with $z_* < \kappa^{-1}$ and $\delta < \kappa^{-1} z_*$. Fix* *an element* $(\hat{\Theta}_-, \hat{\Theta}_+) \in \hat{\mathcal{Z}}^S$. *In the case* $\Delta_{\mathfrak{p}} = 1$, *fix a* $\hat{u} = 0$ *point in* $\mathbb{R} \times I_*$; *and in the case* $\Delta_{\mathfrak{p}} = 2$, *fix a pair of* $\hat{u} = 0$ *points in* $\mathbb{R} \times I_*$. *There exists a unique pair* $(\varphi^{\mathfrak{p}0}, \varsigma^{\mathfrak{p}0})$ *that are described by Parts 2 and 3 and are such that the* $\Psi_{\mathfrak{p}}$-*image of the surface given by (2.9) obeys the conditions stated for* $C_{\mathfrak{p}0}$ *in Part 1*.

The proof of Proposition 2.2 in the case when all $\mathfrak{p} \in \Lambda$ versions of $\Delta_{\mathfrak{p}}$ are zero is given separately in Section 4 because it has fewer components than the proof for the $\Delta_{\mathfrak{p}} > 1$ cases. The proof for general case is given in Section 6. The next section introduces certain analytic tools that are used in Section 4. Section 5 introduces some additional tools to handle the general case.

## 3. Cauchy-Riemann equations on $\mathbb{R} \times \mathcal{X}$

The $\Psi_{\mathfrak{p}}$-image of a graph of the form

$$(x, \hat{u}) \to (x, \hat{u}, \hat{\phi} = \varphi(x, \hat{u}), h = \varsigma(x, \hat{u}))$$

(3.1)

defines a J-holomorphic subvariety if and only if the pair of functions $(\varphi, \varsigma)$ satisfy a certain non-linear Cauchy-Riemann equation as functions of the coordinates $(x, \hat{u})$. The purpose of this section is to describe this equation and to supply various tools that will be used subsequently to construct desired solutions.

### a) Almost complex structures on $\mathbb{R} \times \mathcal{X}$

Suppose that $J_{\mathfrak{p}}$ is a given almost complex structure on $\mathbb{R} \times \mathcal{H}^+_{\mathfrak{p}*}$ with the property that the $\mathbb{R} \times \mathcal{H}^+_{\mathfrak{p}*}$ part of any surface from Proposition II.3.2's moduli space $\mathcal{M}_\Sigma$ and Proposition II.3.4's moduli space $\mathcal{M}_{\mathfrak{p}0}$ are $J_{\mathfrak{p}}$-holomorphic. Use $\Psi_{\mathfrak{p}}$ to view $J_{\mathfrak{p}}$ as an almost complex structure on $\mathbb{R} \times \mathcal{X}$.

The $\Psi_{\mathfrak{p}}$-inverse images of the surfaces from $\mathcal{M}_\Sigma$ and those from $\mathcal{M}_{\mathfrak{p}0}$ are the constant $(x, \hat{u})$ slices of $\mathbb{R} \times \mathcal{X}$. This understood, the fact that they are $J_{\mathfrak{p}}$-holomorphic has the following implication: The $J_{\mathfrak{p}}$ version of $T^{1,0}(\mathbb{R} \times \mathcal{X})$ must contain a form that can be written as

$$dx + i q_0 d\hat{u},$$

(3.2)

where $q_0$ is a $\mathbb{C}$-valued function with strictly positive real part. A second linearly independent 1-form for $J_{\mathfrak{p}}$'s version of $T^{1,0}(\mathbb{R} \times \mathcal{X})$ can always be written as

$$d\hat{\phi} + iq_1 dh + i\mu_0(dx - i\bar{q}_0 d\hat{u}),$$

(3.3)

where $q_1$ is a $\mathbb{C}$-valued function with strictly positive real part, and $\mu_0$ is a $\mathbb{C}$-valued function.

**Lemma 3.1**: *Suppose that $J_p$ obeys the $\mathbb{R} \times \mathcal{H}^+_{p*}$ versions of the first five bullets in Part 1 of Section 1c. Then $q_0, q_1$ and $\mu_0$ depend only on the coordinates $\hat{u}$ and $h$; and both $q_1$ and $\mu_0$ are real valued.*

*Proof of Lemma 3.1*: Whether real or not, the functions $q_0$, $q_1$ and $\mu_0$ are invariant with respect to constant translations along the $\mathbb{R}$ and $\mathbb{R}/(2\pi\mathbb{Z})$ factors in $\mathbb{R} \times X$ if $J_0$ is invariant with respect to the respective translations along the $\mathbb{R}$ and $\mathbb{R}/(2\pi\mathbb{Z})$ factors of $\mathbb{R} \times \mathcal{H}^+_{p*}$.

To see about the imaginary parts of $q_1$ and $\mu_0$ define the adjoint action of $J_p$ on the cotangent bundle by by the following rule: Let $\langle \, , \, \rangle$ denote the pairing between covectors and vectors and let $\mathfrak{e}$ and $w$ denote respective covectors and vectors over the same base point. Then $\langle J_p^T \mathfrak{e}, w \rangle = \langle \mathfrak{e}, J_p w \rangle$. Note that $J_p^T$ acts on (3.2) and on (3.3) as multiplication by $i$. With this in mind, use the identity $J_p \partial_s = v$ with the second bullet in (1.30) to see that $J_p$ acts as multiplication by $-i$ on $\partial_x + i(v^{-1}\partial_{\hat{u}} - \alpha^{-1}\sqrt{6} x \cos\theta \, \partial_{\hat{\phi}} - \varpi \partial_x)$. This vector is therefore sent to zero when paired via $\langle \, , \, \rangle$ with (3.2). Such is the case if and only if $(1 - i\varpi) - q_0 v^{-1} = 0$ so $q_0 = v(1 - i\varpi)$. This same vector is also sent to zero when paired via $\langle \, , \, \rangle$ with (3.3), and this happens if and only if $2\mu_0 = \alpha^{-1}\sqrt{6} x \cos\theta$. To see about $q_1$, use the fourth bullet in (1.30) and Equation (II.3.9) to see that $\Psi_{p*} J_p \partial_h$ is proportional to $\partial_\phi$, and so it follows from the third bulled of (1.30) that $J_p \partial_h$ is proportional to $\partial_{\hat{\phi}}$. Now use (3.3) to see that $q_1$ is real if and only if $\langle dh, J_p \partial_h \rangle = 0$.

Assume in what follows that $J_p$ obeys the assumptions of Lemma 3.1. With it understood that $q_1$ and $\mu_0$ are real, then (3.2) and (3.3) imply the following: A submanifold in $\mathbb{R} \times \mathcal{H}^+_{p*}$ given by the $\Psi_p$ image of the surface given by (3.1) is $J_p$-holomorphic if and only if the functions $(\varphi, \varsigma)$ that appear in (3.1) obey a system of Cauchy-Riemann equations that can be written as

$$a_1 \partial_x \varphi - \partial_{\hat{u}} \varsigma = 0 \quad \text{and} \quad a_2 \partial_x \varsigma + \partial_{\hat{u}} \varphi + b = 0,$$

(3.4)

where $a_1$, $a_2$, and $b$ constitute a set of $\mathbb{R}$-valued functions with $a_1$ and $a_2$ strictly positive. Their respective values at any given point $(x, \hat{u})$ are obtained from an eponymous set of the functions of the variables (û, h) by setting h = ς(x, û). This eponymous set is $\mathfrak{a}_1 = \frac{q_{0R}}{q_1}$, $\mathfrak{a}_2 = q_{0R}q_1$, and $\mathfrak{b} = -q_{0R}\mu_0$ where $q_{0R}$ denotes the real part of $q_0$.

Keep in mind that the expressions that appear on the left hand side of the two equations in (3.4) are only defined in the case that the map $\Psi_\mathfrak{p}$ is defined on the graph in (3.1). This is to say that the graph must define a surface in $\mathbb{R} \times \mathcal{X}$. This requirement constitutes an implicit constraint on the absolute value of ς at any given point (x, û). In particular, an assertion in the subsequent discussions that a given pair (φ, ς) solves (3.4) in all cases implies that ς obeys this implicit constraint.

## b) Maps from $\mathbb{R} \times I_*$ to $\mathbb{R}^2$ and linear operators

The central concern for the rest of this section are first order, linear operators on $C^\infty(\mathbb{R} \times I_*; \mathbb{R}^2)$ that are described next. Let D denote the operator in question and let $(\varphi', \varsigma'): \mathbb{R} \times I_* \to \mathbb{R}^2$ denote a given map. The respective first and second components of the map $D(\varphi', \varsigma')$ are

$$\mathfrak{a}_1 \partial_x \varphi' - \partial_{\hat{u}} \varsigma' + \mathfrak{b}_1 \varsigma' \quad and \quad \mathfrak{a}_2 \partial_x \varsigma' + \partial_{\hat{u}} \varphi' + \mathfrak{b}_2 \varsigma' \;,$$

(3.5)

where $\mathfrak{a}_1, \mathfrak{a}_2, \mathfrak{b}_1$ and $\mathfrak{b}_2$ are smooth functions of (x, û) with the following four properties:

- $\mathfrak{a}_1$ *and* $\mathfrak{a}_2$ *are everywhere positive and they both have uniform limits as* $x \to \pm\infty$ *to positive functions of û.*
- *The function on $\mathbb{R}$ given by the rule* $x \to \sup_{\{x\} \times I_*}(|\partial_x \mathfrak{a}_1| + |\partial_x \mathfrak{a}_2| + |\mathfrak{b}_1|)$ *limits uniformly as* $|x| \to \infty$ *with limit zero.*
- $\mathfrak{b}_2$ *limits uniformly as* $x \to -\infty$ *to a function of û. By the same token, $\mathfrak{b}_2$ limits uniformly as* $x \to \infty$ *to a function of û.*
- *The respective integral over $I_*$ of the $x \to \infty$ and $x \to -\infty$ limits of $\mathfrak{b}_2$ are non-zero and have the same sign.*

(3.6)

Of interest in what follows is a Fredholm incarnation of the operator D given by (3.5) and (3.6) whose domain and range are certain Hilbert spaces of maps from $\mathbb{R} \times I_*$ to $\mathbb{R}^2$ that is characterized as follows: The domain Hilbert space for D is the $L^2_1$ completion of the subspace of smooth maps from $\mathbb{R} \times I_*$ to $\mathbb{R}^2$ whose elements are as follows: A given pair (φ', ς') is in this subspace if and only if the following conditions are met:

- *The pair has compact support on $\mathbb{R} \times I_*$.*
- *The function φ' vanishes on the $|\hat{u}| = R + \frac{1}{2}\ln z_*$ boundaries of $\mathbb{R} \times I_*$.*

(3.7)

The square of the $L^2_1$-norm in question assigns to $(\varphi', \varsigma')$ the value

$$\int_{\mathbb{R}\times I_*} (|(\partial_x\varphi', \partial_x\varsigma')|^2 + |(\partial_{\hat{u}}\varphi', \partial_{\hat{u}}\varsigma')|^2 + |(\varphi', \varsigma')|^2)$$

(3.8)

Here, $|\cdot|^2$ denotes the Eucidean inner product on $\mathbb{R}^2$. The Hilbert space so defined is denoted in what follows by $\mathbb{H}$. The range space for this Fredholm version of D is the $L^2$ Hilbert space completion of the space of compactly supported elements in $C^\infty(\mathbb{R}\times I_*; \mathbb{R}^2)$; this is the completion that is defined using the inner product on $\mathbb{R}\times I_*$ whose square is the integral of $|(\varphi', \varsigma')|^2 = |\varphi'|^2 + |\varsigma'|^2$. This $L^2$ Hilbert space is denoted by $\mathbb{L}$.

The next proposition asserts the central fact about this Fredholm version of D.

**Proposition 3.1**: *The operator* D *as described by (3.5) and (3.6) with domain space* $\mathbb{H}$ *and range space* $\mathbb{L}$ *is Fredholm with index* 0 *and trivial kernel.*

Section 3c proves that D is Fredholm and Section 3d computes the index of D and proves that the kernel is trivial. The remainder of this subsection describes the relevant examples.

Let $\mathfrak{h} = (\varphi, \varsigma)\colon \mathbb{R} \times I_* \to \mathbb{R}^2$ denote a pair of functions with a graph given by (3.1) that lies in $\mathbb{R} \times \mathcal{X}$. The pair $\mathfrak{h}$ and the almost complex structure $J_\mathfrak{p}$ can be used to define a version of (3.5), this denoted by $D_\mathfrak{h}$. The definition is as follows: Let $(\varphi', \varsigma')$ denote a given, bounded map from $\mathbb{R} \times I_*$ to $\mathbb{R}^2$. Take $t \in \mathbb{R}$ near zero and write the expressions on the right hand sides of the two equations in (3.4) using pair $(\varphi + t\varphi', \varsigma + t\varsigma')$ in lieu of $(\varphi, \varsigma)$. View the result as a map from a neighborhood of 0 in $\mathbb{R}$ to $C^\infty(\mathbb{R}\times I_*; \mathbb{R}^2)$. The derivative of this map at $t = 0$ is $D_\mathfrak{h}(\varphi', \varsigma')$.

A more explicit description of $D_\mathfrak{h}$ is given momentarily. To this end, recall that the functions $(a_1, a_2, b)$ that appear in (3.4) are obtained from an eponymous set of functions, $(a_1, a_2, b)$, of the coordinates $(\hat{u}, h)$ for $\mathcal{X}$. Let $(a_{1h}, a_{2h}, b_h)$ denote the functions on $\mathbb{R} \times I_*$ whose respective values at any given point $(x, \hat{u}) \in \mathbb{R} \times I_*$ are those of the partial derivative with respect to h at the point $(x, \hat{u}, h = \varsigma(x, \hat{u}))$ of $(a_1, a_2, b)$. The respective first and second components of the $\mathbb{R}^2$-valued function $D_\mathfrak{h}(\varphi', \varsigma')$ can be written in terms of these partial derivatives as

- $a_1\partial_x\varphi' - \partial_{\hat{u}}\varsigma' + (a_{1h}\partial_x\varphi)\varsigma'$
- $a_2\partial_x\varsigma' + \partial_{\hat{u}}\varphi' + (a_{2h}\partial_x\varsigma + b_{2h})\varsigma'$

(3.9)

This observedly has the form depicted in (3.5).

What follows gives sufficient conditions on $\varsigma$ for (3.6) to hold.

**Lemma 3.2**: *Suppose that $\varsigma \in C^\infty(\mathbb{R} \times I_*)$ is a function with the following properties:*

- *The norm of $\varsigma$ at any given $(x, \hat{u}) \in \mathbb{R} \times I_*$ is such that point $(\hat{u}, \hat{\phi}, h = \varsigma(x,\hat{u}))$ is in $\mathcal{X}$.*
- *$|\varsigma|$ has respective limits as $x \to \pm\infty$; and both limits are less than $x_0 + 4e^{-2R}$.*
- *The function on $\mathbb{R}$ given by the rule $x \to \sup_{\{x\} \times I_*}(|\partial_x \varsigma| + |\partial_{\hat{u}} \varsigma|)$ limits uniformly to zero as $|x| \to \infty$.*

*Use $\varsigma$ with a given bounded function $\varphi$ to define the operator in (3.7). Then the corresponding version of $(\mathfrak{a}_1, \mathfrak{a}_2, \mathfrak{b}_1, \mathfrak{b}_2)$ obey the conditions in (3.6). In particular, this occurs if $\varsigma$ comes from a pair whose corresponding graph in $\mathbb{R} \times \mathcal{X}$ is the $\Psi_p$-inverse image of a surface in $\mathbb{R} \times \mathcal{H}^+_{p*}$ that is $J_p$-holomorphic where $|s| \gg 1$ and also obeys the conditions in the first and second bullets of (2.9).*

*Proof of Lemma 3.2*: The condition in the first bullet is required for defining $D_\mathfrak{h}$. Granted that $D_\mathfrak{h}$ is well defined, the fact that $\mathfrak{a}_1$ and $\mathfrak{a}_2$ in (3.5) are positive follows from (3.9) since the functions $a_1$ and $a_2$ that appear in (3.4) are positive. The conditions in the second and third bullets imply that the $x \to \pm\infty$ limits of $\varsigma$ exist and both are independent of $\hat{u}$. Moreover, the bound on these limits given by the second bullet imply that these respective values for h with a given value for $(\hat{u}, \hat{\phi})$ define a point in $\mathcal{X}$. It follows from the top bullet in (1.30) that the functions $(a_1, a_2, b)$ are independent of x where $|x| \gg 1$ and so depend only on the coordinates $\hat{u}$ and h where $|x| \gg 1$. This and the fact that $\varsigma$ limits to a constant implies that $\mathfrak{a}_1$ and $\mathfrak{a}_2$ have uniform $x \to \pm\infty$ limits that are positive functions of $\hat{u}$. This also implies that $\mathfrak{b}_2$ has uniform $x \to \pm\infty$ limits that are functions of $\hat{u}$. The condition stated in the second bullet of (3.6) follows from the second bullet of Lemma 3.1 via the chain rule.

To prove the fourth bullet in (3.6), use the third bullet of Lemma 3.1 to see that the $x \to \pm\infty$ limits of $\mathfrak{b}_2$ are those of $b_\mathfrak{h}$, and so given by the function $\hat{u} \to b_{2\mathfrak{h}}(\hat{u}, h_\pm)$, where $h_\pm$ are the corresponding $x \to \pm\infty$ limits of $\varsigma$. Meanwhile, the value of $\hat{h}$ in (1.27) is constant along any given integral curve of $v$ in $\mathcal{H}^+_p$ and therefore $\{b_{2\mathfrak{h}}(\hat{u}, h_+)\}_{\hat{u} \in I_*}$ and $\{b_2(\hat{u}, h_-)\}_{\hat{u} \in I_*}$ are the values of $b_\mathfrak{h}$ along the $\Psi_p$-inverse image of integral curves of $v$. To say more, let $\gamma$ denote the $\hat{u} \in I_*$ part of an integral curve of $v$ in $\mathcal{H}^+_{p*}$. The constant value of h on $\gamma$ and the $\phi$ coordinate of the $\hat{u} = 0$ point on $\gamma$ determines $\gamma$. If $\gamma$ is parametrized by the coordinate u on $\mathcal{H}^+_p$, then the coordinate $\phi$ on $\gamma$ changes via the rule in (II.2.5). With $\gamma$ parametrized by $\hat{u}$, the change in $\phi$ along $\gamma$ is given by

$$\frac{d\phi}{d\hat{u}} = -\frac{\sqrt{6}x\cos\theta}{f(1-3\cos^2\theta)} \frac{\partial u}{\partial \hat{u}} ,$$

(3.10)

where f is the function of u given in (1.4). Here, it is understood that $\theta$ is determined by u given the constant value of $\hat{h}$ on $\gamma$; and so it is determined by $\hat{u}$ and the constant value of h. This understood, if follows that $\gamma$'s version of $\mathfrak{b}_2$ is given by

$$\mathfrak{b}_2 = \frac{\partial^2 \phi}{\partial h \partial \hat{u}} = -\frac{\partial}{\partial h}\left(\frac{\sqrt{6} x \cos\theta}{f(1-3\cos^2\theta)} \frac{\partial u}{\partial \hat{u}}\right).$$

(3.11)

Let $\Delta\phi = \phi(\gamma|_{\hat{u}=R+\frac{1}{2}\ln z_*}) - \phi(\gamma|_{\hat{u}=-R-\frac{1}{2}\ln z_*})$. This number depends on $\gamma$ and so defines a function of the parameter h. The integral of $\mathfrak{b}_2$ over the domain $I_*$ is $\frac{d}{dh}(\Delta\phi)$. To compute the latter, use (3.10) to write

$$\Delta\phi = -\sqrt{6} \int_{-R-\frac{1}{2}\ln z_*}^{R+\frac{1}{2}\ln z_*} \frac{x\cos\theta}{f(1-3\cos^2\theta)} du.$$

(3.12)

To compute the h derivative of (3.12), introduce $\theta_0$ to denote the value of $\theta$ at the u = 0 point along $\gamma$. This is determined by h by solving $h|_\gamma = (x_0 + 4e^{-2R})\cos\theta_0 \sin^2\theta_0$ with the constraint that the solution $\theta_0$ is such that $1 - 3\cos^2\theta_0 > 0$. Lemma II.2.2 guarantees a unique solution. Meanwhile $\theta$ along $\gamma$ is determined at any given value of u by $\theta_0$ via the rule $f(u)\cos\theta \sin^2\theta = f(0)\cos\theta_0 \sin^2\theta_0$. This understood, it follows from (3.12) using the chain rule that

$$\frac{d}{dh}(\Delta\phi) = -\sqrt{6} \int_{-R-\frac{1}{2}\ln z_*}^{R+\frac{1}{2}\ln z_*} \frac{x(1+3\cos^2\theta)}{f^2(1-3\cos^2\theta)^3} du.$$

(3.13)

As can be seen, the expression on the right hand side is negative in all cases.

Consider now the final assertion of the lemma that concerns the case where $\varsigma$ comes from a pair $(\varphi, \varsigma)$ whose large $|x|$ values define, via $\Psi_p$ and the graph in (3.1), a surface in $\mathbb{R} \times \mathcal{H}^+_p$ that is J-holomorphic and has the asserted large $|s|$ behavior. The convergence condition implies that $\varsigma$ limits uniformly as $x \to \infty$ to a constant, and likewise as $x \to -\infty$. By the same token, the function $\hat{u} \to \varphi(x, \hat{u})$ also limits uniformly as $x \to \infty$ to a function of $\hat{u}$, and likewise as $x \to -\infty$. Granted these uniform limits, and given that $(\varphi, \varsigma)$ obey (3.4) at large $|x|$, the standard elliptic regularity theorems of the sort that can be found in Chapter 6 of Morrey's book [M] will prove that $\partial_x \varsigma$ and $\partial_{\hat{u}} \varsigma$ converge uniformly to zero as $|x| \to \infty$.

### c) Proof of Proposition 3.1: The Fredholm assertion

The proof that D is Fredholm has five steps.

Step 1: Let $\|\cdot\|$ denote the norm for $\mathbb{L}$. An operator such as D from $\mathbb{H}$ to $\mathbb{L}$ has finite dimensional kernel and closed range if there exists $c \geq 1$ such that the following holds for all elements $\eta$ in $\mathbb{H}$:

- $\|D\eta\|^2 \geq c^{-1} \|d\eta\|^2 - c \|\eta\|^2$.
- *If $\eta$ has support only where $|x| > c$, then* $\|D\eta\|^2 \geq c^{-1} \|\eta\|^2$.

(3.14)

A standard argument using the Rellich lemma uses (3.14) to deduce that D has closed range and finite dimensional kernel. The cokernel of D is isomorphic to the kernel of a certain formal $L^2$ adjoint which is also a bounded operator from $\mathbb{H}$ to $\mathbb{L}$. The upcoming Step 5, explains why (3.14) holds for this formal adjoint, and so the cokernel of D is finite dimensional.

Step 2: To prove what is asserted by the top bullet in (3.14), multiply the square of the left most expression in (3.5) by $\mathfrak{a}_1^{-1}$ and the square of the right most by $\mathfrak{a}_2^{-1}$. Integrate the resulting expressions over $\mathbb{R} \times I_*$. Use $\mathrm{M}(\eta)$ to denote the result of this integration. This number $\mathrm{M}(\eta)$ is relevant because $\mathrm{M}(\eta) > c_0^{-1} \|D\eta\|^2$. This understood, the bound that is asserted in the first bullet of the lemma is obtained with the help of an integration by parts to eliminate the term $(\partial_x \varphi' \partial_{\hat{u}} \varsigma' - \partial_u \varphi' \partial_x \varsigma')$ that appears in the integrand that defines $\mathrm{M}(\eta)$. There are no boundary terms from the integration by parts because of the second bullet in (3.7). With this term absent, the desired bound follows directly using the triangle inequality.

Step 3: To see about the second bullet in (3.14), introduce $\mathfrak{a}_{1\text{-}}$, $\mathfrak{a}_{2\text{-}}$ and $\mathfrak{b}_{2\text{-}}$ to denote the respective $x \to -\infty$ limits of $\mathfrak{a}_1$, $\mathfrak{a}_2$ and $\mathfrak{b}_2$. Let $Q^-$ denote the quadratic function on $C^\infty(I_*; \mathbb{R}^2)$ that is given by the rule

$$(\varphi', \varsigma') \to \int_{I_*} (\mathfrak{a}_{1\text{-}}^{-1} |\partial_{\hat{u}} \varsigma'|^2 + \mathfrak{a}_{2\text{-}}^{-1} |\partial_{\hat{u}} \varphi' + \mathfrak{b}_{2\text{-}} \varsigma'|^2).$$

(3.15)

Restrict this form to the subspace of pairs $(\varphi', \varsigma')$ with $\varphi' = 0$ at the boundary of the interval. On this restricted domain, the function $Q^-$ is such that

$$Q^-(\varphi', \varsigma') \geq c_0 R^{-2} \int_{I_*} (|\varphi'|^2 + |\varsigma'|^2).$$

(3.16)

Indeed, this follows directly given that there is no pair $(\varphi', \varsigma')$ with $\varphi = 0$ on the boundary of $I_*$ and with $Q(\varphi', \varsigma') = 0$. To see why no such pair exists, note first that if $Q(\varphi', \varsigma') = 0$, then $\varsigma'$ is constant and

$$\varphi'(\hat{u}) = \varsigma' \int_{-R-\frac{1}{2}\ln z_*}^{\hat{u}} \mathfrak{b}_{2-} \,.$$

(3.17)

Since the integral over $I_*$ of $\mathfrak{b}_{2-}$ is non-zero, the right hand side is not zero at $\hat{u} = R + \frac{1}{2}\ln z_*$ unless both $\varsigma'$ and $\varphi'$ is zero.

There is the analogous quadratic function, $Q^+$, on $C^\infty(I_*; \mathbb{R}^2)$ that is defined by the $x \to \infty$ limits of $\mathfrak{a}_1, \mathfrak{a}_2$ and $\mathfrak{b}_2$. The latter also dominates what is written on the right hand side of (3.16) when $\varphi'$ is zero where $|\hat{u}| = R + \frac{1}{2}\ln z_*$.

Hold on to the $Q^-$ and $Q^+$ versions of (3.16) for use momentarily.

Step 4: Suppose that $x_1 > 1$ and that $\eta = (\varphi', \varsigma')$ has compact support that lies where $x < -x_1$. Integrate by parts as instructed in Step 3, but now write the resulting expression for $M(\eta)$ as

$$\int_{x<-x_1} Q^-(\varphi', \varsigma') + \int_{\mathbb{R}\times I_*} (\mathfrak{a}_1 |\partial_x \varphi'|^2 + \mathfrak{a}_2 |\partial_x \varsigma'|^2) + 2 \int_{\mathbb{R}\times I_*} \mathfrak{b}_{2-} \varsigma' \partial_x \varsigma' + \mathfrak{e} \,,$$

(3.18)

where $|\mathfrak{e}| \leq \Delta (\|d\eta\|^2 + \|\eta\|^2)$ with $\Delta$ such that $\lim_{x_1 \to \infty} \Delta = 0$. Integrate by parts on the right most integral in (3.18) to see that it is zero. Meanwhile, the left most integral in (3.18) is no less than $c_0^{-1} R^{-2} (\|\varphi'\|^2 + \|\varsigma'\|^2)$. Thus, what is written in (3.18) is greater than

$$\|\partial_x \eta\|^2 + c_0^{-1} R^{-2} \|\eta'\|^2 \quad \text{if } x_1 > c.$$

(3.19)

This last bound implies what is asserted by the second bullet in (3.14) for the case when $\eta$ is supported where $x < -x_1$. But for notation, the same argument using $Q^+$ proves the second bullet of (3.14) for the case when $\eta$ is supported where $x > x_1$.

Step 5: Up to a sign, the formal adjoint in question is defined by using integration by parts to rewrite inner products with $D\eta$ using the inner product on $\mathbb{L}$. To be explicit, the operator sends any given $\eta = (\varphi^\#, \varsigma^\#)$ to the element in $\mathbb{L}$ with respective components

- $-\mathfrak{a}_1 \partial_x \varphi^\# - \partial_{\hat{u}} \varsigma^\# - (\partial_x \mathfrak{a}_1) \varphi^\#$
- $-\mathfrak{a}_2 \partial_x \varsigma^\# + \partial_{\hat{u}} \varphi^\# + \mathfrak{b}_1 \varphi^\# + (\mathfrak{b}_2 - \partial_x \mathfrak{a}_2) \varsigma^\#$

(3.20)

Use $D^\#$ to denote the operator in (3.20). What follows explains why $D^\#$ obeys the assertions made by the two bullets in (3.14).

An argument much like that used in Step 2 proves that $\|D^\#\eta\|^2 \geq c^{-1}\|d\eta\|^2 - c\|\eta\|^2$. Meanwhile, the assumptions that $|\partial_x \mathfrak{a}_1|$ and $|\partial_x \mathfrak{a}_2|$ limit to zero as $|x| \to \infty$ imply that $D^\#$ has the same form as that of D at large $|x|$ but for the sign changes in the derivative terms. As a consequence, the argument that proves the second bullet of (3.14) for D proves it for $D^\#$ as well.

**d) Proof of Proposition 3.1: The kernel and cokernel**

This subsection computes the Fredholm index of D and then its kernel dimension. Both are found equal to zero. The cokernel dimension is therefore zero as well. These computations are done in five steps.

*Step 1*: This step computes the index of D. This is done by deforming D to an operator whose index is readily computable. The discussion that follows concerns the case when the integral that is described in the fourth bullet of (3.6) is negative. A very much analogous discussion holds when the integral in question is positive.

The deformation is through a family of operators from $\mathbb{H}$ to $\mathbb{L}$ that all have the same schematic form as D. The family is parametrized by $[0, 1]$. Fix $r \in [0, 1]$ and the member parametrized by r sends any given $\eta = (\varphi', \varsigma') \in \mathbb{H}$ to the element in $\mathbb{L}$ whose respective components are

- $((1-r)\mathfrak{a}_1 + r)\partial_x \varphi' - \partial_{\hat{u}} \varsigma' + (1-r)\mathfrak{b}_1 \varsigma'$
- $((1-r)\mathfrak{a}_2 + r)\partial_x \varsigma' + \partial_{\hat{u}} \varphi' + ((1-r)\mathfrak{b}_2 - r)\varsigma'$ .

(3.21)

The $r = 0$ member is D and the $r = 1$ member sends $\eta$ to the element in $\mathbb{L}$ with respective components

$$\partial_x \varphi' - \partial_{\hat{u}} \varsigma' \quad and \quad \partial_x \varsigma' + \partial_{\hat{u}} \varphi' - \varsigma' .$$

(3.22)

The index of D is the same as this $r = 1$ version. To see that the latter has index equal to 0, suppose first that $\eta$ is such that what is written in (3.22) vanishes. Then $\varphi'$ obeys the second order equation $\partial_{\hat{u}}^2 \varphi' + \partial_x^2 \varphi' - \partial_{\hat{u}} \varphi' = 0$. Keeping in mind that $\varphi' = 0$ where $|\hat{u}| = R + \frac{1}{2}\ln z_*$, and that $|\varphi'|^2$ is integrable, the maximum principle demands that $\varphi'$ vanish identically, This the case, then $\varsigma'$ must be constant, and hence zero because $|\varsigma'|^2$ is also integrable.

As noted in Step 5 of the preceding subsection, the cokernel of the operator defined by (3.22) is isomorphic to the kernel of latter's version of what is depicted in

(3.20). This is the operator that sends any given η to the element in $\mathbb{L}$ with respective components

$$-\partial_x \varphi' - \partial_{\hat{u}} \varsigma' \quad \text{and} \quad -\partial_x \varsigma' + \partial_{\hat{u}} \varphi' - \varsigma'$$
(3.23)

Granted this form for $D^{\#}$, the same maximum principle argument applies to prove that it lacks a non-trivial square integrable kernel.

<u>Step 2</u>: This step proves that the kernel of D is trivial given a certain claim whose proof occupies the remaing steps. To this end, suppose that $\eta = (\varphi', \varsigma') \in \mathbb{H}$ is such that $D\eta = 0$. Let $\Gamma \subset \mathbb{R} \times I_*$ denote the locus where $\varsigma' = 0$. As explained in the subsequent steps, this set is non-empty; and it is either $\mathbb{R} \times I_*$ or it has the structure of a graph with the following properties:

- *The interiors of the edges are the components of the locus in $\Gamma$ where $d\varsigma' \neq 0$, and each vertex is a critical point of the map to $\mathbb{R}^2$ defined by $(\varphi', \varsigma')$.*
- *Each edge is a $C^1$-embedded, closed interval*
- *Each vertex has but a finite number of incident edges. No pair of distinct incident edges are tangent at any given vertex.*
- *Each interior vertex has an even number of incident edges; this number is at least 4.*
- *Each boundary vertex has an odd number of incident edges.*
- *Each edge is oriented by the restriction of $d\varphi$; and this is the orientation that is induced on the edge by viewing it as a boundary component of the $\varsigma' < 0$ locus.*

(3.24)

These last facts are not compatible with the fact that $\varphi' = 0$ where $|\hat{u}| = R + \frac{1}{2} \ln z_*$ and has limit zero as $|x| \to \infty$ unless $\Gamma = \mathbb{R} \times I_*$, in which case $\varsigma'$ is everywhere zero and thus so is $\varphi'$. To see why $\Gamma$ can not be a graph, suppose to the contrary that $\Gamma$ is described by (3.24). Let $U \subset \mathbb{R} \times I_*$ denote a component of the complement of $\Gamma$. Bullets 2-5 of (3.24) imply that $\partial U$ is piecewise smooth, and so any given differential form can be integrated between points on $\partial U$. Meanwhile, either $\varsigma' > 0$ in U or $\varsigma' < 0$. In either case, the final bullet implies that $d\varphi'$ is positive on the smooth part of $\partial U$ given a suitable orientation. As a consequence, $\varphi'$ increases monotonically along $\partial U$. This is not possible for it precludes an end point of any component of $\partial U$ where $|\hat{u}| = R + \frac{1}{2} \ln z_*$, and it precludes a non-compact component of $\partial U$, and it precludes a component with no boundary. The fact that $\Gamma \neq \emptyset$ precludes the case $U = \mathbb{R} \times I_*$.

<u>Step 3</u>: This step explains why (3.24) describes $\Gamma$ given that $\Gamma \neq \emptyset$ and $\Gamma \neq \mathbb{R} \times I_*$. To this end, let $\Gamma' \subset \Gamma$ denote the subset where $d\varsigma' \neq 0$. This is a smooth, 1-dimensional submanifold in $\mathbb{R} \times I_*$. It follows from (3.5) that $d\varphi´ > 0$ on the tangent line of $\Gamma´$ if the latter is oriented so that $d\varsigma´$ points towards the side where $\varsigma´ > 0$. It also follows that $d\varphi´ = 0$ at the points in $\Gamma-\Gamma´$. Thus $\Gamma-\Gamma´$ is a subset of the set of singular points of the map $\eta$ to $\mathbb{R}^2$. Let $p \in \Gamma-\Gamma´$ and set $\lambda = (\varphi´ - \varphi´(p)) + i\varsigma´$. This $\mathbb{C}$-valued function vanishes at p. In addition, (3.5) when written for $\lambda$ has the form $\overline{\partial}\lambda + \nu\lambda + \mu\overline{\lambda} = 0$ where $\nu$ and $\mu$ are smooth $\mathbb{C}$-valued functions. Here, $\overline{\partial} = \frac{1}{2}(\partial_x + i\partial_{\hat{u}})$.

To exploit this equation for $\lambda$, introduce $w = x - x(p) + i(\hat{u} - \hat{u}(p))$, this being a $\mathbb{C}$-valued coordinate function for $\mathbb{R} \times I_*$. It follows from the equation for $\lambda$ using Taylor's theorem with remainder that $\lambda$ near p must have the form $\lambda = mw^q + e$ where $m$ is a non-zero complex number, $q \geq 2$ is an integer and $|e| \leq c_0|w|^{q+1}$. Note that the unique continuation principle implies that q is finite. This depiction of $\lambda$ implies what is asserted by (3.24) about the interior vertices of $\Gamma$. The argument for the boundary vertices is very much the same after using the Schwarz reflection trick from Theorem 24 in [Ah] to view any given boundary point as an interior point of a domain to which $(\varphi´, \varsigma´)$ extend so as to solve a corresponding extension of (3.5).

Step 4: This step constitutes a digression that is needed to explain why $\Gamma \neq \emptyset$. To start, let $(\mathfrak{a}_{1*}, \mathfrak{a}_{2*}, \mathfrak{b}_*)$ denote either the $x \to \infty$ or $x \to -\infty$ limit of $(\mathfrak{a}_1, \mathfrak{a}_2, \mathfrak{b}_2)$. Introduce the operator $L: C^\infty(I_*; \mathbb{R}^2) \to C^\infty(I_*; \mathbb{R}^2)$ that is defined so as to send $\eta = (\iota, \lambda)$ to

$$L\eta = (-\partial_{\hat{u}}\lambda, \partial_{\hat{u}}\iota + \mathfrak{b}_*\lambda)$$

(3.25)

The relevant domain for L is the subspace in $C^\infty(I_*; \mathbb{R}^2)$ that consists of the pairs $(\sigma, \lambda)$ with $\iota = 0$ at the boundary points of $I_*$. A pair $(\iota, \lambda)$ in this domain is said to be a *weighted* eigenfunction for L if

$$-\partial_{\hat{u}}\lambda = E\mathfrak{a}_{1*}\iota \quad and \quad \partial_{\hat{u}}\iota + \mathfrak{b}_*\lambda = E\mathfrak{a}_{2*}\lambda,$$

(3.26)

with $E \in \mathbb{R}$. The number E is said to be a weighted eigenvalue. Straightforward variations of standard arguments show the following: The set of weighted eigenvectors is discrete, has no accumulation points and is unbounded in both directions. What is said in Step 3 implies that 0 is not a weighted eigenvector. Moreover, at most a finite number of weighted eigenvectors that share the same weighted eigenvalue. Third, if $\eta = (\iota, \lambda)$ and $\eta´ = (\iota´, \lambda´)$ are weighted eigenvectors with different weighted eigenvalues, then

$$\int_{I_*} (\mathfrak{a}_{1*}\iota\iota´ + \mathfrak{a}_{2*}\lambda\lambda´) = 0.$$

(3.27)

Finally, the $L^2$ completion of the domain of L is spanned by the set of weighted eigenvectors.

Let $\vartheta$ denote a minimal spanning set of weighted eigenvectors, here chosen so that if $\eta = (\iota, \lambda) \in \vartheta$, then

$$\int_{I_*} (\mathfrak{a}_{1*} \iota^2 + \mathfrak{a}_{2*} \lambda^2) = 1.$$

(3.28)

Step 5: To see that $\Gamma \neq \emptyset$, first integrate the right most equation in (3.5) on each slice of the form $\{x\} \times I_*$ to obtain

$$\partial_x \left( \int_{\{x\} \times I_*} \mathfrak{a}_2 \varsigma' \right) - \int_{\{x\} \times I_*} (\partial_x \mathfrak{a}_2) \varsigma' = - \int_{\{x\} \times I_*} \mathfrak{b}_2 \varsigma'.$$

(3.29)

To exploit this identity, suppose now that $\varsigma'$ is nowhere zero. No generality is lost by assuming that $\varsigma' > 0$. If the integrals of the $x \to \pm\infty$ limits of $\mathfrak{b}_2$ are negative, then (3.29) is used at points where $x \gg 1$. If the integrals of the $x \to \pm\infty$ limts of $\mathfrak{b}_2$ are positive, then (3.29) is used at points where $x \ll -1$. Except for cosmetics, the argument for the latter case is identical to that for the former. Granted this, only the case where the integrals of these limits of $\mathfrak{b}_2$ are negative is considered in what follows.

To make something of (3.29), use arguments much like those in Section 2.3 of [HT] to see that $(\varphi', \varsigma')$ can be written for $x \gg 1$ as

$$(\varphi', \varsigma') = c(e^{Ex} (\iota, \lambda) + \mathfrak{e})$$

(3.30)

where the notation is as follows: First, $c \in (0, \infty)$. Second, $(\iota, \lambda) \in \vartheta$ is an element with negative, weighted eigenvalue, this being E. Third, $\mathfrak{e}$ is such that the function $x \to e^{Ex} |\mathfrak{e}|$ has limit zero as $x \to \infty$. Use (3.26) to see that if $(\sigma, \lambda)$ is a weighted eigenvector, then $\lambda$ has transversal zero locus. This understood, it follows from (3.30) that $\varsigma'$ is positive where $x \gg 1$ if and only if $(\iota, \lambda)$ is such that $\lambda > 0$ at all interior points of $I_*$.

Granted the preceding, it follows from (3.26) and (3.30) that the weighted eigenvector that appears in (3.30) has $\lambda \geq 0$. It also follows from these equations that

$$\int_{\{x\} \times I_*} \mathfrak{b}_2 \varsigma' < -r |E| (1 - c^{-1}) \int_{\{x\} \times I_*} \mathfrak{a}_2 \varsigma' \quad \text{where } x \gg 1.$$

(3.31)

Use this last bound in (3.29) with the fact that $|\partial_x \mathfrak{a}_2| \to 0$ as $x \to \infty$ to deduce that

$$\partial_x ( \int_{\{x\} \times I_*} \mathfrak{a}_2 \varsigma') > c^{-1} \int_{\{x\} \times I_*} \mathfrak{a}_2 \varsigma' \quad \textit{for } x \gg 1$$

(3.32)

This inequality can not hold if $|\varsigma'|^2$ is integrable. As a consequence, the assumption that $\varsigma' > 0$ is untenable.

**e) The Banach spaces $\mathbb{H}_*$ and $\mathbb{L}_*$**

The norm that defines the Hilbert space $\mathbb{H}$ does not control the supremum norm of its elements. This being the case, the inverse function theorem that is used in what follows employs a slightly stronger norm. The author learned the latter from Morrey's book [M]. The definition requires the choice of a positive number that is less than $\frac{1}{100}$. Use $\upsilon$ in what follows to denote this number. This Banach space is denoted by $\mathbb{H}_*$. It is the closure of the space of pairs that obey (3.7) using a norm that is the sum of the $L^2_1$ norm used for $\mathbb{H}$ and a norm that is defined momentarily. The extra term in the norm for $\mathbb{H}_*$ is the square root of the function that assigns to a given pair $\eta' = (\varphi', \varsigma')$ the number

$$\sup_{(x, \hat{u}) \in \mathbb{R} \times I_*} \sup_{\rho \in (0,1)} \rho^{-\upsilon} \int_{\mathrm{dist}(\cdot, (x, \hat{u})) < \rho} |d\eta'|^2 \ .$$

(3.33)

Here, $\mathrm{dist}(\cdot, \cdot)$ denotes the Euclidean distance function. The norm on $\mathbb{H}_*$ is denoted by $\|\cdot\|_{\mathbb{H}*}$. The lemma below in part justifies the introduction of this space.

**Lemma 3.3**: *Elements in $\mathbb{H}_*$ are Hölder continuous with exponent $\frac{1}{2} \upsilon$ and the inclusion map from $\mathbb{H}_*$ into the corresponding Hölder Banach space is continuous. In particular, there exists a constant $\kappa > 1$ that depends only on $\upsilon$ and has the following significance: If $\mathfrak{f} \in \mathbb{H}_*$, then $|\mathfrak{f}| \leq \kappa \|\mathfrak{f}\|_{\mathbb{H}*}$. In addition, $\lim_{|x| \to \infty} |\mathfrak{f}|$ exists and it is zero; thus, elements in $\mathbb{H}_*$ have pointwise uniform limit zero as $|x| \to \infty$.*

*Proof of Lemma 3.3*: These assertions follow directly from Theorem 3.5.2 in Morrey's book [M].

A corresponding $L^2$ version of $\mathbb{H}_*$ is defined to be the closure of the space of compactly supported elements in $C^\infty(\mathbb{R} \times I_*; \mathbb{R}^2)$ using the norm given by the sum of the $L^2$ norm and that defined by replacing $d\eta'$ in (3.33) by $\eta'$. This last Banach space is denoted in what follows by $\mathbb{L}_*$. The norm on $\mathbb{L}_*$ when needed is denoted by $\|\cdot\|_{\mathbb{L}*}$.

**Lemma 3.4**: *An operator* D: $\mathbb{H} \to \mathbb{L}$ *of the sort described by (3.5) and (3.7) maps* $\mathbb{H}_*$ *to* $\mathbb{L}_*$; *and its inverse restricts to* $\mathbb{L}_*$ *so as to define a bounded linear operator from* $\mathbb{L}_*$ *to* $\mathbb{H}_*$.

*Proof of Lemma 3.4*: Given Proposition 3.1, the assertion follows from Theorems 3.5.2 and 5.4.1 of Morrey's book [M].

## 4. Proof of Proposition 2.2 when $\Delta_p = 0$

The first four subsections prove that there exists at least one pair $(\varphi^{p0}, \varsigma^{p0})$ that satisfies the requirements of Proposition 2.2. By way of a look ahead, the existence proof uses an open/closed argument for a certain 1-parameter family of $|x| \to \infty$ asymptotic conditions and $|\hat{u}| = R + \frac{1}{2} \ln z_*$ boundary conditions for (3.4). The parameter space is the interval [0, 1]; the parameter {1} boundary conditions are those required by Proposition 2.2. Meanwhile the parameter {0} case is designed so as to have an obvious solution. Use $\mathcal{I}$ to denote the subset of parameter values in [0, 1] for which (3.4) has a solution with the corresponding asymptotic conditions and boundary conditions. The set $\mathcal{I}$ is proved to be both open and closed. This being the case, and as $\{0\} \in \mathcal{I}$, so $\mathcal{I} = [0, 1]$ and there is at least one pair $(\varphi^{p0}, \varsigma^{p0})$ that satisfies the requirements of Proposition 2.2.

The final subsection proves that this is the only pair of functions that satisfies all of Proposition 2.2's criteria. This uniqueness proof uses a non-linear version of the argument that is used Section 3d to prove that the operator D in Proposition 3.1 has trivial kernel.

### a) The 1-parameter family

The definition of the family of asymptotic/boundary conditions has three parts.

*Part 1*: Reintroduce $\gamma_{p-}$ and $\gamma_{p+}$ from Part 1 of Section 2d. By way of a reminder, these are the respective segments of integral curves of $v$ in $\mathcal{H}^+_{p*}$ that come from $\Theta_-$ and $\Theta_+$. Each parameter value $\tau \in [0, 1]$ also labels a segment of an integral curve of $v$ that crosses $\mathcal{H}^+_{p*}$. The corresponding segment is denoted $\gamma_\tau$. The upcoming definition uses $\varphi^{S0}_+$ and $\varphi^{S0}_-$ to denote the respective $\hat{u} = R + \frac{1}{2} \ln z_*$ and $\hat{u} = -R - \frac{1}{2} \ln z_*$ versions of the function $\varphi^{S0}(\cdot, z_*)$ that appears in (2.11). The segment $\gamma_\tau$ is the unique integral curve of $v$ in $\mathcal{H}^+_{p*}$ that obeys the following three constraints:

- *The segment $\gamma_\tau$ starts on the surface where* $e^{-2(u+R)}(1 - 3\cos^2\theta) = z_*$ *and it ends on the surface where* $e^{2(u-R)}(1 - 3\cos^2\theta) = z_*$.

- *If $\tau \in (0, 1)$, then the $\phi$ coordinate of $\gamma_\tau$ at its starting point is $\phi^{SO}_-(x = \frac{2\tau - 1}{\tau(1-\tau)}), z_*)$, and the $\phi$ coordinate of $\gamma_\tau$ at its ending point is $\phi^{SO}_+(x = \frac{2\tau - 1}{\tau(1-\tau)}, z_*)$.*
- *The segment $\gamma_{\tau=0}$, is $\gamma_{p-}$ and the segment $\gamma_{\tau=1}$ is $\gamma_{p+}$.*

(4.1)

Lemma II.2.2 supplies the desired segment $\gamma_\tau$.

The next lemma addressed the continuity and differentiability of the family $\{\gamma_\tau\}_{\tau \in [0,1]}$. This lemma views each integral curve from the family $\{\gamma_\tau\}_{\tau \in [0,1]}$ as a map from [0, 1] into $\mathcal{H}^+_{p*}$ that pulls *a* back as a constant multiple of the Euclidean differential.

**Lemma 4.1**: *The assignment of the point $\gamma_\tau(\sigma)$ to any given pair $(\tau, \sigma) \in [0,1] \times [0,1]$ defines a smooth map from $[0,1] \times [0,1]$ into $\mathcal{H}^+_{p*}$.*

*Proof of Lemma 4.1*: By construction, the map is continous on $[0,1) \times [0,1]$ and smooth on $(0,1) \times [0,1]$. It follows from the fourth bullet of Proposition 2.1 using the Chain rule that the map is smooth on $[0,1) \times [0,1]$. By the same token, if the map is continous up to along the $\{1\} \times [0,1]$ boundary, then it is also smooth up to and along $\{1\} \times [0,1]$. To see about continuity along this boundary, remark that $\lim_{\tau=1} \gamma_\tau$ exists; and this limit is a segment of an integral curve of $v$ that crosses $\mathcal{H}^+_{p*}$. Let $\gamma$ denote this limit. The issue is whether $\gamma$ is $\gamma_{p+}$. As explained next, such is the case because $\mathfrak{m}_p = 0$. To prove this, note that $\gamma_{p-}$ and $\gamma_{p+}$ concatenate with the [0,1] parametrized paths $\tau \to \gamma_\tau(0)$ and $\tau \to \gamma_\tau(1)$ to define a piecewise smooth, closed 1-cycle in $\mathcal{H}^+_{p*}$. The $\mathfrak{m}_p = 0$ condition implies that this 1-cycle is null-homotopic. Let $\iota$ denote this 1-cycle. Meanwhile, the paths $\gamma_{p-}$, $\gamma$ and the [0,1] parametrized paths $\tau \to \gamma_\tau(0)$ and $\tau \to \gamma_\tau(1)$ also concatenate to define a closed 1-cycle in $\mathcal{H}^+_{p*}$. Use $\iota'$ to denote the latter. The 1-cycle $\iota'$ is also null-homotopic as it bounds the surface given by the closure of the image of the map from $[0, 1) \times [0,1]$ that sends $(\tau,\sigma)$ to $\gamma_\tau(\sigma)$. Keeping this in mind, remark that $\gamma$ and $\gamma_{p+}$ have the same endpoints, and so the change, $\Delta\phi$, of the angle $\phi$ along $\gamma$ must differ from that along $\gamma_{p+}$ by an integer multiple of $2\pi$. This integer is zero if and only if $\iota$ and $\iota'$ are homotopic. This being the case, it follows from Lemma 2.1 that $\gamma = \gamma_{p+}$.

*Part 2*: Reintroduce the nonincreasing function, $\chi: \mathbb{R} \to [0, 1]$ which has value 1 on $(-\infty, 0]$ and value 0 on $[1, \infty)$. Given $\tau \in [0, 1]$, introduce $x_\tau$ to denote $\frac{2\tau - 1}{\tau(1-\tau)}$ and then define the function $\hat{x}_\tau: \mathbb{R} \to (-\infty, x_\tau)$ by the rule

$$x \to \hat{x}_\tau(x) = x\chi(x - x_\tau + 3) + x_\tau(1 - \chi(x - x_\tau + 3)).$$

(4.2)

This definition is such that $\hat{x}_\tau = x$ for $x < x_\tau - 3$ and $\hat{x}_\tau = x_\tau$ for $x > x_\tau + 1$. The derivative of this function $\hat{x}_\tau$ is non-negative and bounded from above by 4. Its derivatives to any given order greater than 1 also enjoy $\tau$-independent bounds.

*Part 3*: What follows are the parameter $\tau \in [0, 1]$ asymptotic/boundary conditions for (3.4).

- $\lim_{x \to -\infty} (\varphi, \varsigma)|_x = (\phi, h)|_{\gamma_{p-}}$ and $\lim_{x \to \infty} (\varphi, \varsigma)|_x = (\phi, h)|_{\gamma_\tau}$.
- $\varphi(\cdot, \hat{u} = -R - \frac{1}{2} \ln z_*) = \varphi^{SO}_-(\hat{x}_\tau(\cdot), z_*)$ and $\varphi(\cdot, \hat{u} = R + \frac{1}{2} \ln z_*) = \varphi^{SO}_+(\hat{x}_\tau(\cdot), z_*)$.

(4.3)

To say more about these conditions, remark that the $\tau = 1$ version of the top bullet in (4.3) reproduces the first two bullets in (2.9), and the $\tau = 1$ version of the bottom bullet in (4.3) reproduces (2.12). Meanwhile, the $\tau = 0$ version of (4.3) demands that

- $\lim_{s \to -\infty} (\varphi, \varsigma)|_x = (\phi, h)|_{\gamma_{p-}}$ and $\lim_{x \to \infty} (\varphi, \varsigma)|_x = (\phi, h)|_{\gamma_{p-}}$.
- $\varphi(x, \hat{u} = -R - \frac{1}{2} \ln z_*) = \phi(\gamma_{p-}|_{\hat{u} = -R - \frac{1}{2} \ln(z_*)})$ and $\varphi(x, \hat{u} = R + \frac{1}{2} \ln z_*) = \phi(\gamma_{p-}|_{\hat{u} = -R - \frac{1}{2} \ln(z_*)})$.

(4.4)

Note in particular that the equations in (3.4) with the boundary conditions in (4.4) are solved by the x-independent given by the pair $(\varphi, \varsigma)$ with $\varphi(\hat{u}) = \phi(\gamma_{p-}|_{\hat{u}})$ and $\varsigma$ the constant function $\varsigma = h(\gamma_{p-})$. This is to say that the $\Psi_p$ image of the corresponding image of (3.1) is the J-holomorphic surface $\mathbb{R} \times \gamma_{p-}$.

**b) Proof that $\mathcal{T}$ is open**

Let $\mathcal{T} \subset [0, 1]$ denote the set of parameters for which (3.4) has a solution that obeys the given parameter's version of (4.3). As noted at the end of the previous subsection, the set $\mathcal{T}$ contains 0; and so it is not empty. This subsection proves that $\mathcal{T}$ is open. The argument for this has four parts.

*Part 1*: For $\tau \in [0, 1]$, define the pair of functions $(\varphi_\tau, \varsigma_\tau)$ on $\mathbb{R} \times I_*$ using the rule

- $\varphi_\tau(x, \hat{u}) = \chi(R + \frac{1}{2} \ln z_* + \hat{u}) \varphi^{SO}_-(\hat{x}_\tau(x), z_*) + \chi(R + \frac{1}{2} \ln z_* - \hat{u}) \varphi^{SO}_+(\hat{x}_\tau(x), z_*)$
  $+ (1 - \chi(R + \frac{1}{2} \ln z_* + \hat{u})) \chi(x) \phi(\gamma_{p-}|_{\hat{u}}) + (1 - \chi(R + \frac{1}{2} \ln z_* - \hat{u}))(\chi(-x) \phi(\gamma_\tau|_{\hat{u}}))$.
- $\varsigma_\tau(x, \hat{u}) = \chi(x) h(\gamma_{p-}) + \chi(-x) h(\gamma_\tau)$.

(4.5)

Here, the notation has $\phi(\gamma|_{\hat{u}})$ with $\gamma = \gamma_{p-}$ or $\gamma = \gamma_\tau$ denoting the lift to $\mathbb{R}$ of the coordinate $\phi$ on $\gamma$'s intersection with the $\Psi_p$ image of the $\hat{u} \in I_*$ slice of $\mathbb{R} \times X$. Note in this regard that γ has transversal intersection with this slice, this being a consequence of what is said by the second bullet in (1.30). The lift $\phi(\gamma|_{\hat{u}})$ is chosen so that its value at the $\hat{u} = -R - \frac{1}{2}\ln z_*$ starting point of γ is that of the function $\phi^{S0}_-$. With this choice, the value of this lift at the $\hat{u} = R + \frac{1}{2}\ln z_*$ endpoint of γ is that of $\phi^{S0}_+$. Note also that the function $h$ is constant on the integral curves of ν in $\mathcal{H}^+_{p*}$; what is written as $h(\gamma)$ in the lower bullet of (4.5) is the constant value of $h$ on γ.

The pair $(\phi_\tau, \varsigma_\tau)$ obeys the parameter τ boundary condition given by (4.3).

*Part 2*: Reintroduce the Banach space $\mathbb{H}_*$ from Section 3e. Let $\mathbb{B}_* \subset \mathbb{H}_*$ denote a small radius ball about the origin, chosen so that elements in $\mathbb{B}_*$ have pointwise norm bounded by $c_0^{-1} x_0$. The norm is chosen so that any given $(\hat{\phi}, \hat{\varsigma})$ in $\mathbb{B}_*$ has $|\hat{\varsigma}| \leq c_0^{-1} x_0$ at all points in $\mathbb{R} \times X$. Lemma 3.3 supplies such a ball. Reintroduce the Banach space $\mathbb{L}_*$ from Section 3e as well.

The rule that follows defines a map $[0, 1] \times \mathbb{B}_*$ to $\mathbb{L}_*$ if $\mathbb{B}_*$ has small radius. The desired map sends any given element $(\tau, (\hat{\phi}, \hat{\varsigma}))$ to the pair of functions in $\mathbb{L}_*$ with respective components

- $a_1 \partial_x(\hat{\phi} + \phi_\tau) - \partial_{\hat{u}}(\hat{\varsigma} + \varsigma_\tau)$.
- $a_2 \partial_x(\hat{\varsigma} + \varsigma_\tau) + \partial_{\hat{u}}(\hat{\phi} + \phi_\tau) + b$.

(4.6)

To say more about the notation, the functions $(a_1, a_2, b)$ are viewed as functions on $\mathbb{R} \times I_*$ that depend implicitly on $\varsigma = \hat{\varsigma} + \varsigma_\tau$. As in (3.4), their values at any given $(x, \hat{u}) \in \mathbb{R} \times I_*$ are obtained from an eponymous set of $\hat{\phi}$-independent functions on $X$ by evaluating the latter at the point $(\hat{u}, h = \varsigma(x, \hat{u}))$. The size constraint on the radius $\mathbb{B}_*$ is needed to guarantee that the $(\hat{\phi} + \phi_\tau, \hat{\varsigma} + \varsigma_\tau)$ version of (3.1) defines a point in $\mathbb{R} \times X$.

The map defined by (4.6) is denoted in what follows by $\mathcal{F}$. This map is designed so that any given $(\tau, (\hat{\phi}, \hat{\varsigma})) \in [0,1] \times \mathbb{B}_*$ version of $(\phi = \hat{\phi} + \phi_\tau, \varsigma = \hat{\varsigma} + \varsigma_\tau)$ obeys (3.4) and the parameter τ version of (4.3) if and only $\mathcal{F}(\tau, (\hat{\phi}, \hat{\varsigma})) = 0$.

*Part 3*: The next lemma summarizes the contents of this step.

**Lemma 4.2**: *Fix $\tau \in \mathcal{I}$ and let $(\phi, \varsigma)$ denote a corresponding solution to (3.4) with boundary values given by the parameter τ version of (4.3). There exists a neighborhood $\mathcal{I}_\tau \subset [0, 1]$ of τ and a continuous map from $\mathcal{I}_\tau$ to $\mathbb{B}_*$ of the following sort: Given $\tau' \in \mathcal{I}_\tau$, use $(\hat{\phi}, \hat{\varsigma}) \in \mathbb{B}_*$ to denote the corresponding element. Then*

- $\mathcal{F}(\tau', (\hat{\varphi}+\varphi_{\tau'}, \hat{\varsigma})) = 0$ *and so* $(\hat{\varphi}+\varphi_{\tau}, \hat{\varsigma}+\varsigma_{\tau})$ *solves (3.4) with the parameter* $\tau$ *asymptotic/boundary conditions from (4.3).*
- *This solution for* $\tau' = \tau$ *is the given pair* $(\varphi, \varsigma)$.

*Proof of Lemma 4.2*: It follows from Lemma 3.3 with the fourth bullet of Proposition 2.1 that $\mathcal{F}$ defines a smooth map from $[0, 1] \times \mathbb{B}_*$ to $\mathbb{L}_*$ and that any $\tau \in [0, 1]$ version of $\mathcal{F}(\tau, \cdot)$ defines a smooth map from $\mathbb{B}_*$ to $\mathbb{L}_*$ whose derivatives to any given order are bounded uniformly as $\tau$ varies in $[0, 1]$. The proof that $\mathcal{F}$ is smooth on $[0, 1] \times \mathbb{B}_*$ invokes the fourth bullet of Proposition 2.1 to establish that the maps $\tau \to \varphi^{S0}_-(\hat{x}_\tau(\cdot), z_*)$ and $\tau \to \varphi^{S0}_+(\hat{x}_\tau(\cdot), z_*)$ from $[0, 1]$ to $C^\infty(\mathbb{R})$ are smooth on the interval $[0, 1]$. It follows from Lemma 3.4 that the differential along the $\mathbb{B}_*$ component of $[0, 1] \times \mathbb{B}_*$ is an isomorphism from $\mathbb{H}_*$ to $\mathbb{L}_*$. These last facts with the inverse function theorem prove the lemma.

*Part 4*: Granted that the solutions given by Lemma 3.4 are smooth, it then follows that the set $\mathcal{I}$ is an open subset of $[0, 1]$. Meanwhile, the fact that these solutions are smooth follows using elliptic regularity arguments of the sort that can be found in Chapter 6 of [M]. Note in this regard that the equations in (3.4) are linear with constant coefficients on the part of $\mathbb{R} \times I_*$ where $|\hat{u}| > R + \frac{1}{2}\ln z_* + \ln\delta$. This being the case, standard boundary regularity arguments for the Laplace equation can be employed to prove that the solutions are smooth along the boundary of $\mathbb{R} \times I_*$.

**c) Proof that $\mathcal{I}$ is closed**

The assertion that the set $\mathcal{I} \subset [0, 1]$ is a closed set is a consequence of the upcoming Lemma 4.3.

Lemma 4.3 uses the following notation: Given a positive integer k and a function $\alpha$ on $\mathbb{R} \times I_*$, the lemma has $\nabla^{(k)}\alpha$ denoting the tensor of k'th order partial derivatives of $\alpha$.

**Lemma 4.3**: *There exists a purely* S-*dependent (or* $\mathcal{K}$-*compatible) constant* $\kappa > 1$ *with the following significance: Define the geometry of* Y *with* $\delta \leq \kappa^{-1}z_*$. *Then the space of solutions to (3.4) with asymptotic and boundary conditions given by versions of (4.3) is sequentially compact in the strong* $C^\infty$ *sense. To elaborate, let* $\{(\tau_n, (\varphi_n, \varsigma_n))\}_{n=1,2,\ldots}$ *denote a sequence such that* $\tau_n \in [0, 1]$ *and such that* $(\varphi_n, \varsigma_n)$ *is a solution to (3.4) with asymptotic/boundary conditions given by the* $\tau_n$-*version of (4.3). There exists* $\tau \in [0,1]$ *and a solution,* $(\varphi, \varsigma)$, *to (3.4) with asymptotic/boundary conditions given by the* $\tau$ *version* *of (4.3). Moreover, there exists a subsequence from the sequence (hence renumbered consecutively) such that $\{\tau_n\}_{n=1,2,...}$ converges to $\tau$, and such that*

$$\lim_{n\to\infty} \sup_{(x,\hat{u})\in\mathbb{R}\times I_*} (|\nabla^{(k)}(\varphi-\varphi_n)| + |\nabla^{(k)}(\varsigma-\varsigma_n)|) = 0$$

*for any given positive integer* k.

*Proof of Lemma 4.3:* Suppose to start that the conditions stated in the subsequent equation are satisfied if $\delta < \kappa z_*$ for some purely S-dependent (or $\mathcal{K}$-compatible) $\kappa > 1$.

- *There exists $\kappa_* > 1$ such that $1 - 3\cos^2\theta > \kappa_*^{-1}$ on the $\Psi_p$ image of any the graph that is defined via (3.1) by any given solution to (3.4) with asymptotic/boundary conditions given by any given $\tau \in [0,1]$ version of (4.3).*
- *Given $\varepsilon > 0$, there exists $\kappa_\varepsilon > 1$ such that the following is true: Let $(\varphi,\varsigma)$ denote a solution to (3.4) with asymptotic/boundary conditions given by some $\tau \in [0,1]$ version of (4.3). Then $(\varphi,\varsigma)|_x$ defines a pair of functions on $I_*$ that differs by less than $\varepsilon$ from its respective $x \to -\infty$ and $x \to \infty$ limit when $x < -\kappa_\varepsilon$ and $x > \kappa_\varepsilon$.*

(4.7)

If (4.7) holds, then by now standard elliptic regularity arguments as can be found in Chapter 6 of [M] prove Lemma 4.3. Note in this regard that the original set of three $\hat{\phi}$-independent functions on $\mathcal{X}$ that are used to define the $(x, \hat{u})$ dependent coefficient functions $(a_1, a_2, b_2)$ in (3.4) have uniformly bounded derivatives to any given order on the $\Psi_p$-inverse of any subset of $\mathbb{R} \times \mathcal{H}^+_{p*}$ where there is a positive lower bound for $1 - 3\cos^2\theta$. Note also that (3.4) is a linear equation with constant coefficients on the $|\hat{u}| > R + \frac{1}{2}\ln z_* + \ln\delta$ part of $\mathbb{R} \times I_*$ when written in using the variables $(x, z = e^{-2(R-\hat{u})})$ on the positive $\hat{u}$ part and $(-x, z = e^{-2(R+\hat{u})})$ on the negative $\hat{u}$ part. This being the case, standard boundary regularity arguments for the Cauchy-Riemann equations can be employed to prove that the solutions are smooth along the boundary of $\mathbb{R} \times I_*$.

Given what was just said, it remains to prove that (4.7) holds. This is done in five steps.

Step 1: A key input is a bound for the integral of the 2-form $w$ over the $\Psi_p$ image in $\mathbb{R} \times \mathcal{H}^+_{p*}$ of the graph of a solution to (3.4) with boundary values given by some parameter $\tau \in [0, 1]$ version of (4.3). The next lemma is used to derive such a bound. It has a second use in a subsequent step.

**Lemma 4.4**: *There exists a purely S-dependent (or $\mathcal{K}$-compatible) constant $\kappa \geq 1$ with the following significance: Fix $c \in (2\delta^2, \frac{4}{3\sqrt{3}}\delta_*^2)$ and $z_{**} \in (\delta^2, z_*]$. Suppose that $(\varphi, \varsigma)$* *solves (3.4), obeys a given $\tau \in [0, 1]$ version of (4.3), and is such that $|\varsigma| \leq c$ on the slice of $\mathbb{R} \times I_*$ where $|\hat{u}| = R + \frac{1}{2}\ln z_{**}$. Then $|\varsigma| \leq c + \kappa(z_* - z_{**})$ where $R + \frac{1}{2}\ln z_{**} \leq |\hat{u}| \leq R + \frac{1}{2}\ln z_*$.*

It is important to keep in mind for the subsequent applications of Lemma 4.4 that the constant $\kappa$ from this lemma depends neither on $\delta_*$ nor on $z_*$ when the latter are small.

*Proof of Lemma 4.4*: The argument that follows establishes the asserted upper bound for the $\hat{u} = R + \frac{1}{2}\ln z_*$ boundary given that $|\varsigma| \leq c$ where $\hat{u} = R + \frac{1}{2}\ln z_{**}$. A completely analogous argument does the trick for the other boundary component of $\mathbb{R} \times I_*$. To start, write $(\varphi, \varsigma)$ in terms of the coordinate x and a coordinate $z = e^{-2(R-\hat{u})}$. The $\mathbb{C}$-valued function $\varphi + i\varsigma$ is a holomorphic function of $x + iz$ where $z \in [z_{**}, z_*]$ and as a consequence, the function $\varsigma$ is annihilated by the operator $\partial_x^2 + \partial_z^2$. Keep this fact in mind. Now define

$$R = 100 \sup_{\mathbb{R}} |\partial_x \varphi^{SO}_+(\cdot, z_*)|$$

(4.8)

This constant is purely S-dependent (or $\mathcal{K}$-compatible). With R in hand, use w to denote x-independent function on $\mathbb{R} \times [z_{**}, z_*]$ given by the rule $z \to w(z) = c + R(z - z_{**})$. This function is also harmonic. Its value where $z = z_{**}$ is greater than that of $|\varsigma|$ and its $x \to \pm\infty$ limits are greater than those of $|\varsigma|$. Meanwhile, its z-derivative where $z = z_*$ is greater than $|\partial_z \varsigma|$ where $z = z_*$ because the Cauchy-Riemann equations identify $\partial_z \varsigma$ with $\partial_x \varphi$; and $\varphi$ where $z = z_*$ is given by $\varphi^{SO}_+$ via (4.3). These various upper bounds with the maximum principle imply that $w \geq |\varsigma|$ on the whole strip $\mathbb{R} \times [z_{**}, z_*]$.

Step 2: This step states and then proves the desired bound on the integral of w.

**Lemma 4.5**: *There exists a purely S-dependent (or $\mathcal{K}$-compatible) $\kappa \geq 1$, and there exists $\kappa_\Theta \geq 1$ that depends only on S and $\Theta_\pm$; and these have the following significance: Suppose that $\tau \in [0,1]$ and that $(\varphi, \varsigma)$ is a solution to the corresponding version of (3.4). Let $C \subset \mathbb{R} \times \mathcal{H}^+_{p*}$ denote the $\Psi_p$-image of the graph of $(\varphi, \varsigma)$. Then*

- $\int_C w \leq \kappa \delta_*^2$ .

- $\int_{C \cap ([s_0, s_0+1] \times \mathcal{H}^+_{p*})} ds \wedge \hat{a} \leq \kappa_\Theta$ *for all $s_0 \in \mathbb{R}$.*

As in the case of Lemma 4.4, the constant κ supplied by Lemma 4.5 depends neither on $\delta$ nor on $z_*$ when the latter are small.

*Proof of Lemma 4.5*: The fact that the integral of $w$ over the $s \ll -1$ part of $C$ is finite is proved in the next paragraph. But for obvious notational changes, the same argument proves the finiteness of the integral over the $s \gg 1$ part of $C$.

Let $\hat{u} \to (\varphi_{\gamma_{p-}}, \varsigma_{\gamma_{p-}})|_{\hat{u}}$ denote the pair whose graph is the $\Psi_p$-inverse image of the J-holomorphic cylinder $\mathbb{R} \times \gamma_{p-}$. Note in this regard that $\varsigma_{\gamma_{p-}}$ is the constant value of the function $h$ on $\gamma_{p-}$. Standard elliptic regularity theorems of the sort found in Chapter 6 of [M] in conjunction with (4.3) and (3.4) prove the following: Fix $\varepsilon > 0$ and there exists $s_\varepsilon \geq 1$ such that any $s \geq s_\varepsilon$ version of $(\varphi,\varsigma)|_s$ has $C^1$ distance at most $\varepsilon$ from $(\varphi_{\gamma_{p-}}, \varsigma_{\gamma_{p-}})$. This fact with (1.6), the fact that $w$ is non-negative on $\wedge^2 TC$ and Stokes' theorem imply the finiteness claim.

What follows next explains why the bound given by Lemma 4.5's first bullet holds. To start, use the just described application of Stokes' theorem to identify the integral of $w$ over $C$ with the sum of the following two expressions:

- $\int_{\gamma_\tau} h\, d\varphi - \int_{\gamma_{p-}} h\, d\varphi + \int_{\mathbb{R} \times \{\hat{u} = R + \frac{1}{2}\ln z_*\}} \varsigma\, d\varphi - \int_{\mathbb{R} \times \{\hat{u} = -R - \frac{1}{2}\ln z_*\}} \varsigma\, d\varphi \, .$

- $\int_{\gamma_\tau} x(1 - 3\cos^2\theta)\, du - \int_{\gamma_{p-}} x(1 - 3\cos^2\theta)\, du \, .$

(4.9)

To bound the left-most two terms in the top bullet, note first that the function $\varsigma$ is constant on the integral curves of $v$ in $\mathcal{H}^+_{p*}$; its value being that of the function $h$ depicted in (1.27). As a consequence, the difference between these two terms can be written as

$$h(\gamma_\tau)\Delta\varphi_{\gamma_\tau} - h(\gamma_{p-})\Delta\varphi_{\gamma_{p-}} \, ,$$

(4.10)

where $\Delta\varphi_{(\cdot)}$ is the change in the coordinate $\varphi$ along the indicated integral curve. It follows from (1.4) that the two values of $h$ are bounded in absolute value by $c_0 x_0$. Meanwhile, the two values of $\Delta\varphi_{(\cdot)}$ are determined by S and the $\Delta_p = 0$ constraint. It follows as a consequence that the two left most terms in (4.9) are bounded by a purely S-dependent (or $\mathcal{K}$-compatible) multiple of $x_0$.

To bound the two right most terms in the top bullet, keep in mind that $|\varsigma| \leq \frac{4}{3\sqrt{3}} \delta_*^2$ in any event. This the case, the two right most terms are bounded in absolute value by $c_0 \delta_*^2$ times the the integral over $\mathbb{R}$ of the function $x \to (|\partial_x \varphi^{S0}_-(x,z_*)| + |\partial_x \varphi^{S0}_+(x,z_*)|)$. The latter integrals are bounded by a purely S dependent (or $\mathcal{K}$-compatible) constant, this by virtue of the fourth bullet in Proposition 2.1.

Turn now to the second bullet in (4.9). Any given closed integral curve of $v$ in the $|\hat{u}| \leq R + \frac{1}{2} \ln z_*$ part of $\mathcal{H}^+_{p*}$ is determined by the change in $\phi$ between its endpoints, thus $\Delta\phi_{(\cdot)}$. This understood, use the fundamental theorem of calculus to see that what is written in the second bullet of (4.9) is no greater than

$$6 \sup_{\tau \in [0,1]} \int_{\Delta\phi_{\gamma_{p-}}}^{\Delta\phi_{\gamma_{p+}}} x \cos(\theta_{\gamma_\tau}(u)) \sin(\theta_{\gamma_\tau}(u)) \left| \frac{d\theta_{\gamma_\tau}(u)}{d(\Delta\phi_{\gamma_\tau})} \right| du$$

(4.11)

With regards to notation, the coordinate u is used as an affine parameter along the integral curves of $v$ and then what is denoted by $\theta_{(\cdot)}(\cdot)$ is the value of the coordinate $\theta$ at the indicated parameter value along the indicated integral curve.

To bound (4.11), note that that $\Delta\phi_{(\cdot)}$ for any given integral curve determines the angle, $\theta_+$, of its $\hat{u} = R + \frac{1}{2} \ln z_*$ endpoint, and vice versa. In particular, it follows from the two equations in the fourth bullet of Lemma II.2.2 that

$$\frac{d(\Delta\phi_{(\cdot)})}{d\theta_+} = 4\sqrt{6} \, z_* (1 - 3\cos^2\theta_+) \sin\theta_+ \int_{[-R-\ln\delta, R+\ln\delta]} \frac{x(u)}{f(u)^2} \frac{1 + 3\cos(\theta(u))}{(1 - 3\cos^2(\theta(u)))^3} du$$

(4.12)

Meanwhile, the bottom equation of the fourth bullet in Lemma II.2.2 implies that

$$\sin(\theta_{(\cdot)}(u)) \frac{d\theta_{(\cdot)}(u)}{d\theta_+} = 2z_* f(u)^{-1} (1 - 3\cos^2\theta_+) \sin\theta_+ \frac{1}{(1 - 3\cos^2(\theta_{(\cdot)}(u)))}$$

(4.13)

for the variation with $\theta_+$ of in the angle $\theta_{(\cdot)}$ at any given value of u a given integral curve segment of $v$. These last two equations imply that what is written in (4.11) is bounded by a purely S-dependent (or $\mathcal{K}$-compatible) multiple of $\delta_*^2$.

Minor cosmetic changes to the arguments from Step 4 of the proof of Proposition II.5.1 in Section II.5b give the bound on the integral of $ds \wedge \hat{a}$. Note in this regard that the integration by parts used in these arguments has no boundary contributions from the $|\hat{u}| = R + \frac{1}{2} \ln z_*$ boundary of C because $\hat{a}$ near the boundary is the 1-form $df$.

Step 3: This step states and then proves a refined bound on $|\varsigma|$ when $(\varphi, \varsigma)$ obey (3.4) with boundary values given by (4.3).

**Lemma 4.6**: *There exists a purely S-dependent (or K-compatible) constant $\kappa \geq 1$ such that if $\delta^2 < \kappa^{-1} z_*$, then the following is true: Let $(\varphi, \varsigma)$ denote a solution to (3.4) that obeys a given $\tau \in [0, 1]$ version of (4.3). Then $|\varsigma|$ is bounded by $\kappa z_*$ where $|\hat{u}| > R + \frac{1}{2} \ln z_* - 8$.*

It is important to keep in mind that the version of $\kappa$ from Lemma 4.6 does not depend on $\delta$ nor on $z_*$.

*Proof of Lemma 4.6*: Suppose that the lemma is false so as to derive some nonsense. Granted this assumption, there exists a sequence $\{D_n, (\hat{\Theta}_{n-}, \hat{\Theta}_{n+}), (\tau_n, (\varphi_n, \varsigma_n))\}_{n=1,2,...}$ of the following sort: Each index n version of $D_n$ is a data set of the form $(\delta_n, x_{0n}, R_n, J_n)$ suitable for defining the geometry of Y and (3.4), and such that $\delta_n < \frac{1}{n} z_*$. It is assumed here that all $n \in \{1, 2, ...,\}$ versions of $D_n$ use the same data from M; in particular, they define the same almost complex structure and pseudo-gradient vector field $\mathfrak{v}$ on the complement of the their attaching handles. The pair $(\hat{\Theta}_{n-}, \hat{\Theta}_{n+})$ is an $\{\mathfrak{m}_p = 0\}_{p \in \Lambda}$ element from $\hat{\mathcal{Z}}^S$. Meanwhile, the pair $(\varphi_n, \varsigma_n)$ obeys the $D_n$ version of (3.4) with the asymptotic/boundary conditions determined via (4.3) by the data $(\Theta_- = \Theta_{n-}, \Theta_+ = \Theta_{n+})$ and $\tau = \tau_n$. In addition, there are points where $|\hat{u}| > R + \frac{1}{2} \ln z_* - 8$ at which $|\varsigma_n| > 4 \kappa_* z_*$ with $\kappa_*$ here denoting the version of $\kappa$ given by Lemma 4.4. Note that any such point lies in the radius $e^{-4} z_*^{1/2}$ coordinate ball centered at one or the other of the critical points from $\mathfrak{p}$. By passing to a subsequence and renumbering, arrange that this occurs for each n in the radius $e^{-4} z_*^{1/2}$ ball about a fixed critical point from $\mathfrak{p}$. The argument that follows discusses the case when the critical point in question has index 1. But for some sign changes, the same argument works for the index 2 critical point. Let p denote the index 1 critical point in question.

Let $\{C_{p0n}\}_{n=1,2,...}$ denote the corresponding sequence of submanifolds. Note that various index n version different data sets to define the geometry of Y and almost complex structure on $\mathbb{R} \times Y$. Even so, the following is true: Fix an integer N, and then all $n > N$ versions of the almost complex structure agree on the part of the radius $\delta_*$ coordinate ball centered at p where the radius is greater than $(N^{-1} z_*)^{1/2}$. This almost complex structure is denoted by J. Let z denote the function $e^{-2(R-|\hat{u}|)}$ and let U denote the part of the radius $\delta_*$ coordinate ball centered on p where the radius is greater than $e^{-100} z_*^{1/2}$ and less than $z_*$. Granted what was just said, each $n \geq c_0$ versions of $C_{p0n}$ intersects $\mathbb{R} \times U$ as a properly embedded, J-holomorphic submanifold. Let $C_n$ denote this part of $C_{p0n}$.

For each n, let $s_n$ denote a value for *s* of a point in $C_n$ that corresponds to a point where $|\varsigma_n| > 4\kappa_* z_*$. Translate $C_n$ along the $\mathbb{R}$ factor of $\mathbb{R} \times U$ by $-s_n$ so that such a point in $C_n$ sits where $s = 0$ in the new submanifold. Let $C_n'$ denote this new submanifold.

Let $s_0 \in \mathbb{R}$. Then then integral of $ds \wedge \hat{a}$ over the $s \in [s_0, s_0+1]$ part of $C_n´$ is bounded by $z_*$ since $\hat{a} = df$ here and $C_n´$ has intersection number 1 or 0 with each constant $(s, f)$ level set with $f \in (1, 1+z_*)$. This understood, use Lemma 4.5 with Proposition II.5.5 to obtain a subsequence of $\{C_n´\}_{n=1,2,...}$ that converges on compact subsets in $\mathbb{R} \times U$ in the manner described by Proposition II.5.5. Let $\vartheta$ denote the resulting set of pairs consisting of an irreducible, J-holomorphic subvariety and positive integer weight.

As explained momentarily, the set $\vartheta$ must contain a pair whose subvariety component sits entirely in the $1 - 3\cos^2\theta = 0$ locus and is therefore the intersection of U with some element from Proposition II.3.3's moduli space $\mathcal{M}_1$. Granted for the moment that $\vartheta$ contains such a subvariety, it then follows from the manner of convergence described in Proposition 5.5 that all large n versions of $C_{p0n}$ must contain a loop that represents a non-zero mulitple of the generator of $H_1(\mathbb{R} \times \mathcal{H}^+_{p*}; \mathbb{Z})$. To elaborate, recall that such a generator can be taken to be any circle in U on which the coordinate s, the distance from the distance from p and the angle $\theta$ are constant. Such a circle is given by pushing a $1 - 3\cos^2\theta = 0$ circle where s and the distance from p are constant to the part of $\mathbb{R} \times \mathcal{H}^+_{p*}$ where $1 - 3\cos^2\theta$ is slightly positive. Meanwhile, a $1 - 3\cos^2\theta = 0$ circle of this sort is a constant radius slice of U's intersection with any submanifold from $\mathcal{M}_1$.

As just noted, if $\vartheta$ contains U's intersection with a submanifold from $\mathcal{M}_1$, then there is a circle in each large n version of $C_{p0n}$ that represents a non-zero multiple of the generator of the first homology of $\mathbb{R} \times \mathcal{H}^+_{p*}$. But this conclusion is nonsense by virtue of the fact that $C_{p0n}$ is diffeomorphic to $\mathbb{R} \times I_*$ and thus is contractible. This nonsense is what is required to prove the lemma.

What follows is the promised explanation for why $\vartheta$ contains U's intersection with a submanifold from $\mathcal{M}_1$. Lemma 4.4 is the key to the argument, for it implies that $|\varsigma_n| > \kappa_* z_*$ where $z = 2\delta_n^2 < \frac{2}{n} z_*$. Keeping this in mind, write z and h in terms of the variables $(r, \theta, \phi)$ where r is the distance to p. By way of reminder, the coordinate $z = r^2(1 - 3\cos^2\theta)$ and $h = r^2 \cos\theta \sin^2\theta$. Thus,

$$h/z = \frac{\cos\theta\sin^2\theta}{1- 3\cos^2\theta} \quad and \quad z^2 + 6\frac{h^2}{\sin^2\theta} = r^4(1+3\cos^4\theta) \,.$$

(4.14)

As a consequence, a $z < \frac{2}{n} z_*$ point where $|\varsigma_n| > \kappa_* z_*$ is a point where

$$1 - 3\cos^2\theta < \tfrac{1}{n}\tfrac{1}{3\sqrt{3}} \quad and \quad r > z_*^{1/2} \,.$$

(4.15)

This implies that $\vartheta$ contains a subvariety with a $1 - 3\cos^2\theta = 0$ point. Such a subvariety can not have points where $1 - 3\cos^2\theta < 0$ as there are no such points in $C_{p0n}$. Thus, it must sit entirely in the $1 - 3\cos^2\theta = 0$ locus and so constitute the intersection between U and a submanifold from $\mathcal{M}_1$.

<u>Step 4</u>: Assume that the parameters $(\delta, x_0, R)$ are such that Lemma 4.6 can be invoked. Let $\kappa_{**}$ denote the version of $\kappa$ from Lemma 4.6. Take $z_*$ so that $100\kappa_{**}z_*$ is less than $10^{-6}\delta_*^2$.

**Lemma 4.7**: *There exists $\kappa \geq 1$ with the following significance: Let $(\varphi, \varsigma)$ denote a solution to (3.4) that obeys a given version of (4.3). Then $1 - 3\cos^2\theta > \kappa^{-1}$ on the $\Psi_p$ image of the graph $(x, \hat{u}) \to (x, \hat{u}, \hat{\phi} = \varphi(x,\hat{u}), h = \varsigma(x,\hat{u}))$ in $\mathbb{R} \times \mathcal{X}$.*

The preceding lemma asserts the condition in the top bullet of (4.7).

***Proof of Lemma 4.7***: It follows from Lemma 4.6 and (4.14) that that $1 - 3\cos^2\theta > c_0^{-1}$ on the $\Psi_p$ image of the $|\hat{u}| = R + \frac{1}{2}\ln z_*$ boundary of the graph where $c_0 \geq 1$ is a purely S-dependent (or $\mathcal{K}$-compatible) constant. Thus, if $1 - 3\cos^2\theta \leq c_0^{-1}$ on the $\Psi_p$ image of the graph, then this must occur in the interior. This requires that $\cos\theta$ take its maximum on the interior. Lemma II.4.8 asserts that this maximum can occur only where $\hat{u} = 0$.

With the preceding understood, suppose that the lemma is false so as to generate some nonsense. Granted this assumption, there exists a sequence $\{(\tau_n, (\varphi_n, \varsigma_n))\}_{n=1,2,...}$ of the following sort: Each index n version of $(\varphi_n, \varsigma_n)$ is a solution to (3.4) with asymptotic/boundary conditions given by the $\tau = \tau_n$ version of (4.3). Furthermore, there is some $u = 0$ point on the $\Psi_p$-image of the graph of $(\varphi_n, \varsigma_n)$ where $1 - 3\cos^2\theta < \frac{1}{n}$. Given what is said in Lemma 4.5, an application of Proposition II.5.5 analogous to that used to prove Lemma 4.6 generates the same sort of conclusion: There is a loop in the $\Psi_p$ image of every large n graph that generates the first homology of $\mathbb{R} \times \mathcal{H}^+_{p*}$. As noted in the proof of Lemma 4.6, this is a nonsensical conclusion.

<u>Step 5</u>: The next lemma asserts the condition in the lower bullet of (4.7).

**Lemma 4.8**: *Given $\varepsilon > 0$, there exists $\kappa_\varepsilon > 1$ with the following significance: Suppose that $\tau \in [0, 1]$ and that $(\varphi, \varsigma)$ is a solution to (3.4) with boundary values given by the parameter $\tau$ version of (4.3). There are $\mathbb{R}$-valued lifts of $\varphi(x,\hat{u}), \phi(\gamma_{p-})$ and $\phi(\gamma_\tau)$ such that*
- $|\varphi(x, \hat{u}) - \phi(\gamma_{p-}|_{\hat{u}})| + |\varsigma(x,\hat{u}) - h(\gamma_{p-})| < \varepsilon$ *where* $x < -\kappa_\varepsilon$.
- $|\varphi(x, \hat{u}) - \phi(\gamma_\tau|_{\hat{u}})| + |\varsigma(x,\hat{u}) - h(\gamma_\tau)| < \varepsilon$ *where* $x > \kappa_\varepsilon$.

***Proof of Lemma 4.8***: Suppose to the contrary that no such $\kappa_\varepsilon$ exists so as to derive some nonsense. If this is the case, then there exists $\varepsilon_* > 0$ and a sequence $\{(\tau_n, (\varphi_n, \varsigma_n))\}_{n=1,2,...}$

with $\tau_n \in [0,1]$ and with $(\varphi_n, \varsigma_n)$ a solution to (3.4) with asymptotic/boundary conditions given by the $\tau_n$ version of (4.3). Moreover, $(\varphi_n, \varsigma_n)$ violates the $\varepsilon = \varepsilon_*$ conclusions of the lemma at all points where x < -n or x > n. This said, no generality is lost by assuming that the sequence violates the conclusions where x < -n.

Construct a new solution to (3.4) by translating $(\varphi_n, \varsigma_n)$ by the constant amount along the $\mathbb{R}$ factor of $\mathbb{R} \times I_*$ so that the resulting pair violates the $\varepsilon = \varepsilon_*$ version of what is asserted by top bullet of the lemma at some point where x = 0. Let $(\varphi_n', \varsigma_n')$ denote this new solution. The conditions given by the first bullet of (4.7) hold for the sequence $\{(\varphi_n', \varsigma_n')\}_{n=1,2,...}$. This being the case, then the standard elliptic regularity arguments (see again Chapter 6 of [M]) prove that there is a subsequence that converges in the $C^\infty$-Fréchet topology on compact subsets of $\mathbb{R} \times I_*$. Let $(\varphi, \varsigma)$ denote the limit. This pair obeys (3.4) and it obeys the $\tau = 0$ version of the condition given by the second bullet in (4.3). This is to say that the function $\varphi$ on the boundary of $\mathbb{R} \times I_*$ is independent of the $\mathbb{R}$-coordinate and its respective values on the two boundaries are those of $\phi(\gamma_p)$ on the relevant boundary of $I_*$.

The pair $(\varphi, \varsigma)$ also satifies the conditions given by the $\tau = 0$ version of the top bullet in (4.3). To see this, remark that the integral of $w$ over the $\Psi_p$-image of the graph given by this $(\varphi, \varsigma)$ version of (3.1) is finite. As the image of the graph is J-holomorphic, it follows using Lemma II.5.6 that any given sufficiently large, constant |s| slice of the $\Psi_p$-image of the graph must be everywhere very close to the $\hat{u} \in I_*$ segment of an integral curve of $v$ in $\mathcal{H}^+_{p*}$. Given the constant $|\hat{u}| = R + \frac{1}{2} \ln z_*$ value for $\varphi$, this segment must be from $\gamma_{p_-}$.

Note next that the translation that defined $\{(\varphi_n', \varsigma_n')\}_{n=1,2,...}$ guarantees that the solution $(\varphi, \varsigma)$ is *not* the solution to (3.4) and with boundary values the $\tau = 0$ version of (4.3) that is given by the x-independent map $\hat{u} \to (\phi(\gamma_{p_-}|_{\hat{u}}), \hat{h}(\gamma_{p_-}))$.

The conclusion of the previous paragraph is nonsensical given the assertion that there is at most one solution to any given version of (3.4) with a given $\tau \in [0, 1]$ asymptotic/boundary conditions from (4.3). This assertion is proved in the next subsection; it is the latter's Lemma 4.9.

**d) Uniqueness**

The next lemma completes the proof of Lemma 4.8. It also proves the uniqueness assertion of Proposition 2.2

**Lemma 4.9**: *Equation (3.4) has at most one solution whose boundary values are described by a given $\tau \in [0, 1]$ version of (4.3).*

***Proof of Lemma 4.9***: Suppose that $\tau \in [0, 1]$ and that $(\varphi^{(0)}, \varsigma^{(0)})$ and $(\varphi^{(1)}, \varsigma^{(1)})$ are two solutions to (3.4) with asymptotic/boundary conditions that are given by the parameter $\tau$ version of (4.3). Introduce $\varphi´$ to denote $\varphi^{(1)} - \varphi^{(0)}$. This function is zero on the boundary of $\mathbb{R} \times I_*$ and it has limit 0 as $|x| \to \infty$ on $\mathbb{R} \times I_*$. Let $\varsigma´ = \varsigma^{(1)} - \varsigma^{(0)}$. The pair $(\varphi´, \varsigma´)$ obeys an equation that can be written as $D_\mathfrak{h}(\varphi´, \varsigma´) = 0$ where $D_\mathfrak{h}$ is decribed by a version of (3.5) and (3.6). Indeed, just such an equation arises by subtracting the $(\varphi^{(0)}, \varsigma^{(0)})$ version of (3.4) from the corresponding $(\varphi^{(1)}, \varsigma^{(1)})$ where it is understood that $\mathbb{R}$-valued lifts of $\varphi^{(0)}$ and $\varphi^{(1)}$ are chosen so their boundary values agree. The functions $\mathfrak{a}_1$ and $\mathfrak{a}_2$ that appear in this version of (3.5) are the functions $(x, \hat{u}) \to a_1(x, \hat{u}, \varsigma^{(1)}|_{(x,\hat{u})})$ and $(x, \hat{u}) \to a_2(x, \hat{u}, \varsigma^{(1)}|_{(x,\hat{u})})$ that appear in the $(\varphi^{(1)}, \varsigma^{(1)})$ version of (3.4). Meanwhile, $\mathfrak{b}_1$ and $\mathfrak{b}_2$ are given by

- $\mathfrak{b}_1(x,\hat{u}) = [\int_0^1 a_{1\mathfrak{h}}(\cdot, \varsigma^{(1)} + s(\varsigma^{(0)} - \varsigma^{(1)})) ds] \partial_x \varphi^{(0)} + [\int_0^1 a_{1\mathfrak{h}}´(\cdot, \varsigma^{(1)} + s(\varsigma^{(0)} - \varsigma^{(1)})) ds] \partial_{\hat{u}} \varphi^{(0)}$

- $\mathfrak{b}_2(x,\hat{u}) = b(x,\hat{u}) + [\int_0^1 a_{2\mathfrak{h}}(\cdot, \varsigma^{(1)} + s(\varsigma^{(0)} - \varsigma^{(1)})) ds] \partial_x \varsigma^{(0)} + [\int_0^1 a_{1\mathfrak{h}}´(\cdot, \varsigma^{(1)} + s(\varsigma^{(0)} - \varsigma^{(1)})) ds] \partial_{\hat{u}} \varsigma^{(0)}$

(4.16)

Given that $(\varphi^{(1)}, \varsigma^{(1)})$ converges uniformly as $|x| \to \infty$, and given that this pair solves (3.4), standard elliptic regularity theorems as in Chapter 6 of [M] prove that the corresponding pair $(\partial_x \varphi^{(1)}, \partial_x \varsigma^{(1)})$ converges uniformly to zero as $|x| \to \infty$. This implies that the version of $(\mathfrak{a}_1, \mathfrak{a}_2, \mathfrak{b}_1, \mathfrak{b}_2)$ just defined obeys the condtions in (3.6). Thus, Proposition 3.1 can be invoked to see that the just defined version of $D_\mathfrak{h}$ has trivial kernel and so $(\varphi´, \varsigma´) = 0$.

## 5. Analytic background for the $\Delta_\mathfrak{p} > 0$ cases

This section prepares some analytic tools that are used in Section 6 to prove Proposition 2.2 when $\Delta_\mathfrak{p} > 0$. The analysis concerns two related issues that owe allegiance to Item b) in the third bullet of (2.9). This third bullet of (2.9) changes the domain of $(\varphi, \varsigma)$ so as to be the complement of either one or two $\hat{u} = 0$ points in $\mathbb{R} \times I_*$. The first issue is of import with regards to the behavior of the pairs $(\varphi, \varsigma)$ that arise in the $\Delta_\mathfrak{p} > 0$ versions of (2.10) near the missing $\hat{u} = 0$ points. Sections 5a and 5b are devoted to this topic. The second issue concerns the versions of (3.5)'s operator D that arise in the $\Delta_\mathfrak{p} > 0$ versions of (2.10). The domain and range spaces for the $\Delta_\mathfrak{p} > 0$ versions of D change to reflect the changed domain for the corresponding pair in (2.10) and the behavior of this pair near the missing $\hat{u} = 0$ point or points. The remaining subsections use what is said in Sections 5a and 5b to first define the new domain and range spaces for D, and then prove an analog of Proposition 3.1.

### a) J-holomorphic ends and the u = 0, 1 - 3cos²θ = 0 locus

This subsection describes the ends of J-holomorphic submanifolds whose constant $s$ slices converge as $s \to \infty$ in an isotopic fashion to one or the other of the curves $\hat{\gamma}_\mathfrak{p}^+$ and $\hat{\gamma}_\mathfrak{p}^-$. This is precisely the sort of end that appears $\Delta_\mathfrak{p} > 0$ cases of (2.9). What is said here concerns specifically the $\hat{\gamma}_p^+$ story as the story for the other curve can be obtained from the one told here by replacing $\theta$ with $\pi - \theta$ and changing some $\pm$ signs at various points. The story here is told in four parts.

*Part 1*: Introduce by way of notation $(s_+, \phi_+)$ as coordinates for $\mathbb{R} \times \mathbb{R}/(2\pi\mathbb{Z})$ so as to distinguish the latter from the eponymous factor in $\mathbb{R} \times \mathcal{H}_p$. A differential operator mapping $C^\infty(\mathbb{R} \times (\mathbb{R}/2\pi\mathbb{Z}); \mathbb{R}^2)$ to itself is defined by the rule that sends a pair (a, b) to the pair with respective first and second components

- $\partial_{s_+} a + (\sigma_0 \frac{3}{2x_0}) \partial_{\phi_+} b + (3\sigma_0 \frac{x_0 + 4e^{-2R}}{x_0^2}) a$.
- $\partial_{s_+} b - (\frac{3}{4x_0\sigma_0}) \partial_{\phi_+} a - (\frac{4}{x_0\sigma_0} e^{-2R}) b$.

(5.1)

Here, $\sigma_0$ is the value at $u = 0$, $\theta = \theta_0$ of the function $\sigma$ that is used for the fifth bullet of Part 1 in in Section 1c. This operator is denoted by $\mathfrak{D}_0$.

Pairs in the kernel of $\mathfrak{D}_0$ describe deformations of the J-holomorphic submanifold $\mathbb{R} \times \hat{\gamma}_p^+$ that are J-holomorphic to first order in the distance from $\mathbb{R} \times \hat{\gamma}_p^+$. In particular, pairs (a, b) with limit zero as $s_+ \to \infty$ describe the ends of J-holomorphic submanifolds whose constant $s$-slices converge as $s_+ \to \infty$ in an isotopic fashion to $\hat{\gamma}_p^+$. More is said about this in Parts 2-4. What follows directly talks about the kernel of $\mathfrak{D}_0$.

The operator depicted in (5.1) has constant coefficients, and so the kernel has a basis whose elements are irreducible representations of the $\mathbb{R}/(2\pi\mathbb{Z})$ action on the space of maps from $\mathbb{R} \times (\mathbb{R}/2\pi\mathbb{Z})$ to $\mathbb{R}^2$ generated by $\partial_{\phi_+}$. Using this Fourier mode decomposition makes an easy task of writing the kernel of (5.1). To say what this leads to, introduce

$$\lambda_1 = 3\sigma_0 \frac{x_0 + 4e^{-2R}}{x_0^2} \quad and \quad \lambda_2 = \frac{4}{x_0\sigma_0} e^{-2R}.$$

(5.2)

A basis for the kernel of $\mathfrak{D}_0$ is given by the $\phi_+$-independent elements

$$\mathfrak{y}_{0+} = (e^{-\lambda_1 s_+}, 0) \quad and \quad \mathfrak{y}_{0-} = (0, e^{\lambda_2 s_+}),$$

(5.3)

and then, for each $n \in \{1, 2, \ldots\}$, elements that have the form

$$\mathfrak{y}_{n+} = e^{-\lambda_{1n} s_+}(\cos n(\phi_+ - \phi_n), r_{1n} \sin n(\phi_+ - \phi_n)) \quad and \quad \mathfrak{y}_{n-} = e^{\lambda_{2n} s_+}(\cos n(\phi_+ - \phi_n), -r_{2n} \sin n(\phi_+ - \phi_n)),$$

(5.4)

where $\lambda_{1n} = \frac{1}{2} ((\lambda_1 + \lambda_2)^2 + \frac{n^2 9}{2x_0^2})^{1/2} + \lambda_1 - \lambda_2)$ and $\lambda_{2n} = \frac{1}{2} ((\lambda_1 + \lambda_2)^2 + \frac{n^2 9}{2x_0^2})^{1/2} + \lambda_2 - \lambda_1)$; and where $r_{1n}$ and $r_{2n}$ are certain specific positive constants. Meanwhile, $\phi_n \in \mathbb{R}/2\pi\mathbb{Z}$ can be any chosen angle.

Let $S_+$ denote the linear span of $\{\mathfrak{y}_{n+}\}_{n=0,1,...}$ versions of the left hand pair in (5.4). Let $S_-$ denote the linear span of $\{\mathfrak{y}_{n-}\}_{n=0,1,...}$ and those given by the various versions of the right hand pair. The elements from $S_+$ limit to zero as $s_+ \to \infty$ and those from $S_-$ limit to zero as $s_+ \to -\infty$.

*Part 2*: The upcoming description of the ends of J-holomorphic submanifolds invokes some geometric constructions that are described next. To start, remark that the restrictions of the coordinate functions $s$ and $\phi$ parametrized $\mathbb{R} \times \hat{\gamma}_p^+$. The resulting functions on $\mathbb{R} \times \hat{\gamma}_p^+$ are denoted by $(s_+, \phi_+)$. These coordinates with an auxilliary set of Euclidean coordinates $(\theta_+, u_+)$ for a small radius disk in $\mathbb{R}^2$ can be used as coordinates for an $\mathbb{R} \times (\mathbb{R}/2\pi\mathbb{Z})$-invariant, tubular neighborhood in $\mathbb{R} \times \mathcal{H}_p$ of $\mathbb{R} \times \hat{\gamma}_p^+$. This parametrization can be chosen so as to have the properties that are listed in the upcoming equation (5.5). The list uses $U_+ \subset \mathcal{H}^+_{p*}$ to denote the constant $s$ slices of this tubular neighborhood, this being an $\mathbb{R}/(2\pi\mathbb{Z})$-invariant tubular neighborhood of $\hat{\gamma}_p^+$. The list also refers to respective $\mathbb{R} \times \mathbb{R}/(2\pi\mathbb{Z})$ actions on the $(s_+, \phi_+, \theta_+, u_+)$ coordinate domain and on $\mathbb{R} \times U_+$. The action on the former are the constant translations of $s_+$ and $\phi_+$; and the action on the latter are the constant translations along the $\mathbb{R}$ factor and the constant translations of the coordinate $\phi$ for the $U_+$ factor. The final piece of new notation is the use of $\theta_*$ to denote the angle with $\cos\theta_* = \frac{1}{\sqrt{3}}$, this being the value of $\theta$ on $\hat{\gamma}_p^+$. Granted this notation, what follows lists the properties of the parametrization:

- *The constant $(s_+, \phi_+)$ disks are J-holomorphic.*
- *The parametrization has $\theta = \theta_* + \theta_+$ and $u = u_+$*
- *The parameterization is equivariant with respect to the respective $\mathbb{R} \times (\mathbb{R}/2\pi\mathbb{Z})$ actions.*
- *The coordinates $(s_+, \phi_+)$ equal $(s, \phi)$ on the $\theta_+ = 0, u_+ = 0$ cylinder.*

(5.5)

A parametrization of this sort can be constructed using Lemma 5.4 from [T1] with a little help from the inverse function theorem to arrange the condition in the second bullet.

*Part 3*: Granted these coordinates, a deformation of $\mathbb{R} \times \hat{\gamma}_p^+$ can be parametrized as a graph via functions $(a, b): \mathbb{R} \times \mathbb{R}/(2\pi\mathbb{Z}) \to \mathbb{R}^2$ as

$$(s_+, \phi_+) \to (s_+, \phi_+, \theta_+ = a(s_+,\phi_+), u_+ = b(s_+,\phi_+))$$

(5.6)

If $\eta = (a, b)$ is defined on a given open set in $\mathbb{R} \times S^1$ and if $|\eta| \leq c_0^{-1}$, then the resulting graph over the given open set defines a J-holomophic surface if and only if $\eta$ obeys a non-linear equation with the schematic form

$$\mathfrak{D}_0\eta + \mathfrak{r}_1 d\eta + \mathfrak{r}_0,$$

(5.7)

where $\mathfrak{r}_1$ is a smooth map from a certain small radius disk about the origin in $\mathbb{R}^2$ to $\text{Hom}(T^*\mathbb{R}^2, \mathbb{R}^2)$; and where $\mathfrak{r}_0$ is a smooth map from this same disk to $\mathbb{R}^2$. These are such that $|\mathfrak{r}_1| \leq c_0|\eta|$ and $|\mathfrak{r}_0| \leq c_0|\eta|^2$.

By way of an example, the J-holomorphic cylinders that comprise Proposition II.3.4's moduli space $\mathcal{M}_{\mathfrak{p}0}$ foliate the $u = 0$ slice of $\mathbb{R} \times \mathcal{H}^+_{\mathfrak{p}0}$. These cylinders are $\phi$-invariant. Each such cylinder has two ends; their constant $s$ slices converge isotopically as $s \to \infty$ to the respective integral curves $\hat{\gamma}^+_\mathfrak{p}$ and $\hat{\gamma}^-_\mathfrak{p}$. The very large $s$ parts of the end whose slices converge to $\hat{\gamma}^+_\mathfrak{p}$ appears as a $\phi_+$-independent solutions to (5.7) that are defined for $s_+ \gg 1$ with pairs $(a, b)$ such that $a > 0$ and $b = 0$. In particular, integrating the $u = 0$ version of the vector field in Equation (II.3.10), or using arguments much like those in Section 2 of [HT] finds that the relevant version of $\eta$ can be written as

$$\eta = \alpha(e^{-\lambda_1 s_+} + e_1, 0) \quad \text{where } \alpha \in (0,\infty) \text{ and } |e_1| \leq c_0|\alpha|e^{-(\lambda_1 + 1/c_0)s_+}.$$

(5.8)

Meanwhile, the large $s$ part of the end of any given submanifold from Proposition II.3.4's moduli space $\mathcal{M}_{\mathfrak{p}+}$ is described by (5.8) with $\alpha < 0$.

By way of a second example, the end in $\mathbb{R} \times \mathcal{H}^+_\mathfrak{p}$ from Proposition II.3.3's moduli spaces $\mathcal{M}_1$ and $\mathcal{M}_2$ whose constant $s$ slices converge as $s \to \infty$ to $\hat{\gamma}^+_\mathfrak{p}$ are described where $s \ll -1$ by a $\phi_+$-invariant solution to (5.7) that is defined where $s_+ \ll -1$ and has the form

$$\eta = \alpha(0, e^{\lambda_2 s_+} + e_2) \quad \text{where } \alpha \in \mathbb{R}-0 \text{ and } |e_2| \leq c_0|\alpha|^2 e^{-(\lambda_2 + 1/c_0)|s_+|}.$$

(5.9)

The $\alpha > 0$ cases describe the end of the submanifolds from $\mathcal{M}_1$ and the $\alpha < 0$ cases describe the end of the submanifolds from $\mathcal{M}_2$.

A third example involves the submanifolds from Proposition II.3.2's moduli space $\mathcal{M}_\Sigma$. Those parameterized as in the second bullet of Proposition II.3.2 by a pair $(x, y)$ with y near 1 can be written using (5.5)-(5.7) using

$$\mathfrak{y} = (\alpha_+ e^{-\lambda_1 s_+} + e_+, \alpha_- e^{\lambda_2 s_+} + e_-)$$

(5.10)

where $\alpha_+ > 0$ and $\alpha_- > 0$. Here, $e_+$ and $e_-$ are both $\phi_+$-invariant. In addition, there norms are such that $|e_+| \leq c_0 |\alpha_+|(|\alpha_+| + |\alpha_-|) e^{-(\lambda_1 + 1/c_0)s_+}$ and $|e_-| \leq c_0 |\alpha_-|(|\alpha_+| + |\alpha_-|) e^{-(\lambda_2 + 1/c_0)s_+}$. Note that this representation is valid only over a domain $I \times (\mathbb{R}/2\pi\mathbb{Z}) \subset \mathbb{R} \times \mathbb{R}/(2\pi\mathbb{Z})$ where I is a bounded interval whose endpoints are determined by $\alpha_+$ and $\alpha_-$. The left endpoint diverges as $\alpha_- \to 0$ and the right endpoint diverges as $\alpha_+ \to 0$. Meanwhile, a surface from $\mathcal{M}_\Sigma$ parametrized by $(x, y)$ with $y \sim 2$ appears as in (5.10) but with $\alpha_- < 0$.

*Part 4*: This part of the subsection directly addresses the issue of describing ends of J-holomorphic submanifolds using the kernel of $\mathfrak{D}_0$. As just noted, any such end whose large $s \gg 1$ slices sit in $U_+$ and converge to $\hat{\gamma}_p^+$ in an isotopic fashion as $s \to \infty$ is described by a solution to (5.7) that is defined where $s_+ \gg 1$ and has $s_+ \to \infty$ limit equal to zero. By the same token, the any such end whose $s \ll -1$ slices sit in $U_+$ and converge to $\hat{\gamma}_p^+$ in an isotopic fashion is described by a solution to (5.7) defined where $s_+ \ll -1$ and converging to 0 as $s_+ \to -\infty$. The following lemma describes all such solutions to (5.7).

**Proposition 5.1**: *There exists $\kappa \geq 1$ with the following significance: Fix $s_* \geq 1$.*
- *Suppose that $\mathfrak{y}$ is a solution to (5.7) with domain $[s_*, \infty) \times \mathbb{R}/(2\pi\mathbb{Z})$ that converges to 0 as $s_+ \to \infty$ and has pointwise norm bounded by $\kappa^{-2}$. There exists $n \geq 1$ such that $\mathfrak{y}$ can be written as*

$$\mathfrak{y} = c_0(e^{-\lambda_1 s_+} + e_1, 0) + c_n \mathfrak{y}_{n+} + \mathfrak{e}_n$$

*with $c_0, c_n \in (-\kappa^{-1}, \kappa^{-1})$, with $e_1$ given by (5.8), and with $|\mathfrak{e}_n| \leq \kappa |c_n|(|c_n| + |c_0|) e^{-(\lambda_{1n} + 1/\kappa)s_+}$. Conversely, given $n \in \{1, 2, ...\}$ and constants $c_0, c_n \in (-\kappa^{-1}, \kappa^{-1})$, there exists a solution to (5.7) that can be written in this way.*
- *Suppose that $\mathfrak{y}$ is a solution to (5.7) with domain $(-\infty, -s_*] \times \mathbb{R}/(2\pi\mathbb{Z})$ that converges to 0 as $s_+ \to -\infty$ and has pointwise norm bounded by $\kappa^{-2}$. There exists $n \geq 1$ such that $\mathfrak{y}$ can be written as*

$$\mathfrak{y} = c_0(0, e^{\lambda_2 s_+} + e_2) + c_n \mathfrak{y}_{n-} + \mathfrak{e}_n$$

*with $c_0, c_n \in (-\kappa^{-1}, \kappa^{-1})$, with $e_2$ given by (5.9), and with $|\mathfrak{e}_n| \leq \kappa |c_n|(|c_n| + |c_0|) e^{-(\lambda_{2n} + 1/\kappa)|s_+|}$. Conversely, given $n \in \{1, 2, ...\}$ and constants $c_0, c_n \in (-\kappa^{-1}, \kappa^{-1})$, there exists a solution to (5.7) that can be written in this way.*

*Moreover, in either case, the derivatives of $\mathfrak{y}$ to any given order are square integrable where $|s_+| > 2s_*$.*

***Proof of Proposition 5.1***: The analysis from Section 2 and specifically Section 2.3 of [HT] can be used but for one added comment to prove that any given solution to (5.7) with $s_+ \to \infty$ limit zero can be written as described. The extra comment concerns the derivation of the distinct of bounds on the norms for $\mathfrak{e}_0$ and $\mathfrak{e}_1$. These bounds are obtained by projection $\mathfrak{y}$ and the expression in (5.7) onto $\mathbb{R}/(2\pi\mathbb{Z})$–invariant subspace of maps from $\mathbb{R} \times \mathbb{R}/(2\pi\mathbb{Z})$ to $\mathbb{R}^2$. To elaborate, this projection is given by the map

$$q \to \Pi q = \frac{1}{2\pi} \int_{\mathbb{R}/(2\pi\mathbb{Z})} q(\cdot, \phi_+) d\phi_+ \;.$$

(5.11)

The use of such a projection is not discussed in Section 2 of [HT]. Even so, the latter's arguments can be applied separately to the $\mathbb{R}/(2\pi\mathbb{Z})$ invariant part of (5.7) and the remainder with what are little more than notational changes to obtained the distinct bounds for the norms of $\mathfrak{e}_0$ and on $\mathfrak{e}_1$.

The proof of the converse assertion in the first bullet is given below in two steps. The proof of the converse assertion in the second bullet is identical but for straightforward notational and cosmetic changes and so is not given. The proof that the derivatives to any given order are square integrable invokes standard elliptic regularity theorems of the sort that can be found in Chapter 6 of [M].

Step 1: Use $\mathbb{H}$ to denote now the $L^2_1$ completion of the space of smooth, $\mathbb{R}^2$-valued functions on $[0, \infty) \times \mathbb{R}/(2\pi\mathbb{Z})$ with compact support and which lie in Part 1's subspace $\mathcal{S}_-$ on the boundary, $\{0\} \times \mathbb{R}/(2\pi\mathbb{Z})$. Let $\mathbb{L}$ denote the $L^2$ completion of the space of smooth, $\mathbb{R}^2$-valued functions on $[0, \infty) \times \mathbb{R}/(2\pi\mathbb{Z})$. It is straightforward task using integration by parts to prove that $\mathfrak{D}_0$ defines a Fredholm operator from $\mathbb{H}$ to $\mathbb{L}$ with trivial kernel and cokernel. As such, it has a bounded inverse. There is also a version here of the Hilbert spaces $\mathbb{H}_*$ and $\mathbb{L}_*$ that are defined in Section 3e. These are defined by completing the respective dense domains for $\mathbb{H}$ and $\mathbb{L}$ using for $\mathbb{H}_*$ the $[0, \infty) \times \mathbb{R}/(2\pi\mathbb{Z})$ analog of (3.33), and using for $\mathbb{L}_*$ the analog that integrates the square of the norm of $\eta´$ rather than that of its derivatives. The operator $\mathfrak{D}_0$ also defines a bounded, linear map from $\mathbb{H}_*$ to $\mathbb{L}_*$. The analog of Lemmas 3.3 and 3.4 holds in this case: The inverse of $\mathfrak{D}_0$ maps $\mathbb{L}_* \subset \mathbb{L}$ to $\mathbb{H}_* \subset \mathbb{H}$ as a bounded operator. As $\mathfrak{D}_0$ commutes with $\Pi$, the inverse also commutes with $\Pi$.

Step 2: Fix $T \geq 100$ and let $\beta_T$ denote the function on $[0, \infty)$ given by $\beta(s-T)$. Fix $n \in \{1, 2, \ldots\}$. Given the $\mathbb{R} \times \mathbb{R}/(2\pi\mathbb{Z})$ version of Lemma 3.3, there exists $c_0 > 1$ and $T \geq 1$ such that following is true: Suppose that $|c_0| + |c_n| \leq c_0^{-2}$. Let $\mathbb{B} \subset \mathbb{H}_*$ denote the ball about the origin of radius $c_0^{-1}$. Reintroduce $\mathfrak{r}_0$ and $\mathfrak{r}_1$ from (5.7). Let $\mathfrak{y}_+ = (e^{-\lambda_1 s_+} + e_1, 0)$. A smooth map from $(\times_2(-c_0^{-1}, c_0^{-1})) \times \mathbb{B}$ to $\mathbb{L}$ is defined so as to send any given $((c_0, c_n), \mathfrak{q})$

$$\mathfrak{D}_0 \mathfrak{q} + \beta_T \mathfrak{r}_1 (c_0 \mathfrak{y}_+ + c_n \mathfrak{y}_{n+} + \mathfrak{q}) d(c_0 \mathfrak{y}_+ + c_n \mathfrak{y}_{n+} + \mathfrak{q}) + \beta_T \mathfrak{r}_0 (c_0 \mathfrak{y}_+ + c_n \mathfrak{y}_{n+} + \mathfrak{q}) \, . \tag{5.12}$$

The differential of this map at $((c_0, 0) \, 0)$ along the $\mathbb{H}_*$ factor is an isomorphism if $T \geq c_0^{-1}$ This being the case, then the inverse function theorem finds $c_0$ and for $T > c_0$, a smooth map $\mathfrak{e}_n \colon \times_2(-c_0^{-1}, c_0^{-1}) \to \mathbb{B}$ such that the triple $((c_0, c_n), \mathfrak{q} = \mathfrak{e}_n(c_0, c_n))$ is mapped to zero by (5.12). Moreover, this element $\mathfrak{e}_n$ is such that $|\mathfrak{e}_n| \leq c_0 |c_n| (|c_0| + |c_n|) e^{-\lambda_1 T}$. The techniques from Section 2.3 in [HT] can be used to see that $\mathfrak{e}_n$ has the asserted norm bound.

**b) The kernel of $\mathfrak{D}_0$ and graphs over $\mathbb{R} \times \mathbb{I}_*$**

Some of the J-holomorphic cylinders given by Proposition 5.1 via (5.5)–(5.7) will intersect $\mathbb{R} \times \mathcal{H}^+_{p*}$ and so intersect the image of $\Psi_p$. This subsection says something about the $\Psi_p$-inverse image of these intersections. Of particular interest are the cylinders where $s$ is unbounded from above. There are five parts to what follows. Lemma 5.6 in Part 5 gives some indication as to why these cylinders are relevant.

*Part 1*: Fix $n \geq 1$ and a pair $(c_0, c_n) \in \mathbb{R}^2-0$ whose absolute value is small enough to apply the first bullet of Proposition 5.1 to obtain a corresponding solution, $\mathfrak{y}$, to (5.7). Write $\mathfrak{y}$ as $(a, b)$ and use the latter in (5.5) and (5.6) to define a J-holomorphic cylinder in the $s \gg 1$ part of $\mathbb{R} \times U_+ \subset \mathbb{R} \times \mathcal{H}_p$. It follows from the second bullet of (5.5) that $c_0$ must be positive for this to occur.

Assume henceforth that $c_0 > 0$ and that $c_n \neq 0$. If this is so, then the large $s$ part of the cylinder in question has algebraic intersection n with the large $s$ parts of submanifolds from $\mathcal{M}_{p0}$ and from submanifolds from Proposition II.3.2's moduli space $\mathcal{M}_\Sigma$ that come very near the $u = 0$ locus in $\mathbb{R} \times \mathcal{H}_p$ at large $s$. Indeed, this last point is a direct consequence of three facts: First, each such subvariety from $\mathcal{M}_{p0}$ appears as in (5.8), and those from $\mathcal{M}_\Sigma$ appear as in (5.10). Second, the $c_0 \mathfrak{y}_{0+} + \mathfrak{e}_0$ contribution to $\mathfrak{y}$ has the form $(c_0 e^{-\lambda_1 s_+} + e_0, 0)$ with $e_0$ being $\phi$-invariant and having the asserted norm bound. Third, the $u_+$ component of $\mathfrak{y}_{1n}$ has $2n$ zeros on each large, constant $s_+$ circle in $\mathbb{R} \times \mathbb{R}/(2\pi\mathbb{Z})$ and each such zero is transverse.

*Part 2*: The large $s$ part of a cylinder in $\mathbb{R} \times \mathcal{H}^+_\mathfrak{p}$ parametrized via (5.5)-(5.7) and the first bullet of Proposition 5.1 with $c_0 > 0$ has $\mathfrak{y} = c_0(e^{-\lambda_1 s_+}, 0) + \mathcal{O}(e^{-(\lambda_1 + 1/c_0)s_+})$ and so looks to leading order like what is written in (5.8). By way of a reminder, the latter depicts a cylinder from Proposition II.3.2's moduli space $\mathcal{M}_{\mathfrak{p}0}$. This being the case, what follows says more about the $\mathcal{M}_{\mathfrak{p}0}$ case of (5.8) in preparation for what is said in the next parts of the subsection about the $c_n \neq 0$ cases.

To set the stage, keep in mind that a cylinder from Proposition II.3.3's moduli space $\mathcal{M}_{\mathfrak{p}0}$ is $\phi$-invariant and invariant with respect to the involution $\theta \to \pi - \theta$. Such a cylinder has two ends; and the constant $s$ slices of these ends converge in an isotopic fashion as $s \to \infty$ to the respective integral curves $\hat{\gamma}^+_\mathfrak{p}$ and $\hat{\gamma}^-_\mathfrak{p}$. The $\Psi_\mathfrak{p}$-inverse images of these cylinders from $\mathcal{M}_{\mathfrak{p}0}$ are the constant x slices of the $\hat{u} = 0$ locus in $\mathbb{R} \times \mathcal{X}$. The association of the value of x to the corresponding cylinder gives an $\mathbb{R}$-equivariant diffeomorphism between $\mathbb{R}$ and $\mathcal{M}_{\mathfrak{p}0}$. This diffeomorphism from $\mathbb{R}$ to $\mathcal{M}_{\mathfrak{p}0}$ sends any given $y \in \mathbb{R}$ to the cylinder in $\mathcal{M}_{\mathfrak{p}0}$ whose $s = y$ slice is the $(u = 0, \theta = 0)$ circle in $\mathcal{H}^+_\mathfrak{p}$.

Fix $y \in \mathbb{R}$ and let $\Sigma_y \in \mathcal{M}_{\mathfrak{p}0}$ denote the corresponding cylinder. The function $s$ on $\Sigma_y$ has one critical value, this the $s = y$ locus. It restricts to the both components of the complement of this locus as a proper map to $(y, \infty)$. The function $\cos\theta$ increases monotonically as a function of $s$ with $s \to \infty$ limit $\frac{1}{\sqrt{3}}$ on one of these components. Meanwhile, $\cos\theta$ decreases monotonically on the other component with $s \to \infty$ limit equal to $-\frac{1}{\sqrt{3}}$. Let $\mathcal{E}_{0,y} \subset \Sigma_y$ denote the former component, this being the end whose constant $s$ slices converge to $\hat{\gamma}^+_\mathfrak{p}$ and is given via (5.5)-(5.7) by using $\mathfrak{y}$ as depicted in (5.8) for a suitable choice of $\alpha$. Denote the $\mathcal{E}_{0,y}$ version of $\alpha$ by $\alpha_y$.

**Lemma 5.2**: *Fix a pair* $x, y \in \mathbb{R}$ *and let* $\Sigma_x$ *and* $\Sigma_y$ *denote the corresponding surfaces from* $\mathcal{M}_{\mathfrak{p}0}$. *Then* $x - y = \frac{1}{\lambda_1} \ln(\frac{\alpha_x}{\alpha_y})$.

*Proof of Lemma 5.2*: This follows from (5.5)-(5.8) given the fact that $\Sigma_x$ is obtained from $\Sigma_y$ by translating the latter by $x - y$ along the $\mathbb{R}$ factor of $\mathbb{R} \times \mathcal{H}^+_\mathfrak{p}$.

*Part 3*: Fix $n \in \{1, 2, \ldots\}$ and $y \in \mathbb{R}$. Introduce by way of notation $\mathcal{E}_{n,y}$ to denote the large $s$ part of a cylinder that is described via (5.5)-(5.7) and the first bullet of Proposition 5.1 with $c_0 = \alpha_y$ and $c_n \neq 0$. The next lemma describes the $\Psi_p$-inverse image of these sorts of cylinders.

**Lemma 5.3**: *There exists a constant $\kappa \geq 1$ with the following significance: Fix a positive integer* n *and* $\phi_n \in \mathbb{R}/(2\pi\mathbb{Z})$ *so as to specify a particular version of* $\mathfrak{y}_{n+}$ *from (5.4). Choose a real number* y *and set* $c_0 = \alpha_y$; *then choose* $c_n \neq 0$ *so as to define* $\mathcal{E}_{n,y} \subset \mathbb{R} \times \mathcal{H}^+_p$ *via (5.5)-(5.7) using* $\mathfrak{y}$ *as in Proposition (5.1). Fix* $(s_+, \phi_+) \in \mathbb{R} \times S^1$ *with* $s_+ \gg 1$.

- *The $\hat{\phi}$ and h coordinates of the $\Psi_p^{-1}$ image of the corresponding point in $\mathcal{E}_{n,y}$ are*

$$\hat{\phi} = (\phi_+ - \phi_n) \text{ and } h = (x_0 + 4e^{-2R}) \tfrac{2}{3\sqrt{3}} (1 - 3\alpha_y^2 e^{-2\lambda_1 s_+} - 6\alpha_y^2 c_n e^{-(\lambda_1 + \lambda_{1n})s_+} \cos(n(\phi_+ - \phi_n)) + \cdots$$

*where the unwritten term has two parts. The $\phi_+$-invariant part is bounded in absolute value by $e^{-(2\lambda_1 + 1/\kappa)s_+}$. The remainder is bounded in absolute value by $e^{-(\lambda_1 + \lambda_{1n} + 1/\kappa)s_+}$.*

- *The $\hat{u}$ coordinate is*

$$\hat{u}(s_+, \phi_+) = \tfrac{3}{2x_0\sigma_0} (((\lambda_1 + \lambda_2)^2 + \tfrac{n^2 9}{2x_0^2})^{1/2} + \lambda_1 + \lambda_2)(1 + \cdots) c_n e^{-\lambda_{1n}' s_+} \sin(n(\phi_+ - \phi_n)) + \cdots$$

*where $\lambda_{1n}' = \lambda_{1n} + 12e^{-2R}(\sigma_0 x_0)^{-1}$ and where the first unwritten factor is bounded in absolute value by $\kappa e^{-2R}$ and the second by $e^{-(\lambda_{1n}' + 1/\kappa)s_+}$. The absolute values of their derivatives are also bounded by these same respective factors.*

- *The x coordinate is*

$$x(s_+, \phi_+) = y + \lambda_1^{-1} c_n e^{-(\lambda_{1n} - \lambda_1)s_+} \cos(n(\phi_+ - \phi_n)) + \cdots$$

*where the unwritten term is bounded in absolute value by $e^{-(\lambda_{1n} - \lambda_1 + 1/\kappa)s_+}$. The absolute value of its derivatives is also bounded by this same factor.*

A proof is given momentarily. What follows directly is a corollary of what is said by the second and third bullets of Lemma 5.3.

**Corollary 5.4**: *Fix* $s \gg 1$ *and there is an open, contractible neighborhood* $V_s \subset \mathbb{R} \times I_*$ *of the point* $(y, 0)$ *with the following significance: The projection to $\mathbb{R} \times I_*$ of the $\Psi_p$-inverse image of where $s > s$ in $\mathcal{E}_{n,y}$ defines a proper, n to 1 covering map onto $V_s - (y, 0)$.*

The rest of this part of the subsection is occupied with the

*Proof of Lemma 5.3*: The claim in the first bullet follows from what is said in Proposition 5.1 and the fact that $\hat{\phi}$ and h are the respective pull-backs of $\phi$ and $f(u)\cos\theta \sin^2\theta$. Use the identification $\theta = \theta_* + a(s_+, \phi_+)$ and $u = b(s_+, \phi_+)$ in the latter function with Taylor's theorem to obtain the given expression for h.

The formulae for û and x are derived in the four steps that follow. The arguments given take $\phi_n = 0$. The assertion in the general case follows directly from this case by applying a constant $\mathbb{R}/(2\pi\mathbb{Z})$ translation.

<u>Step 1</u>: The $u = 0$ slice $\mathcal{E}_{n,y}$ is parameterized by $s_+$ via the rule

$$\theta(s_+) = \theta_* + \alpha_y (e^{-\lambda_1 s_+} + e_0 \pm c_n e^{-\lambda_{1n} s_+}) + \mathfrak{e}_n ,$$

(5.13)

with $e_0$ and $\mathfrak{e}_n$ given by Proposition 5.1. The plus sign occurs at an angle $\phi_+ = 0 + \mathfrak{e}_+$ and the minus sign occurs at $\phi = \pi + \mathfrak{e}_-$ where $|\mathfrak{e}_+|$ and $|\mathfrak{e}_-|$ are both bounded by $c_0 e^{-(\lambda_{1n} + 1/c_0) s_+}$. From the vantage of $\mathbb{R} \times \mathcal{X}$, the $\phi \sim 0$ intersection locus correspond to points with

$$\hat{u} = 0, \quad \hat{\phi} = 0 + \mathfrak{e}_+ \quad and \quad h = (x_0 + 4 e^{-2R}) \tfrac{2}{3\sqrt{3}} (1 - 3\alpha_y^2 e^{-2\lambda_1 s_+}) + \cdots ,$$

(5.14)

where the unwritten term in the expression for h is bounded in absolute value by $e^{-(2\lambda_1 + 1/c_0) s_+}$. The $u = 0$ locus in $\mathcal{E}_{n,y}$ with $\phi \sim \pi$ has $\hat{\phi}$ coordinate $\pi + \mathfrak{e}_-$ and h coordinate also given by (5.14).

<u>Step 2</u>: It follows from the definition given in (1.29) that the x and û coordinates of the point in $\mathcal{E}_{n,y}$ can be determined from (5.5)–(5.7) by integrating the vector field given in Equation (II.3.10) starting at the point $(s = s_+, \theta = \theta_* + a(s_+, \phi_+), b(s_+, \phi_+))$. The values of x and û of this point on $\mathcal{E}_{n,y}$ are the respective s-coordinates and u coordinates of the point on the relevant integral curve where $\theta = \tfrac{\pi}{2}$.

To see what results, let $\tau \to (s(\tau), \theta(\tau), u(\tau))$ denote for the moment a certain parametrization of this integral curve. Take $\tau = 0$ to be the starting point. As $\theta$ increases along the curve, Equation (II.3.9) implies that $|u|$ decreases from its initially small value as $\tau$ increases. Now, define the parametrization of the curve by $\tau$ so that Taylor's theorem applied to (II.3.9) writes the $\tau$-derivative of u as

$$\tfrac{du}{d\tau} = -2\sqrt{3} e^{-2R} u (1 + \cdots) \cos\theta \sin\theta ,$$

(5.15)

where the unwritten term has absolute value bounded by $c_0 e^{-2\lambda_{1n} s}$. A second application of Taylor's theorem writes

$$\tfrac{d\theta}{d\tau} = \tfrac{1}{2\sqrt{2}} (x_0 + 4 e^{-2R} + \cdots)(1 - 3\cos^2\theta)$$

(5.16)

where the unwritten terms are also bounded by $c_0 e^{-2\lambda_1 s}$. Given the very small $\tau = 0$ value for $\theta - \theta_*$ at the start, it follows from (5.16) that the value of $\tau$ where $\theta = \frac{\pi}{2}$ is given by

$$\tau = \frac{3}{x_0} \sigma_0 s_+ + \cdots$$

(5.17)

where the unwritten term has absolute value bounded by $c_0$. Granted this, use Taylor's theorem to approximate $\cos\theta \sin\theta$ in (5.15) by $\cos\theta_* \sin\theta_{0\backslash *} = \frac{2}{\sqrt{3}}$. Integration produces the formula given for $\hat{u}$ in the second bullet of the lemma.

Step 3: To get an expression for the x-coordinate of a given $u = 0$ point on $\mathcal{E}_{n,y}$, note that the value of $\theta$ on the end $\mathcal{E}_y$ in $\mathcal{M}_{p0}$'s cylinder $\Sigma_y$ is described at large $s_+$ by the $c_n = 0$ version of (5.13). For any given $x \in \mathbb{R}$, use Lemma 5.2 to see that the value of $\theta$ on the end $\mathcal{E}_x$ in the corresponding $\Sigma_x$ is described at large $s_+$ by the version of (5.13) that sets $c_n = 0$ and replaces $s_+$ by $s_+ - (x - y)$. Granted this last observation, use a first order Taylor's approximation to see that the value of $x$ on $\mathcal{E}_{n,y}$ at a given very large $s_+$ and where $u = 0$ is obtained by solving

$$\alpha_y (e^{-\lambda_1 s_+} + e_1)(1 + \lambda_1 (x - y) \cdots) = \alpha_y ((e^{-\lambda_1 s_+} + e_1) \pm c_n e^{-\lambda_{1n} s_+}) + \cdots$$

(5.18)

where the unwritten term on the left hand side involve higher powers of $(x - y)$ and a term with absolute value bounded by $e^{-s_+/c_0} |x - y|$. Meanwhile, the unwritten term on the right hand side has absolute value bounded by $e^{-(\lambda_{1n} + 1/c_0)s_+}$. This last equation implies that the x coordinate of a given $s_+ \gg 1$ point on the $u = 0$ locus in $\mathcal{E}_{n,y}$ is given

$$x(s_+) - y = \pm \lambda_1^{-1} c_n e^{-(\lambda_{1n} - \lambda_1)s_+} + \cdots$$

(5.19)

where the unwritten term has absolute value bounded by $e^{-(\lambda - \lambda_1 + 1/c_0)s}$.

Step 4: Granted that $|u|$ decreases from its initially small value, it also follows from (II.3.9) that the value of $x$ is very nearly the $\mathbb{R}$-parameter of the $\theta = \frac{\pi}{2}$ point on the unique $\Sigma_{(\cdot)}$ surface that contains $(s_+, \theta = \theta_* + a(s_+, \phi_+))$. Given this observation, what is said the preceding steps imply directly the formula for $x$ in the third bullet of the lemma.

*Part 4*: This part of the subsection concerns specifically the case where $n = 1$. The discussion here concerns the normal bundle to the large $s$ part of the surface $\mathcal{E}_{1,y}$ when viewed using (5.5)-(5.7) and Proposition 5.1, and when viewed via $\Psi_{\mathfrak{p}}$ as a submanifold in $\mathbb{R} \times \mathcal{X}$.

To start, use the almost complex structure J and the 2-form $\hat{\omega} = ds \wedge \hat{a} + w$ to define the Riemannian metric $\hat{\omega}(\cdot, J(\cdot))$. Let $N \to \mathcal{E}_{1,y}$ denote the normal bundle to the submanifold $\mathcal{E}_{1,y}$, this being the orthogonal complement in $T(\mathbb{R} \times \mathcal{H}_{\mathfrak{p}})$ of $T\mathcal{E}_{1,y}$. View the large $s$ part of $\mathcal{E}_{1,y}$ using (5.5)-(5.7) and Proposition 5.1 to see that pairing with the 1-forms $(d\theta_+, du_+)$ define an isomorphism between N and the product $\mathbb{R}^2$ bundle. Meanwhile, this same part of $\mathcal{E}_{1,y}$ can be viewed as the $\Psi_{\mathfrak{p}}$ image of a surface in $\mathbb{R} \times \mathcal{X}$; and Corollary 5.4 implies that the 1-forms $(d\hat{\phi}, dh)$ also define an isomorphism between N and the product $\mathbb{R}^2$ bundle.

These two product structures are related in the following way: Let $\mathfrak{y}$ denote a map from the $s_+ \gg 1$ part of $\mathbb{R} \times S^1$ to $\mathbb{R}^2$. Use the product structure defined by $(d\theta_+, du_+)$ to view $\mathfrak{y}$ as a section of N over this part of $\mathcal{E}_{1,y}$. Meanwhile, use the n = 1 version of Corollary 5.4 to view this part of $\mathcal{E}_{1,y}$ as the $\Psi_{\mathfrak{p}}$-image of a graph of the sort depicted in (3.1) with the domain of the relevant version of the pair $(\varphi, \varsigma)$ being the complement of $(y, 0)$ in an $\mathbb{R} \times I_*$ neighborhood of $(y, 0)$. With $\mathcal{E}_{1,y}$ viewed this way, then the image via $(d\hat{\phi}, dh)$ of the section defined by $\mathfrak{y}$ defines a map, $\eta$, from the domain of $(\varphi, \varsigma)$ to $\mathbb{R}^2$. The maps $\mathfrak{y}$ and $\eta$ are related via a rule given by

$$\mathfrak{y}|_{(s_+, \phi_+)} = \mathrm{u} \cdot (\eta|_{(x(s_+, \phi_+), \hat{u}(s_+, \phi_+))})$$

(5.20)

where u is a smooth map from the large $s_+$ part of $\mathbb{R} \times S^1$ to $Gl(2; \mathbb{R})$ with positive determinant. Note that the latter component of $Gl(2; \mathbb{R})$ deformation retracts on to the SO(2) subgroup, and so the restriction of u to any given constant, large $s_+$ circle in $\mathbb{R} \times S^1$ has an integer degree that is independent of the chosen value for $s_+$.

**Lemma 5.5**: *The map* u *just defined has degree 1*.

***Proof of Lemma 5.5***: The constant map from the domain of $(\varphi, \varsigma)$ to $\mathbb{R}^2$ given by the element $(1, 0)$ corresponds via $(d\hat{\phi}, dh)$ to a section of N over $\mathcal{E}_{1,y}$, this being the orthogonal projection to N of the vector field $\partial_{\hat{\phi}}$. The latter generate the deformations of $\mathcal{E}_{1,y}$ that are given by the constant rotations of the $\hat{\phi}$ coordinate. Granted this, use Lemma 5.3 with (5.5) and (5.7) to see that these deformations are generated along the large $s_+$ part of the graph in (5.5) by the section of N that is defined by the orthogonal projection of the vector field

$$re^{-\lambda_{1n}s_+}(-\sin(\phi_+-\phi_n)\,\partial_{\theta_+} + \tfrac{3}{2x_0\sigma_0}((\lambda_1+\lambda_2)^2 + \tfrac{9}{2x_0^2})^{1/2} + \lambda_1 + \lambda_2)\cos(\phi_+-\phi_n)\,\partial_{u_+}) + \cdots$$

(5.21)

where unwritten terms are bounded in absolute value by $c_0 r e^{-(\lambda_{1n}+1/c_0)s_+}$. This last vector rotates once, counter clockwise in $\mathbb{R}^2$ as $\phi_+$ changes from 0 to $2\pi$. This implies that U has degree 1 as claimed.

*Part 5*: The following lemma gives some hint as to the relevance of the cylinders that are described by the first bullet of Proposition 5.1.

**Lemma 5.6**: *Let* $C \subset \mathbb{R} \times \mathcal{H}^+_{p*}$ *denote a properly embedded, J-holomorphic submanifold, and let* $\mathcal{E} \subset C$ *denote an end where s is unbounded from above and whose constant s slices converge in an isotopic fashion to* $\hat{\gamma}^+_p$ *as* $s \to \infty$. *Then the* $s \gg 1$ *part of* $\mathcal{E}$ *can be parametrized via (5.5)-(5.7) by a map of the sort that is described by the first bullet in Proposition 5.1.*

*Proof of Lemma 5.6*: There exists $s_* > 1$ such that the $s \geq s_*$ part of $\mathcal{E}$ is a proper submanifold with boundary in $[s_*, \infty) \times U_+$. Use (5.5)-(5.6) to view this submanifold using the coordinates $(s_+, \phi_+, \theta_+, u_+)$. The function s restricts to the $s \geq s_*$ part of $\mathcal{E}$ as a proper function with no critical points. Granted that this is so, it follows that the projection to the $(\theta_+, u_+) = (0, 0)$ cylinder restricts to the large $s_+$ part of $\mathcal{E}$ as a covering map. This covering map must have degree 1 because the constant s slices of $\mathcal{E}$ are isotopic to $\hat{\gamma}^+_p$. This understood, the large $s_+$ part of $\mathcal{E}$ has intersection number 1 with any given sufficiently large $s_+$ fiber of the projection to the $(\theta_+, u_+) = (0, 0)$ cylinder. This implies that the large $s_+$ part of $\mathcal{E}$ can be written as the graph of a map from the large $s_+$ part of $\mathbb{R} \times \mathbb{R}/(2\pi\mathbb{Z})$ to $\mathbb{R}^2$ that is described by the first bullet of Proposition 5.1.

### c) Fredholm operators

This subsection introduces some new Fredholm domain and range spaces for certain operators of the sort that are described by (3.5) and (3.6). The upcoming Proposition 5.7 supplies the analog of Proposition 3.1 for the new Fredholm incarnations these operators.

To set the stage, let $Q \in \mathbb{R} \times I_*$ denote either the complement of a single $\hat{u} = 0$ point or two $\hat{u} = 0$ points. Suppose that $\mathfrak{h} = (\varphi, \varsigma)$ maps the complement of Q in $\mathbb{R} \times I_*$ to $\mathbb{R}^2$ so as to define a graph in $\mathbb{R} \times X$ via (3.1). Let $C_\mathfrak{h}$ denote the $\Psi_\mathfrak{p}$-image of this graph. Assume in what follows that the large $|s|$ part of $C_\mathfrak{h}$ is J-holomorphic; that it obeys the first bullet in (2.9); and that it obeys the $\Delta_\mathfrak{p} = 1$ or $\Delta_{\mathfrak{p}=2} = 2$ bullets of (2.9). The pair $\mathfrak{h}$ has an associated version of the operator that is depicted in (3.9), this denoted by $D_\mathfrak{h}$. In what follows, D is used to denote an operator that is given by (3.5) and (3.6) with the extra condition

$$D = D_\mathfrak{h} \text{ on the complement of a compact set in } (\mathbb{R} \times I_*) - Q.$$
(5.22)

Operators of this sort play a central role in the upcoming proof of the $\Delta_\mathfrak{p} > 0$ version of Proposition 2.2. Part 1 of what follows defines the new domain and range spaces. This first part of the subsection ends with Proposition 5.7. The subsequent parts of the subsection supply the proof of Proposition 5.7.

*Part 1*: Let $N \to C_\mathfrak{h}$ denote the normal bundle, this being the orthogonal complement to $TC_\mathfrak{h}$ in $T(\mathbb{R} \times \mathcal{H}_\mathfrak{p})$ with orthogonality defined by the metric $\hat{\omega}(\cdot, J(\cdot))$. Here again, $\hat{\omega} = ds \wedge \hat{a} + w$. Identify $C_\mathfrak{h}$ with its inverse image via $\Psi_\mathfrak{p}$ in $\mathbb{R} \times X$. Having done so, use $(d\hat{\phi}, dh)$ to write a section of N as a map from $(\mathbb{R} \times I_*) - Q$ to $\mathbb{R}^2$. Granted this identification, a map from $(\mathbb{R} \times I_*) - Q$ to $\mathbb{R}^2$ can be viewed as a section of N over $C_\mathfrak{h}$. A map with compact support defines a section of N with compact support, and vice-versa. The afore-mentioned Riemannian metric defines a fiber metric for N and an associated metric compatible, covariant derivative for sections of N. It also defines a Riemannian metric on $TC_\mathfrak{h}$ and thus an area form. Use the fiber metric on N and $TC_\mathfrak{h}$, the covariant derivative on N, and integration with respect to this area form to define the $L^2_1$ inner product on the space of sections of N with compact support. This $L^2_1$ inner product gives an inner product on the space of compactly supported maps from $(\mathbb{R} \times I_*) - Q$ to $\mathbb{R}^2$. Use $\mathbb{H}$ to denote the completion using this $L^2_1$ inner product of the subspace of compactly supported maps from $(\mathbb{R} \times I_*) - Q$ to $\mathbb{R}^2$ whose first component is zero along the boundary of $\mathbb{R} \times I_*$. This Hilbert space $\mathbb{H}$ will be the domain space for the desired Fredholm incarnation of D.

The range Hilbert space for the new incarnation of D is a certain $L^2$ inner product space. To set the stage for the definition of this inner product, introduce $e^0$ to denote the denote the pull-back of the $\mathbb{C}$-valued 1-form in (3.2) via $\Psi_\mathfrak{p}^{-1}$. The pull back of the latter to the graph of $\mathfrak{h}$ defines a section, $e^\mathfrak{h}$, of $T^*_\mathbb{C} C_\mathfrak{h}$ and this section defines a polarization of $T^*_\mathbb{C} C_\mathfrak{h}$ since $e^\mathfrak{h} \wedge \bar{e}^\mathfrak{h} \neq 0$. Define $T^{0,1}C_\mathfrak{h}$ to be the span of $\bar{e}^\mathfrak{h}$. By way of comparison, let $\mathcal{E} \subset C_\mathfrak{h}$ denote a J-holomorphic end whose constant $s \gg 1$ slices are circles, this an end whose large $s$ slices converge in $\mathcal{H}_\mathfrak{p}$ to either $\hat{\gamma}_\mathfrak{p}^-$ or $\hat{\gamma}_\mathfrak{p}^+$. As J defines a complex structure on $\mathcal{E}$, so it defines a polarization of $T^*_{\mathbb{C}}\mathcal{E}$ as $T^{1,0}\mathcal{E} \oplus T^{0,1}\mathcal{E}$. This polarization is the same as that given by $\{e^{\flat}, \bar{e}^{\flat}\}$.

Use the metric defined by J and $\hat{\omega}$ to define a Hermitian metric on $T^*_{\mathbb{C}}(\mathbb{R} \times \mathcal{H}^+_{\mathfrak{p}})$; and use the latter to define the norm of $e^{\flat}$. This norm is denoted by $|e^{\flat}|$. Let $|\cdot|_N$ denote the fiber norm described above on N. Reintroduce the functions $a_1$ and $a_2$ from (3.4). Use these functions and $|e^{\flat}|$ to define a norm on the space of maps from $(\mathbb{R} \times I_*) - Q$ to $\mathbb{R}^2$ as follows: Let $\zeta = (\iota, o)$ denote a given map. View the map $(a_1^{-1}\iota, a_2^{-1}o)$ as a section of N over $C_{\mathfrak{h}}$. Set the norm of $\zeta$ to be $|(a_1^{-1}\iota, a_2^{-1}o)|_N |e^{\flat}|$. Use this pointwise norm and integration with respect to the area form on $C_{\mathfrak{h}}$ to define an $L^2$-inner product on the space of compactly supported maps from $(\mathbb{R} \times I_*) - Q$ to $\mathbb{R}^2$. The resulting Hilbert space is denoted by $\mathbb{L}$. This space $\mathbb{L}$ is the new range Hilbert space

By way of an explanation, the trivialization of the normal bundle of $C_{\mathfrak{h}}$ given by the 1-forms $(d\hat{\phi}, dh)$ identifies the latter with the span of the vector fields $\{\partial_{\hat{\phi}}, \partial_h\}$. The almost complex structure J preserves this span, and so endows N with the structure of a complex line bundle. The $\mathbb{C}$-valued 1-form $(a_1^{-1}\iota + i a_2^{-1}o)\bar{e}^{\flat}$ defines a section of $N \otimes T^{0,1}C_{\mathfrak{h}}$. The norm of this section as defined using the induced Hermitian metric is the norm defined above for $\zeta$.

Keep in mind for what follows that the norms that define $\mathbb{H}$ and $\mathbb{L}$ depend on the chosen pair $\mathfrak{h}$. Even so, the spaces $\mathbb{H}$ and $\mathbb{L}$ do not depend on $\mathfrak{h}$. This is so because the respective norms defined as defined by pairs $\mathfrak{h}$ and $\mathfrak{h}'$ are commensurate.

**Proposition 5.7**: *Suppose that D is described by (3.5), (3.6) and (5.22). Then D extends as a Fredholm operator from $\mathbb{H}$ to $\mathbb{L}$ with index $\Delta_{\mathfrak{p}}$ and trivial cokernel.*

*Proof of Proposition 5.7*: The proof is contained in the subsequent parts of this subsection. Part 2 explains why D is Fredholm, Part 3 computes the index, and Part 4 proves that the cokernel is trivial.

*Part 2*: Use $\|\cdot\|_{\mathbb{L}}$ to denote the $L^2$ norm that defined $\mathbb{L}$. Meanwhile, use $\|\cdot\|$ to define the $L^2$ norm on sections of N and on sections of $N \otimes T^*C_{\mathfrak{h}}$. The covariant derivative on sections of N is denoted by $\nabla$. The operator D has closed range and finite dimensional kernel if and only if there exists $c \geq 1$ such that if $\eta \in \mathbb{H}$, then

- $\|D\eta\|_{\mathbb{L}}^2 \geq c^{-1}\|\nabla\eta\|^2 - c\|\eta\|^2$.
- *If $\eta$ has support only where $|\Psi_{\mathfrak{p}}^*s| > c$, then $\|D\eta\|_{\mathbb{L}}^2 \geq c^{-1}\|\eta\|^2$.*

(5.23)

As in the case with Proposition 3.1, the finite dimensionality of the cokernel follows if the formal, $L^2$-adjoint of D also obeys (5.23). Here, the $L^2$ norm is that used to define $\mathbb{L}$. The proof that this is so differs only in notation for the proof that (5.23) holds for D and so will not be given.

To see about (5.23), it is sufficient to restrict attention to two sorts of compactly supported sections of N. With $s_1 > 1$ fixed, the first sort are those with no support where where $s > 4s_1$ on an end $\mathcal{E} \subset C_\mathfrak{h}$ whose constant $s$ slices are circles. The second sort are the sections with support only in the $s > 2s_1$ portion of such an end. The arguments in Section 3c establish the existence of an $s_1$-dependent constant $c$ that makes (5.23) true for all sections of the first sort. The proof that (5.23) for the sections with support where $s > 2s_1$ on an end $\mathcal{E}$ as just describe has three steps. These steps consider the case where the constant $s$ slices of $\mathcal{E}$ converge as $s \to \infty$ to $\hat{\gamma}_\mathfrak{p}^+$. The argument for the other case is identical but for some sign changes.

<u>Step 1</u>: Let $y \in Q$ denote the point that corresponds to the end $\mathcal{E}$. Use what is said in Proposition 5.1, Lemma 5.3 and Corollary 5.4 to view the large $s$ part of the end of $\mathcal{E}$ via (5.5)–(5.7) with $\mathfrak{y}$ in Proposition 5.1 defined using $n = 1$ and $c_0 = \alpha_y$, and with an appropriate choice for $c_1 \in \mathbb{R}-0$ and $\phi_1 \in \mathbb{R}/(2\pi\mathbb{Z})$. Take $s_1$ so that the $s \geq s_1$ part of $\mathcal{E}$ appears in this way. Nothing is lost by assuming that $D = D_\mathfrak{h}$ on the $\Psi_\mathfrak{p}$-inverse image of this part of $\mathcal{E}$.

Use $(d\theta_+, du_+)$ to identify N over the $s \geq s_1$ part of $\mathcal{E}$ with the product bundle and so write a section of N with support on the $s > 2s_1$ part of $\mathcal{E}$ as a map from the large $s_+$ part of $\mathbb{R} \times \mathbb{R}/(2\pi\mathbb{Z})$ to $\mathbb{R}^2$. Let $\mathfrak{r}$ denote such a map, but viewed as a section over $\mathcal{E}$ of N. Multiply this section by $ds_+$ and use the parametrization of $\mathcal{E}$ by $(s_+, \phi_+)$ and the complex line bundle structure on N defined by J to write the latter as a section of $N \otimes_\mathbb{C} (T^*_\mathbb{C}\mathcal{E})$. Use $(\mathfrak{r})_{0,1}$ to denote the $N \otimes T^{0,1}\mathcal{E}$ part of this section of $N \otimes T^{0,1}\mathcal{E}$.

Let $\eta$ denote a map from $(\mathbb{R} \times I_*)-Q$ to $\mathbb{R}^2$ with support only on the part of the domain that parametrizes the $s > 2s_1$ part of $\mathcal{E}$. Write the two components of $D_\mathfrak{h}\eta$ as $(\mathfrak{t}, \mathfrak{o})$ and then view $(a_1^{-1}\mathfrak{t}, a_2^{-1}\mathfrak{o})$ as a section of N over $\mathcal{E}$. With N viewed as a complex line bundle, multiply the latter by $\bar{e}^\mathfrak{h}$ to define a section of $N \otimes T_\mathbb{C}\mathcal{E}$. Denote this last section by $(D_\mathfrak{h}\eta)_{0,1}$.

<u>Step 2</u>: The lemma below is used in Step 3 to write $D_\mathfrak{h}$ near the point $(y, 0)$ in terms of the operator $\mathfrak{D}_0$ from (5.1).

**Lemma 5.8**: *There exists a first order differential operator, $\mathfrak{d}$, on the space of maps from $[s_1, \infty) \times \mathbb{R}/(2\pi\mathbb{Z})$ to $\mathbb{R}^2$, a map $v: [s_1, \infty) \times \mathbb{R}/(2\pi\mathbb{Z}) \to \mathbb{C}-0$, and $\kappa > 1$ with the following properties:*

- *The coefficients of $\mathfrak{d}$ bounded in absolute value by $\kappa\, e^{-s/\kappa}$*
- *The norms of both $v$ and $v^{-1}$ are bounded by $\kappa$.*
- *Let $\eta \in C^{\infty}((\mathbb{R} \times I_*)-Q; \mathbb{R}^2)$ with support only on the part of the domain that parametrizes the $s > 2s_1$ part of $\mathcal{E}$. Let $u$ denote $\mathcal{E}$'s version of the map to $Gl(2;\mathbb{R})$ that appears in Lemma 5.4. Then $((\mathfrak{D}_0 + \mathfrak{d})(u\eta))_{0,1} = v(D_{\mathfrak{h}}\eta)_{0,1}$*

*Proof of Lemma 5.8*: Introduce $U \subset (\mathbb{R} \times I_*)-Q$ to denote the domain that parametrizes the $s > s_1$ part of $\mathcal{E}$. Let $\eta = (\varphi', \varsigma')$ denote a map from $U$ to $\mathbb{R}^2$ that is annihilated by $D_{\mathfrak{h}}$. Let $U' \subset U$ denote an open set with compact closure. For t near zero in $\mathbb{R}$, the $\Psi_p$-image of the graph of the map $(\varphi + t\varphi', \varsigma + t\varsigma')$ defines a deformation of $\Psi_p(U') \subset \mathcal{E}$ that is J-holomorphic to first order in t. Let $(a_{\mathfrak{h}}, b_{\mathfrak{h}})$ denote the map to $\mathbb{R}^2$ from the $s \gg 1$ portion of $\mathbb{R} \times \mathbb{R}/(2\pi\mathbb{Z})$ whose graph parametrizes $\mathcal{E}$ via (5.6). Write $u\eta$ as $(a', b')$. The pair given by $(a_{\mathfrak{h}} + ta', b_{\mathfrak{h}} + tb')$ defines via (5.6) a deformation of $\Psi_p(U')$ that is J-holomorphic to first order in t if and only if $(a', b')$ obeys an equation of the form $(\mathfrak{D}_0 + \mathfrak{d}_{\mathcal{E}})(a', b') = 0$ where $\mathfrak{d}_{\mathcal{E}}$ is a certain first order differential operator whose coefficients are bounded by $c_0 e^{-s/c_0}$. It follows from this that there exists a map $v$ from $[s_1, \infty) \times \mathbb{R}/(2\pi\mathbb{Z}) \to \mathbb{C}-0$ such that the assertion given by the third bullet of the lemma holds using $\mathfrak{d} = \mathfrak{d}_{\mathcal{E}}$ and for any smooth map $\eta$ with support on U. The uniform bounds on $v$ and $v^{-1}$ can be derived using the chain rule from the formulae in Lemma 5.3.

Step 3: Granted what is said in Lemma 5.8, it is sufficient to prove that there exists $c_0 \geq 1$ such that

$$\|\mathfrak{D}_0\mathfrak{y}\|_{L^2} \geq c_0^{-1} (\|d\mathfrak{y}\|_{L^2} + \|\mathfrak{y}\|_{L^2})$$

(5.24)

for all maps $\mathfrak{y}$ with compact support on the $s_+ > 1$ part of $\mathbb{R} \times \mathbb{R}/(2\pi\mathbb{Z})$. That this is so follows from the fact that the symmetric operator

$$(a, b) \to (\partial_\phi b + 2\frac{x_0 + 4e^{-2R}}{x_0} a, -\partial_\phi a - \tfrac{1}{3} e^{-2R} b)$$

(5.25)

on $C^{\infty}(S^1; \mathbb{R}^2)$ has trivial kernel.

*Part 3*: This part of the subsection computes the Fredholm index of D. The computation has three steps.

Step 1: This step first defines from $C_\flat$ a closed manifold with empty boundary that is diffeomorphic to the complement of a point in $\mathbb{R} \times S^1$. This manifold is denoted by Z in what follows. An oriented $\mathbb{R}^2$ bundle is then defined over Z and D is shown to extend over Z as an operator that acts on sections of this bundle.

The manifold Z is defined by a suitable identification of the boundary components. To set the stage, note that any given version of D that obeys (3.5), (3.6) and (5.22) can be continuously deformed through a 1-parameter family of operators obeying (3.5), (3.6) and (5.22) to $D_\flat$. This deformation won't change the index. Granted that such is the case, assume that $D = D_\flat$.

Use (3.1) and the map $\Psi_\mathfrak{p}$ to identify $C_\flat$ with $(\mathbb{R} \times I_*) - Q$. Fix $\varepsilon \in (0, \frac{1}{8})$ and introduce $I'$ to denote the interval $[-R - \frac{1}{2}\ln((1+\varepsilon)z_*), R + \frac{1}{2}\ln((1+\varepsilon)z_*)]$. Granted the aforementioned identification, extend $C_\flat$ as $(\mathbb{R} \times I') - Q$. Having done so, introduce the function $t_+ = e^{-2(R-\hat{u})}$ where $\hat{u} > R + \ln\delta$ on $I'$ and use the pair $(x, t_+)$ to parametrize the part of $(\mathbb{R} \times I') - Q$ where $\hat{u} \in (R + \frac{1}{2}\ln((1-\varepsilon)z_*), R + \frac{1}{2}\ln((1+\varepsilon)z_*)]$. Likewise introduce the function $t_- = -e^{-2(R+\hat{u})}$ and use $(x, t_-)$ to parametrize the $\hat{u} < -R - \frac{1}{2}\ln((1-\varepsilon)z_*)$ portion of the domain $(\mathbb{R} \times I') - Q$. Use these coordinates to identify the $t_+ \in [(1-\varepsilon)z_*, (1+\varepsilon)z_*]$ portion of $(\mathbb{R} \times I') - Q$ with the portion of $(\mathbb{R} \times I_*) - Q$ where $t_- \in [-(1+\varepsilon)z_*, -(1-\varepsilon)z_*]$ using the rule $t_- = -2z_* + t_+$. The slice of Z where $t_+ = z_*$ and so $t_- = -z_*$ in Z is said in what follows to be the $z_*$-*locus*. The complement of this $z_*$-locus in Z is the interior of $(\mathbb{R} \times I') - Q$.

Define an oriented, $\mathbb{R}^2$ bundle over Z as follows: The bundle is obtained from the product $\mathbb{R}^2$ bundle over $(\mathbb{R} \times I') - Q$ by identifying the point $((x, t_+); (\zeta_1, \zeta_2))$ with the point $((x, t_- = -2z_* + t_+), (\zeta_2, -\zeta_1))$. Use $N_Z$ to denote this $\mathbb{R}^2$ bundle.

As explained next, the operator $D_\flat$ extends over the whole of Z as a differential operator on the space of sections of E. This is because $D_\flat$ when written in terms of the coordinates $(x, t_+)$ on the $\hat{u} \in (R + \ln\delta, R + \frac{1}{2}\ln z_*]$ part of $(\mathbb{R} \times I') - Q$ is the standard Cauchy-Riemann operator; and this is also the case for $D_\flat$ when written in terms of the coordinates $(x, t_-)$ on the $\hat{u} = [-R - \frac{1}{2}\ln z_*, -R - \ln\delta)$ part of $(\mathbb{R} \times I') - Q$. This extension of $D_\flat$ to Z is denoted in what follows by $\mathfrak{D}_Z$.

Step 2: This step defines a 1-parameter family of 'matching conditions' for sections of $N_Z$ with discontinuity on the $z_*$-locus in Z. The family is parametrized by the interval $[0, 1]$. A given parameter value is denoted by $\tau$.

Fix $\tau \in [0, 1]$ and suppose that $(\varphi', \varsigma')$ is a map from $\mathbb{R} \times (I_*-(y,0))$ to $\mathbb{R}^2$. This map is said to satisfy the $\tau$-matching condition when the following is true: Let $(\varphi'_-, \varsigma'_-)$ and $(\varphi'_+, \varsigma'_+)$ denote the respective $\mathbb{R}^2$-valued functions that are defined by $(\varphi', \varsigma')$ where the coordinate $t_- \in [-z_*, -\tfrac{1}{2} z_*]$ and where the coordinate $t_+ \in [\tfrac{1}{2} z_*, z_*]$. Then

$$\varphi'_+|_{t_+=z_*} = -\tau \varsigma'_- \quad and \quad \varphi'_-|_{t_-=-z_*} = \tau \varsigma'_+|_{t_+=z_*}$$

(5.26)

For each $\tau \in [0, 1]$, define the Hilbert space $\mathbb{H}_\tau$ by copying the definition of the Hilbert space $\mathbb{H}$ in Part 3 of this subsection but with the $|\hat{u}| = R + \tfrac{1}{2} \ln z_*$ boundary conditions used in Part 1 replaced by those in (5.26). Define the Hilbert space $\mathbb{L}$ as in Part 1. The $\tau = 0$ version of $\mathbb{H}_\tau$ is the Hilbert space $\mathbb{H}$. The $\tau = 1$ version is a Hilbert space of sections of $N_Z$. Meanwhile, $\mathbb{L}$ can be viewed as the closure of the space of sections of $N_Z$ with respect to the $L^2$ norm that is defined as in Part 1. Thus, $\mathbb{H}_1$ and $\mathbb{L}$ can be viewed as respective $L^2_1$ and $L^2$ Hilbert spaces of sections of $N_Z$.

**Lemma 5.9**: *For each $\tau \in [0, 1]$, the operator $D_\mathfrak{h}$ defines a Fredholm map from $\mathbb{H}_\tau$ to $\mathbb{L}$.*

*Proof of Lemma 5.9*: The conditions in (5.23) must be established for $D_\mathfrak{h}$ on the dense domain of smooth, compactly supported maps from $(\mathbb{R} \times I')-Q$ that obey (5.26). By the same token, these same conditions must be established for the formal $L^2$ adjoint. The argument that proves the analog of (5.23) for the formal $L^2$ adjoint of a given $\tau \in [0, 1]$ version of $D_\mathfrak{h}$ is identical to that just given but for cosmetics. Note in this regard that this adjoint has dense domain given by the compactly supported maps from $(\mathbb{R} \times I')-Q$ that obeys (5.26).

The new issues with regards to (5.23) for $D_\mathfrak{h}$ do not concern $D_\mathfrak{h}$ near the points in $Q$; they concern only the part of the argument that comes from Section 3c. This understood, consider the top line in (5.23). The top line is established in Section 3c using an integration by parts with the observation that the boundary terms are separately zero. The same integration by parts for $\tau \neq 0$ now yields respective $t_+ = z_*$ and $t_- = -z_*$ boundary terms that are not identically zero, but are opposite in sign. As a consequence, these terms add to zero and so make no contribution.

Consider next the lower line in (5.23). The key issue is whether (3.16) holds with $\tau \neq 0$. If this is so, then the argument used in Step 4 of Section 3c can be used here with only notational modifications to establish the desired result. To see about (3.16), use (3.17) to see that an element $(\varphi', \varsigma')$ in the kernel of $Q^-$ must be such that $\varsigma'$ is constant and

$$\varphi'|_{t_+ = z_*} = \varphi'|_{t_- = -z_*} + \varsigma' c$$

(5.27)

where $c \neq 0$ is the integral of $\mathfrak{b}_{2-}$ over $I_*$. As can be seen using (3.13), this constant $c$ is less than $-x_0 \delta^{-2} R$ and so significantly less than -2. What with (5.26), this requires that

$$-\tau \varsigma' = (\tau + c) \varsigma'.$$

(5.28)

and so $c = -2\tau$. Thus, (5.28) can not hold.

Step 3: The family of Hilbert spaces $\{\mathbb{H}_\tau\}_{\tau \in [0,1]}$ defines a smooth, Hilbert space bundle $\mathcal{H} \to [0, 1]$. Indeed, a non-isometric isomorphism from $\mathbb{H}_\tau$ to $\mathbb{H}_0$ can be defined as follows: Fix a compactly supported function $\mu: [\frac{1}{2} z_*, z_*] \to [0, 1]$ that is equal to 1 near $z_*$. Let $(\varphi', \varsigma')$ denote a given element in $\mathbb{H}_\tau$. Let $(\iota', o')$ denote the image of this element in $\mathbb{H}_0$. Then $(\iota', o') = (\varphi', \varsigma')$ except where $t_+ > \frac{1}{2} z_*$ and where $t_- < -\frac{1}{2} z_*$. The pair $(\iota', o')$ where $t_+ > \frac{1}{2} z_*$ is

$$(\iota', o')|_{(x,t_+)} = (\varphi_+', \varsigma_+')|_{(x,t_+)} + \tau \mu(t_+)(\varsigma'|_{(x,t_- = -2z_* + t_+)}, 0).$$

(5.29)

A similar formula defines $(\iota', o')$ where $t_- < -\frac{1}{2} z_*$.

The family of Fredholm operators $\{D_\mathfrak{h}: \mathbb{H}_\tau \to \mathbb{L}\}_{\tau \in [0,1]}$ defines a smooth section of the Fredholm homomorphisms from $\mathcal{H}$ to the product Hilbert space bundle $[0, 1] \times \mathbb{L}$. This being the case, all members of this family have the same Fredholm index. In particular, the Fredholm index of $D_\mathfrak{h}$ needed for Proposition 5.6 is that of $\mathfrak{D}_Z$.

Given the latter observation, the arguments used in Sections 3d and 4b,d of the article [T2] can be applied with only cosmetic changes to see that the Fredholm index of $\mathfrak{D}_Z$ on $\mathbb{H}_1$ is equal to $\Delta_\mathfrak{p}$.

*Part 4*: This part of the subsection explains why D has trivial cokernel. This will follow with a proof that the kernel of D has dimension $\Delta_\mathfrak{p}$. The proof that such is the case has five steps.

Step 1: With N viewed as the product bundle over the complement of Q in $\mathbb{R}^2$, the operator D has the schematic form that is depicted in (3.5). Suppose that $(\varphi', \varsigma')$ is in the kernel of D. Let $\Gamma$ denote the locus in $(\mathbb{R} \times I_*)-Q$ where $\varsigma' = 0$. The argument used in Step 2 of the proof of Section 3d can be repeated here to see that $\Gamma$ is described by (3.24) if it is not empty and not all of $(\mathbb{R} \times I_*)-Q$.

The argument used in Step 5 of Section 3d can be repeated to see that $\Gamma$ is not empty. This understood, assume that $\Gamma$ is not all of $(\mathbb{R} \times I_*)-Q$. The argument from this same step in Section 3d proves somewhat more about $\Gamma$. It proves; in particular, that $\Gamma$ has a non-zero, even number of edges with the following property: Either x is unbounded on the edge, or the edge has an end point on the boundary of $\mathbb{R} \times I_*$, or the closure of the edge in $\mathbb{R} \times I_*$ is a point of Q.

Step 2: Fix $T \gg 1$ and define a closed, rectangular path in the interior $\mathbb{R} \times I_*$ with sides parallel to the axis such that the constant û edges obey $|\hat{u}| = R + \frac{1}{2} \ln z_* - \frac{1}{T}$ and the constant x edges obey $|x| = T$. Orient this path so that a circumnavigation in the positive direction travels in the positive û direction on the $x = T$ edge. Use $R_T$ to denote this oriented rectangular path. If T is sufficiently large, then the restriction of $(\varphi', \varsigma')$ to $R_T$ defines a nowhere zero map from $R_T$ to $\mathbb{R}^2$. Indeed, this can be seen for the constant û edges by using the fact that $(\varphi', \varsigma')$ obeys the Cauchy-Riemann equation where $\hat{u} > R + \ln\delta$ when written as function of $(x, t_+ = e^{-2(R-\hat{u})})$; and that it also obeys these equations where $\hat{u} < -R - \ln\delta$ when written as functions of $(x, t_- = -e^{-2(R+\hat{u})})$. Meanwhile, arguments much like those used to prove Proposition 2.4 in [HT] prove that there are no zeros of $(\varphi', \varsigma')$ where $|x| \gg 1$. Granted what was just said, the pair $(\varphi', \varsigma')$ defines a map from $R_T$ to $\mathbb{R}^2-\{0\}$ for all T sufficiently large. Each such large T map has a degree; they are all the same. As explained next, this degree is negative. To see this, note that the degree is equal to the intersection number between the image of $R_T$ and any given outward directed ray in $\mathbb{R}^2$; for example, the positive x-axis. The path $R_T$ intersects the positive x-axis where $\varsigma' = 0$ and $\varphi' > 0$. These are all edges of $\Gamma$, and it follows from Step 1 that this set is non-empty when T is large. Meanwhile, (3.5) implies directly that each intersection point between the image of $R_T$ and the positive x-axis has negative local intersection number.

Step 3: Suppose that q is a zero of $(\varphi', \varsigma')$ in the interior of $(\mathbb{R} \times I_*)-Q$. It follows from (3.5) that there are no zeros of $(\varphi', \varsigma')$ save q in some small radius disk centered at q, and that $(\varphi', \varsigma')$ has positive degree as a map from the boundary of this disk to $\mathbb{R}^2-\{0\}$.

Step 4: Let $U \subset \mathbb{R}^2 \times I_*$ denote a very small radius disk centered at $(y, 0) \in Q$ with the radius such that $D = D_\mathfrak{h}$ on U and such that the graph of $(\varphi', \varsigma')$ over $U-(y,0)$ maps via $\Psi_p$ to the very large s part of the corresponding end. Denote the latter by $\mathcal{E}$. It follows from the upcoming (5.29) that are no zeros in $(\varphi', \varsigma')$ in $U-(y,0)$ if U has small radius. Granted that U has such a small radius, the pair $(\varphi', \varsigma')$ defines a map from the boundary of any concentric disk in U to $\mathbb{R}^2-\{0\}$. This map has a degree, this denoted by $n_y$. This degree is negative if Q has a single point. Indeed, this follows from what is said in Steps 2 and 3. By the same token, if Q has two points, then the sum of the degrees of the maps at the two points is negative.

The pair $(\varphi', \varsigma')$ on $U-\{0\}$ defines a section of the bundle N over the large $s$ part of $\mathcal{E}$. This section can be written with respect to the product structure for N on $\mathcal{E}$ given by the basis $\{d\theta, du\}$. The resulting map to $\mathbb{R}^2-\{0\}$ from the large $s$ part of $\mathcal{E}$ is denoted in what follows by $(a, b)$. It follows from Lemma 5.5 that the pair $(a, b)$ define a map with degree $n_y+1$.

Step 5: As noted in Lemma 5.8, the operator $D_\hbar$ on $\mathcal{E}_0$ can be decomposed as the sum $D_\hbar = \mathfrak{D}_0 + \mathfrak{d}$ where $\mathfrak{D}_0$ is given in (5.1) and where $\mathfrak{d}$ is a first order operator whose symbol and zero'th order terms have norm bounded by $c_0 e^{-s/c_0}$. The arguments for Proposition 2.4 in [HT] prove that any given element in the kernel of $D_\hbar$ at large $s$ on $\mathcal{E}$ appear as follows for some $n \in \{0, 1, 2, \ldots\}$:

$$e^{-\lambda_{1n} s}(\cos n(\phi - \phi_n), r_{1n} \sin n(\phi - \phi_n)) + \mathfrak{e}_n ,$$

(5.30)

where each $n \geq 1$ version of $\lambda_{1n}, r_{1n}$ and $\phi_n$ are as defined in (5.4), and where $\lambda_{10} = \lambda_1$. Meanwhile $\mathfrak{e}_n$ is such that $|\mathfrak{e}_n| \leq e^{-(\lambda_{1n}+1/c_0)s}$. Note that each version of (5.30) defines a map from any constant and sufficiently large $s$ circle in $\mathbb{R} \times S^1$ to $\mathbb{R}^2-\{0\}$. The $n = 0$ version has degree zero and all $n \geq 1$ versions have positive degree, this being n.

Suppose now that Q has a single point. Given that the $n_y + 1 \leq 0$ and $n_y < 0$, the pair $(a, b)$ defined in Step 4 has non-positive degree, it follows that it has degree zero. This must be true for any such pair arising from the kernel of $D_\hbar$. If $(a', b')$ is a second such pair, then a linear combination of the latter with $(a, b)$ can be found so that the result defines an $n > 0$ version of (5.30). This is impossible if the linear combination is not identically zero. The preceding conclusion implies that the kernel of $D_\hbar$ has dimension 1.

Suppose next that Q has two points. Denote these points as $(y, 0)$ and $(y', 0)$. Given that the degree in (5.30) is non-negative, and given that $n_y + n_{y'} < 0$, it follows that only the cases $(n_y = 0, n_{y'} = -1)$, $(n_{y'} = 0, n_y = -1)$ and $(n_y = -1, n_{y'} = -1)$ can occur. The argument from the preceding paragraph can be repeated to see that the kernel of $D_\hbar$ can not contain two linearly independent elements which are such that both have $n_y = 0$ or both have $n_{y'} = 0$. This constraint is satisifed only if the kernel has dimension 2.

**d) The Banach spaces $\mathbb{H}_*$ and $\mathbb{L}_*$**

There is an analog in the context of Proposition 5.7 of Section 3e's Banach spaces $\mathbb{H}_*$ and $\mathbb{L}_*$. To set the stage for the definitions, first reintroduce the notation used in (5.22). Use $\text{dist}(\cdot,\cdot)$ to denote the distance function on $C_\mathfrak{h}$ that is induce by the metric on $\mathbb{R} \times \mathcal{H}^+_{p*}$. Let $\mu$ denote a smooth, non-increasing function on $[0, \infty)$ with value 1 on $[0, \tfrac{1}{2}]$ and value 0 on $[1, \infty)$. Given $\rho > 0$ and $(x, \hat{u}) \in \mathbb{R} \times (I_*-(y,0))$, use $\mu_{\rho,(x,\hat{u})}$ to denote the function $\mu(\rho^{-1}\text{dist}(\cdot,(x,\hat{u})))$. As in (3.33), fix $\upsilon \in (0, \tfrac{1}{100})$. The norm that defines $\mathbb{H}_*$ is the sum of two terms. The first is the norm for $\mathbb{H}$, and the second is the square root of the function that assigns to a given smooth map in $\mathbb{H}$ the value

$$\sup_{(x,\hat{u}) \in \mathbb{R} \times (I_*-(y,0))} \sup_{\rho \in (0,1)} \rho^{-\upsilon} \| \mu_{\rho,(x,\hat{u})} \nabla \eta' \|^2 .$$

(5.31)

This norm is denoted by $\|\cdot\|_{\mathbb{H}*}$. The Banach space $\mathbb{H}_*$ is the completion of the set of smooth, compactly supported elements in $\mathbb{H}$ using this norm.

The Banach space $\mathbb{L}_*$ is the completion of the space of smooth, compactly supported sections of $\mathbb{L}$ using the norm that is the sum of the norm $\|\cdot\|_\mathbb{L}$ with the norm whose square is the function that is given by replacing $\nabla \eta'$ in (5.31) by $\eta'$ and by replacing the norm $\|\cdot\|$ by $\|\cdot\|_\mathbb{L}$.

The following lemma states the analogs of Lemmas 3.3 and 3.4 for these new versions of $\mathbb{H}_*$ and $\mathbb{L}_*$. Given Proposition 5.7, the assertions also follow directly Theorem 3.5.2 in [M].

**Lemma 5.10**: *Define $\mathbb{H}_*$ and $\mathbb{L}_*$ as above.*

- *Elements in $\mathbb{H}_*$ are Hölder continuous with exponent $\tfrac{1}{2}\upsilon$ and the inclusion map from $\mathbb{H}_*$ into the corresponding Hölder Banach space is continuous. In particular, there exists a constant $\kappa > 1$ that depends only on $\upsilon$ and has the following significance: If $\mathfrak{f} \in \mathbb{H}_*$, then $|\mathfrak{f}| \leq \kappa \|\mathfrak{f}\|_{\mathbb{H}*}$.*
- *If $(x, \hat{u})$ is any given point in $(\mathbb{R} \times I_*)-Q$, then $\lim_{\text{dist}(\cdot,(x,\hat{u})) \to \infty} |\mathfrak{f}|$ exists and it is zero; thus, elements in $\mathbb{H}_*$ have pointwise uniform limit zero as $s \to \infty$ on $C_\mathfrak{h}$.*
- *Any operator $D$ given by (3.5), (3.6) and (5.22) maps $\mathbb{H}_*$ to $\mathbb{L}_*$; and its inverse restricts to $\mathbb{L}_*$ so as to define a bounded linear operator from $\mathbb{L}_*$ to $\mathbb{H}_*$.*

As was the case for $\mathbb{H}$ and $\mathbb{L}$, the norms that define the Banach spaces $\mathbb{H}_*$ and $\mathbb{L}_*$ depend on the chosen pair $\mathfrak{h}$ but the spaces do not.

## 6. Proof of Proposition 2.2: The $\Delta_{\mathfrak{p}} > 0$ case

The proof in the case when some $\mathfrak{p} \in \Lambda$ versions of $\Delta_{\mathfrak{p}}$ are 1 or 2 follows much the same path as that given in the preceding section for when $\Delta_{\mathfrak{p}}$ is zero. In particular, an open/closed argument is again used for a certain $[0,1]$-parametrized family of non-linear, elliptic, first order equations for a map from the complement of either one or two $\hat{u} = 0$ points in $\mathbb{R} \times I_*$ to $\mathbb{R}/(2\pi\mathbb{Z}) \times (-\frac{4}{3\sqrt{3}}\delta_*^2, \frac{4}{3\sqrt{3}}\delta_*^2)$. As in Section 4, the $\tau = 0$ member of this family is explicit, and the $\tau = 1$ member is the desired pair $(\varphi^{\mathfrak{p}0}, \varsigma^{\mathfrak{p}0})$. The substantive differences are consequences of two related facts: The first is that the domain of $(\varphi^{\mathfrak{p}0}, \varsigma^{\mathfrak{p}0})$ is now the complement in $\mathbb{R} \times I_*$ of the afore-mentioned $\hat{u} = 0$ point or points. The second stems from Item b) of the third bullet in (2.9); the latter prescribes the behavior of $(\varphi, \varsigma)$ near the missing set in $\mathbb{R} \times I_*$. This prescribed behavior makes for a more complicated $\tau = 0$ member of the family. The new domain and the prescribed asymptotics requires versions of Proposition 5.7's the operator D.

The arguments that follow discuss only the case when $\Delta_{\mathfrak{p}} = 1$ and $\mathfrak{m}_{\mathfrak{p}} = -1$ because the $\Delta_{\mathfrak{p}} = \mathfrak{m}_{\mathfrak{p}} = 1$ arguments and those when $\Delta_{\mathfrak{p}} = 2$ are identical but for cosmetic changes.

### a) An approximation to $(\varphi^{\mathfrak{p}0}, \varsigma^{\mathfrak{p}0})$

This subsection constructs an $\mathbb{R}$-parametrized family of maps from the complement of a $\hat{u} = 0$ point in $\mathbb{R} \times I_*$ to the space $(\mathbb{R}/2\pi\mathbb{Z}) \times (-\frac{4}{3\sqrt{3}}\delta_*^2, \frac{4}{3\sqrt{3}}\delta_*^2)$ such that each member defines a graph in $\mathbb{R} \times X$ whose $\Psi_{\mathfrak{p}}$-image has $C_{\mathfrak{p}0}$'s large $|s|$ asymptotics and $C_{\mathfrak{p}0}$'s behavior near the $|\hat{u}| = R + \frac{1}{2}\ln z_*$ boundaries of $\mathbb{R} \times I_*$. The family is parametrized by the $\mathbb{R}$ coordinate of the missing $\hat{u} = 0$ point. The $\mathbb{R}$ coordinate of this point is denoted by y. There are three parts to the construction.

Looking ahead, Section 6b explain how any one of these approximations can be used as the starting solution for a $[0,1]$-parametrized family of equations whose parameter 1 solution is the desired $(\varphi^{\mathfrak{p}0}, \varsigma^{\mathfrak{p}0})$. The construction of this $[0,1]$-parametrized family requires the analytic tools that are supplied by Section 5.

*Part 1*: Construct the 1-parameter family of arcs $\{\gamma_\tau\}_{\tau \in [0,1]}$ as done in Part 1 of Section 4a. These are described in (4.1). The $\mathfrak{m}_{\mathfrak{p}} = 1$ condition implies that $\gamma_{\tau=1} \neq \gamma_{\mathfrak{p}+}$. Even so, complete the constructions of Section 4 with the family $\{\gamma_\tau\}_{\tau \in [0,1]}$ to obtain a map from $\mathbb{R} \times I_*$ to $(\mathbb{R}/2\pi\mathbb{Z}) \times (-\frac{4}{3\sqrt{3}}\delta_*^2, \frac{4}{3\sqrt{3}}\delta_*^2)$ whose graph in $\mathbb{R} \times X$ has J-holomorphic image via $\Psi_{\mathfrak{p}}$. Use $(\varphi_-, \varsigma_-)$ to denote $\mathbb{R}$-valued functions that define this map.

Let $\gamma_{\mathfrak{p}-}{}'$ denote the integral curve of $v$ in the $|\hat{u}| \leq R + \frac{1}{2}\ln z_*$ part of $\mathcal{H}^+_{\mathfrak{p}*}$ with the following properties: It starts where $\hat{u} = -R - \frac{1}{2}\ln z_*$ at the same $\phi$ angle as $\gamma_{\mathfrak{p}-}$ and it ends where $\hat{u} = R + \frac{1}{2}\ln z_*$ at the same $\phi$ angle as $\gamma_{\mathfrak{p}-}$. Let $\Delta\phi$ and $\Delta\phi'$ denote the respective angle changes along $\gamma_{\mathfrak{p}-}$ and $\gamma_{\mathfrak{p}-}{}'$. These are given by the integral that appears in the second bullet of Lemma 2.1. Require that $\Delta\phi' = \Delta\phi - 2\pi$. Redo the construction in Part 1 of Section 4a starting with the arc $\gamma_{p-}'$ at $\tau = 0$. Use $\{\gamma_\tau'\}_{\tau \in [0,1]}$ to denote the resulting family. This family is such that $\gamma_{\tau=1}' = \gamma_{p+}$. Corresponding arcs $\gamma_\tau$ and $\gamma_\tau'$ have the same $\phi \in \mathbb{R}/2\pi\mathbb{Z}$ values where $\hat{u} = -R - \frac{1}{2}\ln z_*$ and also where $\hat{u} = R + \frac{1}{2}\ln z_*$. Even so, the corresponding versions of $\Delta\phi$ and $\Delta\phi'$ differ by $-2\pi$.

Redo the constructions of Section 4 with this second family $\{\gamma_\tau'\}_{\tau \in [0,1]}$ to obtain a second graph in $\mathbb{R} \times X$ whose image via $\Psi_p$ is J-holomorphic. Use $(\varphi_+, \varsigma_+)$ to denote the map from $\mathbb{R} \times I_*$ to $\mathbb{R} \times (-\frac{4}{3\sqrt{3}}\delta_*^2, \frac{4}{3\sqrt{3}}\delta_*^2)$ that define this second graph. If necessary, add $2\pi$ times an integer to $\varphi_+$ so that $\varphi_-$ and $\varphi_+$ agree where $\hat{u} = -R - \frac{1}{2}\ln z_*$. As a consequence, $\varphi_+ - \varphi_- = -2\pi$ at $\hat{u} = R + \frac{1}{2}\ln z_*$

*Part 2*: Choose a smooth, non-decreasing map w: $[-1, 1] \to [0, 1]$ that is equal to zero on the interval $[-1, -\frac{1}{8}]$, equal to one on $[\frac{1}{8}, 1]$, and is such that $w(-s) = 1 - w(s)$. Fix $\varepsilon > 0$ and introduce functions $w^-_{\varepsilon,y}$ and $w^+_{\varepsilon,y}$ mapping $\mathbb{R}$ to $[0, 1]$ by the rules

$$w^-_{\varepsilon,y}(x) \to w(\tfrac{1}{\varepsilon}(x-y)) \quad and \quad w^+_{\varepsilon,y}(x) = w(-\tfrac{1}{\varepsilon}(x-y))$$

(6.1)

These $\hat{u}$-independent functions are used to define a graph over $\mathbb{R} \times I_*$ that is smooth in the complement of the part of the $\hat{u} = 0$ locus where $|x - y| < \varepsilon$. The graph is defined by the pair of respective $\mathbb{R}/2\pi\mathbb{Z}$ and $\mathbb{R}$ valued functions $(\varphi_{\varepsilon,y,0}, \varsigma_{\varepsilon,y,0})$ that are given by the two rules that follow.

- *Where $\hat{u} \leq 0$*: $(\varphi_{\varepsilon,y,0} = \varphi_- + w^-_{\varepsilon,y}(\varphi_+ - \varphi_-), \varsigma_{\varepsilon,y,0} = \varsigma_- + w^-_{\varepsilon,y}(\varsigma_+ - \varsigma_-))$ .
- *Where $\hat{u} \geq 0$*: $(\varphi_{\varepsilon,y,0} = \varphi_+ + w^+_{\varepsilon,y}(\varphi_- - \varphi_+ - 2\pi), \varsigma_{\varepsilon,y,0} = \varsigma_+ + w^+_{\varepsilon,y}(\varsigma_- - \varsigma_+))$ .

(6.2)

What follows are two key properties of the pair $(\varphi_{\varepsilon,y,0}, \varsigma_{\varepsilon,y,0})$.

- *The function $\varphi_{\varepsilon,y,0}$ where $|\hat{u}| = R + \frac{1}{2}\ln z_*$ is the function $\varphi^{S0}(\cdot, z_*)$ in (2.12).*
- *The $|x - y| > \varepsilon$ parts of the graph of $(\varphi_{\varepsilon,y,0}, \varsigma_{\varepsilon,y,0})$ defines via $\Psi_p$ a J-holomorphic submanifold with boundary in $\mathbb{R} \times \mathcal{H}^+_{p*}$ whose constant $|s|$ slices converge as $s \to -\infty$ to the arc $\gamma_{p-}$ and as $s \to \infty$ to the arc $\gamma_{p+}$.*

(6.3)

Introduce C' to denote the image via $\Psi_p$ of the graph of the pair $(\varphi_{\varepsilon,y,0}, \varsigma_{\varepsilon,y,0})$.

*Part 3*: Set $c_0 = \alpha_y$, fix $c_1 \in \mathbb{R} - 0$ but small and fix $\phi_1 \in \mathbb{R}/(2\pi\mathbb{Z})$ and use the resulting n = 1 version of $\mathfrak{y}$ from Proposition 5.1 to define via (5.5)-(5.7) a J-holomorphic cylinder in the $s \gg 1$ part of $\mathbb{R} \times \mathcal{H}^+_{p*}$. Use $\mathcal{E} \subset \mathbb{R} \times \mathcal{H}^+_{p*}$ to denote this cylinder. The picture supplied by Lemma 5.3 has the following implication: There exists $s_1 > s_y$ such that the $\Psi_p$-inverse image of the complement of a compact set in the $s \geq s_1$ part of $\mathcal{E}$ can be written as a graph over the complement of $(y, 0)$ in a small radius disk in $\mathbb{R} \times I_*$ about $(y, 0)$. This is to say that this part of $\Psi_p^{-1}(\mathcal{E})$ can be written as the graph

$$(x, \hat{u}) \to (x, \hat{u}, \hat{\phi} = \varphi_y(x, \hat{u}), h = \varsigma_y(x, \hat{u})),$$

(6.4)

where $(\varphi_y, \varsigma_y)$ is a smooth map from the complement of $(y, 0)$ in such a small radius disk to $\mathbb{R} \times (-\frac{4}{3\sqrt{3}} \delta_*^2, \frac{4}{3\sqrt{3}} \delta_*^2)$. Use $\rho_y$ to denote the radius of this disk.

What follow are two additional observations that follow directly from Lemma 5.3: First, the function $(x, \hat{u}) \to \varsigma_y(x, \hat{u})$ limits uniformly to $\frac{2}{3\sqrt{3}} (x_0 + 4e^{-2R})$ as $|x-y|^2 + |\hat{u}|^2$ limits to zero. The second concerns the map $\varphi_y$ on circles where $|x-y|^2 + |\hat{u}|^2$ is constant. Fix any $r_* \in (0, \rho_y)$ and define the pair $(\varphi_{\varepsilon,y,0}, \varsigma_{\varepsilon,y,0})$ using (6.8) with $\varepsilon \leq r_*$. Then the restriction of $\varphi_y$ to the circle $|x-y|^2 + |\hat{u}|^2 = r_*^2$ defines a map from $S^1$ to $S^1$ that is homotopic to the restriction of $\varphi_{\varepsilon,y,0}$.

Let $r$ now denote the radial coordinate on the disk of radius $\rho_y$ in $\mathbb{R} \times I_*$ centered on $(y, 0)$. Reintroduce the function $w$ and set $w_y$ to be the function on this same disk given by $w(2\rho_y^{-1} r - 1)$. This function is 1 where $r \geq \frac{5}{8} \rho_y$ and it is zero where $r \leq \frac{3}{8} \rho_y$.

Fix $\varepsilon < \frac{1}{8} \rho_y$ and use (6.8) to define the pair $(\varphi_{\varepsilon,y,0}, \varsigma_{\varepsilon,y,0})$. The function $\varphi_{\varepsilon,y,0}$ can be written on the $r \in (\varepsilon, \rho_y)$ part of this disk as $\varphi_{\varepsilon,y,0} = \varphi_y + \varphi'_{\varepsilon,y,0}$ where $\varphi'_{\varepsilon,y,0}$ is an $\mathbb{R}$-valued function on this part of the disk.

The pairs $(\varphi_{\varepsilon,y,0}, \varsigma_{\varepsilon,y,0})$ and $(\varphi_y, \varsigma_y)$ with the function $w_y$ are used next to define functions $(\varphi_{\varepsilon,y}, \varsigma_{\varepsilon,y})$ on the complement of $(y, 0)$ in $\mathbb{R} \times I_*$. These are given by $(\varphi_{\varepsilon,y,0}, \varsigma_{\varepsilon,y,0})$ on the complement of the radius $\rho_y$ disk centered at $(y, 0)$; and given inside this disk by

$$(\varphi_{\varepsilon,y} = \varphi_y + w_y \varphi'_{\varepsilon,y,0}, \varsigma_{\varepsilon,y} = \varsigma_y + w_y \varsigma_{\varepsilon,y,0}).$$

(6.5)

This pair is such that the $(\varphi, \varsigma) = (\varphi_{\varepsilon,y}, \varsigma_{\varepsilon,y})$ version of the following conditions are obeyed:

- *The function $\varphi$ where $|\hat{u}| = R + \frac{1}{2} \ln z_*$ is the function $\varphi^{SO}(\cdot, z_*)$ that appears in (2.12).*
- *The graph $(x, \hat{u}) \to (x, \hat{u}, \hat{\phi} = \varphi(x, \hat{u}), h = \varsigma(x, \hat{u}))$ lies in $\mathbb{R} \times X$; and as a consequence, this graph is in the domain of the map $\Psi_p$.*
- *There exists a purely S-dependent (or K-compatible) constant $\kappa_c \geq 1$ with the following significance: Assume that $z_* < \kappa_c^{-1}$ and that $\delta^2 < \kappa_c^{-1} z_*$. The $\Psi_p$-image of the $|\hat{u}| \geq R + \frac{1}{2} \ln z_* - 6$ part of the graph is J-holomorphic where $1 - 3\cos^2\theta \leq \kappa_c^{-1}$.*

- *There is a constant $\kappa_{cc} > 1$ with the following property: The $\Psi_p$-image of the graph is J-holomorphic where either $|s| \geq \kappa_{cc}$ or $1 - 3\cos^2\theta \leq \kappa_{cc}^{-1}$.*
- *The $\Psi_p$ image of the $|x - y| \leq \varepsilon$ part of this graph is J-holomorphic.*
- *Each constant $s \ll -1$ slice of the $\Psi_p$-image of this graph consists of a single arc that is isotopic rel boundary in $\mathcal{H}^+_{p*}$ to $\gamma_{p-}$. The corresponding family of such arcs converges as an isotopy rel boundary to $\gamma_{p-}$ as $s \to -\infty$.*
- *Each constant $s \gg 1$ slice of the $\Psi_p$-image of this graph consist of two components.*
  a) *The first is an arc that is isotopic rel boundary in $\mathcal{H}^+_{p*}$ to $\gamma_{p-}$. The corresponding family of such arcs converges as an isotopy rel boundary to $\gamma_{p+}$ as $s \to -\infty$.*
  b) *The second is an embedded circle that is isotopic in $\mathcal{H}^+_{p*}$ to $\hat{\gamma}^+_p$. The corresponding family of such circles converges as an isotopy rel boundary to $\hat{\gamma}^+_p$ as $s \to \infty$.*

(6.6)

By way of explanation for the third bullet, Lemma 4.6 supplies a purely S-dependent (or $\mathcal{K}$-compatible) version of $\kappa_c$ such that the $\Psi_p$-image of the $|\hat{u}| \geq R + \frac{1}{2}\ln z_* - 8$ part of the graph of $(\varphi_{\varepsilon,y}, \varsigma_{\varepsilon,y})$ is J-holomorphic where $1 - 3\cos^2\theta < \kappa_c^{-1}$. This fact is used momentarily. The constant $\kappa_{cc}$ for the fourth bullet is supplied by Lemma 4.7. The remaining bullets follow directly from the definition of $(\varphi_{\varepsilon,y}, \varsigma_{\varepsilon,y})$.

Fix $\varepsilon \in (0, \frac{1}{8})$ so that (6.6) holds.

### b) Deformations to $(\varphi^{p0}, \varsigma^{p0})$

This subsection studies a $[0, 1]$-parametrized family of equations for a map from $\mathbb{R} \times (I_* - (y, 0))$ that obeys (6.6). The initial equation is satisfied by $(\varphi_{\varepsilon,y}, \varsigma_{\varepsilon,y})$ and a solution to the final equation can serve as $(\varphi^{p0}, \varsigma^{p0})$ since the $\Psi_p$ image of its graph in $\mathbb{R} \times X$ is J-holomorphic.

To define these equations, use $(\varphi_{\varepsilon,y}, \varsigma_{\varepsilon,y})$ for $(\varphi, \varsigma)$ in the left hand side of (3.4) and write the resulting pair of functions on $\mathbb{R} \times (I_* - (y, 0))$ as $(\mathfrak{g}_1, \mathfrak{g}_2)$. These have compact support in $\mathbb{R} \times (I_* - (y, 0))$, a consequence of the third, fourth and fifth bullets of (6.6).

Use Lemma 4.6 to find a purely S-dependent (or $\mathcal{K}$-compatible) constant $r_c \geq 1$ so that the $\Psi_p$ image of the $|\hat{u}| \geq R + \frac{1}{2}\ln z_* - 8$ part of the graph of $(\varphi_{\varepsilon,y}, \varsigma_{\varepsilon,y})$ is J-holomorphic where $1 - 3\cos^2\theta \leq r_c^{-1}$. Use $\Psi_p$ to view the angle $\theta$ as a function on $\mathbb{R} \times X$. Having done so, let $\chi_c$ denote function on $\mathbb{R} \times X$ given by

$$\chi_c = 1 - \chi(2(R + \tfrac{1}{2}\ln z_* - 7 - |\hat{u}|)) \chi(4r_c(1 - 3\cos^2\theta) - 3) .$$

(6.7)

This function is equal to 1 where the $\Psi_p$-image of the graph of $(\varphi_{\varepsilon,y}, \varsigma_{\varepsilon,y})$ is not J-holomorphic, and it is equal to zero on the part of this graph that is described in the third bullet of (6.6).

Use the fourth bullet of (6.6) to find $r > 1$ such that the $\Psi_p$-image of the graph of $(\varphi_{\varepsilon,y}, \varsigma_{\varepsilon,y})$ is J-holomorphic where $1 - 3\cos^2\theta < \frac{1}{r}$. Let $\chi_{cc}: \mathbb{R} \times X \to [0,1]$ denote the function $\chi(2 - 2r(1 - 3\cos^2\theta))$. This function equals 1 where $1 - 3\cos^2\theta > \frac{1}{r}$; and it vanishes where $(1 - 3\cos^2\theta) \leq \frac{1}{2r}$.

The $\tau \in [0, 1]$ member of the family of equations asks for a pair $(\varphi, \varsigma)$ that obeys (6.12) and is such that

$$a_1 \partial_x \varphi - \partial_{\hat{u}} \varsigma - (1-\tau)(\chi_c \chi_{cc})|_{h=\varsigma} \, \mathfrak{g}_1 = 0 \quad \text{and} \quad a_2 \partial_x \varsigma + \partial_{\hat{u}} \varphi + b - (1-\tau)(\chi_c \chi_{cc})|_{h=\varsigma} \, \mathfrak{g}_2 = 0. \tag{6.8}$$

Here, as in (3.4), what is written as $a_1$, $a_2$ and b are function on $\mathbb{R} \times (I_* - (y,0))$ that are obtained from the eponymous set of functions of the variables $(\hat{u}, h)$ by setting $h = \varsigma$.

An open/closed strategy is used in what follows to construct a smoothly parametrized family $\{(\varphi_\tau, \varsigma_\tau)\}_{\tau \in [0,1]}$ such that each $\tau \in [0,1]$ member obeys (6.6) and (6.8) and with the $\tau = 0$ member given by $(\varphi_{\varepsilon,y}, \varsigma_{\varepsilon,y})$. The image via $\Psi_p$ of the graph in $\mathbb{R} \times X$ of $(\varphi_{\tau=1}, \varsigma_{\tau=1})$ is J-holomorphic since the $(\mathfrak{g}_1, \mathfrak{g}_2)$ terms in (6.8) are absent when $\tau = 1$. This being the case, this $\tau = 1$ member of the family serves for the desired $(\varphi^{p0}, \varsigma^{p0})$. To set up the open/closed argument, use $\mathcal{I}$ to denote the subset of points $\tau \in [0,1]$ for which the corresponding version of (6.8) has a solution. Since $\tau = 0$ is in $\mathcal{I}$, this set is not empty. Part 1 of this subsection explains why $\mathcal{I}$ is open. The remaining four parts explain why $\mathcal{I}$ is closed. Given that $\mathcal{I}$ is not empty, and both open and closed, this set can only be $[0,1]$. Section 6c completes the proof of the $\mathfrak{m}_p = -1$ version of Proposition 2.2 by explaining why there is but a single $(\varphi^{p0}, \varsigma^{p0})$ with the desired properties.

*Part 1*: This part of the subsection proves that $\mathcal{I}$ is open. To this end, suppose that $\tau \in \mathcal{I}$ and let $\mathfrak{h} = (\varphi_\tau, \varsigma_\tau)$ denote a corresponding pair that obeys (6.6) and (6.8). It follows from Lemma 5.9 that there exists a ball $\mathbb{B}_* \subset \mathbb{H}_*$ about the origin with two essential properties. To state them, fix for the moment $(\varphi', \varsigma') \in \mathbb{B}_*$ and use $(\varphi, \varsigma)$ to denote $(\varphi_\tau + \varphi', \varsigma_\tau + \varsigma')$. Here is the first property: The graph of $(\varphi, \varsigma)$ is in $\mathbb{R} \times X$. The second property is that the assignment to any given $(\varphi', \varsigma') \in \mathbb{B}_*$ of the corresponding $(\varphi, \varsigma) = (\varphi_\tau + \varphi', \varsigma_\tau + \varsigma')$ version of the expressions on the left hand hand side of (6.8) defines a smooth map from $\mathbb{B}_*$ to $\mathbb{L}_*$. Let $\mathcal{F}_\mathbb{B}$ denote this map.

Let $I \subset [0, 1]$ denote an open neighborhood of $\tau$, and define a map $\mathcal{F}: I \times \mathbb{B}_* \to \mathbb{L}_*$ by the rule

$$(\tau', \eta) \to \mathcal{F}(\tau', \eta) = \mathcal{F}_{\mathbb{B}}(\eta) - (\tau - \tau')(\chi_c \chi_{cc})|_{h=\varsigma} (\mathfrak{g}_1, \mathfrak{g}_2).$$
(6.9)

Fix $\tau' \in I$. The differential of $\mathcal{F}$ at $(\tau, 0)$ along the $\mathbb{B}$ factor of its domain is an operator D that obeys (3.5), (3.6) and (5.22). Lemma 5.10 asserts that D maps $\mathbb{H}_*$ surjectively to $\mathbb{L}_*$ and so the differential of $\mathcal{F}$ at any such $(\tau', 0)$ point is an isomorphism. This with the inverse function theorem supplies a smooth map, $q$, from a neighborhood $I' \subset I$ of $\tau$ to $\mathbb{B}_*$ such that $\mathcal{F}(\tau', q(\tau')) = 0$. This being the case, any given $\tau' \in I'$ version of the pair $(\varphi, \varsigma)$ = $(\varphi_\tau, \varsigma_\tau) + q(\tau')$ obeys the $\tau'$ version of (6.8). Note in this regard that $(\varphi, \varsigma)$ is smooth, a fact that can be proved in a straightforward fashion using standard elliptic regularity techniques, for example, those in Chapter 6 of [M].

Granted that $(\varphi, \varsigma)$ is smooth, and granted that the pair $(\mathfrak{g}_1, \mathfrak{g}_2)$ has compact support $\mathbb{R} \times (I_* - (y, 0))$, it follows that $(\varphi, \varsigma)$ is described by (6.6). Thus $I' \subset \mathcal{I}$ and so $\mathcal{I}$ is open.

*Part 2*: This part of the subsection outlines the proof that $\mathcal{I}$ is closed. The proof starts with a lemma which describes a compact set in $\mathbb{R} \times \mathcal{H}^+_{p*}$ with the following significance: If $\tau \in [0, 1]$ and if $(\varphi, \varsigma)$ obeys (6.6) and (6.8), then the $\Psi_p$-image of the graph of $(\varphi, \varsigma)$ is J-holomorphic on the complement of this set. The proof then derives $\tau$ and $(\varphi, \varsigma)$ independent bounds for the integral of $w$ over such a graph, and for the integral of $ds \wedge \hat{a}$ over any subset of the graph where $s$ is bounded. These integrals bounds are used with Proposition II.5.5 to control the part of the graph that lies in the complement of the $\Psi_p$-inverse image of the aforementioned compact set. The resulting control over this part of the graph is used in conjunction with some standard elliptic regularity tools to obtain $\tau$- independent pointwise bounds for the derivatives to any given order for $(\varphi, \varsigma)$.

Granted all of this preliminary work, the proof proceeds as follows: Fix a point, $\tau_0$, in the closure of $\mathcal{I}$. A sequence $\{\tau_n, (\varphi_n, \varsigma_n)\}_{n=1,2,\ldots}$ is chosen with $\{\tau_n\}_{n=1,2,\ldots} \subset \mathcal{I}$ converging to $\tau_0$ and with any given $n \in \{1, 2, \ldots\}$ version of $(\varphi_n, \varsigma_n)$ obeying (6.6) and solving the $\tau = \tau_n$ version of (6.8). The control described in the preceding paragraph is used to obtain a subsequence of $\{(\varphi_n, \varsigma_n)\}_{n=1,2,\ldots}$ that converges to a pair that obeys (6.6) and the $\tau = \tau_0$ version of (6.8).

The details of the arguments proving $\mathcal{I}$ is closed occupy the remaining Parts 3-6 of this subsection.

*Part 3*: The second lemma gives the needed integral bounds for $w$ and $ds \wedge \hat{a}$. To set the stage, reintroduce $\kappa_c$ from the third bullet of (6.6), and let $\mathcal{W} \subset \mathbb{R} \times \mathcal{H}^+_{p*}$ denote the set of points from the $|\hat{u}| > R + \frac{1}{2} \ln z_* - 6$ part of $\mathbb{R} \times X$ where $1 - 3\cos^2\theta < \frac{1}{2} \kappa_c^{-1}$.

Reintroduce $\kappa_{cc}$ from the fourth bullet of (6.6) and let $\mathcal{W}_* \subset \mathbb{R} \times \mathcal{H}^+_p$ denote the subset where both $1 - 3\cos^2\theta \leq \kappa_{cc}^{-1}$ and $|s| \geq \kappa_{cc}$.

**Lemma 6.1**: *There exists a purely $S$-dependent (or $\mathcal{K}$-compatible) constant $\kappa > 1$ with the following significance: Fix $\tau \in [0,1]$ and suppose that $\mathfrak{h} = (\varphi, \varsigma)$ is a pair that obeys (6.6) and (6.8). Use $C_\mathfrak{h}$ to denote the $\Psi_p$-image of the graph of $\mathfrak{h}$. Let $I \subset \mathbb{R}$ denote an interval of length 1. Then*

$$\int_{(C_\mathfrak{h} \cap \mathcal{W}) \cap (I \times \mathcal{H}^+_{p*})} (ds \wedge \hat{a} + w) \leq \kappa.$$

*There is a $(\varphi, \varsigma)$ and $\tau$-independent constant $\kappa_* > 1$ such that*

$$\int_{C_\mathfrak{h} \cap \mathcal{W}_*} w \leq \kappa_* \quad \text{and} \quad \int_{(C_\mathfrak{h} \cap \mathcal{W}_*) \cap (I \times \mathcal{H}^+_{p*})} ds \wedge \hat{a} \leq \kappa_*.$$

*Proof of Lemma 6.1*: The proof that follows is given in seven steps.

<u>Step 1</u>: This step, Step 2 and Step 3 consider the integral that involves $C_\mathfrak{h} \cap \mathcal{W}$. To this end, remark that this integral has two regions of support; the part of the graph of $(\varphi, \varsigma)$ where the function $\hat{u}$ is greater than $R + \frac{1}{2}\ln z_* - 6$ and where $\hat{u} < -R - \frac{1}{2}\ln z_* + 6$. The argument that follows considers the former region as the argument for the other is identical but for some sign changes.

Reintroduce $r_c$ from (6.7) and set $K = \chi(4r_c(1 - 3\cos^2\theta) - 1)$, here viewed using $\Psi_p$ as a function on $\mathbb{R} \times X$. The latter is equal to 1 where $1 - 3\cos^2\theta < \frac{1}{4} r_c^{-1}$ and it is equal to zero where $1 - 3\cos^2\theta > \frac{1}{2} r_c^{-1}$. Let $s_0$ denote the midpoint of the interval $I$ and introduce $L$ to denote the function $\chi(2|s - s_0| - 2)\chi(2(R + \frac{1}{2}\ln z_* - 6 - \hat{u}))$. This function is equal to 1 where both $s \in I$ and $\hat{u} \geq R + \frac{1}{2}\ln z_* - 6$ and it is equal to zero where $|s - s_0| > \frac{3}{2}$ or where $\hat{u}$ is less than $R + \frac{1}{2}\ln z_* - 6.5$.

<u>Step 2</u>: Consider the integral of $-\sqrt{6} d(L^2 K^2 h d\hat{\phi})$ over the graph of $(\varphi, \varsigma)$. Note that the integrand is supported on the part where the $\Psi_p$ image is $J$-holomorphic. Moreover, the integrand is equal to $\Psi_p^* w$ on $(C_\mathfrak{h} \cap \mathcal{W}) \cap (I \times \mathcal{H}^+_{p*})$; this being a consequence of the formula for $w$ in (1.6). To say more about the integrand, it proves useful to introduce the coordinate $v = e^{-2(R-\hat{u})}$ for the $\hat{u} > R + \ln\delta$ part of $\mathbb{R} \times (I_* - (y, 0))$. Using $(x, v)$ now to parametrize the graph, the 2-form $-\sqrt{6} d(L^2 K^2 h d\hat{\phi})$ appears as

$$\sqrt{6} \{L^2(K^2 + 2\varsigma K \partial_h K)(\partial_x \varphi \partial_v \varsigma - \partial_x \varsigma \partial_v \varphi) + 2\sqrt{6} KL\varsigma(\partial_v(KL)\partial_x \varphi - \partial_x(KL)\partial_v \varphi)\} dx \wedge dv.$$

(6.10)

Granted that the $\Psi_p$-image of the graph of $(\varphi,\varsigma)$ is J-holomorphic on the support of the form $-\sqrt{6}\,d(L^2K^2h\,d\hat{\phi})$, and given the properties of J in Section 1c, so $(\varphi,\varsigma)$ obey the Cauchy-Riemann equations on the support of this form when viewed as functions of $(x,v)$. This is to say that $\partial_x\varphi - \partial_v\varsigma = 0$ and $\partial_x\varsigma + \partial_v\varphi$. Thus, the 2-form in (6.10) is

$$\sqrt{6}\,\{L^2(K^2 + 2\varsigma K\,\partial_h K)(|\partial_x\varphi|^2 + |\partial_v\varphi|^2) + 2\sqrt{6}\,KL\varsigma(\partial_v(KL)\,\partial_x\varphi - \partial_x(KL)\,\partial_v\varphi)\}\,dx \wedge dv\,.$$
(6.11)

What follows is a key observation: The function $h\partial_h K$ is non-negative; this being a consequence (1.27) and the definition of $\Psi_p$. Granted that this is so, the function that multiplies the form $dx \wedge dv$ in (6.11) is no less than

$$\sqrt{6}\,\{\tfrac{1}{2}L^2K^2(|\partial_x\varphi|^2 + |\partial_v\varphi|^2) - c\delta_*^{\,4}(|\partial_x(KL)|^2 + |\partial_v(KL)|^2)\}$$
(6.12)

where $c \geq 1$ is a purely S-dependent (or $\mathcal{K}$-compatible) constant. Given the definitions of K and L, it follows that this function is no less than $-cz_*^{-2}$ with $c \geq 1$ being another purely S-dependent (or $\mathcal{K}$-compatible) constant. This means that the integral of $-\sqrt{6}\,d(L^2K^2h\,d\hat{\phi})$ over the complement of the part of the $\Psi_p$-inverse image of $(C_\mathfrak{h} - \mathcal{W}) \cap (I \times \mathcal{H}^+_{p*})$ is bounded from below by $-cz_*^{-1}$ where $c$ is again purely S-dependent (or $\mathcal{K}$-compatible).

Step 3: Stoke's theorem equates the integral of $-\sqrt{6}\,d(L^2K^2h\,d\hat{\phi})$ with the line integral

$$\int_{\mathbb{R}} (K^2L^2\varsigma\,d\varphi)\,|_{\hat{u}=R+\frac{1}{2}\ln z_*}\,.$$
(6.13)

Given the boundary condition in the first bullet of (6.6) and given the fourth bullet of Proposition 2.1, this integral is no greater than a purely S-dependent (or $\mathcal{K}$-compatible) constant.

Lemma 6.1's bound for the integral of $w$ over $(C_\mathfrak{h} \cap \mathcal{W}) \cap (I \times \mathcal{H}^+_p)$ follows directly from the conclusions of the preceding paragraph and from the conclusions of Step 3 because $w$ is a non-negative multiple of $dx \wedge d\hat{u}$ on the J-holomorphic part of $C_\mathfrak{h}$ which contains $C_\mathfrak{h} \cap \mathcal{W}$.

The asserted bound for the integral of $ds \wedge \hat{a}$ over $(C_\mathfrak{h} \cap \mathcal{W}) \cap (I \times \mathcal{H}^+_p)$ follows from the fact that the latter form when pulled back by $\Psi_p$ and written in terms of the coordinates $(x, v)$ is $dx \wedge dv$. This the case, its integral is bounded by $z_*$.

<u>Step 4</u>: The set $\Psi_p^{-1}(\mathcal{W}_*)$ is a compact set in $\mathbb{R} \times \mathcal{X}$ and so the image in $\mathbb{R} \times I_*$ of $\Psi_p^{-1}(\mathcal{W}_*)$ via the projection is compact. It follows from this that there exist $d \geq 1$ whose significance is explained next. To set the background, introduce W to denote the portion of $\mathbb{R} \times (I_* - (y, 0))$ where the following conditions are met:

$$|x| \leq d, \quad |x - y| + |\hat{u}| \geq d^{-1} \quad and \quad |\hat{u}| \leq R + \ln\delta - \tfrac{1}{8}.$$

(6.14)

Let $\tau \in [0,1]$ and suppose that $\mathfrak{h} = (\varphi, \varsigma)$ obeys (6.6) and (6.8). Then the $\Psi_p$-image of the graph of $(\varphi, \varsigma)$ over the complement of W lies in the complement of $\mathcal{W}$. This last observation implies the existence of $c \geq 1$ such that

- *The functions $a_1$ and $a_2$ that appear in (6.6) are bounded from below on W by $c^{-1}$ and bounded from above on W by $c$. Likewise b is bounded on W by $c$.*
- *The metric on the $(x, \hat{u}) \in W$ part of the graph of $(\varphi, \varsigma)$ coming from the Euclidean metric from $\mathbb{R} \times I_* \times \mathbb{R}/(2\pi\mathbb{Z}) \times (-\tfrac{4}{3\sqrt{3}}\delta_*^2, \tfrac{4}{3\sqrt{3}}\delta_*^2)$ pushes forward via $\Psi_p$ to a metric on $C_\mathfrak{h} \cap \mathcal{W}$ that is bounded respectively above and below $c$ and $c^{-1}$ times the metric that comes from $\mathbb{R} \times \mathcal{H}^+_{p*}$.*

(6.15)

The second bullet in particular implies that the Euclidean inner product on $\mathbb{R} \times \mathcal{X}$ can be used when deriving an upper bound for $-w$ and for $-ds \wedge \hat{a}$.

<u>Step 5</u>: Use (1.6) and (1.30) to write

- $\Psi_p^* w = \sqrt{6}\, d\hat{\phi} \wedge dh + \alpha_w d\hat{u} \wedge dh$,
- $\Psi_p^*(ds \wedge \hat{a}) = \upsilon_0 dx \wedge d\hat{u} + x^2 \Gamma d\hat{\phi} \wedge dh + \alpha_x dx \wedge d\hat{\phi} + \alpha_{\hat{u}} d\hat{u} \wedge dh$,

(6.16)

where $\upsilon_0$ and $\Gamma$ are positive and such that the following is true: There exists $c_0 \geq 1$ such that $\upsilon_0 > c_0^{-1}$ and $|\Gamma| + |\alpha_w| + |\alpha_x| + |\alpha_y| \leq c_0$ on $\Psi_p^{-1}(\mathcal{W})$.

<u>Step 6</u>: Suppose that $\tau \in [0,1]$ and that $\mathfrak{h} = (\varphi, \varsigma)$ obeys (6.12) and (6.13). Use the coordinates $(x, \hat{u}) \in \mathbb{R} \times (I_* - (y, 0))$ for the graph of $(\varphi, \varsigma)$ to parametrize $\Psi_p^{-1}(C_\mathfrak{h})$. This coordinate map pulls back the form $d\hat{\phi} \wedge dh$ as

$$(\partial_x \varphi \partial_{\hat{u}} \varsigma - \partial_x \varsigma \partial_{\hat{u}} \varphi)\, dx \wedge d\hat{u}.$$

(6.17)

Use (6.8) to see that the latter expression can be written as

$$(a_1|\partial_x\varphi|^2 + a_2|\partial_x\varsigma|^2 + (1-\tau)\chi(\partial_x\varphi g_1 + \partial_x\varsigma g_2))\, dx \wedge d\hat{u}\ .$$
(6.18)

Since $dx \wedge d\hat{u}$ gives the proper orientation for $C_\mathfrak{h}$, what is written above with (6.6) and (6.16) imply that

$$w|_{\wedge^2 TC_\mathfrak{h}} \geq c_1^{-1}(|d\varphi|^2 + |d\varsigma|^2)\, dx \wedge d\hat{u} - c_1 dx \wedge d\hat{u} \quad and \quad (ds \wedge \hat{a})|_{\wedge^2 TC_\mathfrak{h}} \geq -c_1 dx \wedge d\hat{u}$$
(6.19)

where $c_1 \geq 1$ enjoys a $\tau$ and $(\varphi,\varsigma)$ independent upper bound on W. These formula supply a $(\varphi,\varsigma)$ and $\tau$ independent lower bound for the respective parts of $C_\mathfrak{h}$ where $w$ and $ds \wedge \hat{a}$ are negative multiples of the area form.

Step 7: To see about the integral of $w$ over the whole of $C_\mathfrak{h}$ remark first that the argument used at the start of Lemma 4.5 has what are purely cosmetic modifications that prove that the integral of $w$ over $C_\mathfrak{h}$ is finite. With (1.6) used to identify $w$ on $\mathbb{R} \times \mathcal{H}^+_{p*}$ as $w = d(x(1 - 3\cos^2\theta\, d\theta) - \sqrt{6}\, d(\hat{h}\, d\phi)$, this same argument justifies an application of Stokes' theorem to write the integral of $w$ over $C_\mathfrak{h}$ as a sum of five terms. The first two are integrals over the arcs $\gamma_{p+}$ and $\gamma_{p-}$. That over $\gamma_{p+}$ is given by Equation (II.5.9) and that over $\gamma_{p-}$ is $(-1)$ times the $\gamma_{p-}$ version of (II.5.9). Steps 2 and 3 in the proof of Proposition II.5.1 bound the total contribution from these two terms by a purely S-dependent (or $\mathcal{K}$ compatible) constant. The third term in the sum is the integral of $-\sqrt{6}\hat{h}\, d\phi$ over $\hat{\gamma}^+_p$. This is $\frac{4\pi\sqrt{2}}{3}(x_0 + 4e^{-2R})$. The last two terms are the integrals of $\hat{h}\, d\phi$ over the two boundary components of $C_\mathfrak{h}$. Up to an overall plus/minus sign, one is the integral of $\varsigma\partial_x\varphi\, dx$ along the $\hat{u} = R + \frac{1}{2}\ln z_*$ boundary of $\mathbb{R} \times I_*$ and the other is the integral of $\varsigma\partial_x\varphi\, dx$ along the $\hat{u} = -R - \frac{1}{2}\ln z_*$ boundary of $\mathbb{R} \times I_*$. Both integrals are bounded by the integral of $\frac{4}{\sqrt{3}}\delta_*^2|\partial_x\varphi^{S0}|$ where $\varphi^{S0}$ is the value at $z_*$ of the function that appears in (2.12). In particular, it follows from what is asserted by the fourth bullet of Proposition 2.1 that both versions of the latter integral are no greater than a purely S-dependent (or $\mathcal{K}$-compatible) constant.

A bound on the integral of $ds \wedge \hat{a}$ over $C_\mathfrak{h} \cap (I \times \mathcal{H}^+_{p*})$ is obtained by mimicking what is done in Step 4 of the proof of Proposition II.5.1. Note in this regard that the integration by parts done in the latter proof has no contributions from the boundary of $C_\mathfrak{h}$ because the 1-form $\hat{a}$ annihilates the tangent space of each level set of $f$ in $M_\delta$.

*Part 4*: This first lemma below supplies $\tau$ and $(\varphi,\varsigma)$ independent, $\mathcal{O}(z_*)$ upper bound for $|\varsigma|$ near the boundary of $\mathbb{R} \times X$. This lemma states what Lemma 4.6 states for the $\Delta_p = 0$ case. The second lemma uses what is said in Lemma 6.2 to obtain $\tau$ and $(\varphi,\varsigma)$ independent, positive lower bounds for the function $1 - 3\cos^2\theta$ on various parts of the corresponding surface $C_\mathfrak{h}$. The latter are the analogs of those given in the $\mathfrak{m}_p = 0$ case by Lemma 4.7.

**Lemma 6.2**: *There exists a purely* S*-dependent (of* $\mathcal{K}$*-compatible) constant* $\kappa \geq 1$ *such that if* $\delta^2 < \kappa^{-1} z_*$, *then the following is true: Suppose that* $\tau \in [0, 1]$ *and that* $(\varphi, \varsigma)$ *is described by (6.6) and obeys (6.8). Then* $|\varsigma|$ *is bounded by* $\kappa z_*$ *where* $|\hat{u}| > R + \frac{1}{2}\ln z_* - 6$.

*Proof of Lemma 6.2*: Except for two modifications, the argument is identical to that given to prove Lemma 4.6. The first modification replaces the appeal to Lemma 4.5 with an appeal to the first inequality of Lemma 6.1. The second modification concerns the hypothetical nonsense loop in each large n version of $C_{p0n}$, this being the loop that must define a non-zero homology class in $\mathbb{R} \times \mathcal{H}^+_{p*}$. Although $C_{p0n}$ in this case does contain a loop that generates the homology of $\mathbb{R} \times \mathcal{H}^+_{p*}$, the hypothetical nonsense loop would sit entirely in either the $u > 0$ or the $u < 0$ part of $C_{p0n}$. Each of these parts is contractible in $\mathbb{R} \times \mathcal{H}^+_{p*}$, so no such loop can exist.

The next lemma gives the promised lower bounds for $1 - 3\cos^2\theta$.

**Lemma 6.3**: *There exists* $\kappa > 1$, *and given* $\varepsilon \in (0, 1]$, *there exists* $\kappa_\varepsilon > 1$; *and these have the following significance: Suppose that* $\tau \in [0,1]$ *and that* $(\varphi, \varsigma)$ *is described by (6.6) and is a solution to (6.8).*
- $\cos\theta < \frac{1}{\sqrt{3}} - \kappa_\varepsilon^{-1}$ *on the* $|u| \geq \varepsilon$ *part of* $C_\mathfrak{h}$.
- $\cos\theta < \frac{1}{\sqrt{3}} - \kappa^{-1}$ *on the* $s < -\kappa$ *part of* $C_\mathfrak{h}$.
- $\cos\theta > -\frac{1}{\sqrt{3}} + \kappa^{-1}$ *on the whole of* $C_\mathfrak{h}$.

*Proof of Lemma 6.3*: The proof has seven steps.

Step 1: This step proves that there exists a $(\varphi, \varsigma)$ and $\tau$ independent $\kappa_\varepsilon \geq 1$ such that $1 - 3\cos^2\theta > \kappa_\varepsilon^{-1}$ on the part of $C_\mathfrak{h}$ where $|u| > \varepsilon$. This implies what is asserted by the first bullet. To start, use Lemma 6.2 to choose a $(\varphi, \varsigma)$ and $\tau$ independent constant $r \geq 1$ such that the $1 - 3\cos^2\theta \leq \frac{1}{r}$ part of $C_\mathfrak{h}$ is J-holomorphic and so that $1 - 3\cos^2\theta > \frac{1}{r}$ on the boundary of $C_\mathfrak{h}$ and on the segments $\gamma_{p+}$ and $\gamma_{p-}$ that are used to describe the arc components of the large $|s|$ slices of $C_\mathfrak{h}$. Granted this, then the proof of Lemma 4.7 can be copied to prove the existences of $\kappa_\varepsilon$ but for two modifications. The first modification replaces the appeal to Lemma 4.5 with an appeal to Lemma 6.1. As in the proof of Lemma 6.2, the second modification concerns the hypothetical non-sense loop. This loop would sit entirely in either the u > 0 part of $C_{p0n}$ or in the u < 0 part, and both parts are contractible. Thus, no such loop can exist.

Step 2: This step and Step 3-6 explain why there exists a $(\varphi, \varsigma)$ and $\tau$ independent constant $\kappa > 1$ such that $1 - 3\cos^2\theta > \kappa^{-1}$ on the $s < -\kappa$ part of $C_\mathfrak{h}$. The existence of such a constant implies what is asserted by the second bullet. Existence is proved by assuming to the contrary that no such constant exists so as to derive nonsense. Granted this assumption, there exists a sequence $\{\tau_n, (\varphi_n, \varsigma_n)\}_{n=1,2,...}$ with the following properties: First, any given $n \in \{1, 2, ...\}$ version of $\tau_n \in [0, 1]$ and $(\varphi_n, \varsigma_n)$ obeys (6.6) plus the $\tau = \tau_n$ version of (6.8). Furthermore, the corresponding $\mathfrak{h} = (\varphi_n, \varsigma_n)$ version of $C_\mathfrak{h}$ has a point where both $1 - 3\cos^2\theta < 1/n$ and $s < -n$. Given what is said by Lemma II.4.8 nothing is lost by taking this point to be a local minimum of $1 - 3\cos^2\theta$ and thus a point where $u = 0$. No generality is lost by assuming that such a point occurs where $\cos\theta \sim \frac{1}{\sqrt{3}}$.

Introduce from Part 2 of Section 5a the tubular neighborhood $U_+ \subset \mathcal{H}^+_{p*}$ of the $\cos\theta = \frac{1}{\sqrt{3}}$ and $u = 0$ integral curve of $v$, this being the loop $\hat{\gamma}^+_p$. Let r be as in Step 1, and use what is said in Step 1 to choose $\varepsilon < \frac{1}{r}$ so that points in $C_\mathfrak{h}$ with $\cos\theta > \frac{1}{\sqrt{3}} - \varepsilon$ are mapped via the projection to $\mathcal{H}^+_{p*}$ to a subset in $U_+$ with compact closure. Use $V \subset U_+$ to denote the subset of points where $1 - 3\cos^2\theta < \varepsilon$.

The ensuing discussion uses the coordinates $(s_+, \phi_+, \theta_+, u_+)$ for $\mathbb{R} \times U_+$ from Section 5b and (5.5). By way of reminder, the $(\theta_+ = 0, u_+ = 0)$ locus is $\mathbb{R} \times \hat{\gamma}^+_p$, and any given constant $s_+$ and $\phi_+$ disk is J-holomorphic.

For each $n \in \{1, 2, ...\}$, choose a point in the $\mathfrak{h} = (\varphi_n, \varsigma_n)$ version of $C_\mathfrak{h}$ with $s_+$ coordinate less than $-n$ and with $\cos\theta > \frac{1}{\sqrt{3}}(1 - 1/n)^{1/2}$. Use $s_n$ to denote the value of $s_+$ at this point; and use $V_n \subset C_\mathfrak{h}$ to denote the component of this chosen point in the $\mathbb{R} \times V$ part of $C_\mathfrak{h}$. Fix $\varepsilon' \in (\frac{1}{4}\varepsilon, \frac{1}{2}\varepsilon)$ so that each index n version of $V_n$ is transversal to the locus where $1 - 3\cos^2\theta = \varepsilon'$. Introduce $V' \subset V$ to denote the $1 - 3\cos^2\theta \leq \varepsilon'$ part and use $V_n'$ to denote the connected component of the $\mathbb{R} \times V'$ part of $V_n$ that contains the chosen point where $s_+ = s_n$ and $1 - 3\cos^2\theta < 1/n$. Let $\partial V_n' \subset V_n$ denote the boundary of $V_n'$. The ensuing discussion here and in Steps 3-5 assume that at least one of the following two conditions hold for an infinite subset of $n \in \{1, 2, ...\}$.

- *$s_+$ is bounded from above on $V_n'$.*
- *There are points on $\partial V_n'$ where $s_+ \geq s_n$.*

(6.20)

The case when neither condition holds when n is sufficiently is large treated in Step 6.

Assume now that one or the other of the conditions in (6.20) hold for all indices n. If the first condition holds, use $s_{n+}$ to denote the maximum value of $s_+$ on $V_n'$. If the first condition fails but the second condition holds, use $s_{n+}$ to denote the minimum value of $s_+$ on the $s \geq s_n$ part of $\partial V_n'$. Meanwhile, $s_+$ is bounded from below on $V_n'$ in any event. If it is the case that $s \leq s_n$ on $\partial V_n'$ set $s_{n-}$ to denote the maximum value of $s_+$ on the $s \leq s_n$ part of $\partial V_n'$. If $s > s_n$ on $\partial V_n'$, then set $s_{n-}$ to be the minimum of $s$ on $V_n'$. Extra arguments are needed when the following occurs:

*Neither $\{s_n - s_{n+}\}_{n=1,2,...}$ nor $\{s_n - s_{n-}\}_{n=1,2,...}$ have convergent subsequences.*
(6.21)

The next step assumes that one or both of these requirements is violated.

Step 3: Assume that (6.21) is violated. Pass to a subsequence (hence renumbered consecutively from 1) so that one of the sequences in question is convergent. For each large n in $\{1, 2, ....\}$, translate $V_n$ by $-s_n$ along the $\mathbb{R}$ factor of $\mathbb{R} \times V_+$ and let $V_{n*}$ denote the resulting J-holomorphic submanifold. This translate is a properly embedded submanifold in $\mathbb{R} \times V$. Moreover, $s_+$ on $\mathcal{V}_n$ takes value zero, and it is bounded either from above or below by some n-independent constant $s_*$.

Use Proposition II.5.5 with Lemma II.5.6 and Lemma 6.1 to obtain a subsequence of $\{V_{n*}\}_{n=1,2,...}$ (hence renumbered consecutively) that converges on compact subsets of $\mathbb{R} \times V_+$ in the manner dictated by Proposition II.5.5. The geometric limit is a closed, J-holomorphic subvariety of $\mathbb{R} \times V_+$ that sits where $1 - 3\cos^2\theta \leq 0$ and contains a point where $1 - 3\cos^2\theta = 0$. Moreover, $s$ is bounded from either above or below on this subvariety. But this is impossible because the bound $1 - 3\cos^2\theta \leq 0$ with a point of equality implies that the limit subvariety is $\mathbb{R} \times \hat{\gamma}_p^+$.

Step 4: Now assume that (6.21) holds. Construct each index n version of $V_{n*}$ as directed in Step 3. In this case, Lemma 6.1 with Proposition II.5.5 and Lemma II.5.6 provide a subsequence of $\{V_{n*}\}_{n=1,2,...}$ (hence relabled consecutively from 1) that converges in the manner dictated by Proposition II.5.5 to $\mathbb{R} \times \hat{\gamma}_p^+$.

Let $\pi_+$ denote the projection map $(s_+, \phi_+, \theta_+, u_+) \to (s_+, \phi_+)$ from $\mathbb{R} \times U_+$ to $\mathbb{R} \times \hat{\gamma}_p^+$. The submanifold $V_n'$ has positive local intersection numbers with the constant $(s_+, \phi_+)$ disks in $\pi(V_n')$. This has the following consequence: Let $D \subset \pi(V_n')$ denote a disk whose inverse image in $V_n'$ is disjoint from the boundary. Then the restriction of $\pi_+$ to $\pi_+^{-1}(D) \cap V_n'$ is a finite to 1, branched cover with purely positive ramification points. This observation has an important consequence that is described momentarily. To set the stage, suppose that $\upsilon \subset \mathbb{R} \times \hat{\gamma}_p^+$ is an embedded, oriented loop with the following three properties: First, $\upsilon$ is -1 times the generator of $H_1(\mathbb{R} \times \hat{\gamma}_p^+; \mathbb{Z})$. Second, $\upsilon \in \pi_+(V_n')$ and $\pi_+^{-1}(\upsilon) \cap V_n{'}$ is disjoint from the boundary of $V_n{'}$. Third, $\upsilon$ does not contain any branch points of $\pi_+$ on $V_n{'}$. Granted these conditions, the projection map $\pi_+: \pi_+^{-1}(\upsilon) \cap V_n{'} \to \upsilon$ must be 1-1 on each component of $\pi_+^{-1}(\upsilon) \cap V_n{'}$. This is proved in the next paragraph.

To prove this is, note that each component of $\pi_+^{-1}(\upsilon) \cap V_n{'}$ comes via the graph $(x, \hat{u}) \to (x, \hat{u}, \hat{\phi} = \varphi_n(x,\hat{u}), h = \varsigma_n(x,\hat{u}))$ of an embedded loop in $\mathbb{R} \times (I_* - (y,0))$ that has positive linking number in $\mathbb{R} \times I_*$ with the point $(y,0)$. As each such loop is embedded, it must have linking number 1 with $(y, 0)$ and so its image via $\Psi_{p*}$ must be -1 times the generator $H_1(\mathbb{R} \times \mathcal{H}^+_{p*}; \mathbb{Z})$. This would not be the case were $\pi_+^{-1}(\upsilon) \cap V_n$ a non-trivial covering map.

Step 5: Each index n version of $V_n{'}$ has strictly positive, and locally constant intersection number with the fibers of $\pi_+$ over $(s_{n-}, s_{n+}) \in \mathbb{R} \times \hat{\gamma}_p^+$. Let $m_n$ denote this intersection number. Granted what was said in the previous step, it must be the case that the $s_+ \in (s_{n-}, s_{n+})$ part of $V_n{'}$ has $m_n$ components and $\pi_+$ restricts to each component so as to map it diffeomorphically onto $(s_{n-}, s_{n+}) \times \hat{\gamma}_p^+$. Fix a component of this part of $V_n{'}$ whose closure has an $s_+ = s_{n+}$ point where $1 - 3\cos^2\theta = \varepsilon'$. Use $A_n$ to denote the chosen component.

Let $\mathcal{V}_n$ denote the translate of $V_n$ by the constant factor $-s_{n+}$ of the $s_+$ coordinate. Let $\mathcal{V}_n{'} \subset \mathcal{V}_n$ denote the corresponding translate of $V_n{'}$ and let $\mathcal{A}_n \subset \mathcal{V}_n{'}$ denote the translate of $A_n$. The subvariety $\mathcal{V}_n$ intersects $(s_{n-} - s_{n+}, 1) \times V_+$ as a J-holomorphic submanifold with an $s_+ = 0$ point in the closure of $\mathcal{A}_n$ where $1 - 3\cos^2\theta = \varepsilon'$.

Use Proposition II.5.5 and Lemma II.5.6 with Lemma 6.1 to find a subsequence of $\{\mathcal{V}_n\}_{n=1,2,\ldots}$ (henceforth renumbered consecutively from 1) that converges on compact subsets of $(-\infty, 1) \times U_+$ to a non-empty, properly embedded J-holomorphic subvariety. Let $\mathcal{V}$ denote this limit and let $\mathcal{V}'$ and $\mathcal{A}$ denote the respective subsets of $\mathcal{V}$ that arise from the corresponding limits of $\{\mathcal{V}_n{'}\}_{n=1,2,\ldots}$ and $\{\mathcal{A}_n\}_{n=1,2,\ldots}$.

The function $1 - 3\cos^2\theta \geq 0$ on $\mathcal{V}$ and so $\mathcal{V}$ cannot contain a fiber of the projection $\pi_+$ for the latter has points where $1 - 3\cos^2\theta < 0$. This implies that $1 - 3\cos^2\theta = \varepsilon'$ at a point where $s_+ = 0$ in the closure of $\mathcal{A}$ and so $1 - 3\cos^2\theta > 0$ on the $(-\infty, -1]$ part of $\mathcal{A}$. Meanwhile, the functions $1 - 3\cos^2\theta$ and u have limit 0 as $s_+$ and thus s limit to $-\infty$ on $\mathcal{A}$. In addition, $\mathcal{A}$ has intersection number 1 with each $s_+ \leq -1$ fiber of $\pi_+$.

Granted these properties, it follows that the $s \ll -1$ part of $\mathcal{A}$ is given via (5.5)-(5.7) by a map $\mathfrak{y}$ of the sort that is described by the second bullet of Proposition 5.1. This last conclusion is nonsense for the following reason: The function $1 - 3\cos^2\theta$ is positive on $\mathcal{A}$, but any given subvariety that comes via the second bullet of Proposition 5.1 has points where $1 - 3\cos^2\theta < 0$ and s is less than any specified value.

<u>Step 6</u>: This step considers the case where neither bullet in (6.20) is satisfied. If this is the case, then the $s_+ < s_{n-}$ part of $V_n'$ will contain the end of $C_\mathfrak{h}$ whose constant $s$ slices converge in an isotopic fashion as $s \to \infty$ to $\hat{\gamma}_\mathfrak{p}^+$. This implies, in particular, that the map $\pi_+$ restricts to the $s < s_{n-}$ part of $V_n'$ as a 1-1 diffeomorphism onto $\mathbb{R} \times \hat{\gamma}_\mathfrak{p}^+$. Let $\mathcal{V}_n$ denote the translate of $V_n$ that adds $-s_{n-}$ to the $s_+$ coordinate of each point. Use Lemma 6.1 with Proposition II.5.5 and Lemma II.5.6 to obtain a subsequence of $\{\mathcal{V}_n\}_{n=1,2,...}$ (hence renumbered consecutively from 1) that converges on compact subsets of $(1, \infty) \times V$ in the manner dictated by Proposition II.5.5. Let $\mathcal{V}$ denote this limit and let $\mathcal{V}' \subset \mathcal{V}$ denote the part that comes as a limit from $\{V_n'\}$. Let $\mathcal{A} \subset \mathcal{V}'$ denote the $s_+ \geq 1$ part of $\mathcal{V}'$. A repeat of the arguments from the second to last paragraph in Step 5 prove that $\mathcal{A}$ is an embedded cylinder with boundary with the following properties: First, $1 - 3\cos^2\theta > 0$ on $\mathcal{A}$, but also $1 - 3\cos^2\theta \leq \frac{1}{2}\varepsilon$ on $\mathcal{A}$. Second, both $\theta_+$ and $u_+$ limit to 0 on $\mathcal{A}$ as $s_+ \to \infty$. Third, the projection $\pi_+$ maps $\mathcal{A}$ diffeomorphically to $[1, \infty) \times \hat{\gamma}_\mathfrak{p}^+$.

What is said in the first bullet of Proposition 5.1 and the first bullet in Lemma 6.4 imply that the constant $\varepsilon$ can be chosen in advance so that $\mathcal{A}$ is the graph of a smooth map $(s_+, \phi_+) \to \mathfrak{y} = (\theta_+ = a_+(s_+, \phi_+), u_+ = b_+(s_+, \phi_+))$ that obeys (5.7) on $[1, \infty) \times \mathbb{R}/(2\pi\mathbb{Z})$ and has the form

$$\mathfrak{y} = c_0(e^{-\lambda_1 s_+} + e_1, 0) + c_1 \mathfrak{y}_{1+} + \mathfrak{e}_1$$

(6.22)

where $c_0 > 0$ and $c_1 > 0$, where $e_1$ is given in (5.8) and $\mathfrak{y}_{1+}$ is some $\phi_1$ version of (5.4). Meanwhile, $\mathfrak{e}_1$ is described in the first bullet of Proposition 5.1.

The missing $\hat{u} = 0$ point in $\mathbb{R} \times I_*$ that defines the domain of $(\varphi_n, \varsigma_n)$ has coordinates $(y, 0)$. Let $\alpha_y$ denote the strictly positive constant from Lemma 5.2 that is assigned to $y$. Let $\mathcal{E}_n$ denote the index $n$ version of the end $C_{\mathfrak{h}=(\varphi_n, \varsigma_n)}$ whose constant $s$ slices converge as $s \to \infty$ to $\hat{\gamma}_\mathfrak{p}^+$. Lemmas 5.3 and 5.6 imply that the $s \gg 1$ part of $\mathcal{E}_n$ is given via (5.5)-(5.6) by a map from Proposition 5.1's first bullet that can be written as

$$\mathfrak{y}_{(n)} = \alpha_y(e^{-\lambda_1 s_+} + e_1, 0) + c_{1n} \mathfrak{y}_{1+} + \mathfrak{e}_{1n}$$

(6.23)

where $c_{1n} \in \mathbb{R} - 0$, where $\mathfrak{y}_{1+}$ is given by some $\phi_1 = \phi_{1n}$ version of (5.4) and where $\mathfrak{e}_{1n}$ obeys the bounds given in Proposition 5.1 for the latter's $\mathfrak{e}_1$. What follows is now a direct consequence of (6.22) and (6.23): Given $\varepsilon_* > 0$, there exists $n_* \geq 1$ such that

$$|\alpha_y - c_0 e^{-\lambda_1 |s_{n-}|}| < \varepsilon_*$$

(6.24)

when $n \geq n_*$. But this is nonsense given that $\alpha_y > 0$ and $\{s_n\}_{n=1,2,...}$ is unbounded from below.

Step 7: This step proves the third bullet of Lemma 6.3. The proof starts by assuming that the assertion is false so as to derive some nonsense. Granted there is no such $\kappa$, there is a sequence $\{\tau_n, (\varphi_n, \varsigma_n)\}_{n=1,2...}$ with the following properties: First, any given $n \in \{1, 2, ...\}$ version of $\tau_n \in [0, 1]$ and $(\varphi_n, \varsigma_n)$ obeys (6.6) plus the $\tau = \tau_n$ version of (6.8). Furthermore, the corresponding $\mathfrak{h} = (\varphi_n, \varsigma_n)$ version of $C_\mathfrak{h}$ has a point where both $\cos\theta < -\frac{1}{\sqrt{3}} + 1/n$. Let $(s_n, p_n) \in \mathbb{R} \times \mathcal{H}^+_{p*}$ denote such a point. Use Proposition II.5.5 and Lemma II.5.6 with Lemma 6.1 to see that the sequence $\{|s_n|\}_{n=1,2,...}$ can not have convergent subsequences. What is said in the second bullet of Lemma 6.3 implies that $\lim_{n \to \infty} s_n = \infty$. This understood, the arguments used in Steps 2-6 can be used with only cosmetic changes to generate the desired nonsense.

*Part 5*: Let $\tau \in [0,1]$ and let $\mathfrak{h} = (\varphi, \varsigma)$ be as described in (6.6) and a solution to (6.8). The subvariety $C_\mathfrak{h}$ has an end whose constant $s$ slices converge in an isotopic fashion as $s \to \infty$ to $\hat{\gamma}^+_p$. The following is a consequence of Lemmas 5.3, 5.10 and Proposition 5.10: There exist constants $s_\mathfrak{h} \geq 1$ and $c_{1\mathfrak{h}} \in \mathbb{R}-0$ and $\phi_{1\mathfrak{h}} \in \mathbb{R}/(2\pi\mathbb{Z})$ such that $\mathcal{E}$ intersects the $s_+ \in [s_\mathfrak{h}, \infty)$ part of $\mathbb{R} \times U_+$ as a smooth, properly embedded submanifold with boundary on the $s_+ = s_\mathfrak{h}$ slice. Furthermore, this intersection is given by the graph of a smooth map as described in the first bullet of Proposition 5.1 with domain $[s_\mathfrak{h}, \infty) \times \mathbb{R}/(2\pi\mathbb{Z})$ that has the form depicted in (6.23) with $c_1 = c_{1\mathfrak{h}}$ and with $\mathfrak{y}_{1+}$ defined using $\phi_1 = \phi_{1\mathfrak{h}}$. The next lemma says something about the constants $s_\mathfrak{h}$ and $c_{1\mathfrak{h}}$.

**Lemma 6.4**: *There exists $\kappa > 1$ with the following significance: Let $\tau \in [0,1]$ and suppose that $\mathfrak{h} = (\varphi, \varsigma)$ is described by (6.6) and that it obeys (6.8). Then the corresponding constant $s_\mathfrak{h}$ can be chosen so that $s_\mathfrak{h} \leq \kappa$. Meanwhile, $|c_{1\mathfrak{h}}| \in [\kappa^{-1}, \kappa]$.*

*Proof of Lemma 6.4*: An upper bound for $s_\mathfrak{h}$ is obtained using what are essentially the same arguments as those in Step 6 of the proof of Lemma 6.3. The salient difference in this case is that the assumption of no uniform upper bound gives a sequence $\{\tau_n, (\varphi_n, \varsigma_n)\}_{n=1,2,...}$ with the property that that the corresponding sequence $\{s_{n-}\}_{n=1,2,...}$ is now unbounded from above instead of from below. This understood, the inequality in (6.24) is replaced by $|\alpha_y - c_0 e^{\lambda_1 s_{n-}}| < \varepsilon_*$ which cannot hold when n is large if $\{s_n\}$ diverges.

The upper and lower bounds on $c_\mathfrak{h}$ follow in a straightforward manner given the apriori bound on $s_\mathfrak{h}$. In fact, the upper bound follows from the constraint that $1 - 3\cos^2\theta$ is positive. The lower bound follows by assuming the contrary and deriving a contradiction from a limit submanifold that is described by (6.22) and obtained with the help of Proposition II.5.5 from a sequence $\{\tau_n, (\varphi_n, \varsigma_n)\}_{n=1,2,\ldots}$ that has $|c_{1n}| < 1/n$

Let $\tau$ and $\mathfrak{h} = (\varphi, \varsigma)$ be as described above. The $s \ll -1$ slices of $C_\mathfrak{h}$ converge in an isotopic fashion to the arc $\gamma_{p-}$ as $s \to -\infty$; and the $s \gg 1$ slices have a component that converges in an isotopic fashion as $s \to \infty$ to the arc $\gamma_{p+}$. The next lemma asserts that the convergence in both cases is suitably $\mathfrak{h}$ and $\tau$ independent. This lemma is the analog of what Lemma 4.8 states for the $\Delta_\mathfrak{p} = 0$ case.

**Lemma 6.5**: *Given $\varepsilon > 0$, there exists $\kappa_\varepsilon > 1$ with the following significance: Suppose that $\tau \in [0, 1]$ and that $(\varphi, \varsigma)$ is described by (6.6) and obeys (6.8). Then*
- $|\varphi(x, \hat{u}) - \phi(\gamma_{p-}|_\hat{u})| + |\varsigma(x, \hat{u}) - h(\gamma_{p-})| < \varepsilon$ *where* $x < -\kappa_\varepsilon$.
- $|\varphi(x, \hat{u}) - \phi(\gamma_\tau|_\hat{u})| + |\varsigma(x, \hat{u}) - h(\gamma_\tau)| < \varepsilon$ *where* $x > \kappa_\varepsilon$.

***Proof of Lemma 6.5***: What with Lemmas 6.3 and 6.4, the proof of Lemma 4.8 can be quoted in an essentially verbatim fashion to prove Lemma 6.5

*Part 6*: This last part of the subsection completes the proof that the set $\mathcal{T}$ is closed. To start, suppose that $\tau \in [0,1]$ and that $\mathfrak{h} = (\varphi, \varsigma)$ is described by (6.6) and obeys (6.8). Lemmas 6.2-6.5 supply a $\tau$ and $\mathfrak{h}$ independent disk $U \subset \mathbb{R} \times I_*$ centered on $(y, 0)$ and a compact set in $\mathbb{R} \times \mathcal{X}$ such that the graph of $(\varphi, \varsigma)$ over $(\mathbb{R} \times I_*) - U$ maps into this compact set and has uniform limits as $x \to \pm\infty$. This implies that the functions $a_1, a_2$ and $b$ that appear in (6.8) have $\mathfrak{h}$ and $\tau$ independent bounds for $(x, \hat{u}) \in (\mathbb{R} \times I_*) - U$, and that $a_1$ and $a_2$ are bounded away from zero by $\mathfrak{h}$ and $\tau$ independent, positive constants. As a consequence, standard elliptic regularity arguments of the sort that can be found in Chapter 6 of [M] can be employed to see that the absolute values of the deriviates of $\mathfrak{h}$ on $(\mathbb{R} \times I_*) - U$ to any given order have $\mathfrak{h}$ and $\tau$ independent bounds. Lemma 6.2 insures that the boundary values also enjoy $\mathfrak{h}$ and $\tau$ independent bounds. Lemma 6.5 insures that these $\tau$ and $\mathfrak{h}$ independent derivative bounds hold uniformly as $x \to \pm\infty$.

Meanwhile, the part of $C_\mathfrak{h}$ that is parametrized via $\Psi_{p*}$ by $U$ maps into the $\mathbb{R} \times U_+$ part of $\mathbb{R} \times \mathcal{H}_{p*}$; and in particular the part where $s_+ \geq \kappa$ with $\kappa$ as in Lemma 6.4. As a consequence, this part of $C_\mathfrak{h}$ can be described using (5.5)-(5.6) by a solution to (5.7) from the first bullet of Proposition 51 that has the form given in (6.23) with $c_{1\mathfrak{h}}$ as described in Lemma 6.5.

Granted all of this, suppose that $\{\tau_n\}_{n=1,2,\ldots} \in \mathcal{I}$ converges to $\tau_0 \in [0, 1]$. For each index n, let $(\varphi_n, \varsigma_n)$ denote a pair described by (6.6) and obeying the $\tau = \tau_n$ version of (6.8). Use the uniform bounds described in the preceding paragraphs for $\{(\varphi_n, \varsigma_n)\}_{n=1,2,\ldots}$ on $(\mathbb{R} \times I_*) - U$ with the Arzoli-Ascoli theorem to obtain a subsequence that converges on $(\mathbb{R} \times I_*) - U$ in the strong $C^\infty$ topology to a pair $(\varphi, \varsigma)$ whose graph over $(\mathbb{R} \times I_*) - U$ maps into $\mathbb{R} \times X$. Meanwhile, use the uniform bounds on the constants $\{c_{1\mathfrak{h}=(\varphi_n,\varsigma_n)}\}_{n=1,2,\ldots}$ to obtain a subsequence as above whose corresponding sequence of constants converges to a non-zero limit. Let $c_{1\mathfrak{h}}$ denote the latter. Granted this convergence, it follows that $(\varphi, \varsigma)$ extends over $U-(y,0)$ to give a solution to the $\tau = \tau_0$ version of (6.8) that is described by (6.6). Thus $\tau_0 \in \mathcal{I}$.

### c) Uniqueness

This subsection completes the proof of Proposition 2.2 by proving that there is only one pair $(\varphi^{p0}, \varsigma^{p0})$ that obeys the conditions imposed by Proposition 2.2. This uniqueness assertion is one consequence of the lemma that follows.

**Lemma 6.6**: *Fix $y \in \mathbb{R}$ and $\tau \in [0, 1]$. There exists exactly one pair $(\varphi, \varsigma)$ with the domain $\mathbb{R} \times (I_* - (y,0))$ that obeys the conditions set forth in (6.6) and (6.8).*

*Proof of Lemma 6.6:* Suppose that $y \in \mathbb{R}$, that $\tau \in [0, 1]$ and that $(\varphi^{(0)}, \varsigma^{(0)})$ and $(\varphi^{(1)}, \varsigma^{(1)})$ are two pair with domain $\mathbb{R} \times (I_* - (y,0))$ that obey (6.6) and (6.8). Write

$$(\varphi^{(1)}, \varsigma^{(1)}) = (\varphi^{(0)} + \varphi', \varsigma^{(0)} + \varsigma')$$

(6.25)

where $\eta = (\varphi', \varsigma')$ is a smooth map from $(\mathbb{R} \times I_*) - (y,0)$ to $\mathbb{R}$ that obeys

- $\varphi' = 0$ *where* $|\hat{u}| = R + \tfrac{1}{2} \ln z_*$.
- $\lim_{|x| \to \infty} (\varphi', \varsigma') = 0$.

(6.26)

This pair obeys an equation of the form $D\eta = 0$ with $D$ described by (3.5) and (3.6) with coefficient functions $\mathfrak{a}_1, \mathfrak{a}_2, \mathfrak{b}_1$ and $\mathfrak{b}_2$ as described in the proof of Lemma 4.9.

Let $C_0$ and $C_1$ denote the respective $\mathfrak{h} = (\varphi^{(0)}, \varsigma^{(0)})$ and $\mathfrak{h} = (\varphi^{(1)}, \varsigma^{(1)})$ versions of $C_\mathfrak{h}$. These have corresponding ends where the constant $s$ slices converge in an isotopic fashion to $\hat{\gamma}_{\mathfrak{p}}^+$ as $s \to \infty$. These ends can be written via (5.5) and (5.6) in the $s_+ \geq c_0$ part of $\mathbb{R} \times U_+$ as graphs over $[s_+, \infty) \times \mathbb{R}/(2\pi\mathbb{Z})$. These respective maps, $\mathfrak{y}_{(0)}$ and $\mathfrak{y}_{(1)}$, are described by (6.23). Write $\mathfrak{y}_{(1)}$ as $\mathfrak{y}_{(0)} + \mathfrak{y}$ where $\mathfrak{y} = (a, b)$ here denoting a smooth map from $[c_0, \infty) \times \mathbb{R}/(2\pi\mathbb{Z})$ to $\mathbb{R}^2$ with limit zero as $|s| \to \infty$. In particular, it follows from (6.23) that $\mathfrak{y}$ can be written as $\mathfrak{y} = c\mathfrak{y}_{1+} + \mathfrak{e}$ where $c \in \mathbb{R}$ is non-zero if $\mathfrak{y}_{(0)} \neq \mathfrak{y}_{(1)}$, where $\mathfrak{y}_{1+}$ is given by (5.4) with $\phi_1$ determined by $\mathfrak{y}_{(0)}$ and $\mathfrak{y}_{(1)}$, and where $|\mathfrak{e}| \leq c_0 |c| e^{-(\lambda_{11}+1/c_0)s_+}$. Note in particular that $\mathfrak{y}$, if not identically zero, it defines a degree 1 map from any constant, $s_+ \gg 1$ circle in $\mathbb{R} \times \mathbb{R}/(2\pi\mathbb{Z})$ to $\mathbb{R}^2 - 0$.

A change of variables relates the pair $\eta = (\varphi', \varsigma')$ at points $(x, \hat{u})$ near $(y, 0)$ to the pair $\mathfrak{y}$. This formula takes the form $\mathfrak{y} = \mathfrak{u}(1 + \mathfrak{e}')\eta$ where $\mathfrak{u}$ is the linear map from Lemma 5.5 and where $|\mathfrak{e}'| < c_0 \, e^{-s_+/c_0}$. With the preceding as background, use the arguments in Steps 1-4 of Part 4 in Section 5c to see that $\eta = (\varphi', \varsigma')$, if not identically zero, defines a map from the boundary of any very small radius disk about the point $(y, 0)$ to $\mathbb{R}^2 - \{0\}$ with negative degree. The relation $\mathfrak{y} = \mathfrak{u}(1 + \mathfrak{e})\eta$ implies that $\mathfrak{y}$ must define a non-positive degree map to $\mathbb{R}^2 - \{0\}$ from any sufficiently large $s_+$ circle in $\mathbb{R} \times \mathbb{R}/(2\pi\mathbb{Z})$; and in particular $\mathfrak{y}$ can not have degree 1.

This paradox is avoided only if $(\varphi^{(0)}, \varsigma^{(0)}) = (\varphi^{(1)}, \varsigma^{(1)})$.

**7) Cobordisms to the ech-HF submanifold moduli space**

Section 2 describes sets of the form $\mathcal{C}_0 = \{C_{S0}, \{C_{\mathfrak{p}0}\}_{\mathfrak{p} \in \Lambda}\}$ with $C_{S0}$ being a surface with boundary in $\mathbb{R} \times M_\delta$ and with each $\mathfrak{p} \in \Lambda$ version of $C_{\mathfrak{p}0}$ being a surface with boundary in $\mathbb{R} \times \mathcal{H}^+_{\mathfrak{p}*}$. The interiors of these surfaces are J-holomorphic. Each looks much like the portion of an ech-HF submanifold in the relevant part of $\mathbb{R} \times Y$. However, these surfaces with boundary do not necessarily fit together so as to define a closed surface in $\mathbb{R} \times Y$. As noted in Part 2 of Section 1f, the sets that are described in Section 2 comprise one boundary of a cobordism with two boundary components, the other being the moduli space of ech-HF submanifolds. This cobordism has an associated proper function mapping it to $[0, 1]$ with the inverse image over 1 being the boundary composed of ech-HF submanifolds. This section first describes and then constructs these cobordism spaces.

Section 2a describes a data set of the form $(z_*, \delta, x_0, R)$ along with an almost complex structure $J_{HF}$ for the construction of any given version of $\mathcal{C}_0$. With $(\delta, x_0, R)$ and $J_{HF}$ specified, an almost complex structure for $\mathbb{R} \times Y$ is then chosen subject to the conditions given in Part 1 of Section 1c. This almost complex structure is again denoted by J.

The definitions in Section 2 and the constructions in the previous sections require the choice of an orbit in $\mathcal{A}_{HF}/\mathbb{R}$ of a given Lipshitz surface. As in Section 2, the latter determine an upper bound on $z_*$ and $\delta$. As explained in Section 1a, the choice of $\delta$ determines an upper bound for $x_0$ and the choice for $x_0$ determines one for R. The required upper bounds for $z_*$ and $\delta$ may need some refinement in order to construct the cobordism space. The refined upper bounds are stated as needed for the various constructions that follow. In Section 2 and in what follows, the upper bounds in question for $z_*$ and $\delta$ can be chosen so as to hold for all subvarieties chosen from a given finite or compact set in $\mathcal{A}_{HF}/\mathbb{R}$. Likewise, if $\mathcal{K} \subset \mathcal{A}_{HF}$ is a given $\mathbb{R}$-invariant, weakly compact set, then the parameters $z_*$ and $\delta$ can be chosen to be $\mathcal{K}$-compatible.

As in the previous sections, S is used to denote a chosen Lipshitz submanifold. Likewise, $(\hat{\Theta}_-, \hat{\Theta}_+) \in \hat{\mathcal{Z}}^S$ is chosen. Use $(\Theta_-, \Theta_+)$ again to denote the corresponding pair in $\mathcal{Z}_{ech,M}$.

**a) The cobordism space**

A point in the cobordism space consists of a pair $(\tau, \mathcal{C})$ where $\tau \in [0, 1]$ and where $\mathcal{C} = \{C_S, \{C_\mathfrak{p}\}_{\mathfrak{p} \in \Lambda}\}$ is a set of $G+1$ submanifolds with boundary in $\mathbb{R} \times Y$. The map to $[0,1]$ sends $(\tau, \mathcal{C})$ to $\tau$. The three parts that follow describe $\mathcal{C}$. Part 1 describes $C_S$, the second part describes the various $\mathfrak{p} \in \Lambda$ versions of $C_\mathfrak{p}$. The third part describes how $\tau$ enters the picture. The notation used in Section 2 is used here also.

*Part 1*: What is denoted by $C_S$ is a properly embedded submanifold with boundary in the $f^{-1}([1+z_*, 2-z_*])$ part of $\mathbb{R} \times M_\delta$ whose interior is J-holomorphic. This submanifold is characterized in part by the four properties that are listed in what follows.

To set the stage for the statement of the first property, introduce the constant $\rho_S$, the disk bundle $N_0 \to S$, the map $\mathfrak{e}_S$ and the other notation from Part 5 of Section 2c. Let $\kappa_0$ denote the constant from Lemma II.6.5 and let $\kappa$ denote the constant from Lemma II.6.6. Introduce $\kappa_S$ to denote $10^6 \kappa \kappa_0$. Let $U \subset \mathbb{R} \times [1, 2] \times \Sigma$ denote the tubular neighborhood of S that is described in Lemma II.6.5.

PROPERTY 1: *View $C_S$ as a submanifold in $\mathbb{R} \times [1+z_*, 2-z_*] \times \Sigma$. As such, $C_S$ lies in* U *and in the image via the exponential map $\mathfrak{e}_S$ of the radius $\kappa_S^{-1} \rho_S^2$ disk subbundle in $N_0$. Moreover, $C_S$ has intersection number 1 with the $\mathfrak{e}_S$ image of each fiber of this disk bundle over the $t \in [1+z_*, 2-z_*]$ part of* S.

The next property writes $C_S$ as the image of a map from the $t \in [1+z_*, 2-z_*]$ part of S that has the form $\mathfrak{e}_S \circ \eta$ where $\eta$ is a section of $N_0$. This upcoming property also refers to the Fredholm operator $D_S$ that is discussed in Parts 2-4 of Section II.6e and depicted in (1.25). The *kernel* of $D_S$ is the vector space of sections in the domain Hilbert space that are annihilated by $D_S$.

PROPERTY 2: *The section $\eta$ is $L^2$ orthogonal to the restriction of each element in the kernel of $D_S$ to the part of S where $t \in [1+z_*, 2-z_*]$.*

The third property speaks to the large $|s|$ behavior of $C_S$. The notation borrows from the fourth bullet of Proposition 2.1. By way of a reminder, let $\hat{\upsilon}_-$ and $\hat{\upsilon}_+$ denote the respective HF-cycles that are used to define $\Theta_-$ and $\Theta_+$. Let q denote a given intersection point from either with $\Sigma$, this a point in $C_- \cap C_+$. The corresponding integral curve of $\mathfrak{v}$ from $\hat{\upsilon}_-$ or $\hat{\upsilon}_+$ appears as $(1,2) \times q$ when writing the $f^{-1}(1,2) \subset M$ as $(1,2) \times \Sigma$. The point q labels a corresponding $s << -1$ or $s >> 1$ end of S, this denoted by $\mathcal{E}_{Sq}$. This end of S is in q's component of $\mathbb{R} \times [1,2] \times (T_- \cap T_+)$. Note also that the functions $(s,t)$ restrict as coordinate functions to $\mathcal{E}_{Sq}$. The normal bundle $N_S$ over $\mathcal{E}_{Sq}$ is identified with the product $\mathbb{R}^2$ bundle in the manner that is described just prior to (2.6) and this identification is used to view a section as a map to $\mathbb{R}^2$. Meanwhile, the exponential map $\mathfrak{e}_S$ over $\mathcal{E}_{Sq}$ is written as in (2.6).

As in Proposition 2.1, $q_*$ is used to denote the point in $T_- \cap T_+$ near q where the corresponding segment of an integral curve of $v$ from $\Theta_-$ or $\Theta_+$ intersects $\Sigma$. This point $q_*$ has distance $c_0\delta$ or less from q and so lies in q's component of $T_- \cap T_+$. Writing $\mathfrak{e}_S$ over $\mathcal{E}_{Sq}$ writes $q_*$ as a section of the normal bundle $N_S$ over $\mathcal{E}_{Sq}$.

PROPERTY 3: *Let $q \in C_- \cap C_+$ denote an intersection point of an integral curve of $\mathfrak{v}$ from either $\hat{\upsilon}_-$ or $\hat{\upsilon}_+$. The section $\eta$ over the $t \in [1+z_*, 2+z_*]$ part of $\mathcal{E}_{Sq}$ converges pointwise as $s \to -\infty$ or $s \to \infty$ to $q_*$.*

The final property views $C_S$ as sitting in $\mathbb{R} \times [1+z_*, 2-z_*] \times \Sigma$. It talks about the behavior of $C_S$ where $t$ is near the end-points of the interval $[1+z_*, 2-z_*]$. The statement of this property uses the notation from Part 1 of Section 1c. In particular, Part 1 of Section 1c uses the coordinates $(\varphi_+, \hat{h}_+)$ for any given component of the region $T_+ \subset \Sigma$, and it uses the coordinate $(\varphi_-, \hat{h}_-)$ for any given component of $T_-$.

PROPERTY 4: *View $C_S$ as a submanifold with boundary in $\mathbb{R} \times [1+z_*, 2-z_*] \times \Sigma$. As such, a neighborhood of its boundary has the following properties:*

- *The $[1+z_*, 1+z_S]$ portion of $C_S$ has G components, with one mapping to each component of $T_+$. A given component of this portion of $C_S$ is the image of a map from $\mathbb{R} \times [z_*, z_S]$ to $\mathbb{R} \times [1+z_*, 1+z_S] \times \mathbb{R}/(2\pi\mathbb{Z}) \times (\frac{4}{3\sqrt{3}}\delta_*^2, \frac{4}{3\sqrt{3}}\delta_*^2)$ that has the form*

$$(x, z) \to (s = x, t = 1+z, \varphi_+ = \varphi^S(x,z), h_+ = \varsigma^S(x,z))\,.$$

- *The $[2-z_S, 2-z_*]$ portion of $C_S$ has G components, with one mapping to each component of $T_-$. A given component of this portion of $C_S$ is the image of a map from $\mathbb{R} \times [z_*, z_S]$ to $\mathbb{R} \times [2-z_S, 2-z_*] \times \mathbb{R}/(2\pi\mathbb{Z}) \times (\frac{4}{3\sqrt{3}}\delta_*^2, \frac{4}{3\sqrt{3}}\delta_*^2)$ of the form*

$$(x, z) \to (s = x, t = 2-z, \varphi_- = \varphi^S(x, z), h_- = \varsigma^S(x, z))\,.$$

In short, these five properties say that $C_S$ looks much like the $\mathbb{R} \times M_\delta$ part of an ech-HF submanifold. Note in particular that these properties are satisfied if $C_S = C_{S0}$ with the latter coming from a set of the sort that is described in Section 2.

*Part 2*: What is denoted by $C_p$ is a properly embedded submanifold with boundary in $\mathbb{R} \times \mathcal{H}^+_{p*}$ with J-holomorphic interior. There are two boundary components, one on the $u > 0$ component of the boundary of $\mathbb{R} \times \mathcal{H}^+_{p*}$ and the other on the $u < 0$ component. The submanifold $C_p$ is diffeomorphic to the complement of $\Delta_p$ interior points of the product of $\mathbb{R}$ with a closed interval. What follows lists two additional properties.

PROPERTY 1: *The large $|s|$ part of $C_p$ is described by (2.9)*

PROPERTY 2: *The submanifold $C_p$ is the $\Psi_p$ image of a graph in the $|\hat{u}| \leq R + \frac{1}{2}\ln z_*$ part of $\mathbb{R} \times \mathcal{H}^+_{p*}$ over a domain in $\mathbb{R} \times I_*$ having the form*

$$(x, \hat{u}) \to (x, \hat{u}, \hat{\phi} = \varphi^p(x,\hat{u}), h = \varsigma^p(x, \hat{u}))$$

*The domain for the functions $(\varphi^p, \varsigma^p)$ is $\mathbb{R} \times I_*$ when $\Delta_p = 0$, it is the complement of a single $\hat{u} = 0$ point when $\Delta_p = 1$, and it is the complement of two $\hat{u} = 0$ points when $\Delta_p = 2$.*

In short, these properties say that $C_p$ looks much like the $\mathbb{R} \times \mathcal{H}^+_{p*}$ part of an ech-HF submanifold.

*Part 3*: The parameter $\tau$ enters the story here. To set the stage, fix $\mathfrak{p} \in \Lambda$ and write the part of $C_S$ in $\mathbb{R} \times \mathcal{H}^+_{p*}$ as in PROPERTY 4 in Part 2 using functions $(\varphi^S, \varsigma^S)$. and write $C_p$ as in (7.3). Meanwhile, reintroduce from (2.5) the functions $(\varphi^{S0}, \varsigma^{S0})$ that are defined by the surface $C_{S0}$. The functions $(\varphi^S, \varsigma^S)$ that define $C_S$ and the pair $(\varphi^p, \varsigma^p)$ that define $C_\mathfrak{p}$ are constrained on the common $t = 1 + z_*$ and $t = 2 - z_*$ boundaries of their domains to obey

- $\varsigma^S - \varsigma^{S0} = \tau(\varsigma^\mathfrak{p} - \varsigma^{S0})$ .
- $\tau(\varphi^S - \varphi^{S0}) = \varphi^\mathfrak{p} - \varphi^{S0}$ .

(7.1)

What follows are two remarks concerning these matching conditions. The first remark concerns the $\tau = 1$ version of (7.1): This version asserts that $C_S$ and $C_\mathfrak{p}$ fit seemlessly together across their common boundary in $\mathbb{R} \times \mathcal{H}^+_{\mathfrak{p}*}$. As a consequence, any given $\tau = 1$ version of $C = C_S \cup (\cup_{\mathfrak{p} \in \Lambda} C_\mathfrak{p})$ is an ech-HF submanifold.

The second remark concerns the $\tau = 0$ case. A set $\mathcal{C}_0 = \{C_{S0}, \{C_{\mathfrak{p}0}\}_{\mathfrak{p} \in \Lambda}\}$ of the sort described in Section 2 obeys all of $\tau = 0$ conditions. Moreover, it follows from Propositions 2.1 and Proposition 2.2 that these are the only sets that obey the $\tau = 0$ conditions.

A set $\mathcal{C} = \{C_S, \{C_\mathfrak{p}\}_{\mathfrak{p} \in \Lambda}\}$ that is described by Parts 1 and 2 above, and obeys a given $\tau \in [0,1]$ version of (7.1) is said to be a $(J, \tau)$-holomorphic submanifold.

**b) The structure of the cobordism space**

Introduce $\mathcal{M}^*$ to denote the set of pairs of the form $(\tau, \mathcal{C})$ with $\tau \in [0,1]$ and with $\mathcal{C}$ being a $(J, \tau)$-holomorphic submanifold. This set is given the topology whereby open neighborhoods of a given element $(\tau, C = \{C_S, \{C_\mathfrak{p}\}_{\mathfrak{p} \in \Lambda}\})$ are generated by sets of the following sort: Fix $\varepsilon > 0$ and a compactly supported 2-form $\upsilon$ on $\mathbb{R} \times Y$. The set in question contains a given $(\tau', C' = \{C'_S, \{C'_\mathfrak{p}\}_{\mathfrak{p} \in \Lambda}\})$ if $|\tau - \tau'| < \varepsilon$ and if the conditions in (1.16) hold with the pair $(C, C')$ replaced by each pair from $\{(C_S, C'_S), \{(C_\mathfrak{p}, C'_\mathfrak{p})\}_{\mathfrak{p} \in \Lambda}\}$.

The map from $\mathcal{M}^*$ to $[0,1]$ defined by the rule $(\tau, \mathcal{C}) \to \tau$ is denoted by $\pi_{\mathfrak{l}}$. A second map, this one from $\mathcal{M}^*$ to a Euclidean space, also enters the story. The latter is denoted by $p$ and its definition follows directly. To start, introduce $\Lambda_* \subset \Lambda$ to denote the subset of $\Delta_\mathfrak{p} \geq 1$ elements. The map $p$ sends $\mathcal{M}^*$ to $\times_{\mathfrak{p} \in \Lambda}(\times_{\Delta_\mathfrak{p}} \mathbb{R})$. To give the rule that defines $p$, write a given element in $\mathcal{M}^*$ as $(\tau, \mathcal{C} = \{C_S, \{C_\mathfrak{p}\}_{\mathfrak{p} \in \Lambda}\})$. Each $\mathfrak{p} \in \Lambda$ version of $C_\mathfrak{p}$ is defined by a pair of functions whose domain is the complement in $\mathbb{R} \times I_*$ of $\Delta_\mathfrak{p}$ with $\hat{u} = 0$. The $\mathbb{R}$ coordinate of these missing $\hat{u} = 0$ points are the $\mathbb{R}$ coordinates of $p(\tau, \mathcal{C})$ in $\mathfrak{p}$'s factor of $\times_{\mathfrak{p} \in \Lambda}(\times_{\Delta_\mathfrak{p}} \mathbb{R})$ with it understood that when $\Delta_\mathfrak{p} = 2$, then the first coordinate in in the corresponding factor $\times_2 \mathbb{R}^2$ corresponds to the end of $C_\mathfrak{p}$ where $\cos(\theta)$ limits to $\frac{1}{\sqrt{3}}$ as $s \to \infty$. The upcoming propositions set $n_* = \sum_{\mathfrak{p} \in \Lambda} \Delta_\mathfrak{p}$ and they $\mathbb{R}^{n_*} = \times_{\mathfrak{p} \in \Lambda}(\times_{\Delta_\mathfrak{p}} \mathbb{R})$.

The following proposition describes the structure $\mathcal{M}^*$.

**Proposition 7.1**: *Fix a Lipshitz submanifold S such that $D_S$ has trivial cokernel. There exists a purely S-dependent constant $\kappa \geq 1$, and with $z_* < \kappa^{-1}$, there exists a constant $\kappa_*$ that depends on $z_*$ but is otherwise purely S-dependent with the following property: Use $\delta < \kappa_*^{-1} z_*$ with a pair $(\hat{\Theta}_-, \hat{\Theta}_+)$ from $\hat{\mathcal{Z}}^S$ to define $\mathcal{M}^*$. Fix $(\tau, C) \in \mathcal{M}^*$.*

- *There exist an integer $n \geq 0$, a neighborhood $U \subset \mathbb{R}^{n_*+n}$ of the origin, an open neighborhood $I \subset [0,1]$ of $\tau$, a smooth map $\mathfrak{f}: I \times U \to \mathbb{R}^n$ that sends $(\tau, 0)$ to the origin, and a topological embedding $\Phi: \mathfrak{f}^{-1}(0) \to \mathcal{M}^*$ onto an open set that sends the pair $(\tau, 0)$ to $C$ and is such that $\pi_I \circ \Phi$ gives the projection from $I \times U$ to $I$.*
- *The subspace of elements $\mathcal{M}^*_{smooth} \subset \mathcal{M}^*$ where $\mathfrak{f}$ is a submersion is open and a smooth $(n_*+1)$-dimensional manifold with boundary. The maps $\pi_I$ and $p$ are smooth on this smooth subset.*
- *The integer $n$ can be taken to be either 0 or 1 at points in $\mathcal{M}^*_{smooth}$; and it can be taken equal to zero at the points in $\mathcal{M}^*_{smooth}$ where $d\pi_I \neq 0$.*
- *An open neighborhood of $\pi_I^{-1}(0)$ in $\mathcal{M}^*_{smooth}$ is mapped by $\pi_I \times p$ diffeomorphically onto an open neighborhood of $\{0\} \times (\times_{p \in \Lambda} (\times_{\Delta_p} \mathbb{R}))$ in $[0,1] \times (\times_{p \in \Lambda} (\times_{\Delta_p} \mathbb{R}))$.*

Proposition 7.1 can be generalized in a straightforward manner to account for variations in the choice of the pair S. The formulation of this more general version is omitted.

The next proposition refers to the notion from Section 2a of a weakly compact subset of Lipshitz submanifolds. The proposition asserts that J, in particular, can be chosen so as to make Proposition 7.1's map $\mathfrak{f}$ everywhere a submersion for a residual set of Lipshitz submanifolds from any given weakly compact subset Lipshitz submanifolds.

To set the stage for the proposition, introduce the notion of a Lipshitz subvariety. The latter is a certain sort of 2-dimensional, $J_{HF}$-holomorphic subvariety in $\mathbb{R} \times [1,2] \times \Sigma$. The definition is identical to that in Section 1g's for a Lipshitz submanifold but for three items. First, the subvariety need not be a submanifold as it is allowed to have a finite number of interior singular points. Second, no irreducible component lies in a constant $(s,t)$ slice of $\mathbb{R} \times [1,2] \times \Sigma$. Third, PROPERTY 8 in Section 1g need not be obeyed. To say more about this last point, note that any given Lipshitz subvariety can be viewed as a pair, $(S, u)$ where S is a smooth complex curve with 2G boundary components and $u$ a $J_{HF}$ holomorphic map from S into $\mathbb{R} \times [1,2] \times \Sigma$ whose image is the subvariety in question. The pair $(S, u)$ is described by the first six bullets in (II.6.2) and the modified version of the seventh bullet of (II.6.2) that requires $u$ to embed the complement of a finite set of interior points. If $u$ is an immersion, there is a holomorphic line bundle over S whose restriction to any given small radius disk is the normal bundle to its $u$-image. In this case an operator $D_S$ that maps sections of the latter to sections of its tensor product with $T^{0,1}S$ which has the form depicted in (1.25). When $u$ is not an embedding, there is an operator that plays the role of $D_S$ and is denoted by $D_S$. This operator is obtained from what is denoted by $D\bar{\partial}$ in the proof of Proposition 3.4 in [L] by restricting the latter to elements of the form $(\xi, Y, 0)$. The latter is Fredholm when viewed as a bounded, linear map from the $L^2_1$ completion of its domain to the $L^2$ completion of its range. Because $u$ is singular at only a finite number of points, all in the interior, this follows directly from what is said in Part 4 of Section II.6e.

**Proposition 7.2**: *Fix a countable set in $(\times_3 (0,1)) \times (1, \infty)$ of possible choices for the data $(z_*, \delta, x_0, R)$ and there is a $C^\infty$-residual set of allowed choices for $J_{HF}$ for which the assertions that follow are true. Choose $J_{HF}$ from this residual set.*
- *Let S denote a Lipshitz subvariety. Then $D_S$ has trivial cokernel.*
- *Let $\mathcal{K}$ denote a given $\mathbb{R}$-invariant, weakly compact subset of Lipshitz submanifolds and there exists a constant $\kappa \geq 1$ that depends only on $\mathcal{K}$ and, given $z_* < \kappa^{-1}$, there exists $\kappa_* > 1$ that depends on $z_*$ and $\mathcal{K}$ with the following significance: Choose a $\mathcal{K}$-compatible data set $(z_*, \delta, x_0, R)$ from the given set with $z_* < \kappa^{-1}$ and $\delta < \kappa_*^{-1} z_*$. Use this data to the geometry of Y.*
    a) *There is a certain residual set of almost complex structures pursuant to the constraints given in Section 1c and there exists a residual subset in $\mathcal{K}$ such that if J is chosen from the former and S from the latter then all $(\hat{\Theta}_-, \hat{\Theta}_+) \in \hat{\mathcal{Z}}^S$ version of $\mathcal{M}^*$ are such that*
    1) *The corresponding $\mathcal{M}^*_{smooth}$ is the whole of $\mathcal{M}^*$ and so $\mathcal{M}^*$ is a smooth $(n_*+1)$-dimensional manifold with boundary; and $\pi_I \times p: \mathcal{M}^* \to [0, 1] \times (\times_{p \in \Lambda} (\times_{\Delta_p} \mathbb{R}))$ is a smooth map.*
    2) *The critical values of $\pi_I$ are in $(0,1)$; and only finitely many of them are critical values of $\pi_I$'s restriction to any given compact set in $\mathcal{M}^*$.*
    b) *If $\mathcal{K}$ is an open set, then the various versions of $\mathcal{M}^*|_{\tau=1}$ as defined by the elements in $\mathcal{K}$ and a given choice for J define a smooth manifold such that the tautological map to $\mathcal{K}$ is smooth.*

Propositions 7.1 and 7.2 are proved in Section 7e.

By way of a parenthetical remark, the tools that are developed in this section can be used to strengthen both the fourth bullet of Proposition 7.1 and Item 2) of the second bullet of Proposition 7.2. With regards to the fourth bullet of Proposition 7.1, the open neighborhood of $\pi^{-1}(0)$ can be taken to be the $\pi_I$-inverse image of an open neighborhood of 0 in [0,1]. The strengthened version of Item 2 of the second bullet of Proposition 7.2 asserts that the set of critical values of $\pi_I$ on the whole of $\mathcal{M}^*$ is finite. The proofs of these strengthened versions do not involve any new technology. Even so, a full presentation is lengthy and so these stronger versions are not proved here.

The final proposition implies that $\pi_{\mathfrak{l}}$ is in all cases a proper map. This proposition refers to the strong $C^\infty$ topology on spaces of sections and maps. This topology is defined as follows: An open neighborhood of a given section or map is indexed by a positive integer and a positive number. Let q denote the given section or map and let $(k, \varepsilon)$ denote a given positive integer and positive number. Elements in the corresponding neighborhood of q have $C^k$ distance less that $\varepsilon$ from q over the whole of q's domain.

The upcoming proposition also refers to the normal bundle and various associated notions for a submanifold that is described by Part 2 of Section 7a, an example being some $\mathfrak{p} \in \Lambda$ version of $C_\mathfrak{p}$. Let C denote the relevant submanifold. The fiber metric on C's normal bundle is defined by the ambient metric on $\mathbb{R} \times \mathcal{H}^+_{\mathfrak{p}*}$ that comes from the chosen almost complex structure J and the compatible 2-form $ds \wedge \hat{a} + w$. The latter metric defines a metric on C and the covariant derivative on sections of the normal bundle and tensor bundles over C. The normal bundle can also be endowed with an exponential map that embeds a constant radius disk subbundle into $\mathbb{R} \times \mathcal{H}^+_{\mathfrak{p}*}$. This exponential map gives the canonical identification between the zero section and C, its differential on the zero section is the identity map, and it maps the disk bundle over the boundary of C to the boundary of $\mathbb{R} \times \mathcal{H}^+_{\mathfrak{p}*}$.

**Proposition 7.3**: *Fix a Lipshitz submanifold* S. *There exists a purely* S-*dependent constant* $\kappa \geq 1$, *and with* $z_* < \kappa^{-1}$, *there exists a constant* $\kappa_*$ *that depends on* $z_*$ *but is otherwise purely* S-*dependent with the following property: Fix* $z_* < \kappa^{-1}$ *and* $\delta < \kappa_*^{-1} z_*$, *and then a pair* $(\hat{\Theta}_-, \hat{\Theta}_+) \in \hat{\mathcal{Z}}^S$ *to define* $\mathcal{M}^*$. *The map* $\mathfrak{p}: \mathcal{M}^* \to \times_{\mathfrak{p} \in \Lambda} (\times_{\Delta_\mathfrak{p}} \mathbb{R})$ *is proper in the following strong sense: Let* $\{(\tau_n, C_n = \{C_{Sn}, \{C_{\mathfrak{p}n}\}_{\mathfrak{p} \in \Lambda}\}, \}_{n=1,2,...} \subset \mathcal{M}^*$ *denote any given sequence with fixed* $\mathfrak{p}$-*image. There is an element* $(\mathcal{C} = \{C_S, \{C_\mathfrak{p}\}_{\mathfrak{p} \in \Lambda}\}, \tau) \in \mathcal{M}^*$ *with the given* $\mathfrak{p}$-*image and a subsequence (hence renumbered consecutively from 1) with the properties listed next.*
- *For each* $n \in \{1, 2, ...\}$, *let* $\eta_n$ *denote the section of the disk bundle* $N_S$ *that defines* $C_{Sn}$. *The resulting sequence of sections* $\{\eta_n\}_{n=1,2,...}$ *converges over the* $t \in [1 + z_*, 2 - z_*]$ *part of* S *in the strong* $C^\infty$ *topology to the section that defines* $C_S$.
- *Fix* $\mathfrak{p} \in \Lambda$. *There exists a sequence of sections* $\{\eta_{\mathfrak{p}n}\}_{n=1,2,...}$ *of the disk subbundle of* $C_\mathfrak{p}$'s *normal bundle that converges to zero in the strong* $C^\infty$-*topology on* $C^\infty(C_\mathfrak{p}; N)$ *and is such that each index* n *version of* $C_{\mathfrak{p}n}$ *is the image of the composition of the exponential map with the corresponding section* $\eta_{\mathfrak{p}n}$.

This proposition is proved in Sections 7d.

### c) Boundary conditions for the Cauchy-Riemann equations

This subsection constitutes a digression to introduce the analytic tools that are needed to handle the boundary matching conditions given by (7.1). By way of background, the pair $(\varphi^S, \varsigma^S)$ that appears in the index 1 critical point version of (7.1) obeys the standard Cauchy-Riemann equations in coordinates $(x, z)$ for $\mathbb{R} \times [z_*, z_S]$ that are defined by the rule $s = x, t = 1 + z$. Meanwhile, the pair $(\varphi^p, \varsigma^p)$ obey these equations on the domain $\mathbb{R} \times [\delta^2, z_*]$ if $z$ is identified with $e^{-2(R-\hat{u})}$. There is are analogous Cauchy-Riemann equations near an index 2 critical point version of (7.1). The Cauchy-Riemann equations here are obeyed by $(\varphi^S, -\varsigma^S)$ using the coordinates $(x, -z)$ for $\mathbb{R} \times [z_*, z_S]$ that are defined by writing $(s, t)$ as $s = x$ and $t = 2 - z$. Mean the pair $(\varphi^p, -\varsigma^p)$ obey the same Cauchy-Riemann equations in terms of coordinates $(x, z)$ on the domain $\mathbb{R} \times [\delta^2, z_*]$ when $z$ is defined by $z = e^{-2(R+\hat{u})}$. These coordinate identifications are used implicitly in the rest of this section and in the subsequent sections. The pair $(\varphi^{S0}, \varsigma^{S0})$ or $(\varphi^{S0}, -\varsigma^{S0})$, as the case may be, obeys the Cauchy-Riemann equations on the domain $\mathbb{R} \times [e^{-8}z_*, z_S]$ and thus on both sides of the $z = z_*$ locus where (7.1) holds. The fact that all of these pairs obey a linear equation near the $z = z_*$ locus explains the focus in this subsection on the coupled, linear boundary value problem that is described next.

The boundary value problem is that for pairs $(\varphi_+, \varsigma_+)$ and $(\varphi_-, \varsigma_-)$ which obey the Cauchy-Riemann equations on the respective domains $\mathbb{R} \times [z_*, z_S]$ and $\mathbb{R} \times [e^{-8}z_*, z_*]$. This is to say that

$$\partial_x \varphi_\pm - \partial_z \varsigma_\pm = 0 \quad and \quad \partial_x \varsigma_\pm + \partial_z \varphi_\pm = 0 .$$

(7.2)

on the relevant domain, and so $\eta_+ = \varphi_+ + i\varsigma_+$ and $\eta_- = \varphi_- + i\varsigma_-$ are holomorphic functions of the complex coordinate $x + iz$. Their boundary values are constrained on the common boundary of their respective domains by a given $\tau \in [0, 1]$ version of

$$\varsigma_+ = \tau \varsigma_- \quad and \quad \tau \varphi_+ = \varphi_- \quad where \ z = z_* .$$

(7.3)

They are also constrained so that

$$\lim_{|x| \to \infty} (|\varsigma_\pm| + |\varphi_\pm|) = 0 .$$

(7.4)

The five parts that follow in this subsection discuss various aspects of this coupled, linear boundary value problem.

*Part 1*: This part of the subsection describes energy bounds that hold for the pairs just described. These are summarized by the next lemma. The pairs of functions that appear in this lemma are not assumed to obey (7.3). The lemma uses the notation $\|\cdot\|$ to denote the $L^2$ norm of a given function with it understood that the integration domain is the domain where the function is defined.

**Lemma 7.4**: *Suppose that $(\varphi_+, \varsigma_+)$ and $(\varphi_-, \varsigma_-)$ are pairs of compactly supported, smooth functions that are defined on the respective domains $\mathbb{R} \times [z_*, z_S)$ and $\mathbb{R} \times (e^{-8}z_*, z_*]$, and that obey (7.3) on the common boundary of their respective domains of definition. Let $\eta_\pm$ denote the $\mathbb{C}$-valued functions $\varphi_\pm + i\varsigma_\pm$. Then $\|\bar{\partial}\eta_+\|^2 + \|\bar{\partial}\eta_-\|^2 = \tfrac{1}{4}\|d\eta_+\|^2 + \tfrac{1}{4}\|d\eta_-\|^2$.*

*Proof of Lemma 7.4*: Integration by parts finds that

$$\|\bar{\partial}\eta_+\|^2 + \|\bar{\partial}\eta_-\|^2 = \tfrac{1}{4}\|d\eta_+\|^2 + \tfrac{1}{4}\|d\eta_-\|^2 + \int_{\mathbb{R}\times z_*} (\varphi_- \partial_x \varsigma_- - \varphi_+ \partial_x \varsigma_+)$$

(7.5)

Use (7.3) to see that the boundary integral is zero.

*Part 2*: This part describes a version of the maximum principle for pairs $(\varphi_\pm, \varsigma_\pm)$ that obey (7.2)–(7.4). Such is the content of the next lemma.

**Lemma 7.5**: *Suppose that $(\varphi_+, \varsigma_+)$ and $(\varphi_-, \varsigma_-)$ are pairs of smooth functions that are defined on the respective domains $\mathbb{R} \times [z_*, z_S]$ and $\mathbb{R} \times [e^{-8}z_*, z_*]$ and obey (7.2)–(7.4). Define functions $\varphi$ and $\varsigma$ on the domain $\mathbb{R} \times [e^{-8}z_*, z_S]$ by the rule:*
- *$\varphi = \varphi_-$ where $z \in [e^{-8}z_*, z_*]$ and $\varphi = \tau\varphi_+$ where $z \in [z_*, z_S]$.*
- *$\varsigma = \tau\varsigma_-$ where $z \in [e^{-8}z_*, z_*]$ and $\varsigma = \varsigma_+$ where $z \in [z_*, z_S]$.*

*If either function is not identically zero, then neither function can have a local maximum or minimum on $\mathbb{R} \times (e^{-8}z_*, z_S)$.*

*Proof of Lemma 7.5*: Consider, for example $\varsigma$. The function $\varsigma_-$ is harmonic as is $\varsigma_+$. Thus, neither can have local extremal values in the interior of their domain of definition. If $\tau = 0$, then the claim follows directly from this last observation because $\varsigma$ is zero where $z \leq z_*$. Suppose next that $\tau > 0$ and suppose that $\varsigma$ takes a local maximum or minimum at a given point $(x_*, z_*)$. The simplest case to consider is that where $\partial_z \varsigma_+ \neq 0$ at this point. Suppose for the sake of argument that $\partial_z \varsigma_+ < 0$. Then the function $z \to \tau\varsigma_-(x_*, z)$ is a decreasing function of $z$ for $z$ near to but slightly less than $z_*$. Thus, $\tau\varsigma_-$ and hence $\varsigma$ will not have a local maximum at $(x_*, z_*)$.

To see about the general case, suppose for the sake of argument that $(x_*, z_*)$ is a local maximum for $\varsigma$. What follows generates some nonsense from this assumption. To set up the notation, introduce $q$ to denote the value of $\varsigma$ at this point. No generality is lost by assuming that $q \neq 0$. Introduce r to denote the Euclidean distance on $\mathbb{R} \times [e^{-8} z_*, z_S]$ from $(x_*, z_*)$. Fix $\varepsilon \in (0, q)$ but very small, chosen in particular so that $q$-$\varepsilon$ is a regular value of both $\varsigma_+$ and $\tau \varsigma_-$. Let U denote the component of the set where $\varsigma \geq q$-$\varepsilon$ that contains the point $(x_*, z_*)$. This set is compact. Introduce $\partial U_+$ and $\partial U_-$ to denote the respective $z \geq z_*$ and $z \leq z_*$ components of the boundary of U. These are smooth arcs. It follows from (7.2) that these arcs are oriented by the respective 1-forms $d\varphi_+$ and $d\varphi_-$. In particular, the respective integrals

$$\int_{\partial U_+} d\varphi_+ \quad and \quad \tau \int_{\partial U_-} d\varphi_-$$

(7.6)

are positive. But this last conclusion is nonsense they sum to zero, a consequence of (7.2) with the fundamental theorem of calculus (Stokes' theorem along an arc).

*Part 3*: This part explains how a bound $|\varsigma_+|$ for $z > \frac{1}{2} z_S$ and on $|\varsigma_-|$ for $z < \frac{1}{2} z_*$ can be used to obtain $\tau$-independent bounds on $\varsigma^p$ at $z = z_*$. Lemma 7.5 supplies such a bound for $|\varsigma_S|$. The following lemma makes a quantative statement.

**Lemma 7.6**: *There exists a $z_*$ and $z_S$ independent constant $\kappa \geq 1$ with the following significance: Suppose that $(\varphi_\pm, \varsigma_\pm)$ are pairs of functions as in Lemma 7.5 that obey (7.2)-(7.4). Fix constants $r_* > 0$ and $r_S > 0$ such that $|\varsigma_+| \leq r_S$ for $z > \frac{1}{2} z_S$ and $|\varsigma_-| \leq r_*$ for $z < \frac{1}{2} z_*$. Then $|\varsigma_-| \leq \kappa (r_* + z_* r_S / z_S)$ where $z = z_*$.*

***Proof of Lemma 7.6:*** Let $\Delta \subset \mathbb{R} \times [e^{-8} z_*, \frac{1}{2} z_*]$ denote a disk of radius $\frac{1}{4} z_*$, and introduce $\Delta'$ to denote the concentric disk with radius $\frac{1}{4} z_*$. Use of the standard Green's function for the Laplacian with a cut-off function that is 1 on $\Delta'$ and zero on the complement of $\Delta$ will find that

$$|\nabla \varsigma_-| \leq c_0 r_* z_*^{-1}$$

(7.7)

at the origin of $\Delta$. It follows as a consequence of (7.2) that $|\nabla \varphi_-| \leq c_0 r_* z_*^{-1}$ where $z \leq \frac{1}{4} z_*$. Much the same argument finds that $|\nabla \varphi_+| \leq c_0 r_S z_S^{-1}$ where $z \geq \frac{1}{4} z_S$. Hold on to these bounds for the moment.

The two pair of functions $(\varphi_\pm', \varsigma_\pm') = (\partial_x \varphi_\pm, \partial_x \varsigma_\pm)$ obey (7.2) and they also obey (7.3). As explained momentarily, Lemma 7.5 can be invoked using the function $\varphi'$ which is defined to be $\partial_x \varphi_-$ where $z \leq z_*$ and $\tau \partial_x \varphi_+$ where $z \geq z_*$. Granted Lemma 7.5, it follows from the conclusions of the preceding paragraph that $\partial_x\varphi$ at $z = z_*$ can have absolute value no greater that $c_0(r_* z_*^{-1} + r_S z_S^{-1})$. This with the Cauchy-Riemann equations imply that $|\partial_z \varsigma_-| \leq c_0(r_* z_*^{-1} + r_S z_S^{-1})$ at $z = z_*$.

Let w denote the function on $[\frac{1}{4} z_*, z_*]$ given by the rule

$$z \to w(z) = r + c_0(r_* z_*^{-1} + r_S z_S^{-1}) z. \tag{7.8}$$

View w as an x-independent function on $\mathbb{R} \times [\frac{1}{4} z_*, z_*]$. As such, it is harmonic, and the maximum principle implies that $w(z) \geq \varsigma_-(z)$ on its domain of definition. This the case, it it follows that $\varsigma_- \leq w$ at all $z \in [\frac{1}{4} z_*, z_*]$ and so $c_0 r$ at $z = z_*$. This gives the asserted upper bound for $\varsigma_-$ where $z = z_*$. The exact same argument with $-\varsigma_-$ replacing $\varsigma_-$ gives the asserted upper bound at $z = z_*$ for $-\varsigma_-$.

Return now to the assertion that Lemma 7.5 can be invoked using $\varphi'$. There is no issue if it is known apriori that $|\partial_x \varsigma_\pm|$ and $|\partial_x \varphi_\pm|$ limit uniformly to zero as $|x| \to \infty$. If this has not been established, the argument proceeds as follows: Let $\beta: \mathbb{R} \to [0, 1]$ denote a smooth function with compact support with integral equal to 1. Given $\varepsilon > 0$, introduce $\beta_\varepsilon: \mathbb{R} \to [0, 1]$ to denote the function given by the rule $x \to \varepsilon^{-1}\beta(\varepsilon^{-1}x)$. Consider now the pair $(\varphi_{\varepsilon\pm}, \varsigma_{\varepsilon\pm})$ given by the mollifying formula

$$(\varphi_{\varepsilon\pm}, \varsigma_{\varepsilon\pm})|_{(x,z)} = \int_\mathbb{R} \beta_\varepsilon(x')(\varphi_\pm, \varsigma_\pm)|_{(x'+x,z)} \, dx'. \tag{7.9}$$

Every $\varepsilon > 0$ version of $(\varphi_{\varepsilon\pm}, \varsigma_{\varepsilon\pm})$ obeys (7.2)–(7.4). This is also the case for their partial derivatives to any given order with respect to x. This understood, define for $\varepsilon > 0$ the mollified function $\varphi_\varepsilon'$ given by $\partial_x \varphi_{\varepsilon-}$ where $z \leq z_*$ and $\tau \partial_x \varphi_{\varepsilon+}$ where $z \geq z_*$. Lemma 7.5 holds for this function. Meanwhile the family $\{\varphi_\varepsilon'\}_{\varepsilon>0}$ converges uniformly as $\varepsilon \to 0$ on compact subsets of the domain $\mathbb{R} \times [e^{-8} z_*, z_S]$ to the function $\partial_x \varphi$. This understood, take the $\varepsilon \to 0$ limit of $|\partial_x \varphi_\varepsilon|$ to see that the assertion of Lemma 7.5 holds for $\partial_x \varphi$.

*Part 4*: This part says something about apriori estimates near the $z = z_*$ locus for solutions to (7.2)–(7.4).

**Lemma 7.7**: *Fix $k \in \{1, 2, \ldots\}$ and there exists a $z_*$ and $z_S$ independent constant $\kappa \geq 1$ with the following significance: Suppose that $(\varphi_\pm, \varsigma_\pm)$ are smooth functions that obey (7.2)–(7.4). Fix $r_* > 0$ and $r_S > 0$ such that $|\varsigma_+| \leq r_S$ for $z > \frac{1}{2} z_S$ and $|\varsigma_-| \leq r_*$ for $z < \frac{1}{2} z_*$. Then the norms of the derivatives of $(\varphi_-, \varsigma_-)$ to order k where $z > \frac{1}{4} z_*$ and those of $(\varphi_+, \varsigma_+)$ to order k where $z < \frac{1}{4} z_S$ are bounded by $\kappa(r_* z_*^{-k} + r_S z_S^{-k})$.*

***Proof of Lemma 7.7***: The argument given in the proof of Lemma 7.6 bounds $|\partial_x \varphi_-|$ where $z \in [\frac{1}{4} z_*, z_*]$ by $c_0 (r_* z_*^{-1} + r_S z_S^{-1})$, and an analogous argument bounds $|\partial_x \varsigma_+|$ where $z \in [z_*, \frac{1}{4} z_S]$. Much the same argument bounds $|\partial_x^k \varphi_-|$ by $c_k (r_* z_*^{-k} + r_S z_S^{-k})$ for $z \in [\frac{1}{4} z_*, z_*]$ and $|\partial_x^k \varsigma_+|$ by this same constant where $z \in [z_*, \frac{1}{4} z_S]$. An iterative bootsrapping argument uses the bound for $|\partial_x^k \varphi_-|$ to obtain the desired bound on $|\partial_x^{k-1} \varsigma_-|$; and it uses the bound for $|\partial_x^k \varsigma_+|$ to obtain the desired bound on $|\partial_x^{k-1} \varphi_+|$. What follows describes how this works for $|\partial_x^{k-1} \varsigma_-|$. The argument for $|\partial_x^{k-1} \varphi_+|$ is identical but for the notation.

Note first that $|\partial_x^{k-1} \varsigma_-|$ is apriori bounded where $z = \frac{1}{4} z_*$ by $c_{k-1}(r_* z_*^{-k+1} + r_S z_S^{-k+1})$. This the case, note next that $\partial_z (\partial_x^{k-1} \varsigma_-) = \partial_x^k \varphi_-$ because of (7.2). The assumed bound on $\partial_x^k \varphi_-$ implies that $|\partial_z \partial_x^{k-1} \varsigma_-| \le c_k (r_* z_*^{-k} + r_S z_S^{-k})$ where $z = z_*$. Granted that such is the case, use the fact that $\partial_x^{k-1} \varsigma_-$ is harmonic to invoke the maximum principle for the harmonic function

$$(x, z) \to (\partial_x^{k-1} \varsigma_-)|_{(x,z)} - c_{k-1}(r_* z_*^{-k+1} + r_S z_S^{-k+1}) - c_k (r_* z_*^{-k} + r_S z_S^{-k}) z .$$
(7.10)

on $\mathbb{R} \times [\frac{1}{4} z_*, z_*]$ to obtain the desired upper bound on $\partial_x^{k-1} \varsigma_-$. The desired lower bound is obtained by this argument by replacing $\varsigma_-$ with $-\varsigma_-$.

Bounds on the partial derivatives to order $k$ in both the variables $x$ and $z$ are obtained via the Cauchy-Riemann equations from those for just the partial derivatives with respect to $x$.

### d) Proof of Proposition 7.3

The argument for the proposition when all $\mathfrak{p} \in \Lambda$ versions of $\Delta_\mathfrak{p}$ are zero is very much like that given in Section 4c for the proof of Lemma 4.3. This version of the argument is given in the first six parts of this subsection. Part 7 adds what is needed to prove the proposition when some $\mathfrak{p} \in \Lambda$ versions of $\Delta_\mathfrak{p}$ are 1 or 2.

*Part 1*: This part derives an upper bound for the integral of the 2-form $w$ over the J-holomorphic submanifolds with boundary from any given element in $\mathcal{M}^*$. The lemma below states such a bound.

**Lemma 7.8**: *There exists a purely* S-*dependent* $\kappa \ge 1$ *with the following significance: Define the space* $\mathcal{M}^*$ *using* $z_* \le \kappa^{-1}$, $\delta < \kappa^{-2} z_*$ *and a* $\{\Delta_\mathfrak{p} = 0\}_{\mathfrak{p} \in \Lambda}$ *pair from* $\hat{\mathcal{Z}}^S$. *Suppose that* $\tau \in [0, 1]$ *and* $C = \{C_S, \{C_\mathfrak{p}\}_{\mathfrak{p} \in \Lambda}\}$ *is a* (J, $\tau$)-*holomorphic submanifold. Let* $I \subset \mathbb{R}$ *denote an interval of length* 1. *Then*

$$\int_{C_S} w + \sum_{\mathfrak{p}\in\Lambda} \int_{C_\mathfrak{p}} w \le \kappa \quad \text{and} \quad \int_{C_S\cap(I\times M_\delta)} ds\wedge d\hat{a} + \sum_{\mathfrak{p}\in\Lambda} \int_{C_\mathfrak{p}\cap(I\times(\mathcal{H}^+_{\mathfrak{p}*}\cap M_\delta))} ds\wedge d\hat{a} \le \kappa.$$

*Meanwhile,* $\sum_{\mathfrak{p}\in\Lambda} \int_{C_\mathfrak{p}\cap(I\times\mathcal{H}^+_{\mathfrak{p}*})} ds\wedge d\hat{a} \le \kappa_*$ *where* $\kappa_*$ *is a* $\mathcal{C}$ *and* $\tau$ *independent constant.*

*Proof of Lemma 7.8*: The argument is much like that used for Lemma 4.5. To start, remark that what are essentially cosmetic modifications to the arguments used at the beginning of the Lemma 4.5's proof can be used to prove that the integrals in question are finite, and that Stokes' theorem in various guises can be used to compute them.

With the preceding in mind, use (1.6) to write the 2-form $w$ as $w = w' + \sum_{\mathfrak{p}\in\Lambda} d\mathfrak{b}_\mathfrak{p}$ where any given $\mathfrak{p}\in\Lambda$ version of the 1-form $\mathfrak{b}_\mathfrak{p}$ has compact support in $\mathcal{H}^+_\mathfrak{p}$ and has the form $\mathfrak{b}_\mathfrak{p} = x(1-3\cos^2\theta)\,du - \mathrm{N}\sqrt{6}\,f\cos\theta\sin^2\theta\,d\phi$ with $\mathrm{N}$ being the function of $u$ given by the rule $u \to \mathrm{N}(u) = \chi(|u|-R+\ln\delta_*)$. The form $w'$ is zero on each $\mathfrak{p}\in\Lambda$ version of $C_\mathfrak{p}$, and its support on $C_S$ is disjoint from the boundary of $C_S$. Given that $w'$ is closed, Stokes' theorem with what is said in Corollary II.2.6 can be used to see that its integral over $C_S$ differs by no more than $c_0\delta$ from its integral over $S$.

Fix $\mathfrak{p}\in\Lambda$ and use Stokes' theorem to write the integral of $d\mathfrak{b}_\mathfrak{p}$ as the sum of the integrals given in (4.9) as defined using $(\varphi,\varsigma) = (\varphi^\mathfrak{p},\varsigma^\mathfrak{p})$. As explained in the proof of Lemma 4.5, the sum of left most two terms in the top bullet of (4.9) and the two terms in the lower bullet of (4.9) are bounded by a purely S-dependent (or $\mathcal{K}$-compatible) constant. The as yet unspoken for terms in the top bullet of (4.9) are

$$\int_{\mathbb{R}\times\{\hat{u}=R+\frac{1}{2}\ln z_*\}} \varsigma^\mathfrak{p}\,d\varphi^\mathfrak{p} - \int_{\mathbb{R}\times\{\hat{u}=-R-\frac{1}{2}\ln z_*\}} \varsigma^\mathfrak{p}\,d\varphi^\mathfrak{p}.$$

(7.11)

Meanwhile, integration by parts identifies the integral of $\mathfrak{b}_\mathfrak{p}$ over $C_S$ with the sum

$$\int_{\gamma'_{\mathfrak{p}+}} \mathrm{N}\hat{h}\,d\phi - \int_{\gamma'_{\mathfrak{p}-}} \mathrm{N}\hat{h}\,d\phi - \int_{\mathbb{R}\times\{\hat{u}=R+\frac{1}{2}\ln z_*\}} \varsigma^S\,d\varphi^S + \int_{\mathbb{R}\times\{\hat{u}=-R-\frac{1}{2}\ln z_*\}} \varsigma^S\,d\varphi^S.$$

(7.12)

where $\gamma'_{\mathfrak{p}+}$ and $\gamma'_{\mathfrak{p}-}$ are the parts of the integral curves of $v$ that extend $\gamma_{\mathfrak{p}+}$ and $\gamma_{\mathfrak{p}-}$, and lie in the $f \ge 1+z_*$ and $f < 2-z_*$ parts of the radius $4\delta_*$ coordinate balls centered on the index 1 and 2 critical points from $\mathfrak{p}$. The left most two integrals in (7.12) are zero as $\phi$ is constant in the integral curve segments in question.

To say more about the integrals in (7.11) and the two rightmost integrals in (7.12) note first that the corresponding $\hat{u} = R+\frac{1}{2}\ln z_*$ integrals in (7.11) and in (7.12) come with opposite signs. This is also the case for the corresponding $\hat{u} = -R-\frac{1}{2}\ln z_*$ integrals. With the preceding understood, use (7.1) to identify the sum of what is in (7.11) and (7.12) with

$$-(1-\tau)\int_{\mathbb{R}\times\{\hat{u}=R+\frac{1}{2}\ln z_*\}}(\varsigma^{S0}d\varphi^S - \varsigma^p d\varsigma^{S0}) + (1-\tau)\int_{\mathbb{R}\times\{\hat{u}=-R-\frac{1}{2}\ln z_*\}}(\varsigma^{S0}d\varphi^S - \varsigma^p d\varphi^{S0})$$

(7.13)

The integrals in (7.13) of the 1-form $\varsigma^p d\varphi^{S0}$ are bounded by $c\delta_*^2$ with $c$ here denoting a purely S-dependent (or $\mathcal{K}$-compatible) constant. This is because $|\varsigma^p| \le \frac{4}{\sqrt{3}}\delta_*^2$ and because of what is said by the fourth bullet in Proposition 2.1. The integrals that involve $\varsigma^{S0} d\varphi^S$ are written using Stokes' theorem as the sum of two integrals with support in the $z \ge z_*$ part of the radius $4\delta_*$ coordinate balls about $\mathfrak{p}$'s critical points. The term with support in the index 1 critical point coordinate ball is

$$-(1-\tau)\int_{z\ge z_*}(\varsigma^{S0} dN \wedge d\varphi^S + N d\varsigma^{S0} \wedge d\varphi^S) \; ;$$

(7.14)

and the term with support in the index 2 critical point coordinate ball has the same form but no minus sign in front. To bound these integrals, it is important to keep in mind two facts. First, the form $w$ on $\wedge^2 TC_S$ written in the radius $4\delta_*$ coordinate ball in terms of the pair $(\varphi^S, \varsigma^S)$ is given by $\sqrt{6}(\partial_x \varphi^S \partial_v \varsigma^S - \partial_v \varphi^S \partial_x \varsigma^S) dx \wedge dv$. Second, the pair $(\varphi^S, \varsigma^S)$ obeys the Cauchy-Riemann equations as functions of $(x, v)$. These two facts imply that $w$ on $\wedge^2 TC_S$ on this part of $C_S$ is $\sqrt{6}(|\partial_x \varphi^S|^2 + |\partial_{\hat{u}} \varphi^S|^2) dx \wedge dv$. This understood, use the triangle inequality with the fourth bullet of Proposition 2.1 to see that what is written in (7.14) is no greater than $\frac{1}{1000}\int_{C_S} w + c$ where $c$ is another purely S-dependent (or $\mathcal{K}$-compatible) constant. This last bound with those derived previously imply that

$$\int_{C_S} w + \sum_{\mathfrak{p}\in\Lambda}\int_{C_\mathfrak{p}} w \le \frac{1}{1000}\int_{C_S} w + c$$

(7.15)

with $c$ being purely S-dependent (or $\mathcal{K}$-compatible). The inequality in (7.15) implies what asserted by the first bullet of Lemma 7.8.

Turn now to the integrals of $ds \wedge \hat{a}$. With the $f \in (1,2)$ part of M written as $(1,2)\times \Sigma$, the form $\hat{a}$ is $dt$ with $t$ the Euclidean coordinate on the $(1,2)$ factor. This being the case, the integral of $ds \wedge \hat{a}$ over $C_S \cap (I \times M_\delta)$ is no greater than $G$ and the integral of $ds \wedge \hat{a}$ over the $I \times M_\delta$ part of any given $\mathfrak{p} \in \Lambda$ version of $C_\mathfrak{p}$ is no greater than $z_*$. Minor cosmetic changes to the arguments from Step 4 of the proof of Proposition II.5.1 in Section II.5b give the bound on the integral of $ds \wedge \hat{a}$ over the whole of $C_\mathfrak{p} \cap (I \times \mathcal{H}^+_{\mathfrak{p}*})$.

Note that the integration by parts that is used in this Step 4 does not lead to boundary terms because $\hat{a}$ annihilates the tangent space to the boundary of $\mathcal{H}^+_{\mathfrak{p}*}$.

*Part 2*: This part supplies an upper bound for the distance in $\mathbb{R} \times [1, 2] \times \Sigma$ from any given $\mathcal{M}^*$ version of $C_S$ and S. Here again, the $f \in (1, 2)$ part of $M_\delta$ is identified with $(1, 2) \times \Sigma$. This part of section also supplies a positive bound for the function $1 - 3\cos^2\theta$ on any given $\mathcal{M}^*$ and $\mathfrak{p} \in \Lambda$ version of $C_\mathfrak{p}$. These bounds are summarized in the respective lemmas that follow.

**Lemma 7.9**: *There is a purely S-dependent (or $\mathcal{K}$-compatible) $\kappa > 1$ and given $z_* < \kappa^{-1}$, there exists $\kappa_* > 1$ that depends only on $z_*$ but is otherwise purely S-dependent (or $\mathcal{K}$ compatible) with the following significance: Define the space $\mathcal{M}^*$ using $z_* < \kappa^{-1}, \delta < \kappa^{-2} z_*$ and a $\{\Delta_\mathfrak{p} = 0\}_{\mathfrak{p} \in \Lambda}$ pair from $\hat{\mathcal{Z}}^S$. Suppose that $(\tau, C = \{C_S, \{C_\mathfrak{p}\}_{\mathfrak{p} \in \Lambda}\}) \in \mathcal{M}^*$. Then $C_S$ is the image via the exponential map $\mathfrak{e}_S$ of a section of the the radius $\frac{1}{100} \kappa_S^{-1} \rho_S^2$ subbundle of $N_0$.*

The next lemma concerns the function $1 - 3\cos^2\theta$ on the various $\mathfrak{p} \in \Lambda$ versions of $C_\mathfrak{p}$.

**Lemma 7.10**: *There is a purely S-dependent (or $\mathcal{K}$-compatible) $\kappa > 1$ and given $z_* < \kappa^{-1}$, there exists $\kappa_* > 1$ that depends only on $z_*$ but is otherwise purely S-dependent (or $\mathcal{K}$ compatible) with the following significance: Define $\mathcal{M}^*$ using $z_* < \kappa^{-1}$ and $\delta < \kappa_*^{-1} z_*$ and a $\{\Delta_\mathfrak{p} = 0\}_{\mathfrak{p} \in \Lambda}$ pair from $\hat{\mathcal{Z}}^S$. There exists $\kappa_{**} > 1$ such that if $(\tau, C = \{C_S, \{C_\mathfrak{p}\}_{\mathfrak{p} \in \Lambda}\}) \in \mathcal{M}^*$ is a given element, then $1 - 3\cos^2\theta$ is greater than $\kappa_{**}^{-1}$ on all $\mathfrak{p} \in \Lambda$ versions of $C_\mathfrak{p}$.*

The proofs of Lemmas 7.9 and 7.10 require a preliminary lemma to supply an apriori bound on the norm of any $\mathfrak{p} \in \Lambda$ version of $|\varsigma^S|$ on the common boundary of $C_S$ and $C_\mathfrak{p}$. The proof of Lemma 7.10 requires in addition a bound for any $\mathfrak{p} \in \Lambda$ version of $|\varsigma^\mathfrak{p}|$ where $z = e^{64} \delta^2$. This is the content of the upcoming Lemma 7.11. Lemma 7.11 refers to the coordinate z that is defined on either component of the $|\hat{u}| \geq R + \ln\delta$ portion of $\mathcal{H}^+_{\mathfrak{p}*}$ by the rule $z = e^{-2(R-|\hat{u}|)}$.

**Lemma 7.11**: *There is a purely S-dependent (or $\mathcal{K}$-compatible) $\kappa > 1$, and given $z_* < \kappa^{-1}$ and $\varepsilon > 0$, there exists $\kappa_\varepsilon > 1$ that depends on $z_*$ and $\varepsilon$ but is otherwise purely S-dependent (or $\mathcal{K}$-compatible) such that the following is true: Define the space $\mathcal{M}^*$ using $z_* < \kappa^{-1}, \delta < \kappa_\varepsilon^{-1} z_*$ and a $\{\Delta_\mathfrak{p} = 0\}_{\mathfrak{p} \in \Lambda}$ pair from $\hat{\mathcal{Z}}^S$. Let $(\tau, C = \{C_S, \{C_\mathfrak{p}\}_{\mathfrak{p} \in \Lambda}\}) \in \mathcal{M}^*$, and let $\mathfrak{p} \in \Lambda$.*

- *The pair $(\varphi^S, \varsigma^S)$ is such that $|\varsigma^S| \le \varepsilon \delta_*^2$ where $z = z_*$, this the common boundary of $C_S$ and $C_p$.*
- *The pair $(\varphi^p, \varsigma^p)$ is such that $|\varsigma^p| \le \varepsilon \delta_*^2$ where $z = e^{64} \delta^2$.*

This lemma is proved in Part 4. Given the lemma, what follows in this Part 2 is the proof of Lemma 7.9. Part 3 contains the proof of Lemma 7.10.

*Proof of Lemma 7.9*: The proof has five steps.

Step 1: Write $t$ on $[1, 1+z_S]$ as $1+z$ with $z \in [0, z_S]$. Likewise write $t$ as $t = 2-z$ on $[2-z_S, 2]$ with $z$ again in $[0, z_S]$. Let $c_S > 1$ denote the version of the constant $\kappa$ that is supplied by Lemma 7.11. Given $\varepsilon > 0$ and $z_* < \min(c_S^{-1} z_S, \varepsilon \delta_*^2)$, choose $\delta$ so that the conclusions of Lemma 7.11 can be invoked. Having made such choices, fix $\mathfrak{p} \in \Lambda$. The function $\varsigma^{S0}$ that appears in $\mathfrak{p}$'s version of (7.1) has norm bounded by $c z_*$ where $c \ge 1$ is purely S-dependent (or $\mathcal{K}$-compatible); this follows from Proposition 2.1. It follows from what is said in Corollary II.2.6 that the $|x| \to \infty$ limits of $|\varsigma^S|$ are bounded by $c\delta$. As a consequence, these limits are bounded by $c\varepsilon \delta_*^2$. Here again, $c$ is purely S-dependent (or $\mathcal{K}$-compatible).

Step 2: The pair $(\varphi^S, \varsigma^S)$ is a solution to the Cauchy-Riemann equations on the domain $\mathbb{R} \times [z_*, z_S]$. As a consequence, $\varsigma^S$ is a harmonic functions on this domain. This is to say that it is annihilated by the Laplacian $(\partial_x^2 + \partial_z^2)$. With the preceding understood, fix $c > 1$ and introduce the function

$$z \to (\tfrac{(z - z_*)}{(z_S - z_*)}) (\rho_S^2 - c\varepsilon \delta_*^2) + c\varepsilon \delta_*^2 \,.$$

(7.16)

This is a harmonic function on $\mathbb{R} \times [z_*, z_S]$; and if $c \ge c$ with $c$ purely S-dependent (or $\mathcal{K}$-compatible), then this function is greater than $|\varsigma^S|$ on the boundaries and at large $|x|$. Granted that such is the case, then the maximum principle demands that this function be greater than $|\varsigma^S|$ on the whole of $\mathbb{R} \times [z_*, z_S]$. Choose $c$ in (7.16) so as to be purely S-dependent (or $\mathcal{K}$-compatible) and so that this last conclusion holds.

Step 3: Suppose now that there is no $\kappa$ as claimed by Lemma 9.1 so as to derive some nonsense. Granted this assumption, there is a sequence $\{(D_n, \eta_n)\}_{n=1,2,\ldots}$ of the following sort: First, $D_n$ is a data set with elements $((\hat{\Theta}_{n-}, \hat{\Theta}_{n+}), z_{*n}, \delta_n, x_{0n}, R_n)$ that is suitable for defining the geometry of $Y$ and $\mathcal{M}^*$, and is such that the conclusions of the top bullet in Lemma 7.11 can be invoked with $\varepsilon = \tfrac{1}{n}$ and $z_{*n} \in (0, \tfrac{1}{n} z_S)$ and $\delta_n \le \tfrac{1}{n^2} z_{*n}$.

Second, $\eta_n$ is a section of the bundle $N_S$ over the $t \in [1+z_{*n}, 2-z_{*n}]$ part of S with the properties given in the list below.

- $\eta_n$ is $L^2$ orthogonal to the $t \in [1+z_{*n}, 2-z_{*n}]$ restriction of the elements in the kernel of the operator $D_S$.
- The norm of $|\eta_n|$ at some point in its domain is greater than $\frac{1}{100}\kappa_S^{-1}\rho_S^2$.
- The $|s| \to \infty$ limit of $|\eta_n|$ is bounded by $\frac{1}{n}\delta_*^2$.
- The composition $\mathfrak{e}_S \circ \eta_n$ has J-holomorphic image.
- The absolute value of the $\eta_n$ analog of the function $\varsigma^S$ is bounded where $z \in [z_{*n}, z_S]$ by the $z = z_{*n}$ and $\varepsilon = \frac{1}{n}$ version of the function in (7.16).

(7.17)

As noted in Part 2 of Section II.6e, the third fourth bullet in (7.17) implies that $\eta = \eta_n$ obeys an equation with the schematic form

$$\bar{\partial}\eta + \mathfrak{r}_1(\eta)\cdot\partial\eta + \mathfrak{r}_0(\eta) = 0$$

(7.18)

where the notation is that in (II.6.10). By way of a reminder, $\bar{\partial}$ signifies the d-bar operator on sections of $N_S$ as defined using the hermitian metric to give the bundle a holomorphic structure. What is written as $\partial$ is the adjoint operator. Meanwhile, the map $\mathfrak{r}_1: N_0 \to N_S \otimes \mathrm{Hom}(T^{1,0}S; T^{0,1}S)$ and the map $\mathfrak{r}_0: N_0 \to N_S \otimes T^{0,1}S$ are smooth, fiber preserving maps that vanish along the zero section.

Step 4: The properties listed in Step 3 together with (7.18) imply via standard elliptic regularity arguments that there is a subsequence of $\{\eta_n\}_{n=1,2,...}$, hence renumbered consecutively from 1, that is described either by CASE 1, CASE 2, or CASE 3 given below.

CASE 1: T*he subsequence converges uniformly on compact domains in S to a non-trivial section of $N_0$ over the whole of S with the following properties:*
  a) *It obeys (7.18) and the conditions in the second and third bullets of (II.6.12).*
  b) *It's norm is no greater than $\kappa_S^{-1}\rho_S^2$*
  c) *It is $L^2$-orthogonal to* kernel($D_S$)

The first two conditions above imply that the section is described by Lemma II.6.6. This understood, the third condition implies that the section is identically zero. Given that the limit section is asserted to be non-trivial, CASE 1 can not describe $\{\eta_n\}_{n=1,2,...}$.

Case 2 below uses terminology from Section 1g and Sections II.6c and II.6e.

CASE 2: *There is a negative point $q \in \underline{S}-S$ and a sequence $\{q_n\}_{n=1,2,...} \subset S_0$ of points that converges in $\underline{S}$ to q such that $|\eta_n|(q_n) > \frac{1}{100}\kappa_S^{-1}\rho_S^2$. Write a neighborhood of q in*

S *as in Part 4 of Section 6e. This done, view* $\eta_n$ *on this part of S as a* $\mathbb{C}$*-valued function as done in Section II.6e. There is an unbounded, decreasing sequence* $\{x_n\} \in (-\infty, x_0]$ *such that the translated sequence with n'th member* $\eta'_n$ *given by* $\eta'_n|_x = \eta_n|_{x-x_n}$ *converges on compact domains of* $\mathbb{R} \times [1, 2]$ *to a non-trivial function with the following properties:*

a) *It is holomorphic.*
b) *It's imaginary part vanishes on* $\mathbb{R} \times \{1\}$.
c) *It's real part vanishes on* $\mathbb{R} \times \{2\}$.
d) *It is bounded.*

Properties a)-c) plus (II.6.18) are incompatible with Property d). This being the case, the subsequence that is given at the end of Step 4 is not described CASE 2.

CASE 3 is the analog of CASE 2 where the point q is a positive point of S−S. The analogous conclusion applies: The subsequence that is given at the end of Step 4 is not described by CASE 3 either.

Step 5: Neither CASE 1 nor CASE 2 nor CASE 3 describe the subsequence given by Step 4. This contradicts what is said at the end of Step 4. The contradiction is avoided if and only if Lemma 7.9 is true.

*Part 3*: This part contains the

***Proof of Lemma 7.10***: The proof has five steps.

Step 1: Let $c_S > 1$ denote the larger of the versions of constant $\kappa$ that are supplied by Lemmas 7.11 and 7.9. Given $\varepsilon > 0$ and $z_* < \max(c_S^{-1} z_S, \varepsilon^2 \delta_*^2)$, choose $\delta < (c_S + \kappa_{\varepsilon^2 \delta_*^2})^{-1} z_*$ so as to invoke the conclusions of Lemmas 7.11 and 7.9 using $\varepsilon^2 \delta_*^2$ in lieu of $\varepsilon$. Having made such a choice, suppose that there existst $\varepsilon_* \in (0, 1)$ with the following property:

*Let* $(\tau, \mathcal{C} = \{C_S, \{C_p\}_{p \in \Lambda}\})$ *denote an element from the resulting version of* $\mathcal{M}^*$. *Fix* $p \in \Lambda$. *Then* $|\varsigma^p| \leq (1 - \varepsilon_*) \frac{4}{3\sqrt{3}} \delta_*^2$ *where* $z = z_*$.

(7.19)

Granted (7.19), it then follows from (4.14) that $1 - 3\cos^2\theta > c_0^{-1} \varepsilon_*$ on the boundary of any $(\{C_S, \{C_p\}_{p \in \Lambda}\}, \tau) \in \mathcal{M}^*$ and $p \in \Lambda$ version of $C_p$. With a bound of this sort in hand, a repeat of the arguments for Lemma 4.7 using Lemma 7.8 in lieu of Lemma 4.5 proves Lemma 7.10. The steps that follow give an existence proof for a suitable $\varepsilon_*$.

What with the top line in (7.1) and the top bullet of Lemma 7.11, the assertion made by (7.19) holds automatically for $\tau > c_0^{-1}\varepsilon^2$. This the case, only small values of $\tau$ are of any concern. Even so, no upper bound for $\tau$ is assumed in the remaining steps.

<u>Step 2</u>: Fix $r \in (0, 1)$ and let $D \subset \mathbb{R} \times ((1+r)z_*, z_S)$ denote a disk of radius $rz_*$. The function $|\varsigma^S|$ is bounded by what is written in (7.16); and because $\varsigma^S$ is harmonic on D, this implies that the norm of $d\varsigma^S$ at the center of D is no greater than $c_0 (r z_*)^{-1} \varepsilon^2 \delta_*^4$. This with the Cauchy-Riemann equations implies that $|d\varphi^S|$ is bounded where $z = (1+r)z_*$ by $c_0 (r z_*)^{-1} \varepsilon^2 \delta_*^4$. Meanwhile $|\varsigma^p| < \frac{4}{3\sqrt{3}} \delta_*^2$. As $\varsigma^p$ is harmonic on the domain $\mathbb{R} \times [\delta^2, z_*]$, it follows that $|d\varsigma^p|$ and thus $|d\varphi^p|$ are bounded by $c_0 z_*^{-1} \delta_*^2$ where $z = \frac{1}{2} z_*$.

<u>Step 3</u>: Define a function $\varphi'$ on $\mathbb{R} \times [e^{-8} z_*, z_S]$ as follows: Set $\varphi' = \partial_x(\varphi^p - \varphi^{S0})$ on the $z \leq z_*$ part of this domain, and set $\varphi' = \tau \partial_x(\varphi^S - \varphi^{S0})$ on the $z \geq z_*$ part. Define next a function w of the coordinate z on $\mathbb{R} \times [\frac{1}{2} z_*, (1+r)z_*]$ by the rule

$$z \to w(z) = z_*^{-1}(r^{-1} \varepsilon^2 \delta_*^4 + z_*^{-1} \delta_*^2 (z - (1+r)z_*)).$$
(7.20)

The function w is harmonic and there exists a constant $c_w \leq c_0$ such that $\varphi' - c_w w$ is less than zero where $z = \frac{1}{2} z_*$ and $z = (1+r)z_*$. As explained in the next step, the maximum principle can be used to infer that $\varphi' < c_w w$ at all points in $\mathbb{R} \times [\frac{1}{2} z_*, (1+r)z_*]$. Granted this, it then follows using Lemma 7.1 that

$$|\partial_x \varphi^p| \leq c_0 z_*^{-1} (r^{-1} \varepsilon^2 \delta_*^4 + r) \quad \text{where } z = z_*.$$
(7.21)

Thus, $|\partial_z \varsigma^p| \leq c_0 z_*^{-1} (r^{-1} \varepsilon^2 \delta_*^4 + r)$ where $z = z_*$. Taking $r = \varepsilon \delta_*^2$ finds $|\partial_x \varsigma^p| \leq c_0 z_*^{-1} \varepsilon \delta^2$.

<u>Step 4</u>: What follows considers the case when both $\lim_{|x| \to \infty} |\partial_x(\varphi^p - \varphi^{S0})|$ and $\lim_{|x| \to \infty} |\partial_x(\varphi^S - \varphi^{S0})|$ are zero. The general case is handled using mollifiers as done in the proof of Lemma 7.6. Granted this assumption about the $|x| \to \infty$ limits, it follows that the function $\varphi' - c_w w$ is negative at large $|x|$. This function is harmonic where $z \neq z_*$ and so it lacks local maxima and local minima in $\mathbb{R} \times (\frac{1}{2} z_*, z_*)$ and in $\mathbb{R} \times (z_*, (1+r)z_*)$.

To see about local maxima or minima where $z = z_*$, note that the function on $\mathbb{R}$ given by

$$x \to q(x) = z_*^{-2} \delta_*^2 x$$
(7.22)

is a conjugate harmonic function for w when viewed as a z-independent function on the domain $\mathbb{R} \times [\frac{1}{2} z_*, (1+r) z_*]$. This is to say that the pair (q, w) obey the Cauchy-Riemann equations. This understood, define the function $\varsigma'$ on $\mathbb{R} \times [\frac{1}{2} z_*, (1+r) z_*]$ by setting it to equal $\tau(\varsigma^p - c_w \tau^{-1} q)$ where $z \leq z_*$, and to equal $\varsigma^S - c_w q$ where $z \geq z_*$. With the pair $(\varphi', \varsigma')$ in hand, repeat the argument in the paragraph preceding (7.6) with $\varphi'$ playing the role of $\varsigma$ and $\varsigma'$ the role of $-\varphi'$ to rule out local extreme points for $\varphi' - c_w w$.

An analogous argument with the signs of w and q reversed rules out local extreme points for $\varphi' + c_w w$.

Step 5: The lower bullet of Lemma 7.11 asserts that $|\varsigma^p| \leq \varepsilon \delta_*^2$ where $z = e^{64} \delta^2$. Meanwhile, the Step 4 finds that $|\partial_z \varsigma^p| \leq c_0 z_*^{-1} \varepsilon \delta_*^2$ where $z = z_*$. With this in mind, let w now denote the x-independent function on $\mathbb{R} \times [2\delta^2, z_*]$ given by

$$z \to w(z) = z_* + z_*^{-1} z.$$

(7.23)

There exists $c_0 \geq 1$ such that the function $\varsigma^p - c_0 \varepsilon \delta_*^2 w$ on $\mathbb{R} \times [2\delta^2, z_*]$ has the following properties: It is negative where $|x| \gg 1$ and where $z = e^{64} \delta^2$. Meanwhile, its z-derivative is negative where $z = z_*$. Granted these facts, a version of the maximum priniciple implies that $\varsigma^p < c_0 \varepsilon \delta_*^2$ where $z = z_*$. The analogous argument using $-\varsigma^p$ in lieu of $\varsigma^p$ proves that $\varsigma^p > -c_0 \varepsilon \delta_*^2$ where $z = z_*$. A suitable choice of $\varepsilon$ establishes what is asserted by (7.19).

*Part 4*: This part contains the

***Proof of Lemma 7.11***: The proof has four steps. The first three steps prove the assertion made by the lower bullet in the lemma

Step 1: To prove the assertion made by the lower bullet, assume the contrary to generate some nonsense. If lower bullet is false, there is a sequence $\{(D_n, (\mathcal{C}_n, \tau_n)\}_{n=1,2,...}$ of the following sort: First, $D_n$ is a set with elements $((\hat{\Theta}_{n-}, \hat{\Theta}_{n+}), z_{*n}, \delta_n, x_{0n}, R_n, J_n)$ that are suitable for defining the geometry of Y and the corresponding version of $\mathcal{M}^*$. This data is such that the pair $(\hat{\Theta}_{n-}, \hat{\Theta}_{n+})$ is an $\{\Delta_p = 0\}_{p \in \Lambda}$ element in $\hat{\mathcal{Z}}^S$. The constant $z_{*n} < \frac{1}{n} z_S$ and the constant $\delta_n < \frac{1}{n} z_{*n}$. Meanwhile, $(\mathcal{C}_n = \{C_{Sn}, \{C_{pn}\}\}, \tau_n)$ is an element in the corresponding version of $\mathcal{M}^*$. In addition, there exists $\mathfrak{p} \in \Lambda$ such that the $C_{pn}$ version of the function $\varsigma^p$ has absolute value greater than $\varepsilon \delta_*^2$ at some point where $z = e^{64} \delta_n^2$.

For each $n \in \{1, 2, ...\}$, fix a point in $C_{pn}$ where $z = e^{64} \delta_n^2$ and where $|\varsigma^p| \geq \varepsilon \delta_*^2$. Let $q_n$ denote a given such point. The point $q_n$ projects to the $\delta = \delta_n$ version of $M_\delta$ and the

image lies in the radius $2\delta_*$ coordinate ball about one or the other of the critical points from $\mathfrak{p}$. It follows from (4.14) that all sufficiently large n versions of $q_n$ project in $M_\delta$ so as to have coordinate radius greater than $\frac{1}{25}\varepsilon^{1/2}\delta_*$. The image of $q_n$ also lies where the function $f$ differs from either 1 or 2 by less than $\frac{1}{n}z_{*n}$.

For each n, let $U_n$ denote the following part of the union of the respective radius $2\delta_*$ coordinate balls centered on the index 1 and index 2 critical points from $\mathfrak{p}$: It is the part where the radius is greater than $\frac{1}{100}\varepsilon^{1/2}\delta_*$ and where the function either $|f-1|$ or $|f-2|$ is less than $\frac{1}{2}z_{*n}$. If n is large, $C_{\mathfrak{p}n} \cap (\mathbb{R} \times U_n)$ is a non-empty, J-holomorphic submanifold; it is non-empty as it contains the point $q_n$.

Pass to a subsequence so that after renumbering consecutively from 1, the projection to $M_\delta$ of the resulting set of points $\{q_n\}_{n=1,2,\ldots}$ lies very near one critical point of $\mathfrak{p}$. For each n, let $C_n \subset (C_{\mathfrak{p}n} \cap (\mathbb{R} \times U_n))$ denote the component that contains $q_n$. Translate $C_n$ by a constant amount along the $\mathbb{R}$ factor in $\mathbb{R} \times U_n$ so that the point $q_n$ is moved to where the $\mathbb{R}$ coordinate is 0. Let $\{C_n'\}_{n=1,2,\ldots}$ denote the resulting, translated sequence of submanifolds.

Step 2: This step reviews some relevant geometry. To start, let $p \in \mathfrak{p}$ denote the critical point that is described at the end of the last step. The critical point p labels an irreducible component of the locus $c_{1+} \cup c_{2-}$ in $\Sigma$. Use c to denote this irreducible component. As noted in Section 1c, the circle c has an annular neighborhood, T, with the following property: The identification between $f^{-1}(1,2) \subset M$ and $(1,2) \times \Sigma$ identifies the $\mathbb{C}$-valued 1-form $d\phi + id\hat{h}$ on the part of $f^{-1}(1,2) \cap M_\delta$ in the radius $2\delta_*$ coordinate ball centered on p with a 1-form on an annular neighborhood of c in T. The latter 1-form is holomorphic of type (1, 0).

Let $\hat{U}$ denote the part of the radius $\frac{3}{2}\delta_*$ coordinate ball centered on p where the radius is greater than $\frac{1}{100}\varepsilon^{1/2}\delta_*$ and where either $|f-1| < e^{-100}\delta_*^2$ or $|f-2| < e^{-100}\delta_*^2$ as the case may be. The map from $\hat{U}$ to $\Sigma$ that is given by the flow along the integral curves of $v$ identifies $\hat{U}$ with $(-e^{-100}\delta_*^2, e^{-100}\delta_*^2) \times U$ where $U \subset T$ can be written as $U_1 - U_2$ where $U_1$ is an annular neighborhood of c and where $U_2 \subset U_1$ is a smaller width annular neighborhood of c.

Let z denote the Euclidean coordinate on $(-e^{-100}\delta_*^2, e^{-100}\delta_*^2)$. The identification just described extends so as to identify $\mathbb{R} \times \hat{U}$ with $\mathbb{R} \times (-e^{-100}\delta_*^2, e^{-100}\delta_*^2) \times U$. Let x denote the Euclidean coordinate on the $\mathbb{R}$ factor. The identification just described identifies $T^{1,0}(\mathbb{R} \times \hat{U})$ with the span of the pair $(dx + idz, d\phi + id\hat{h})$. This integrable complex structure is observedly compatible with the symplectic form $dx \wedge dz + d\phi \wedge d\hat{h}$.

Step 3:   Each large n version of $U_n$ lies in $\hat{U}$; it appears in the coordinates just defined as the subset where $z \in (-\frac{1}{2} z_{*n}, \frac{1}{2} z_{*n})$. Given such large n, set $\alpha_n = 2 z_{*n}^{-1} e^{-100} \delta_*^2$ and define a diffeomorphism $\psi \colon \mathbb{R} \times U_n \to \mathbb{R} \times \hat{U}$ by the rule

$$\psi_n(x, z, \phi, \hat{h}) = (\alpha_n x, \alpha_n z, \phi, \hat{h}).$$

(7.24)

This diffeomorphism is J-holomorphic.

It follows as a consequence that $\psi_n(C_n')$ is a properly embedded, J-holomorphic submanifold in $\mathbb{R} \times \hat{U}$. This submanifold sits entirely where $1 - 3\cos^2\theta > 0$, this being the locus where $z = 0$. Even so the point $\psi_n(q_n)$ has distance bounded by $c_0 \frac{1}{n}$ from this locus. Moreover, $\psi_n(q_n)$ has distance at least $c_0^{-1}$ from the boundary of the closure of $\mathbb{R} \times \hat{U}$ inside $\mathbb{R} \times Y$. Note as well that $\psi_n(q_n)$ sits on the $x = 0$ locus. Given that the 2-form $w$ on $\hat{U}$ appears as $-2\sqrt{6} \, d\hat{h} \wedge d\phi$, it follows from Lemma 7.8 that $\psi_n^* w = w$. As a consequence, if $I \subset \mathbb{R}$ is any unit length interval, there is an I and n-independent bound on the integral of $w$ over $\psi_n(C_n') \cap (I \times \hat{U})$. Meanwhile, the integral of $dx \wedge dz$ over $\psi_n(C_n') \cap (I \times \hat{U})$ is no greater than $2e^{-100}\delta_*^2$. Indeed, this follows from the fact that the corresponding $C_{pn}$ is given as a graph over the $|\hat{u}| \leq R + \frac{1}{2} \ln z_*$ part of $\mathbb{R} \times \mathcal{H}_{p+}$. To elaborate, $dx \wedge dz$ on $I \times \hat{U}$ is the pull-back of via the projection map of its name sake on $I \times (-e^{-100}\delta_*^2, e^{-100}\delta_*^2)$. Meanwhile, the projection to this product restricts to $\psi_n(C_n') \cap (I \times \hat{U})$ to define a 1-1 map into $I \times (-e^{-100}\delta_*^2, e^{-100}\delta_*^2)$ because $C_{pn}$ is a graph. This implies that the integral of the 2-form $dx \wedge dt$ over $\psi_n(C_n') \cap (I \times \hat{U})$ can not be greater than its integral over the whole of $I \times (-e^{-100}\delta_*^2, e^{-100}\delta_*^2)$, which is $2e^{-100}\delta_*^2$.

Granted these last observations, invoke Proposition II.5.5 using the sequence $\{\psi_n(C_n')\}_{n \gg 1}$ to obtain a subsequence that converges on compact subsets of $\mathbb{R} \times \hat{U}$ in the manner dictated by Proposition II.5. to a weighted J-holomorphic subvariety in $\mathbb{R} \times \hat{U}$. Let $\vartheta$ denote this subvariety. Given what is said about $\psi_n(q_n)$, the set $\vartheta$ must contain a pair whose subvariety component is the $x = 0$, $z = 0$ locus in $\mathbb{R} \times \hat{U}$.

The latter conclusion constitutes the required nonsense because the existence of such a pair in $\vartheta$ has the same implications as its existence in the analogous version of $\vartheta$ given by Step 1:  There is a circle in each large n version of $C_{pn}$ whose image via the projection to $\mathcal{H}^+_{p*}$ defines a non-zero generator of the latter's first homology. This nonsense proves what is asserted by the lower bullet of Lemma 7.11.

Step 4:   With $\varepsilon$ chosen, fix $z_*$ and $\delta$ to invoke the lower bullet of Lemma 7.12 as a guarantee that $|\varsigma^p| < \frac{1}{4} \varepsilon \delta_*^2$ where $z = 2\delta^2$. Require in addition that $z_* < \frac{1}{4} \varepsilon z_S \delta_*^2$ and that $\delta$ is chosen less that $c_0^{-1} e^{-100} \varepsilon \delta_*^2$. Define the x-independent function $w$ on $\mathbb{R} \times [2\delta^2, z_S]$ by the rule

$$z \to \tfrac{1}{2} \varepsilon \delta_*^2 + z\, z_S^{-1} \rho_S^2 .$$

(7.25)

Meanwhile, let $\varsigma$ denote the function on $\mathbb{R} \times [2\delta^2, z_S]$ that is given by $\tau \varsigma^p$ - w where $z \le z_*$ and given by $\varsigma^S + (1 - \tau) \varsigma^{S0}$ - w where $z > z_*$. It follows from (7.1) that this function is continuous across the $z = z_*$ locus. Meanwhile, the lower bullet of Lemma 7.11 and Property 1 in Part 1 of Section 7a with Corollary II.6 imply that $\varsigma \le 0$ where $z = e^{64}\delta^2$, where $z = z_S$ and where $|x| \gg 1$. As $\varsigma$ is a harmonic function where $z \ne z_*$, it has no local maxima where $z \ne z_*$. As explained momentarily, it has no local maxima where $z = z_*$. Granted that such is the case, it follows that $\varsigma^S \le \varepsilon \delta_*^2$ where $z = z_*$.

To see why $\varsigma$ has no local maxima where $z = z_*$, define the function $\varphi$ on the domain $\mathbb{R} \times [2\delta^2, z_S]$ as follows: Set $\varphi$ equal to $\varphi^p$ - $x z_S^{-1} \rho_S^2$ where $z \le z_*$ and set $\varphi$ to equal $\tau \varphi^S + (1-\tau)\varphi^{S0}$ - $x z_S^{-1} \rho_S^2$ where $z \ge z_*$. Note that the pair $(x z_S^{-1} \rho_S^2, w)$ obey the Cauchy-Riemann equations. This understood, repeat the argument in the paragraph preceding (7.6) using the just defined version of $(\varphi, \varsigma)$ to rule out local maxima for $\varsigma$ on the $z = z_*$ locus.

To prove that $\varsigma^S \ge -\varepsilon \delta_*^2$ where $z = z_*$, repeat the preceding argument with the sign of w reversed and with $x z_S \rho_S^2$ added to $\varphi^p$ and to $\tau \varphi^S + (1-\tau)\varphi^{S0}$ rather than subtracted when defining $\varphi$.

*Part 5*: This part proves that $\mathcal{M}^*$ is compact in a topology that is slightly weaker than the one defined at the outset of Section 7b. This topology is defined as follows: The open neighborhoods of a given element $(\tau, \mathcal{C} = \{C_S, \{C_\mathfrak{p}\}_{\mathfrak{p} \in \Lambda}\})$ are generated by sets that are labeled by data of the form $(\varepsilon, \upsilon, V, I)$ where $\varepsilon$ is a positive number, $\upsilon$ is a smooth, compactly supported 2-form on $\mathbb{R} \times Y$ and $V \subset \mathbb{R} \times Y$ is an open set with compact closure. Meanwhile, $I \subset [0, 1]$ is an open neighborhood of $\tau$. The corresponding open set in $\mathcal{M}^*$ consists of pairs of the form $(\tau', \mathcal{C}' = \{C_S', \{C_\mathfrak{p}'\}_{\mathfrak{p} \in \Lambda}\})$ with $\tau' \in I$ and with $\mathcal{C}'$ obeying

- $\sup_{z \in C \cap V} \text{dist}(z, C' \cap V) + \sup_{z \in C' \cap V} \text{dist}(z, C \cap V) < \varepsilon.$
- $|\int_C \upsilon - \int_{C'} \upsilon| < \varepsilon.$

(7.26)

Part 6 to come proves convergence in Section 7b's topology on $\mathcal{M}^*$ and convergence in the strong $C^\infty$ topology as asserted by the two bullets of Proposition 7.3.

The arguments that follow assume that the data that defines the geometry of Y and then $\mathcal{M}^*$ is such that the conclusion of Lemma 7.11 with a given, small choice for $\varepsilon$. The arguments require $\varepsilon$ to be less than an S-independent, positive number. The subsequent arguments also assume that the defining data is such that the conclusions of Lemmas 7.9 and 7.10 hold.

To set up the arguments, suppose that $\{(\tau_n, \mathcal{C}_n = \{C_{Sn}, \{C_{\mathfrak{p}n}\}_{\mathfrak{p} \in \Lambda}\})\}_{n=1,2,\ldots} \in \mathcal{M}^*$ is a given sequence. No generality is lost by taking a sequence for which the corresponding sequence $\{\tau_n\}_{n=1,2,\ldots}$ converges. Use $\tau \in [0, 1]$ to denote the limit. The six steps that follow prove that this sequence has a subsequence that converges to some $(J, \tau)$-holomorphic submanifold.

<u>Step 1</u>: Fix $\mathfrak{p} \in \Lambda$ and fix attention on one or the other of the $|\hat{u}| \geq R + \frac{1}{2} \ln z_* - 8$ parts of $\mathcal{H}^+_{\mathfrak{p}*}$. Use the coordinates $(x, z, \phi, \hat{h})$ here. For each $n \in \{1, 2, \ldots\}$, write the relevant part of $C_{Sn}$ as a graph using functions $(\varphi^S_n, \varsigma^S_n)$. Meanwhile, write the relevant part of $C_{\mathfrak{p}n}$ as a graph using functions $(\varphi^\mathfrak{p}_n, \varsigma^\mathfrak{p}_n)$.

Invoke Lemma 7.7 to conclude the following: There sequence $\{(\tau_n, \mathcal{C}_n)\}_{n=1,2,\ldots}$ has a subsequence (hence renumbered consecutively from 1) with the following property: The sequence of pairs $\{(\varphi^S_n, \varsigma^S_n)\}_{n=1,2,\ldots}$ and the sequence of pairs $\{(\varphi^\mathfrak{p}_n, \varsigma^\mathfrak{p}_n)\}_{n=1,2,\ldots}$ converge in the $C^\infty$ topology on compact subsets of $\mathbb{R} \times [\frac{1}{4} z_*, \frac{1}{4} z_S]$ to respective pairs denoted by $(\varphi^S, \varsigma^S)$ and $(\varphi^\mathfrak{p}, \varsigma^\mathfrak{p})$. The latter obey (7.1).

In addition, the surface defined where $z \in [z_*, \frac{1}{4} z_S]$ in $\mathbb{R} \times \mathcal{H}_\mathfrak{p}$ by the graph of the pair $(\varphi^S, \varsigma^S)$ when viewed as a surface in $\mathbb{R} \times (1, 2) \times \Sigma$ lies in the radius $\kappa_S^{-1} \rho_S^2$ tubular neighborhood of S. Meanwhile, the $|\hat{u}| \in [R + \frac{1}{2} \ln z_* - \ln 2, R + \frac{1}{2} \ln z_*]$ part of the graph of the pair $(\varphi^\mathfrak{p}, \varsigma^\mathfrak{p})$ defines via $\Psi_\mathfrak{p}$ a surface in portion of $\mathbb{R} \times \mathcal{H}^+_{\mathfrak{p}*}$.

<u>Step 2</u>: Fix $\mathfrak{p} \in \Lambda$. For each $n \in \{1, 2, \ldots\}$, use $(\varphi^\mathfrak{p}_n, \varsigma^\mathfrak{p}_n)$ to denote the pair that defined $C_{\mathfrak{p}n}$ as a graph over $\mathbb{R} \times I_*$. It follows from Lemma 7.10 that the corresponding positive integer sequence of coefficient functions that appear in the various $(\varphi, \varsigma) \in \{(\varphi^\mathfrak{p}_n, \varsigma^\mathfrak{p}_n)\}_{n=1,2,\ldots}$ versions of (3.4) are obtained from the restriction of their eponymous brethren on $\mathcal{X}$ to a compact set. The partial derivatives of the latter set of functions to any given order have uniformly bounded absolute values on such a compact set.

Granted these last observations, and granted the conclusions from Step 1, standard elliptic regularity arguments in Chapter 6 of [M] can be applied to find a subsequence of $\{(\varphi^\mathfrak{p}_n, \varsigma^\mathfrak{p}_n)\}_{n=1,2,\ldots}$, hence renumbered consecutively from 1, that converges in the $C^\infty$ topology on compact subsets of $\mathbb{R} \times I_*$ to a pair, $(\varphi^\mathfrak{p}, \varsigma^\mathfrak{p})$, that obeys (3.4). Given Lemma 7.10, it follows that the graph of this pair is in the domain of the map $\Psi_\mathfrak{p}$ and so defines a J-holomorphic surface $\mathbb{R} \times \mathcal{H}^+_{\mathfrak{p}*}$. Let $C_\mathfrak{p}$ denote the latter surface.

<u>Step 3</u>: Each n ∈ {1, 2, …} version of $C_{Sn}$ is defined by a section, $\eta_n$, of the bundle $N_S$ over the $t \in [1+z_*, 2-z_*]$ part of S. As noted by Lemma 7.9, each $C_{Sn}$ is the image via $\mathfrak{e}_S$ of a section, $\eta_n$, of the radius $\frac{1}{100} \kappa_S^{-1} \rho_S^2$ subbundle in $N_0$.

The maps $\mathfrak{r}_0$ and $\mathfrak{r}_1$ that appear in (7.18) enjoy uniform bounds on their derivatives to any given order on the $\frac{1}{50} \kappa_S^{-1} \rho_S^2$ subbundle in $N_0$. As each $\eta \in \{\eta_n\}_{n=1,2,…}$ obeys (7.18), it follows using what is said in Step 1 with the aforementioned elliptic regularity arguments from Chapter 6 in [M] that the sequence $\{\eta_n\}_{n=1,2,…}$ has a subsequence (now renumbered consecutively from 1) that converges in $C^\infty(S; N_0)$ on compact subsets of the $t \in [1 + \frac{1}{2} z_*, 2 - \frac{1}{2} z_*]$ part of S to a smooth section of $N_0$ with norm no greater than $\frac{1}{100} \kappa_S^{-1} \rho_S^2$. Let $\eta$ denote this section.

It follows from Lemma 7.9 that the composition $\mathfrak{e}_S \circ \eta$ defines a J-holomorphic surface in $\mathbb{R} \times [1 + \frac{1}{2} \ln z_*, 2 - \frac{1}{2} \ln z_*]$ that lies in the radius $\kappa_S^{-1} \rho_S$ tubular neighborhood of S. Let $C_S$ denote this surface. Meanwhile, Step 1 has the following additional implication: The set $\mathcal{C} = \{C_S, \{C_p\}_{p \in \Lambda}\}$ obey the parameter $\tau$ version of the matching conditions given by (7.1).

<u>Step 4</u>: It follows from Steps 1-3 that $(\tau, \mathcal{C})$ satisfies all of the requirements for membership in $\mathcal{M}^*$ except perhaps for the conditions on the $|s| \to \infty$ limits of the various surfaces that comprise $\mathcal{C}$. This step with Steps 5 and 6 prove that this last requirement is met. The assertion that any given sequence in $\mathcal{M}^*$ has a subsequence that limits in the manner described by Steps 1-3 to an element in $\mathcal{M}^*$ verifies the claim that $\mathcal{M}^*$ is compact in the topology that is defined by (7.26).

To start the story on the large $|s|$ behavior, let $\eta$ again denote the section of $N_0$ over the $t \in [1+z_*, 2-z_*]$ part of S that defines $C_S$. The manner of convergence of the sequence $\{\eta_n\}_{n=1,2,…}$ to $\eta$, and the fact that $C_S$ and each submanifold from $\{C_{Sn}\}_{n=1,2,…}$ are J-holomorphic implies that the integral of $w$ over $C_S$ is finite, and this limit is no larger than the lim-sup of the sequence whose n'th component is the integral of $w$ over $C_{Sn}$. As noted previously, Lemma 7.8 provides an upper bound for this lim-sup.

What with Lemma II.5.6, these last observations imply that each very large $|s|$ slice of $C_S$ is very close to a union of curves in $[1+z_*, 2-z_*] \times \Sigma$ whose members projects to points in $\Sigma$. Moreover, given that $\eta$ lies in the radius $\frac{1}{100} \kappa_S^{-1} \rho_S^2$ subbundle of $N_0$, the following must be true: If $s \gg 1$, then the corresponding set of components of the large $s$ slice of $\Sigma$ has genus($\Sigma$) elements, and this set enjoys a 1-1 correspondence with the points in $C_{1+} \cap C_{2-}$ that define the $s \to \infty$ limit of S. This is such that each component of a given constant $s$ slice of $C_S$ lies in the radius $c_0 \sqrt{\varepsilon} \delta_*$ disk neighborhood of the corresponding point in $C_{1+} \cap C_{2-}$. Indeed, this last bound on the radius follows from Lemma 7.11 because the latter implies that any integral curve segment in question has endpoints in the radius $\sqrt{\varepsilon}\delta_*$ coordinate ball centered on some index 1 and index 2 critical point of $f$. Use the correspondence just described to label the components of the large, constant $s$ slices of $C|_s$. Introduce now $d(s)$ to denote the diameter of the image in $\Sigma$ via projection from $[1+z_*, 2-z_*] \times \Sigma$ of the component of $C_S|_s$ with a given label. Then $\lim_{s\to\infty} d(s) = 0$.

There is the analogous story for the constant, $s \ll -1$ slices. Each has genus($\Sigma$) components, and the components enjoy an analogous 1-1 correspondence with the points in $\Sigma$ that define the $s \to -\infty$ limit of S.

<u>Step 5</u>: Fix $\mathfrak{p} \in \Lambda$. This step considers the large $s$ slices of $C_\mathfrak{p}$. Lemma 7.8 supplies an upper bound for the sequence whose n'th member is the integral of $w$ over $C_{\mathfrak{p}n}$. The fact that the sequence $\{(\varphi^\mathfrak{p}_n, \varsigma^\mathfrak{p}_n)\}_{n=1,2,\ldots}$ converges smoothly on compact domains in $\mathbb{R} \times I_*$ to $(\varphi^\mathfrak{p}, \varsigma^\mathfrak{p})$ implies that the integral of $w$ over $C_\mathfrak{p}$ is finite. As a consequence, each very large $|s|$ slice of $C_\mathfrak{p}$ must be every close to the segment of some integral curve of $v$ in the part of $\mathcal{H}^+_{\mathfrak{p}*}$ where $\hat{u} \in I_*$. When viewed via $\Psi_\mathfrak{p}$ in terms of the functions $(\varphi^\mathfrak{p}, \varsigma^\mathfrak{p})$, this means the following: If $|x|$ is very large, then the function $\hat{u} \to \varsigma^\mathfrak{p}(x, \hat{u})$ on $I_*$ is nearly constant. In particular, this implies the following: Given $r > 0$, and then given $s$ sufficiently large, there exists a segment of an integral curve of $v$ in the $\hat{u} \in I_*$ part of $\mathcal{H}^\mathfrak{p}_+$ with the following property: Let $\gamma$ denote the segment. Then

$$\sup_{q \in C^\mathfrak{p}|_s} \text{dist}(q, \gamma) + \sup_{q \in \gamma} \text{dist}(C^\mathfrak{p}|_s, \gamma) < r .$$

(7.27)

Fix a critical point in $\mathfrak{p}$ so as to consider the corresponding version of (7.1). Fix some $x \in \mathbb{R}$ with $x \gg 1$ so that the corresponding constant $s = s(x)$ slices of $C_S$ and $C_\mathfrak{p}$ are very near respective segments of integral curves of $v$. Given that $|\varphi^S - \varphi^{S0}| < c_0 \sqrt{\varepsilon}\delta_*$ at the given value of x on the common boundary of $C_S$ and $C_\mathfrak{p}$, the lower bullet in (7.1) asserts that $|\varphi^\mathfrak{p} - \varphi^{S0}| \leq c_0 \sqrt{\varepsilon}\delta_*$ at the given value of x on the boundary of $I_*$. Let $\gamma^x \subset \mathcal{H}^+_{\mathfrak{p}*}$ denote the $\hat{u} \in I_*$ part of an integral curve of $v$ that lies very close to the given constant $s(x)$ slice of $C_\mathfrak{p}$. The values of the angle $\phi$ on the respective boundary points of $\gamma^x$ and $\gamma_{\mathfrak{p}+}$ differ by at most $c_0 \sqrt{\varepsilon}\delta_*$. This has the following consequence: Let $\Delta\phi^x$ denote the change in the angle $\phi$ along $\gamma^x$ and let $\Delta\phi^+$ denote the corresponding angle change along $\gamma_{\mathfrak{p}+}$. Then there exists $m \in \mathbb{Z}$ such that

$$|\Delta\phi^x - \Delta\phi^+ + 2\pi m| \leq c_0 \sqrt{\varepsilon}\delta_* .$$

(7.28)

To see that $m = 0$ in (7.28) if $\varepsilon < c_0^{-1}$, remark that the same constant $s$ slice of any given large n version of $C_{\mathfrak{p}n}$ must be every where close to $\gamma^x$ also. This understood, define a closed curve in $C_{\mathfrak{p}n}$ as follows: This constant x slice of $C_{\mathfrak{p}n}$ with the $x' > x$ parts of the boundary of $C_{pn}$ and the arc $\gamma_{p+}$ concatenate to define a closed loop in $\mathcal{H}^+_{p*}$. This loop is null-homologous; it is homologous to the image via the projection from $\mathbb{R} \times \mathcal{H}^+_{p*}$ of a loop in $C_{pn}$ and the latter space is contractible. Such a loop is null-homologous if and only if the integer m that appears in (7.28) is zero when $\varepsilon < c_0^{-1}$. Note in particular that this bound on $\varepsilon$ is purely S-dependent.

As usual, there is an analogous description of the $s << -1$ part of $C_p$.

Step 6: Let $\mathcal{L} \subset [0, 1)$ denote the following set: A number $D \in [0, 1)$ lies in $\mathcal{L}$ if the following is true: Fix any $\mathfrak{p} \in \Lambda$. Then for any given, but sufficiently large $s \in \mathbb{R}$, the slice $C_p|_s$ has distance less than D from $\gamma_{p+}$.

What follows momentarily proves that $0 \in \mathcal{L}$. There is an analogous, $s \to -\infty$ version of $\mathcal{L}$, and the analogous argument proves that the latter also contains 0. These last facts with (7.1) imply that the constant $s$ slices of $C_S$ and each $\mathfrak{p} \in \Lambda$ version of $C_p$ converge as $s \to \pm\infty$ to the segments of integral curves of $v$ that are defined by $(\hat{\Theta}_-, \hat{\Theta}_+)$. Given what is said by Steps 1-3, this means that $(\tau, \mathcal{C}) \in \mathcal{M}^*$.

To prove that $0 \in \mathcal{L}$, it is enough to prove the following: If $D \in \mathcal{L}$, then so is $\frac{1}{2}D$. To see that such is the case, let $T_D \subset \mathcal{H}^+_{p*}$ denote the $\hat{u} \in I_*$ part of the radius D disk neighborhood of $\gamma_{p+}$. Let $\gamma \subset T_{2D}$ denote an integral curve of $v$. It follows as a consequence of (II.2.4) and (4.12) that

$$|h(\gamma_{p+}) - h(\gamma)| < c_0 x_0 R^{-1} D .$$

(7.29)

This last bound implies that all sufficiently large x values of $\varsigma^p$ are such that

$$|\varsigma^p|_x - h(\gamma_{p+})| \le c_0 x_0 R^{-1} D .$$

(7.30)

Given (7.1), this last point implies that $|\varsigma^S|_x - h(\gamma_{p+})| \le c_0 x_0 R^{-1} D$ for all sufficiently large x. Use this with (II.2.6) to deduce that any given $\mathfrak{p} \in \Lambda$ version of $|(\varphi^S - \varphi^{S0})|_x|$ is bounded by $c_0 x_0 R^{-1} D$ for all sufficiently large x on the common boundary of $C_S$ and $C_p$. What with (7.2), this implies that

$$|(\varphi^p - \varphi^{S0})|_x| \le c_0 x_0 R^{-1} D \quad \text{for all large } x.$$

(7.31)

The preceding inequality has the following implication: If $s$ is sufficiently large, then there is a segment of an integral curve of $v$ in the $\hat{u} \in I_*$ part of $\mathcal{H}^p_+$ that is described by some $r \le c_0 x_0 R^{-1} D$ version of (7.27). This implies imparticular that $\frac{1}{2} D \in \mathcal{L}$.

As usual, the $s \to -\infty$ limit of the constant $s$ slices of $C|_s$ behave in an analogous fashion.

*Part 6*: This step completes the proof of the $\{\Delta_{\mathfrak{p}} = 0\}_{\mathfrak{p} \in \Lambda}$ case of Proposition 7.3 by verifying that the sequence $\{(\tau_n, \mathcal{C}_n = \{C_{Sn}, \{C_{\mathfrak{p}n}\}_{\mathfrak{p} \in \Lambda}\})\}_{n=1,2,...}$ has a subsequence that converges in the appropriate manner.

The upcoming Lemma 7.12 plays a central role in this argument. Lemma 7.12 assumes that the data chosen to define the geometry of Y and $\mathcal{M}^*$ is such that the conclusions of Lemma 7.11 holds with a constant $\varepsilon$ chosen less than a certain purely S-dependent constant. The data is also such that the conclusions of Lemmas 7.9 and 7.10 hold.

Let $(\Theta_-, \Theta_+) \in \mathcal{Z}_{\text{ech},M} \times \mathcal{Z}_{\text{ech},M}$ denote the pair that lies under $\{\hat{\Theta}_-, \hat{\Theta}_+\}$. Since $\Delta_{\mathfrak{p}}$ is zero for all $\mathfrak{p} \in \Lambda$, both lack integral curves of $v$ from $\{\hat{\gamma}_{\mathfrak{p}}^+, \hat{\gamma}_{\mathfrak{p}}^-\}_{\mathfrak{p} \in \Lambda}$. Lemma 7.12 and the subsequent arguments view the element $\Theta_-$ as a set of segments of integral curves of $v$ written as $\{\Gamma_{S-}, \{\gamma_{\mathfrak{p}-}\}_{\mathfrak{p} \in \Lambda}\}$ where $\Gamma_{S-}$ denotes the union of the parts of the curves in $\Theta_-$ where $f \in [1 + z_*, 2 - z_*]$. The set $\Theta_+$ is likewise written as $\{\Gamma_{S+}, \{\gamma_{\mathfrak{p}+}\}_{\mathfrak{p} \in \Lambda}\}$.

**Lemma 7.12**: *Given $r > 0$, there exists $s_r > 1$ with the following significance: Suppose that $(\tau, \mathcal{C} = \{C_S, \{C_{\mathfrak{p}}\}_{\mathfrak{p} \in \Lambda}\})$ is a given element from $\mathcal{M}^*$. If $s < -s_r$, then*

$$\sup\nolimits_{q \in (C_S \cup (\cup_{\mathfrak{p} \in \Lambda} C_{\mathfrak{p}}))|_s} \text{dist}(q, (\Gamma_{S-} \cup (\cup_{\mathfrak{p} \in \Lambda} \gamma_{\mathfrak{p}-}))) + \sup\nolimits_{q \in (\Gamma_{S-} \cup (\cup_{\mathfrak{p} \in \Lambda} \gamma_{\mathfrak{p}-}))} (C_S \cup (\cup_{\mathfrak{p} \in \Lambda} C_{\mathfrak{p}}))|_s q) < r.$$

*If $s > s_r$, then*

$$\sup\nolimits_{q \in (C_S \cup (\cup_{\mathfrak{p} \in \Lambda} C_{\mathfrak{p}}))|_s} \text{dist}(q, (\Gamma_{S-} \cup (\cup_{\mathfrak{p} \in \Lambda} \gamma_{\mathfrak{p}-}))) + \sup\nolimits_{q \in (\Gamma_{S-} \cup (\cup_{\mathfrak{p} \in \Lambda} \gamma_{\mathfrak{p}-}))} (C_S \cup (\cup_{\mathfrak{p} \in \Lambda} C_{\mathfrak{p}}))|_s q) < r.$$

Lemma 7.12 is proved momentarily.

Lemma 7.12 with what is said in Parts 1-5 imply that $\mathcal{M}^*$ is compact in the topology that is defined in Section 7b. Lemma 7.12 has the following additional implication: Given the lemma, standard elliptic regularity arguments using (7.8) and the various $\mathfrak{p} \in \Lambda$ versions of (3.4) can be applied to prove that the topology as defined in Section 7a on $\mathcal{M}^*$ is the same as the strong $C^\infty$ topology. To elaborate, these tools can be used to bootstrap from the uniform $L^\infty$ convergence that is asserted by Lemma 7.12 at large $s$ to prove strong convergence of the sort asserted by Proposition 7.3 but with respect to some Hölder topology with exponent $\upsilon > 0$. The tools are used again with (7.8) and (3.4) to prove the strong $C^{1+\upsilon}$ convergence, then again to prove strong $C^{2+\upsilon}$ convergence, and so on. The sorts of tools needed can be found in Chapter 6 of [M].

This equivalence between the topology from Section 7a and the strong $C^\infty$ topology implies what is asserted by the two bullets of Proposition 7.3.

***Proof of Lemma 7.12***:  Suppose that the lemma is false so as to generate nonsense. Granted that such is the case, there exists $r > 0$ and a sequence $\{(\tau_n, \mathcal{C}_n)\}_{n=1,2,...}$ of the following sort: Each $n \in \{1, 2, ...\}$ version of $(\tau_n, \mathcal{C}_n)$ is an element in $\mathcal{M}^*$ and such that the condition stated by Lemma 7.12 fails to hold for $s = s_n \geq n$ with some fixed choice of $r$. Given the conclusions of the preceding Parts 1-5 of this subsection, no generality is lost by requiring that the sequence $\{(\tau_n, \mathcal{C}_n)\}_{n=1,2,...}$ converge in the topology on $\mathcal{M}^*$ from Part 5 to element $(\mathcal{C}, \tau) \in \mathcal{M}^*$. The derivation of nonsense from the existence of such a sequence has two steps.

    <u>Step 1</u>:  Let $\hat{\upsilon}_+$ denote the HF-cycle that is used to define $\Theta_+$. Write $f^{-1}(1, 2) \subset M$ as $(1, 2) \times \Sigma$ so as to identify this set of integral curves with a set given by $(1, 2) \times \hat{q}$ with $\hat{q}$ denoting a certain set of G distinct elements from $C_{1+} \cap C_{2-}$. Let $S_-$ denote the submanifold $\mathbb{R} \times [1, 2] \times \hat{q}_-$ in $\mathbb{R} \times [1, 2] \times \Sigma$. This $S_-$ is a Lipshitz submanifold.

    Introduce $\mathcal{M}^*_-$ to denote the $S_-$ version of the space $\mathcal{M}^*$ as defined using $\hat{\Theta}_-$ for both the $s \to -\infty$ and $s \to \infty$ limit conditions on its constituent elements. For any given $\tau \in [0, 1]$, the set $\{\mathbb{R} \times \Gamma_{S_-}, \{\mathbb{R} \times \gamma_{\mathfrak{p}_-}\}_{\mathfrak{p} \in \Lambda}\}$ is a $(J, \tau)$-holomorphic submanifold. Thus $\mathcal{M}^*_-$ is non-empty. These are the only elements in $\mathcal{M}^*_-$. To prove this, let $(\tau, \mathcal{C} = \{C_S, \{C_\mathfrak{p}\}_{\mathfrak{p} \in \Lambda}\})$ denote a given element $\mathcal{M}^*_-$. Compute the sum whose constituent terms are the integral of $w$ over $C_S$ and the integral of $w$ over the respective $\mathfrak{p} \in \Lambda$ versions of $C_\mathfrak{p}$. If this sum is zero, then $C = \{\mathbb{R} \times \Gamma_{S_-}, \{\mathbb{R} \times \gamma_{\mathfrak{p}_-}\}_{\mathfrak{p} \in \Lambda}\}$. The fact that this sum of integrals is indeed zero is proved in the next paragraph.

    To see that the sum of integrals is zero, remark that the form $w$ is exact on a neighborhood in Y of the union of the sets $\{\mathcal{H}^+_{\mathfrak{p}*}\}_{\mathfrak{p} \in \Lambda}$ with a uniform radius tubular neighborhood of the arcs that comprise $\Gamma_{S_-}$. This being the case, integration by parts as in the proof of Lemma 7.8 writes the sum of the integrals of $w$ over the constituents of $\mathcal{C}$ as a sum of three terms:  The first is the integral of the anti-derivative 1-form over the union of the curves that comprise $\Theta_-$. The second is minus the integral of this same 1-form over the same union of curves. The third is itself a sum, this sum indexed by $\Lambda$ with any given term given by (7.13). As both $\varphi^{s0}$ and $\varsigma^{s0}$ are constant in each $\mathfrak{p} \in \Lambda$ version of (7.13), each $\mathfrak{p} \in \Lambda$ contribution to the third term is zero.

    <u>Step 2</u>:  For each $n \in \{1, 2, ...\}$, translate each element in $\mathcal{C}_n$ by $-s_n$ along the $\mathbb{R}$ factor in either $\mathbb{R} \times [1, 2] \times \Sigma$ or the appropriate $\mathfrak{p} \in \Lambda$ version of $\mathbb{R} \times \mathcal{H}^+_{\mathfrak{p}*}$. Let $\mathcal{C}_n'$ denote the corresponding set. Minor modifications of the arguments given in Parts 1-5 prove that the sequence $\{(\tau_n', \mathcal{C}_n')\}_{n=1,2,...}$ converges in the topology given by (7.26) to some element in the space $\mathcal{M}^*_-$. Given the assumptions about the initial sequence, this

element can not be $(\tau, \{\mathbb{R} \times \Gamma_{S_-}, \{\mathbb{R} \times \gamma_{p_-}\}_{p \in \Lambda}\})$ as at least one of the submanifolds of each $n \in \{1, 2, \ldots\}$ version of $\mathcal{C}_n'$ contains a point that has distance at least $r$ from a point in the corresponding submanifold from the set $\{\mathbb{R} \times \Gamma_{S_-}, \{\mathbb{R} \times \gamma_{p_-}\}_{p \in \Lambda}\}$. But this is nonsense given what is said in Step 1.

*Part 7*: This part adds what is needed to the arguments in Parts 1-6 so as to prove Proposition 7.3 when some $\mathfrak{p} \in \Lambda$ versions of $\Delta_\mathfrak{p}$ are 1 or 2. The arguments given in Parts 1-6 can be seen as having three components. The first component is summarized by Lemma 7.8; this component gives bounds on the integrals of $w$ and $ds \wedge \hat{a}$. The second component is summarized by Lemma 7.11; this component controls the behavior the relevant submanifolds where they intersect $\mathbb{R} \times M_\delta$. The final component controls the behavior of the remaining portions of the constituent subvarieties. This component is comprises Lemmas 7.9 and 7.10 and Parts 5 and 6. The three steps that follow speak to these three components in the case when some $\mathfrak{p} \in \Lambda$ versions of $\Delta_\mathfrak{p}$ are 1 or 2.

<u>Step 1</u>: The analog of Lemma 7.8 in the general case makes the same assertion as the original with the $\{\Delta_\mathfrak{p} = 0\}_{\mathfrak{p} \in \Lambda}$ condition omitted. The proof copies the arguments that prove Lemma 6.1 to control the integral of $w$ on the relevant parts of any given $(\tau, \mathcal{C})$ and $\mathfrak{p} \in \Lambda$ version of $C_\mathfrak{p}$ and it copies what is said in Part 1's proof of Lemma 7.8 to control the integral of $w$ over the $(\tau, \mathcal{C})$ version of $C_S$. What is said in Part 1 with regards to (7.13) hold what ever the value of $\Delta_\mathfrak{p}$.

<u>Step 2</u>: The analog of Lemma 7.11 in the general case is identical to its namesake but for the absense of the $\{\Delta_\mathfrak{p} = 0\}_{\mathfrak{p} \in \Lambda}$ assumption. There is but one change in the proof. The arguments for lemma when $\Delta_\mathfrak{p} = 0$ derive nonsense from a certain apriori assumption about the sequence any given $\mathfrak{p} \in \Lambda - \Lambda_*$ version of the sequence $\{C_{\mathfrak{p}n}\}_{n=1,2,\ldots}$ as they find a loop in the all large n versions of $C_{\mathfrak{p}n}$ that generates the first homology of $\mathbb{R} \times \mathcal{H}^+_{\mathfrak{p}*}$. This nonsense comes from the fact that $C_{\mathfrak{p}n}$ is contractible if $\mathfrak{p} \in \Lambda - \Lambda_*$. The space $C_{\mathfrak{p}n}$ is not contractible if $\Delta_\mathfrak{p} \neq 0$, but even so the existence of the corresponding loop in $C_{\mathfrak{p}n}$ is nonsense: The loop in question sits in a component of the $|u| > 0$ part of $C_{\mathfrak{p}n}$, and each such component is contractible.

<u>Step 3</u>: The analog of Lemma 7.9 in the case when some $\mathfrak{p} \in \Lambda$ versions of $\Delta_\mathfrak{p}$ are non-zero is identical to its namesake but for the absence of the $\{\Delta_\mathfrak{p} = 0\}_{\mathfrak{p} \in \Lambda}$ assumption. The proof of the analog is identical to that given for Lemma 7.9 with it understood that Lemma 7.11 holds when various $\mathfrak{p} \in \Lambda$ versions of $\Delta_\mathfrak{p} > 0$.

The analog of Lemma 7.10 replaces the latter with two lemmas. The first one is much like Lemma 6.3.

**Lemma 7.13**: *There is a purely S-dependent (or K-compatible) $\kappa > 1$ and given $z_* < \kappa^{-1}$, there exists $\kappa_* > 1$ that depends only on $z_*$ but is otherwise purely S-dependent (or K compatible); and these have the following significance: Define $\mathcal{M}^*$ using $z_* < \kappa^{-1}$, using $\delta < \kappa_*^{-1} z_*$ and using a pair from $\hat{\mathcal{Z}}^S$. There exists there exists $\kappa_{**} > 1$ and given $\varepsilon \in (0, 1]$, there exists $\kappa_\varepsilon > 1$ such that any given $(\tau, C = \{C_S, \{C_p\}_{p \in \Lambda}\}) \in \mathcal{M}^*$ has the properties listed next.*

- *If $\Delta_p = 0$, then $1 - 3\cos^2\theta > \kappa_{**}^{-1}$.*
- *If $\Delta_p > 0$, then*
  a) *$1 - 3\cos^2\theta > \kappa_{**}^{-1}$ on the $s < -\kappa$ part of $C_p$.*
  b) *$1 - 3\cos^2\theta > \kappa_\varepsilon^{-1}$ on the $|u| \geq \varepsilon$ part of $C_p$.*
- *If $\Delta_p = 1$ and $\mathfrak{m}_p = -1$, then $\cos\theta > -\frac{1}{\sqrt{3}} + \kappa_{**}^{-1}$ on the whole of $C_p$*
- *If $\Delta_p = 1$ and $\mathfrak{m}_p = 1$, then $\cos\theta < \frac{1}{\sqrt{3}} - \kappa_{**}^{-1}$ on the whole of $C_p$*

*Proof of Lemma 7.13*: Given that the assertions of Lemma 7.11 are true when some $\mathfrak{p} \in \Lambda$ versions of $\Delta_p > 0$, the proof of Lemma 7.13 is obtained by using the arguments for the proof of Lemma 7.10 to prove what is asserted in the first bullet; and by using the arguments for the proof of Lemma 6.3 or its ($\Delta_p = 1$, $\mathfrak{m}_p = 1$) or $\Delta_p = 2$ incarnations to prove the assertions of the second and third bullets.

To set the stage for the second lemma, suppose for the sake of argument that $\mathfrak{p} \in \Lambda$ has $\Delta_p = 1$ and $\mathfrak{m}_p = -1$. There is, in this case, constants $s_p \geq 1$ and $c_{1p} \in \mathbb{R} - 0$ and $\phi_{1p} \in \mathbb{R}/(2\pi\mathbb{Z})$ with the following significance: The complement of a certain compact set of the $s \geq 1$ part of $C_p$ has two components. One is a strip diffeomorphic to $[s_p, \infty) \times I_*$ whose image via the projection to $\mathcal{H}^+_{p*}$ is very close to $\gamma_{p+}$. The other is a cylinder whose image via this projection lies in the tubular neighborhood $U_+$ of $\hat{\gamma}_p^+$ that is described in Section 5a. Let $\mathcal{E}$ denote this end of $C_p$. Reintroduce the coordinates $(s_+, \phi_+, \theta_+, u_+)$ for $\mathbb{R} \times U_+$ as defined in (5.5). Then $\mathcal{E}$ sits in the $s_+ \in [s_p, \infty)$ part of $\mathbb{R} \times U_+$ as a smooth, properly embedded submanifold with boundary on the $s_+ = s_p$ slice. Furthermore, this intersection is given by the graph of a smooth map as described in the first bullet of Proposition 5.1 with domain $[s_p, \infty) \times \mathbb{R}/(2\pi\mathbb{Z})$ that has the form depicted in (6.23) with $c_1 = c_\mathcal{E}$ and with $\mathfrak{y}_{1+}$ defined using $\phi_1 = \phi_\mathcal{E}$. There is a completely analogous picture when $(\Delta_p = 1, \mathfrak{m}_p = 1)$ and when $\Delta_p = 2$.

**Lemma 7.14**: *There is a purely S-dependent (or K-compatible) $\kappa > 1$ and given $z_* < \kappa^{-1}$, there exists $\kappa_* > 1$ that depends only on $z_*$ but is otherwise purely S-dependent (or K* *compatible); and these have the following significance: Define $\mathcal{M}^*$ using $z_* < \kappa^{-1}$, using $\delta < \kappa_*^{-1} z_*$ and using a pair from $\hat{\mathcal{Z}}^S$. There exists $\kappa_{**} \geq 1$ such that if $(\tau, \{C_S, \{C_\mathfrak{p}\}_{\mathfrak{p} \in \Lambda}\})$ is from $\mathcal{M}^*$ and if $\mathfrak{p} \in \Lambda$ is such that $\Delta_\mathfrak{p} > 0$, then $s_\mathfrak{p} < \kappa_{**}$ and $\kappa_{**}^{-1} \leq |c_\mathcal{E}| \leq \kappa_{**}$ for each end $\mathcal{E} \subset C_\mathfrak{p}$ whose constant $s >> 1$ slices converge as $s \to \infty$ to either $\hat{\gamma}_\mathfrak{p}^+$ or $\hat{\gamma}_\mathfrak{p}^-$.*

***Proof of Lemma 7.14***: But for notation, the argument is the same as that used to prove Lemma 6.4.

With Lemmas 7.13 and 7.14 in hand, the remaining arguments for the third component differ only cosmetically from the arguments given in Parts 5 and 6 of this section and so no more will be said.

### e) The structure of $\mathcal{M}^*$

There are four parts to this subsection; Part 4 contains the proofs of Propositions 7.1 and 7.2. Parts 1-3 set up the necessary machinery.

*Part 1*: Let $(\tau, \mathcal{C} = \{C_S, \{C_\mathfrak{p}\}_{\mathfrak{p} \in \Lambda}\})$ denote a given element of $\mathcal{M}^*$. This step associates a certain Fredholm operator to $(\tau, \mathcal{C})$ whose cokernel provides a specific version of $\mathbb{R}^n$ for use in Proposition 7.1. This operator is also used in the proof of Proposition 7.2. The operator is denoted in what follows by $\mathcal{D}_\mathcal{C}$.

A dense domain for $\mathcal{D}_C$ consists of a direct sum of function spaces with the first summand labeled by S and the others labeled by the set $\Lambda$. An element in the summand labeled by S is, among other things, a smooth and compactly supported section of the normal bundle $N_S$ on the $t \in [1+z_*, 2-z_*]$ part of S. These sections are further constrained to be $L^2$ orthogonal to the restriction of the $L^2$ kernel of (1.25)'s operator $D_S$ to the $t \in [1+z_*, 2-z_*]$ part of S. Meanwhile, any given $\mathfrak{p} \in \Lambda$ labeled summand consists of a certain sorts of compactly supported, smooth sections of the normal bundle to $C_\mathfrak{p}$. The respective elements in the S-summand and any given $\mathfrak{p} \in \Lambda$ summand are further constrained on the common boundary of $C_\mathfrak{p}$ and $C_S$. The next paragraph describes this constraint.

Write the intersection of $C_S$ with $\mathcal{H}^+_\mathfrak{p}$ as a graph in the manner of PROPERTY 4 in Part 1 of Section 1a so as to identify a section of the normal bundle of $C_S$ on either boundary component as a map from $\mathbb{R}$ to $\mathbb{R}^2$. Use x for the $\mathbb{R}$ coordinate. The coordinate z for the $[z_*, z_S]$ factor is $z = t - 1$ for the index 1 critical point side of $\mathcal{H}^+_\mathfrak{p}$ and $z = 2 - t$ for the index 2 critical point side. With a section of S given, write the components of the corresponding map as $(x, z) \to (\varphi^{S'}, \varsigma^{S'})|_{(x,z)}$. View $C_\mathfrak{p}$ via $\Psi_\mathfrak{p}$ as a graph in $\mathbb{R} \times \mathcal{X}$ in the manner of PROPERTY 2 in Part 2 of Section 7a; and use the 1-forms $(d\hat{\phi}, dh)$ to identify the normal bundle with the product $\mathbb{R}^2$ bundle. Having done so, a section of this normal bundle becomes a map to $\mathbb{R}^2$ from the domain in $\mathbb{R} \times I_*$ of the pair $(\varphi^\mathfrak{p}, \varsigma^\mathfrak{p})$. Introduce the coordinate $z = e^{-2(R-\hat{u})}$ for the $\hat{u} \geq R + \ln\delta$ part of the domain and $z = e^{-2(R+\hat{u})}$ for the $\hat{u} \leq -R - \ln\delta$ part. Granted this notation, the following constraint holds where $z = z_*$ in either case:

- $\varsigma^{S\prime} = \tau \varsigma^{\mathfrak{p}\prime}$.
- $\tau \varphi^{S\prime} = \varphi^{\mathfrak{p}\prime}$.

(7.32)

What follows describes the operator $\mathcal{D}_\mathcal{C}$. To start, $\mathcal{D}_\mathcal{C}$ acts diagonally with respect to the labeling of the summands of the domain and range spaces. To define its action on the S-labeled summand, write $C_S$ as in Section 7a in terms of a section of the normal bundle over the $t \in [1+z_*, 2-z_*]$ part of S. Let $\eta$ denote this section. The action of $\mathcal{D}_\mathcal{C}$ on the S-labeled summand is that of an operator that is denoted in what follows by $D_\eta$. Let $\eta\prime$ denote a section of $N_S$ over the part of S where $t \in [1+z_*, 2-z_*]$. Then operator $D_\eta$ sends $\eta\prime$ to

$$\bar{\partial}\eta\prime + \mathfrak{r}_1(\eta) \cdot \partial\eta\prime + (\nabla_{\eta\prime} \mathfrak{r}_1)|_\eta \partial\eta + \nabla_{\eta\prime} \mathfrak{r}_{0\eta})|_\eta$$

(7.33)

where the notation uses $\nabla_{\eta\prime}$ to denote the directional derivative along the fiber of $N_S$ in the direction given by $\eta\prime$. The terms that involve $\nabla_{\eta\prime}$ are zero'th order and $\mathbb{R}$-linear. By way of an example, the operator in (7.18) is the $\eta = 0$ version of (7.33).

To define the action of $\mathcal{D}_\mathcal{C}$ on the remaining summands, fix $\mathfrak{p} \in \Lambda$ and let $(\varphi^\mathfrak{p}, \varsigma^\mathfrak{p})$ denote the pair that defines $C_\mathfrak{p}$. Use $D_\mathfrak{p}$ to denote the $(\varphi = \varphi^\mathfrak{p}, \varsigma = \varsigma^\mathfrak{p})$ version of the operator that is depicted in (3.6). With the elements in the $\mathfrak{p}$-labeled summand viewed as maps from the relevant domain in $\mathbb{R} \times I_*$ to $\mathbb{R}^2$, the action of $\mathcal{D}_\mathcal{C}$ on the $\mathfrak{p}$-labeled summand is given by $D_\mathfrak{p}$.

To say slightly more about the operators from the set $\{D_\eta, \{D_\mathfrak{p}\}_{\mathfrak{p}\in\Lambda}\}$, keep in mind that $C_S$ and each $\mathfrak{p} \in \Lambda$ version of $C_\mathfrak{p}$ is J-holomorphic. This being the case, the normal bundle of each can be viewed as a complex line bundle. Meanwhile, the Riemannian metric defined by J and the compatible 2-form $ds \wedge \hat{a} + w$ endow these submanifolds and their normal bundles with holomorphic structures. This understood, let N denote the normal bundle to a given $C \in \{C_S, \{C_\mathfrak{p}\}_{\mathfrak{p}\in\Lambda}\}$ but viewed now as a complex line bundle. Writing $C_S$ in terms of the section $\eta$ identifies the $t \in [1+z_*, 2-z_*]$ part of the normal bundle $N_S$ with the $C_S$ version of N. The S-labeled part of the domain of $\mathcal{D}_\mathcal{C}$ can be viewed as a section of this version of N. Do so and $D_\eta$ appears as a first order differential operator that maps $C^\infty(C_S; N)$ to $C^\infty(C_S; N \otimes T^{0,1}C_S)$ which has the schematic form of the operator in (1.25).  Meanwhile, each $\mathfrak{p} \in \Lambda$ version of $D_\mathfrak{p}$ maps the space $C^\infty(C_\mathfrak{p}; N)$ to $C^\infty(C_\mathfrak{p}; N \otimes T^{0,1}C_\mathfrak{p})$ and it also can be written in the form depicted by (1.25).

The Banach space domain for $\mathcal{D}_\mathcal{C}$ is obtained by completing the dense domain described above using the norm whose square is the sum of the squares of certain $L^2_1$ norms on the constituent summands. The $L^2_1$ norm on the S-labeled summand is that defined by the covariant derivative and metric on S. With $C_S$ viewed as the image of $\mathfrak{e}_S$ of the section $\eta$, this is norm is equivalent to the one defined by the induced Riemannian metric on $TC_S$ and the induced metric and covariant derivative for the normal bundle of $C_S$. The $L^2_1$ norm on any given $\mathfrak{p}$-labeled summand is defined using the induced Riemannian metric for $TC_\mathfrak{p}$ and the induced metric and covariant derivative for the normal bundle to $C_\mathfrak{p}$.

The operator $\mathcal{D}_\mathcal{C}$ maps the Banach space just define to the Banach space that is defined by completing the respective space of compactly supported sections of each $C \in \{C_S, \{C_\mathfrak{p}\}_{\mathfrak{p}\in\Lambda}\}$ version of $C^\infty(C; N \otimes T^{0,1}C)$ using the induced $L^2$ inner products. In the case $C = C_S$, this is the same as completing a space of compactly supported sections of the appropriate bundle over the $t \in [1+z_*, 2-z_*]$ part of S using the $L^2$ norm on S.

The respective domain and range Banach spaces for $\mathcal{D}_\mathcal{C}$ are denoted in what follows by $\mathbb{H}_\tau$ and $\mathbb{L}$. The norm on $\mathbb{H}_\tau$ is denoted by $\|\cdot\|_1$, and the $L^2$ norm on either $\mathbb{H}_\tau$ or $\mathbb{L}$ is denoted by $\|\cdot\|$. The operator $\mathcal{D}_\mathcal{C}$ defines a bounded operator from $\mathbb{H}_\tau$ to $\mathbb{L}$.

Respective domain and range Banach spaces with slightly stronger norms are also needed in what follows. The domain version is denoted in what follows by $\mathbb{H}_{\tau*}$. The norm that defines $\mathbb{H}_{\tau*}$ is the sum of the $L^2_1$ norm defined above and the norm whose square is given on each $\mathfrak{p} \in \Lambda$ labeled summand by the $C_\mathfrak{p}$ analog of (5.31), and is given on the S-labeled summand by the analog of (5.31) for sections of the normal bundle of S. Lemma 5.10 and has an analog for S; they assert that elements in $\mathbb{H}_{\tau*}$ are Hölder continuous and the associated map from $\mathbb{H}_{\tau*}$ to the relevant Hölder space is continuous. The strengthened range Hilbert space is denoted by $\mathbb{L}_*$; The square of the norm that defines this space is given on the various summands by replacing in (5.31) the length of the covariant derivative of a given section with that of the section itself. The operator $\mathcal{D}_\mathcal{C}$ defines also a bounded map from $\mathbb{H}_{\tau*}$ to $\mathbb{L}_*$

*Part 2*:  The lemma that follows states what is need concerning the operator $\mathcal{D}_\mathcal{C}$.

**Lemma 7.15**: *There is a purely S-dependent (or $\mathcal{K}$-compatible) $\kappa \geq 1$ such that if $\mathcal{M}^*$ is defined using $\rho_S \leq \kappa^{-1}$, $z_* \leq \kappa^{-1}$, $\delta < \kappa^{-2}z_*$, and any pair from $\hat{\mathcal{Z}}^S$, then the following is true: Use any given $(\tau, \mathcal{C}) \in \mathcal{M}^*$ to define the operator $\mathcal{D}_\mathcal{C}$ and the spaces $\mathbb{H}_\tau$ and $\mathbb{L}$. The operator $\mathcal{D}_\mathcal{C}$ is a Fredholm operator from $\mathbb{H}_\tau$ to $\mathbb{L}$ with index equal to $\sum_{\mathfrak{p}\in\Lambda}\Delta_\mathfrak{p}$. Moreover,* *elements in the kernel of $\mathcal{D}_\mathcal{C}$ are $C^\infty$ elements in $\mathbb{H}_{\tau*}$; and elements in the cokernel of $\mathcal{D}_\mathcal{C}$ are represented by smooth elements in $\mathbb{L}_*$.*

The remainder of this Part 2 contains the

***Proof of Lemma 7.15***: The proof of this lemma has five steps.

Step 1: Just as in Part 4 of Section II.6e and in the proofs of Propositions 3.1 and 5.7, the assertion that $\mathcal{D}_\mathcal{C}$ has closed range and finite dimensional kernel follows if there exists a constant $c > 1$ such that the following two conditions hold:

- $\|\mathcal{D}_\mathcal{C}\mathfrak{h}\|^2 \geq c^{-1}\|\mathfrak{h}\|_1^2 - c\|\mathfrak{h}\|^2$ *for any* $\mathfrak{h} \in \mathbb{H}_\tau$.
- *There exists $s_1 > 1$ such that if $\mathfrak{h} \in \mathbb{H}$ has S-summand with support where $|s| > s_1$ and each $\Lambda$-labled summand has support where $|s| \geq s_1$, then $\|\mathcal{D}_\mathcal{C}\mathfrak{h}\|^2 \geq c^{-1}\|\mathfrak{h}\|^2$.*

(7.34)

There is a purely S-dependent (or $\mathcal{K}$-compatible) constant $c > 1$ such that if both $\rho_S < c^{-1}$ and $\delta < c^{-1}$, then the conditions in (7.34) hold for those $\mathfrak{h} \in \mathbb{H}_\tau$ with the following property: The closure of the support of $\mathfrak{h}$ is disjoint from the boundary of the submanifolds from $\mathcal{C}$. Indeed, if a section of $N_S$ has compact support on the interior of $C_S$ then this follows from what is said in Parts 3 and 4 of Section II.6e given that what is denoted in (7.33) by $\mathfrak{r}_1$ obeys $|\mathfrak{r}_1(\eta)| \leq c_0^{-1}|\eta|$. If a section has support in some $\mathfrak{p} \in \Lambda$ version of $C_\mathfrak{p}$, then this follows from Propositions 3.1 and 5.7.

To see about the top bullet in (7.34) when the support of the summands of $\mathfrak{h}$ intersect the boundaries of the defining domains, write $\mathcal{C}$ as $\{C_S, \{C_\mathfrak{p}\}_{\mathfrak{p}\in\Lambda}\}$ and fix attention on a given $\mathfrak{p} \in \Lambda$. Suppose that $\mathfrak{h} \in \mathbb{H}_\tau$ is in the dense domain, and write its S-summand near the intersection of $\mathbb{R} \times M_\delta$ with $\mathbb{R} \times \mathcal{H}^+_{\mathfrak{p}*}$ as a pair of functions on the domain $\mathbb{R} \times [z_*, z_S]$ as done in (7.32). These are denoted by $(\varphi^{S\prime}, \varsigma^{S\prime})$. Likewise, write an element in the $\mathfrak{p}$-summand as a pair $(\varphi^{\mathfrak{p}\prime}, \varsigma^{\mathfrak{p}\prime})$, these being functions on $\mathbb{R} \times [\delta^2, z_*]$. Both $D_\eta$ and $D_\mathfrak{p}$ appear here as the Cauchy-Riemann operator that acts on a given $(\varphi\prime, \varsigma\prime)$ to give $(\partial_x\varphi\prime - \partial_z\varsigma\prime, \partial_x\varsigma\prime + \partial_z\varphi\prime)$. Given that (7.32) holds, so does Lemma 7.4. The conclusions of this lemma as applied for all $\mathfrak{p} \in \Lambda$ prove the top bullet in (7.34).

The argument for the second bullet in (7.34) when the closure of the support of a given element in $\mathbb{H}_\tau$ intersects the boundaries of the submanifolds from $\mathcal{C}$ occupies Steps 2 and 3 of what follows. These steps focus on the case where the respective summands in a given element from the dense domain of $\mathbb{H}_\tau$ have support where $s \ll -1$ on $C_S$ and each $\mathfrak{p} \in \Lambda$ version of $C_\mathfrak{p}$. With regards to $C_\mathfrak{p}$, no generality is lost by assuming that there exists a purely S-dependent (or $\mathcal{K}$-compatible) constant $c \geq 1$ such that the support of the $C_\mathfrak{p}$ summand lies where $1 - 3\cos^2\theta > c^{-1}$. This is because of what was said at the outset: Equation (7.34) holds for sections which are supported on closed sets in $C_\mathfrak{p}$'s interior. But for notation, the argument works when the summands have support where $s \gg 1$ on the various submanifolds from $\mathcal{C}$.

<u>Step</u> 2: As in Part 4 of Section II.6e and the proofs of Propositions 3.1 and 5.7, the second bullet in (7.34) is obeyed if a certain non-negative definite, quadratic form has trivial kernel. The quadratic form is defined on a certain space of sections of the normal bundle in Y to the union of the curves in $\Theta_-$. This step defines the relevant space of sections and the quadratic form. The quadratic form in question is denoted by $\mathcal{Q}$ in what follows and the domain by $\mathbb{V}_\tau$.

To give the definitions, first decompose the curves from $\Theta_-$ that are not contained in a single $\mathfrak{p} \in \Lambda$ version of $C_\mathfrak{p}$ so as to define a set of segments of the form $\{\Gamma_S, \{\gamma_{\mathfrak{p}-}\}_{\mathfrak{p}\in\Lambda}\}$ where $\Gamma_S \subset M_\delta$ is the union of the G segments that comprise the intersection of the union of the curves from $\Theta_-$ with the $f \in [1+z_*, 2-z_*]$ part of $M_\delta$.

View the $f \in [1, 2]$ part of $M_\delta$ as a subset of $[1, 2] \times \Sigma$. Doing so identifies $\Gamma_S$ as $[1, 2] \times \Lambda^S$ where $\Lambda^S \subset \Sigma$ is a set of N distinct points. Fix a holomorphic coordinate centered on each point in $\Lambda^S$. Doing so identifies a section of the normal bundle of $\Gamma_S$ with a pair of real functions on the interval $[1+z_*, 2-z_*]$, these corresponding to the real and imaginary parts of the complex coordinate. A section of the normal bundle in Y to a segment from $\Gamma_S$ that concatenates with a given $\mathfrak{p} \in \Lambda$ version of $\gamma_{\mathfrak{p}-}$ appears at a shared endpoint as a pair of functions on $[z_*, z_S]$, these denoted by $(\varphi^{S'}, \varsigma^{S'})$.

Fix $\mathfrak{p} \in \Lambda$ and use the map $\Psi_\mathfrak{p}$ to view $\gamma_{\mathfrak{p}-}$ as the locus in $\mathbb{R} \times \mathcal{X}$ where the coordinates are such that $x = 0$ and $h = \hat{h}(\gamma_{\mathfrak{p}-})$. Doing so identifies any given section of the normal bundle to $\gamma_{\mathfrak{p}-}$ with a pair of functions on $I_*$. Such a pair is written as $(\varphi^{\mathfrak{p}'}, \varsigma^{\mathfrak{p}'})$. Near the boundary of $I_*$ these can be written as a pair of functions on the interval $[\delta^2, z_*]$ by writing z on this integral as $z = e^{-2(R-|\hat{u}|)}$.

Given this notation, the domain, $\mathbb{V}_\tau$, for $\mathcal{L}$ consists of the direct sum of G + 1 function spaces, the first labeled by S and the others labeled by $\Lambda$. The S-labeled summand consists of a suitably constrained set of sections of the normal bundle to $\Gamma_S$. Meanwhile, any given $\mathfrak{p} \in \Lambda$ labeled summand consists of a suitably constrained set of pairs of functions on $I_*$. The constraints on the components of the various summands are the boundary constraints given by the various $\mathfrak{p} \in \Lambda$ versions of (7.32).

The quadratic form $\mathcal{Q}$ is the sum of quadratic forms that are defined on the various summands of $\mathbb{V}_\tau$. To define the contribution to this sum from the S-labeled summand, view an element of the latter as a map from $[1+z_*, 2-z_*]$ to $\mathbb{R}^2$ suitably constrained on its boundary. Let $t \to \eta'(t)$ denote such a map. The value of the S-labled quadratic form on $\eta$ is

$$\mathcal{Q}_S(\eta') = \int_{1+z_*}^{2-z_*} |\tfrac{d}{dt}\eta'|^2 dt \ .$$

(7.35)

Let $\mathfrak{p}$ denote a given element in $\Lambda$. The contribution to $\mathcal{Q}$ from the summand labled by $\mathfrak{p}$ is given by the expression in (3.15) with it understood that the coefficients $\mathfrak{a}_{1-}$, $\mathfrak{a}_{2-}$ and $\mathfrak{b}_{2-}$ are defined using $\gamma_{\mathfrak{p}-}$. In particular, $\mathfrak{b}_{2-}$ is the $\gamma = \gamma_{\mathfrak{p}-}$ version of what appears on the right hand side of (3.11).

<u>Step 3</u>: Part 4 of Section II.6e and Section 3c prove that $\mathcal{Q}$ is positive definite in the case $\tau = 0$. This understood, assume $\tau > 0$ in what follows. To see that $\mathcal{Q}$ is postive definite when $\tau > 0$, suppose that $\mathfrak{h} \in \mathbb{V}_\tau$ and $\mathcal{Q}(\mathfrak{h}) = 0$. Fix $\mathfrak{p} \in \Lambda$ and write the corresponding component of $\mathfrak{h}$ as $(\varphi^{\mathfrak{p}'}, \varsigma^{\mathfrak{p}'})$. It follows from (3.15) that $\varsigma^{\mathfrak{p}'}$ is constant along $I_*$, and it follows from (3.17) that $\varphi^{\mathfrak{p}'}$ increases along $I_*$. To elaborate, let $\varphi^{\mathfrak{p}'}_-$ denote the $\hat{u} = -R - \tfrac{1}{2} \ln z_*$ value of $\varphi^{\mathfrak{p}'}$ and let $\varphi^{\mathfrak{p}'}_+$ denote the value at $\hat{u} = R + \tfrac{1}{2} \ln z_*$. Then (3.15) and (3.13) imply that

$$\varphi^{\mathfrak{p}'}_+ - \varphi^{\mathfrak{p}'}_- = c_\mathfrak{p} x_0^{-1} R \, \varsigma^{\mathfrak{p}'} ,$$

(7.36)

where $c_\mathfrak{p} > c_0^{-1}$. Let $\eta'$ denote the component of $\mathfrak{h}$ in the S-summand of $\mathbb{V}_\tau$. It follows from (7.35) that the latter is constant.

Let $\gamma \in \Gamma_S$ denote the segment whose starting point is the $\hat{u} = R + \tfrac{1}{2} \ln z_*$ endpoint of $\gamma_{\mathfrak{p}-}$. Write $\eta'$ near the starting point of $\gamma$ as $(\varphi^\gamma_-, \varsigma^\gamma_-)$, so as to see the implications of what was just said. This pair of functions is independent of the parameter $z$ on the interval $[z_*, z_S]$. The constraint in (7.32) demands that

$$\varsigma^\gamma_- = \tau \varsigma^{\mathfrak{p}'} \quad and \quad \varphi^\gamma_- = \tau^{-1} \varphi^{\mathfrak{p}'}_+$$

(7.37)

The end point of $\gamma$ is the $\hat{u} = -R - \tfrac{1}{2} \ln z_*$ boundary of some $\mathfrak{q} \in \Lambda$ version of $\gamma_\mathfrak{q}$. Introduce the analogous pair $(\varphi^\gamma_+, \varsigma^\gamma_+)$. It follows from what is said in the proof of Proposition II.2.7 that these are determined by $(\varphi^\gamma_-, \varsigma^\gamma_-)$ via a formula of the form

$$\varsigma^\gamma_+ = -a_\gamma \varsigma^\gamma_- - b_\gamma \varphi^\gamma_- \quad and \quad \varphi^\gamma_+ = c_\gamma \varsigma^\gamma_- + d_\gamma \varphi^\gamma_- ,$$

(7.38)

where the coefficients here are such that $b_\gamma \neq 0$ and such that

$$\begin{pmatrix} a_\gamma & b_\gamma \\ c_\gamma & d_\gamma \end{pmatrix}$$

(7.39)

has determinant -1. Meanwhile, it follows from (7.32) that the pair $(\varphi^{q\prime}_-, \varsigma^{q\prime}_-)$ is given by

$$\varsigma^{q\prime}_- = -\tau^{-1}\varsigma^\gamma_+ \quad and \quad \varphi^{q\prime}_- = \tau\varphi^\gamma_+ .$$

(7.40)

The simplest case to analyze is that where $q = p$ and so the endpoint of $\gamma$ is the starting point of $\gamma_p$. If this is the case, then (7.36), (7.37), (7.38) and (7.40) require that

$$\varsigma^{p\prime} = a_\gamma \varsigma^{p\prime} + b_\gamma \tau^{-2}(\varphi^p_- + c_p x_0^{-1} R\varsigma^{p\prime}) \quad and \quad \varphi^{p\prime}_- = c_\gamma \tau^2 \varsigma^{p\prime} + d_\gamma (\varphi^p_- + c_p x_0^{-1} R\varsigma^{p\prime}).$$

(7.41)

It follows from (7.41) that $\mathfrak{h} \neq 0$ if and only if the matrix

$$\begin{pmatrix} +a_\gamma + \tau^{-2} b_\gamma c_p x_0^{-1} R & \tau^{-2} b_\gamma \\ \tau^2 c_\gamma + d_\gamma c_p x_0^{-1} R & d_\gamma \end{pmatrix}$$

(7.42)

has +1 as an eigenvalue. But its determinant is -1 and its trace is larger than $c_0^{-1} x_0^{-1} R$, so one eigenvalue has norm greater than $c_0^{-1} x_0 R$ and the other has norm less than $c_0 x_0 R^{-1}$.

The argument in the general case is much the same. In this case, $\mathfrak{h} \neq 0$ if and only if the product of some $k \in \{1, \ldots, G\}$ versions of (7.42) have +1 as an eigenvalue. The details of the linear algebra are straightforward and so left to the reader.

Step 4: Granted that $\mathcal{D}_c$ has closed range, it follows that its cokernel is isomorphic to the kernel of its adjoint, this a bounded operator from $\mathbb{L}$ to $\mathbb{H}_\tau$. The kernel of the latter is isomorphic to the kernel of the formal, $L^2$-adjoint of $\mathcal{D}_c$. The proof that such is the case amounts to a standard application of linear, elliptic regularity arguments as applied on the interiors of the domains of the summands of $\mathbb{L}$ and an appeal to Lemma 7.7 to deal with the boundary conditions for these domains. This formal $L^2$ adjoint is denoted by $\mathcal{D}_c^\#$. The latter is elliptic, first order and $\mathbb{R}$-linear, with leading order symbol given by the adjoint of the $\bar{\partial}$ symbol. Its kernel obeys (7.32). The kernel of $\mathcal{D}_c^\#$ is finite dimensional if the $\mathcal{D}_c^\#$ version of (7.34) holds. The proof that such is the case is, but for notation, identical to the proof just given for the $\mathcal{D}_c$ version.

To see about the index of $\mathcal{D}_c$, fix $\tau' \in [0, 1]$. The operator $\mathcal{D}_c$ also defines a bounded map from $\mathbb{H}_{\tau'}$ to $\mathbb{L}$, and the argument given above that the $\mathcal{D}_c : \mathbb{H}_\tau \to \mathbb{L}$ is Fredholm can be repeated using now the $\tau'$ version of (7.32) to see that $\mathcal{D}_c$ is a Fredholm map from $\mathbb{H}_{\tau'}$ to $\mathbb{L}$ also. Meanwhile, the family $\{\mathbb{H}_{\tau'}\}_{\tau' \in [0,1]}$ defines a smooth, Banach space bundle over [0, 1]. The proof that this is so mimics what is said in Step 3 of Part 3 in Section 5c. Granted this smoothness, and granted that the set of Fredholm maps between Banach spaces is an open subset in the Banach space of bounded maps, it follows that the index of $\mathcal{D}_\mathcal{C}$ as a map from any given $\tau'$ version of $\mathbb{H}_\tau$ is independent of $\tau'$. In particular, the index of the $\mathbb{H}_\tau$ version is the same as its index on the $\mathbb{H}_{\tau'=0}$ version.

The $\tau' = 0$ matching conditions do not couple the various summands that comprise $\mathbb{H}_\tau$; they give separate conditions on each summand. As a consequence, the index of $\mathcal{D}_\mathcal{C}$ on $\mathbb{H}_\tau$ is the sum of the indices of the operators on the various summands with certain boundary conditions. The operator on the S-summand is $D_\eta$ with the boundary conditions $\varsigma^{S'} = 0$. The latter has index zero because it is given up to a term with small norm by the operator discussed in Section II.6e on a restricted domain, the orthogonal complement of its kernel. Meanwhile, any given $\mathfrak{p} \in \Lambda$ version of $D_\mathfrak{p}$ with the boundary condition $\varphi^{\mathfrak{p}'} = 0$ is of the form that is considered by Propositions 3.1 and 5.7; and the latter have index $\Delta_\mathfrak{p}$. Thus, the index of the $\mathbb{H}_{\tau'=0}$ version of $\mathcal{D}_\mathcal{C}$ is $\sum_{\mathfrak{p} \in \Lambda} \Delta_\mathfrak{p}$.

Step 5: The assertion that elements in the kernel of $\mathcal{D}_\mathcal{C}$ are smooth and in $\mathbb{H}_{\tau*}$ follows using standard linear elliptic estimates for smoothness in the interiors of the submanifolds from $\mathcal{C}$. Lemma 7.7 gives smoothness and the $\mathbb{H}_{\tau*}$ norm bound near the boundaries. The same argument as applied to the formal $L^2$-adjoint $\mathcal{D}_\mathcal{C}^\#$ proves that the elements in the kernel of the latter operator are smooth and are in $\mathbb{L}_*$. These last remarks imply what is asserted by the last sentence in the statement of Lemma 7.15.

*Part 3*: This part of the subsection starts with a lemma that is used subsequently to say more about the behavior at large $|s|$ of the various submanifolds from any given $(\tau, \mathcal{C})$ from $\mathcal{M}^*$. The lemma is also invoked to prove both Propositions 7.1 and 7.2.

To set the notation for this lemma, a domain in $C_S \cup (\cup_{\mathfrak{p} \in \Lambda} C_\mathfrak{p})$ is said to be *semi-bounded* when the coordinate $s$ is unbounded from above and bounded from below, or else unbounded from below and bounded from above. When U is a semi-bounded domain, use $\mathbb{H}_{\tau;U,\mathrm{loc}}$ to denote the vector space whose elements are as follows: An element in $\mathfrak{k}$ consists of a set $(\mathfrak{k}_S, \{\mathfrak{k}_\mathfrak{p}\}_{\mathfrak{p} \in \Lambda})$ where $\mathfrak{k}_S$ is a locally $L^2_1$ section of the normal bundle of $C_S \cap U$, and where any given $\mathfrak{p} \in \Lambda$ version of $\mathfrak{k}_\mathfrak{p}$ is a section of the normal bundle of $C_\mathfrak{p} \cap U$. Moreover, with the normal bundle of $C_S$ identified with the bundle $N_S$ over the $t \in [1+z_*, 2-z_*]$ part of S as in Part 1, then the $\mathfrak{k}_S$ and various $\mathfrak{p} \in \Lambda$ versions of $\mathfrak{k}_\mathfrak{p}$ obey (7.32) on the boundaries of their domain of definition.

Fix $(\mathcal{C}, \tau) \in \mathcal{M}^*$. Suppose U is a semi-bounded domain and that $\mathfrak{h} = (\mathfrak{h}_S, \{\mathfrak{h}_\mathfrak{p}\}_{\mathfrak{p} \in \Lambda})$ $\in \mathbb{H}_{\tau;U,\mathrm{loc}}$ and that it obeys an equation of the form

$$\mathcal{D}_\mathcal{C}\mathfrak{h} + \mathfrak{t}[\mathfrak{h}] + \mathfrak{z} = 0$$

(7.43)

where $\mathfrak{z}$ and $\mathfrak{t}[\cdot]$ are as follows: First, $\mathfrak{z}$ is a fixed, smooth element in $\mathbb{L}$. Meanwhile, $\mathfrak{t}[\cdot]$ = $(\mathfrak{t}_S, \{\mathfrak{t}_\mathfrak{p}\}_{\mathfrak{p} \in \Lambda})$ defines a map from a subspace of uniformly pointwise bounded elements in $\mathbb{H}_{\tau;U;loc}$ to locally square integrable sections of the appropriate bundles over $C_S$ and the various $\mathfrak{p} \in \Lambda$ versions of $C_\mathfrak{p}$. The element $\mathfrak{t}[\mathfrak{k}]$ depends explicitly on the domain coordinates and implicitly on the chosen $\mathfrak{k} \in \mathbb{H}_{\tau;U;loc}$. It has smooth dependence on the domain coordinates. Note that $\mathfrak{t}[\mathfrak{k}]$ need not be a local function of the entries of $\mathfrak{k}$. The map $\mathfrak{k} \to \mathfrak{t}[\mathfrak{k}]$ must have a certain additional requirement whose statement uses the following notation: Let W denote an open set with closure in U. Use $\|\cdot\|_W$ to denote the $L^2$ norm on W, and $\|\cdot\|_{W;1}$ to denote the $L^2_1$ norm on W. Finally, use $\|\cdot\|_{W;\infty}$ to denote the $L^\infty$ norm on W. The additional requirement involves parameters $n \geq 1$ *and* $r > 0$. What follows is the additional requirement:

*Suppose that* $W \subset U$ *is a given semi-bounded subdomain; and suppose that* $\mathfrak{h}$ *and its first derivatives are square integrable on* W. *Suppose in addition that* $\|\mathfrak{h}\|_{W;\infty} \leq n^{-1}$. *Then* $\|\mathfrak{t}[\mathfrak{h}]\|_W^2 \leq (r + n\|\mathfrak{h}\|_{W;\infty}) \|\mathfrak{h}\|_{W;1}^2$ .

(7.44)

Assume that $\mathfrak{z}$ and $\mathfrak{t}$ are as just described.

**Lemma 7.16**: *There exists* $\kappa > n$ *that such that if* $r < \kappa^{-1}$, *then the following is true: Suppose that* $\mathfrak{h}$ *is defined on a semi-bounded domain where it obeys (7.43). Suppose in addition that the entries of* $\mathfrak{h}$ *have absolute value bounded by* $\kappa^{-1}$, *and that their first derivatives are bounded. Then* $\mathfrak{h}$ *and its first derivatives are square integrable on some semi-bounded subdomain in its original domain of definition.*

*Proof of Lemma 7.16*: Write the S-labeled component of $\mathfrak{h}$ as $\eta'$ and each $\mathfrak{p} \in \Lambda$ labeled component as $(\varphi^{\mathfrak{p}'}, \varsigma^{\mathfrak{p}'})$. This done, fix $\varepsilon > 0$ and define $\mathfrak{h}_\varepsilon$ as follows: Its S-labeled component is $e^{-\varepsilon|s|}\eta'$ and any given $\mathfrak{p} \in \Lambda$ labeled component is $e^{-\varepsilon|s|}(\varphi^{\mathfrak{p}'}, \varsigma^{\mathfrak{p}'})$. Note that $\mathfrak{h}'$ obeys (7.32). The assertion that $\mathfrak{h}$ is pointwise bounded implies that $\mathfrak{h}_\varepsilon$ is square integrable on a semi-bounded domain. The derivatives of $\mathfrak{h}_\varepsilon$ are likewise square integral on some semi-bounded domain.

Let $m$ denote the supremum norm of $\mathfrak{h}$. Gvien (7.44), the fact that $\mathfrak{h}$ obeys (7.43) has the following implication: The restriction of $\mathfrak{h}_\varepsilon$ to a suitable semi-bounded subdomain $W \subset U$ obeys

$$\|\mathcal{D}_\mathcal{C}\mathfrak{h}_\varepsilon\|_W^2 \leq c_0(r + nm)\|\mathfrak{h}_\varepsilon\|_{W;1}^2 + c_0\|\mathfrak{z}\|^2$$

(7.45)

Note that all of the integrals in (7.45) are finite. Reintroduce $c$ and $s_1$ from (7.34), and suppose that $W' \subset W$ is a semi-bounded sudomain such that $|s| > 100 s_1$ at each point, and such that each point in $W'$ has distance at least 1 from some point in $W$. Then (7.34), (7.44) and (7.45) imply that

$$c^{-1} \|\mathfrak{h}_\varepsilon\|_{W;1}^2 \leq (c_0 (r + nm) - \tfrac{1}{100} c^{-4}) \|\mathfrak{h}_\varepsilon\|_{W;1}^2 + \mathfrak{w}(\mathfrak{h}) + c_0 \|\mathfrak{z}\|^2. \tag{7.46}$$

where $\mathfrak{w}(\mathfrak{h})$ is finite and $\varepsilon$-independent. If $r \leq c_0^{-1} c^{-4}$ and $m < c_0^{-1} n^{-1} c^{-4}$, then (7.46) supplies an $\varepsilon$-independent bound for $\|\mathfrak{h}_\varepsilon\|_{W;1}^2$. The existence of such a bound implies that $\mathfrak{h}$ is square integrable on an unbounded domain.

The next lemma states one consequence of Lemma 7.16. To set the stage for this lemma, decompose $\Theta_-$ into segments so as to define a set $\{\Gamma_S, \{\gamma_{\mathfrak{p}_-}\}_{\mathfrak{p} \in \Lambda}, \{\mathfrak{o}_\mathfrak{p}\}_{\mathfrak{p} \in \Lambda}\}$ where $\Gamma_S$ denotes a disjoint union of segments of integral curves of $v$ in the $f \in [1 + z_*, 2 - z_*]$ part of $M_\delta$, and where each $\mathfrak{p} \in \Lambda$ version of $\mathfrak{o}_\mathfrak{p}$ is the subset of integral curves from $\{\hat{\gamma}_\mathfrak{p}^+ \ \hat{\gamma}_\mathfrak{p}^-\}_{\mathfrak{p} \in \Lambda}$. If $T \gg 1$, then $(-\infty, -T] \times \Gamma_S$ can be written as a section of the normal bundle of $S$ over the portion where both $t \in [1 + z_*, 2 - z_*]$ and $s \leq -T$. Use $\eta_\Theta$ to denote this section. Meanwhile, each $\mathfrak{p} \in \Lambda$, version of $(-\infty, -T] \times \gamma_{\mathfrak{p}_-}$ can be written via $\Psi_\mathfrak{p}$ as the graph of a pair of functions, these denoted by $(\varphi^\mathfrak{p}_\Theta, \varsigma^\mathfrak{p}_\Theta)$. The function $\varsigma^\mathfrak{p}_\Theta$ is the constant value of $h$ on $\gamma_{\mathfrak{p}_-}$ and $\varphi^\mathfrak{p}_\Theta$ is the function of $\hat{u}$ that gives the $\phi$-angle along $\gamma_{\mathfrak{p}_-}$. An analogous $\eta_\Theta$ and set $\{(\varphi^\mathfrak{p}_\Theta, \varsigma^\mathfrak{p}_\Theta)\}_{\mathfrak{p} \in \Lambda}$ are defined using $\Theta_+$ with it understood that these are defined respectively where both $t \in [1 + z_*, 2 - z_*]$ and $s \gg 1$ on $S$ and where $x \gg 1$ on each $\mathfrak{p} \in \Lambda$ version of $\mathbb{R} \times I_*$.

Fix next a smooth function on $(\tfrac{1}{2} z_*, \tfrac{3}{2} z_*)$ with compact support and which is equal to 1 where $z = z_*$. Use $\beta$ here to denote the chosen function. Fix $\mathfrak{p} \in \Lambda$ and introduce the functions $(\varphi^{S0}, \varsigma^{S0})$ on $\mathbb{R}$ that appear in (7.1). This data is used to modify $\eta_\Theta$ and each $\mathfrak{p} \in \Lambda$ version of $(\varphi^\mathfrak{p}_\Theta, \varsigma^\mathfrak{p}_\Theta)$ near the common boundary of their domains of definition. The modification requires writing $\eta_\Theta$ where $t$ differs from 1 or 2 by less than $z_S$ as a pair of constant functions, $(\varphi^S_\Theta, \varsigma^S_\Theta)$, on $\mathbb{R} \times [z_*, z_S)$. Likewise, view $(\varphi^\mathfrak{p}_\Theta, \varsigma^\mathfrak{p}_\Theta)$ near the boundary of $\mathbb{R} \times I_*$ as constant functions on $\mathbb{R} \times [\delta^2, z_*]$. The modification is given by the replacements

- $(\varphi^S_\Theta, \varsigma^S_\Theta) \to (1 - \beta)(\varphi^S_\Theta, \varsigma^S_\Theta) + \beta(\varphi^{S0}, \varsigma^{S0})$.
- $(\varphi^\mathfrak{p}_\Theta, \varsigma^\mathfrak{p}_\Theta) \to (1 - \beta)(\varphi^\mathfrak{p}_\Theta, \varsigma^\mathfrak{p}_\Theta) + \beta(\varphi^{S0}, \varsigma^{S0})$.

(7.47)

Let $\hat{\eta}_\Theta$ and $\{(\hat{\varphi}^\mathfrak{p}_\Theta, \hat{\varsigma}^\mathfrak{p}_\Theta)\}_{\mathfrak{p} \in \Lambda}$ denote the modified versions. These obey (7.1).

To make contact with Lemma 7.16, fix $(\tau, \mathcal{C}) \in \mathcal{M}^*$ and let $\eta$ and $\{(\varphi^{\mathfrak{p}}, \varsigma^{\mathfrak{p}})\}_{\mathfrak{p} \in \Lambda}$ denote the defining data for $C$, the former as described in Part 1 of Section 1a and each $\mathfrak{p} \in \Lambda$ version of $(\varphi^{\mathfrak{p}}, \varsigma^{\mathfrak{p}})$ as described in Part 2 of Section 1a. Fix $x_* > 1$ so as to be larger than the absolute value of the $\mathbb{R}$ coordinates of the deleted $\hat{u} = 0$ points from $\mathbb{R} \times I_*$ that define the domain of the $\Delta_{\mathfrak{p}} > 0$ versions of $(\varphi^{\mathfrak{p}}, \varsigma^{\mathfrak{p}})$.

Define $\mathfrak{h}_{\mathcal{C}} = (\mathfrak{h}_{\mathcal{C},S}, \{\mathfrak{h}_{\mathcal{C},\mathfrak{p}}\}_{\mathfrak{p} \in \Lambda})$ with $\mathfrak{h}_{\mathcal{C},S}$ the section of $N_S$ over the $t \in [1+z_*, 2-z_*]$ part of S given by $\eta - \hat{\eta}_\Theta$, and with any given $\mathfrak{p} \in \Lambda$ version of $\mathfrak{h}_{\mathfrak{p}}$ being the map from the domain of $(\varphi^{\mathfrak{p}}, \varsigma^{\mathfrak{p}})$ given by $(x, \hat{u}) \to \chi(x_* - |x| + 2) (\varphi^{\mathfrak{p}} - \hat{\varphi}^{\mathfrak{p}}{}_\Theta, \varsigma^{\mathfrak{p}} - \hat{\varsigma}^{\mathfrak{p}}{}_\Theta)|_{(x,\hat{u})}$. Viewed $\mathfrak{h}_{\mathcal{C},S}$ as a section of the normal bundle of $C_S$ and view each $\mathfrak{p} \in \Lambda$ version of $\mathfrak{h}_{\mathcal{C},\mathfrak{p}}$ as a section of the normal bundle of $C_{\mathfrak{p}}$. This done, then $\mathfrak{h}_{\mathcal{C}}$ has support on two disjoint, semi-bounded domains, one where $s \gg 1$ and the other where $s \ll -1$. It obeys (7.32) on each. Moreover, given $r > 0$, there exists $s_r > 1$ such that $|\mathfrak{h}_{\mathcal{C}}| \leq r$ where $|s| \geq s_r$.

This $\mathfrak{h}_{\mathcal{C}}$ also obeys a version of (7.43) on an $s \gg 1$ semi-bounded domain, and also on an $s \ll -1$ semi-bounded domain. This is a consequence of (II.6.10) and (3.4). The function $\mathfrak{h} \to \mathfrak{t}[\mathfrak{h}]$ in this case is a local function of $\mathfrak{h}$ and its derivatives, this is to say that its value at any given point is determined by the point in question and the value of the relevant component of $\mathfrak{h}$ and its derivatives at this same point. In particular the relevant version of $\mathfrak{t}[\cdot]$ has the following schematic form

$$\mathfrak{t}[\mathfrak{h}] = \mathfrak{t}_1[\mathfrak{h}] \cdot \nabla \mathfrak{h} + \mathfrak{t}_0[\mathfrak{h}],$$

(7.48)

where any given component of $\mathfrak{t}_1$ and $\mathfrak{t}_0$ is a non-linear, local, smooth function of the corresponding component of $\mathfrak{h}$ and the relevant domain coordinates. These functions are such that

- $|\mathfrak{t}_1(\mathfrak{h})| \leq c_0 |\mathfrak{h}|$ and $|\mathfrak{t}_0(\mathfrak{h})| \leq c_0 |\mathfrak{h}|^2$.
- *The first derivatives of $\mathfrak{t}_0(\cdot)$ at $\mathfrak{h} = 0$ with respect to variations in $\mathfrak{h}$ and the domain variables are zero; and the analogous first derivatives of $\mathfrak{t}_1$ and the second derivatives of $\mathfrak{t}_0$ are bounded by $c_0$ where $|\mathfrak{h}| \leq c_0^{-1}$.*
- *In general, the derivatives of $\mathfrak{t}_0$ and $\mathfrak{t}_1$ with respect to variations in $\mathfrak{h}$ and the domain variables to any given order are bounded if $|\mathfrak{h}| \leq c_0^{-1}$*

(7.49)

In addition, the first derivatives of $\mathfrak{t}_0(\cdot)$ at $\mathfrak{h} = 0$ are zero; and the first derivatives of $\mathfrak{t}_1$ and the second derivatives of $\mathfrak{t}_0$ are bounded by $c_0$ where $|\mathfrak{h}| \leq c_0^{-1}$.

Granted this background, what follows gives the first of the promised applications of Lemma 7.16.

**Lemma 7.17**: *Let $(\tau, \mathcal{C}) \in \mathcal{M}^*$ and define $\mathfrak{h}_{\mathcal{C}}$ as above using either $\Theta_-$ or $\Theta_+$. Then $\mathfrak{h}_{\mathcal{C}}$ and its derivatives to any given order are square integrable on a semi-bounded domain.*

***Proof of Lemma 7.17***: The fact that $\mathfrak{h}_{\mathcal{C}}$ is square integrable follows from Lemma 7.17. The assertion concerning its derivatives follows using standard elliptic regularity techniques with Lemma 7.7 as applied to the version of (7.43) that uses (7.48) for $\mathfrak{t}$. (As before, the elliptic regularity techniques from Chapter 6 in [M] will do the trick.)

*Part 4*: This part of the subsection contains the proofs of Propositions 7.1 and 7.2. These are taken in order.

***Proof of Proposition 7.1***: The proof has four steps.

Step 1: This step sets the ground work for an application of the implicit function theorem. To this end, fix $(\tau, \mathcal{C}) \in \mathcal{M}^*$ and let $(\tau', \mathcal{C}') \in \mathcal{M}_\tau$ denote an element in a neighborhood of $(\tau, \mathcal{C})$ with the neighborhood chosen so that $\tau'$ is very close to $\tau$ and so that each point in any submanifold from $\mathcal{C}'$ is very close to the corresponding submanifold from $\mathcal{C}$ and vice-versa. Let $\{\eta, \{(\varphi^\mathfrak{p}, \varsigma^\mathfrak{p})\}_{\mathfrak{p} \in \Lambda}\}$ denote the data that defines $\mathcal{C}$ with $\eta$ as described in Part 1 of Section 7a and with each $\mathfrak{p} \in \Lambda$ version of $(\varphi^\mathfrak{p}, \varsigma^\mathfrak{p})$ as described in Part 2 of Section 7a. A corresponding data set is used to define $(\tau', \mathcal{C}')$, but the latter is written now as $\{\eta + \eta', \{(\varphi^\mathfrak{p} + \varphi^{\mathfrak{p}'}, \varsigma^\mathfrak{p} + \varsigma^{\mathfrak{p}'})\}_{\mathfrak{p} \in \Lambda}\}$. It is a consequence of the final assertion of Lemma 7.18 that $\eta'$ when viewed as a section over $C_S$ of the latter's normal bundle is a smooth, $L^2_1$ section. Meanwhile, Lemma 7.18 with the final assertion of Proposition 5.1 imply that each $\mathfrak{p} \in \Lambda$ version of $(\varphi^{\mathfrak{p}'}, \varsigma^{\mathfrak{p}'})$ defines a smooth, $L^2_1$ section of the normal bundle of $C_\mathfrak{p}$ in $\mathbb{R} \times \mathcal{H}^+_{\mathfrak{p}*}$.

Use $\mathfrak{h}$ to to denote the set $\{\eta', (\varphi^{\mathfrak{p}'}, \varsigma^{\mathfrak{p}'})\}_{\mathfrak{p} \in \Lambda}$. This version of $\mathfrak{h}$ does not obey (7.32) unless $\tau = \tau'$; but if it did, then it would define an element in $\mathbb{H}_{\tau*}$, this a consequence of Lemma 7.18 and Proposition 5.1. In any event, $\mathfrak{h}$ obeys a version of (7.43) with $\mathfrak{t}[\cdot]$ given by (7.48) with $\mathfrak{t}_0$ and $\mathfrak{t}_1$ described by (7.49).

The failure of (7.32) is rectified in the next step by encorporating $\tau' - \tau$ in a new definition of $\mathfrak{h}$ and compensating with a corresponding $\tau' - \tau$ dependent term added to $\mathfrak{t}[\cdot]$. The resulting version of $\mathfrak{t}$ is not a local function of $\mathfrak{h}$.

Step 2: This step supplies the new definitions of $\mathfrak{h}$ and $\mathfrak{t}$. To start, fix $\mathfrak{p} \in \Lambda$ and write $\eta$ near the common boundary of $C_\mathfrak{p}$ and $C_S$ as functions $(\varphi^S, \varsigma^S)$ on $\mathbb{R} \times [z_*, z_S]$ in the manner of (7.1). Define $\eta'$ as above and write the latter near this same boundary as a pair of functions $(\varphi^{S'}, \varsigma^{S'})$ on $\mathbb{R} \times [z_*, z_S]$. View $(\varphi^\mathfrak{p}, \varsigma^\mathfrak{p})$ and also the pair $(\varphi^{\mathfrak{p}'}, \varsigma^{\mathfrak{p}'})$ near this boundary as functions on $\mathbb{R} \times [\delta^2, z_*]$. In this guise, the primed pairs obey the following where $z = z_*$:

$$\varsigma^{S\prime} = \tau \varsigma^{p\prime} + (\tau\prime - \tau)(\varsigma^p + \varsigma^{p\prime} - \varsigma^{S0}) \quad and \quad \varphi^{p\prime} = \tau \varphi^{S\prime} + (\tau\prime - \tau)(\varphi^S + \varphi^{S\prime} - \varphi^{S0})$$
(7.51)

To obtain something that obeys (7.32), reintroduce the function β from the proof of Lemma 7.17 and (7.47). Define pairs $(\hat{\varphi}^{S\prime}, \hat{\varsigma}^{S\prime})$ and $(\hat{\varphi}^{p\prime}, \hat{\varsigma}^{p\prime})$ as follows:

- $\hat{\varphi}^{S\prime} = \varphi^{S\prime}$ and $\hat{\varsigma}^{S\prime} = \varsigma^{S\prime} - (\tau\prime - \tau)\beta(\varsigma^p + \varsigma^{p\prime} - \varsigma^{S0})|_{z_* - (z-z_*)}$ where $z \geq z_*$.
- $\hat{\varphi}^{p\prime} = \varphi^{p\prime} - (\tau\prime - \tau)\beta(\varphi^S + \varphi^{S\prime} - \varphi^{S0})|_{z_* - (z-z_*)}$ where $z \leq z_*$.

(7.52)

Extend the latter as η´ and $(\varphi^{p\prime}, \varsigma^{p\prime})$ over the rest of their respective domains. Use these to define the promised new version of $\mathfrak{h}$. The entries of this new version obey (7.32) and it follows from Lemma 7.17 that this new version is in $\mathbb{H}_{\tau*}$.

This new version also obeys a version of (7.43), but with a non-local version of $\mathfrak{t}[\cdot]$ and also with a $\mathfrak{z} \neq 0$ term proportional to $(\tau\prime - \tau)$. To elaborate, the non-local version of $\mathfrak{t}$ is obtained from the version from Step 1 by adding a term that is supported near the boundaries of $C_S$ and $\cup_{\mathfrak{p} \in \Lambda} C_\mathfrak{p}$. Fix a given $\mathfrak{p} \in \Lambda$. The portion of this term that lies where $z > z_*$ is the pair of functions with respective left and right hand components

$$(\tau\prime - \tau) \partial_z(\beta \, \varsigma^{p\prime}|_{z_* - (z-z_*)}) \quad and \quad -(\tau\prime - \tau) \beta (\partial_x \varsigma^{p\prime})|_{z_* - (z-z_*)} \, .$$
(7.53)

The part that lies where $z < z\prime$ has components

$$-(\tau\prime - \tau) \beta (\partial_x \varphi^{S\prime})|_{z_* - (z-z_*)} \quad and \quad (\tau\prime - \tau) \partial_z(\beta \, \varphi^{S\prime}|_{z_* - (z-z_*)}) \, .$$
(7.54)

Meanwhile, the $\mathfrak{z}$ term has support near these same boundaries where it is given by replacing $\varsigma^{p\prime}$ with $(\varsigma^p - \varsigma^{S0})$ in (7.53) and $\varphi^{S\prime}$ with $(\varphi^S - \varphi^{S0})$ in (7.54).

Step 3: Given the preceding definition of $\mathfrak{t}$, it follows that the left hand side of the corresponding version of (7.43) defines a smooth map to $\mathbb{L}_*$ from the product of a ball about the origin in $\mathbb{H}_{\tau*}$ with an interval centered on $\tau \in [0, 1]$. Use $\mathfrak{F}$ to denote this map. The inverse function theorem finds a ball, $\mathcal{B}$, about the origin in the kernel of $\mathcal{D}_C$, an interval, $I$, centered on τ in [0, 1], a smooth map, $\mathfrak{b}: \mathcal{B} \times I \to \mathbb{H}_{\tau*}$, and these such that

- $\mathfrak{b}(0, \tau) = 0$ and $\nabla^{\mathbb{H}} \mathfrak{b}|_{(0,\tau)} = 0$.
- *Let* $(\mathfrak{h}_0, \tau\prime) \in \mathcal{B} \times I$. *Then* $\mathfrak{h}_0 + \mathfrak{b}(\mathfrak{h}_0, \tau\prime) \in \mathbb{B}$ *and* $(1 - \prod)\mathfrak{F}(\mathfrak{h}_0 + \mathfrak{b}(\mathfrak{h}_0, \tau\prime), \tau\prime) = 0$.
- *Let* $\mathbb{B}\prime \subset \mathbb{B}$ *denote the concentric, half radius ball. Suppose that* $(\mathfrak{h}, \tau\prime) \in \mathbb{B}\prime \times I$ *and suppose that* $(1 - \prod)\mathfrak{F}(\mathfrak{h}\prime, \tau\prime) = 0$, *then* $\mathfrak{h} = \mathfrak{h}_0 + \mathfrak{b}(\mathfrak{h}_0, \tau\prime)$ *with* $\mathfrak{h}_0 \in \mathcal{B}$.

(7.55)

A map, $\mathfrak{f}$, from $\mathcal{B} \times \mathrm{I}$ to the cokernel of $\mathcal{D}_C$ is defined now by the rule

$$(\mathfrak{h}_0, \tau') \to \mathfrak{f}((\mathfrak{h}_0, \tau')) = \Pi \mathfrak{F}(\mathfrak{h}_0 + \mathfrak{b}(\mathfrak{h}_0, \tau'), \tau') \,.$$

(7.56)

The map from $\mathcal{B} \times \mathrm{I}$ to $\mathbb{B} \times \mathrm{I}$ given by $(\mathfrak{h}_0, \tau') \to (\mathfrak{h}_0 + \mathfrak{b}(\mathfrak{h}_0, \tau'), \tau',)$ embeds $\mathfrak{f}^{-1}(0) \subset \mathcal{B}$ homeomorphically onto an open set in $\mathbb{B} \times \mathrm{I}$ of solutions to (7.43) that contains $\mathbb{B}' \times \mathrm{I}$.

This map $\mathfrak{f}$ is the map required by Proposition 7.1, and the embedding just described gives the homeomorphism $\Phi$.

Step 4: The claim made by the second bullet of Proposition 7.1 that $\mathcal{M}^*$ is smooth where $\mathfrak{f}$ is a submernsion and the claim that $\pi_\mathrm{I}$ is a smooth on this same set are standard consequences of the inverse function theorem as used in Step 3. The proof that $p$ is continuous follows from Proposition 5.1 and Lemma 5.2.

To see about the derivatives of $p$, let $\mathfrak{p} \in \Lambda$ be such that $\Delta_\mathfrak{p} > 0$. For the sake of argument, suppose that $\mathcal{E} \subset C_\mathfrak{p}$ is an end where the $s \gg 1$ portion maps to the tubular neighborhood $U_+ \subset \mathcal{H}_\mathfrak{p}$ of the integral curve $\hat{\gamma}_\mathfrak{p}^+$ via the projection from $\mathbb{R} \times \mathcal{H}^+_{\mathfrak{p}*}$. Use the coordinates $(s_+, \phi_+, \theta_+, u_+)$ for $\mathbb{R} \times U_+$ and parametrize the $s \gg 1$ part of $\mathcal{E}$ as in Proposition 5.1. It is a consequence of Proposition 5.1 that the operator $D_\mathfrak{p}$ on $\mathcal{E}$ appears as an operator on the space of $\mathbb{R}^2$ valued functions of the coordinates $(s_+, \phi_+)$. An essentially verbatim repeat of the proof of the first bullet of Lemma 5.8 proves that the latter has the form $\mathfrak{D}_0 + \mathfrak{d}$ where $\mathfrak{d}$ is a a first order differential operator on the space of maps from to $\mathbb{R}^2$ whose coefficients are bounded in absolute value by $c_0 e^{-s/c_0}$.

Now let $\mathfrak{h}_0 \in \mathcal{B} \subset \ker(\mathcal{D}_C)$ denote a given element, and let $r$ denote its norm. The $C_\mathfrak{p}$ component of $\mathfrak{h}_0$ appears using this parametrization as a square integrable map, $\mathfrak{y}_\mathfrak{p}$, from the very large $s_+$ part of $\mathbb{R} \times \mathbb{R}/(2\pi\mathbb{Z})$ to $\mathbb{R}^2$ that obeys an equation that has the schematic form $\mathfrak{D}_0 \mathfrak{y}_\mathfrak{p} + \mathfrak{d} \mathfrak{y}_\mathfrak{p} = 0$. Granted that this is so, then the techniques used in Section 2.3 in [HT] can be employed to prove that $\mathfrak{y}_\mathfrak{p}$ can be written as in

$$\mathfrak{y}_\mathfrak{p} = r((c\, e^{-\lambda_1 s_+}, 0) + \mathfrak{e})$$

(7.57)

where $c \in \mathbb{R}$ and where $|\mathfrak{e}| \leq c_0\, e^{-(\lambda_1 + 1/c_0) s_+}$. Meanwhile, (7.43), (7.48) and (7.49) and the techniques from Section 2.3 in [HT] can be employed in a straightforward manner to see that $\mathfrak{b} = \mathfrak{b}(\mathfrak{h}_0, \tau')$ from (7.55) can be written on $\mathcal{E}$ as a square integrable map from the very large $s_+$ part of $\mathbb{R} \times \mathbb{R}/(2\pi\mathbb{Z})$ to $\mathbb{R}^2$ that is bounded by $c_0 (r^2 + \tau'^2) e^{-\lambda_1 s_+}$. Granted Lemma 5.2, this last observation implies that $p$ is a $C^1$ map on $\mathcal{M}^*_{\mathrm{smooth}}$. The proof that $p$ has derivatives to any given order has a similar flavor and is omitted.

The assertions of the third bullet follows from (7.56) because the latter equation depicts a version of $\mathfrak{f}$ with the property that $|\mathfrak{f}(\mathfrak{h}_0)| \leq c_0 |\mathfrak{h}_0|^2$. It follows as a consequence, that $\mathfrak{f}$ is a submersion at $(0,\tau)$ only if $\dim(\text{cokernel}(\mathcal{D}_\mathcal{C}) \leq 1$; and that if $\mathfrak{f}$ is a submersion at $(0,\tau) \in \mathcal{B} \times I$ and $d\pi_I \neq 0$ at $(0,\tau)$, then $\text{cokernel}(\mathcal{D}_\mathcal{C}) = 0$.

The proof of the fourth bullet starts with Propositions 2.1 and 2.2 for they jointly asserted that $\pi_I^{-1}(0)$ is mapped diffeomorphically by the map $p$ to $\times_{p \in \Lambda}(\times_{\Delta_p} \mathbb{R})$ given bounds for $z_*$ and $\delta$ that are purely S-dependent (or $\mathcal{K}$-compatible). Meanwhile, Propositions 2.1, 3.1 and 5.7 imply that the $\mathcal{D}_\mathcal{C}$ along $\pi_I^{-1}(0)$ has trivial cokernel. This fact with (7.56) implies what is asserted by the fourth bullet.

*Proof of Proposition 7.2* The assertion made by the first bullet is proved in Section 3 of [L]. To prove Item a) of the second bullet, note first that given Propositions 3.1 and 5.7, it is enough to consider the case for the part of $\mathcal{M}^*$ where $\tau > 0$. With this restriction on $\tau$ understood, the allowed variations of the almost complex structure on the portion of $\mathbb{R} \times Y$ in the $f^{-1}([1+\delta_*^2, 2-\delta_*^2])$ part of $\mathbb{R} \times M_\delta$ constitute a sufficiently large set for applying standard Smale-Sard arguments as done in Section 3 of [L]. In particular, straightforward modifications to the arguments from this same section of [L] prove what is asserted by this item.

By way of a parenthetical remark, note that the variations in J that are allowed on any $\cup_{p \in \Lambda} \mathbb{R} \times \mathcal{H}_p$ may not form a set that is large enough to invoke the Smale-Sard theorem. This is because the almost complex structures here are constrained to be invariant with respect to both the group of constant translations along the $\mathbb{R}$ factor and the group of constant rotations of the angle $\phi$. The set of allowed variations of J on the $f \in (1+\delta_*^2, 2-\delta_*^2)$ part of $\mathbb{R} \times M_\delta$ is sufficiently large precisely because Lipshitz allows almost complex structures on $\mathbb{R} \times [1, 2] \times \Sigma$ that depend on the coordinate $t \in [1, 2]$. In fact, the set of $t$-independent almost complex structures is likely too small for the applications in Section 3 of [L]. The assertion made by Item b) of the second bullet is proved using the Smale-Sard theorem using the aforementioned arguments from [L]. The details are also straightforward and also omitted.

### 8) Counting ech-HF submanifolds

The section starts with an existence assertion for ech-HF submanifolds, and then a sort of uniqueness assertion. The existence result is stated as Proposition 8.1 and the uniqueness result is stated by Proposition 8.2. These two propositions are in Section 8a. Sections 8b-f explains, among other things, how to count the ech-HF submanifolds that are provided by Proposition 8.1.

### a) Existence and uniqueness of ech-HF subvarieties

Fix a countable set in $(\times_3 (0,1)) \times (1, \infty)$ of possible choices for the data $(z_*, \delta, x_0, R)$ and then choose $J_{HF}$ from Proposition 7.2's residual set.

To set the stage for the upcoming Propositions 8.1 and 8.2, fix a finite or weakly compact subset $\mathcal{K} \subset \mathcal{A}_{HF}$. With $\mathcal{K}$ in hand, select the data set $(z_*, \delta, x_0, R)$ and the almost complex structure J as described by Propositions 7.1-7.3 so that their conclusions can be assumed.

Proposition 8.1 assumes implicitly that a submanifold S has been chosen from $\mathcal{K}$ and choice of $(\hat{\Theta}_-, \hat{\Theta}_+) \in \hat{\mathcal{Z}}^S$ has been made so as to define the corresponding version of $\mathcal{M}^*$. By way of a reminder, $\mathcal{M}^*$ is a smooth manifold with boundary and $\pi_1 \times p$ is a smooth, proper map. Supposing that $y \in \times_{p \in \Lambda} (\times_{\Delta_p} \mathbb{R})$ is a regular value of the map $p$, use $\mathcal{M}^*_y$ to denote $p^{-1}(y)$. The latter is a smooth, 1-dimensional manifold with boundary. The fourth bullet of Proposition 7.1 asserts that there is one and only one component with a boundary point on $\pi_1^{-1}(0)$. Meanwhile, it follows from what is said in Proposition 7.2 that the differential of $\pi_1$ at each point in $\pi_1^{-1}(1)$ is surjective, and it follows from Propositions 7.2 and 7.3 that there are at most a finite set of points in $\pi_1^{-1}(1)$.

The manifold $\mathcal{M}^*_y$ is orientable because this is the case for any 1-manifold. Orientations for the components of $\mathcal{M}^*_y$ are defined by requiring that $\pi_1$ be orientation preserving where it is increasing and orientation reversing where it is decreasing.

Assign to any $\pi_1^{-1}(1)$ point in $\mathcal{M}^*_y$ the weight +1 if the differential of $\pi_1$ at the point is orientation preserving, and assign -1 if not.

**Proposition 8.1**: *The set $\pi_1^{-1}(1)$ in $\mathcal{M}^*_y$ is non-empty. Moreover, the sum of the ±1 weights of these elements is equal to* 1.

Given that the $\pi_1^{-1}(1)$ points in $\mathcal{M}^*_y$ are ech-HF submanifolds, this proposition supplies an existence theorem for ech-HF submanifolds.

*Proof of Proposition 8.1*: If $z_*$ and $\delta$ are chosen small, then Proposition 7.3 asserts that $\mathcal{M}^*_y$ is a compact, oriented 1-manifold with boundary. Each component has either 0 or 2 boundary points. There is one component with a boundary point where $\pi_1 = 0$ and the latter must have a second boundary point, thus where $\pi_1 = 1$. As $\pi_1$ is in no case greater than 1 and as its differential is non-zero at this point, so this point has weight +1. There is a finite set of other components. Those with boundary points must have both boundary points where $\pi_1 = 1$. As the differential of $\pi_1$ is non-zero at both points, one must have weight +1 and the other weight -1. Granted this accounting, then the sum of the weights of the elements in $\mathcal{M}^*_y$ is equal to 1.

Proposition 8.2 uses $\mathcal{K}$ to denote the following subset of $\mathcal{A}_{HF}$: A Lipshitz submanifold S is in $\mathcal{K}$ if and only if the operator $D_S$ has Fredholm index no greater than 1. The quotient space $\mathcal{K}/\mathbb{R}$ is finite, this a consequence of Lemma 5.4 and Corollary 7.2 in [L]. To say more about notation, suppose that $(\hat{\Theta}', \hat{\Theta})$ is a chosen pair from $\times_2 \hat{\mathcal{Z}}_{ech,M}$. Given a data set $(\delta, x_0, R)$ to define the geometry of Y, and given an almost complex structure, J subject to the constraints in Part 1 of Section 1c, the proposition refers to the space $\mathcal{M}_1(\hat{\Theta}', \hat{\Theta})$ defined in Part 2 of Section 1c. Any given element in $\mathcal{M}_1(\hat{\Theta}', \hat{\Theta})$ will have a union of components that comprise an ech-HF submanifold. This ech-HF submanifold part is either $\mathbb{R}$-invariant or not. If the ech-HF submanifold is $\mathbb{R}$-invariant, then there must be a single component from some $\mathfrak{p} \in \Lambda$ version of Proposition II.3.4's moduli spaces $\mathcal{M}_{\mathfrak{p}-}$ and $\mathcal{M}_{\mathfrak{p}+}$. There can also be $\mathbb{R}$-invariant cylinder components from the set $\cup_{\mathfrak{p} \in \Lambda}\{\mathbb{R} \times \hat{\gamma}_\mathfrak{p}^+, \mathbb{R} \times \hat{\gamma}_\mathfrak{p}^-\}$. If the ech-HF submanifold is not $\mathbb{R}$ invariant, then it sits in some $(\hat{\Theta}_-, \hat{\Theta}_+) \in \times_2 \hat{\mathcal{Z}}_{ech,M}$ version of $\mathcal{M}_1(\hat{\Theta}_-, \hat{\Theta}_+)$ with $(\hat{\Theta}_-, \hat{\Theta}_+)$ as described in (2.2).

**Proposition 8.2**: *Fix the data* $(z_*, \delta, x_0, R)$ *and J as described by Propositions 7.1-7.3 with a suitably large choice for their respective versions of* $\kappa$ *and* $\kappa_*$, *and with any choice of Lipshitz submanifold from* $\mathcal{K}$. *Fix* $(\hat{\Theta}_-, \hat{\Theta}_+) \in \times_2 \hat{\mathcal{Z}}_{ech,M}$ *that obey (2.2) and are such that* $\mathcal{M}_1(\hat{\Theta}_-, \hat{\Theta}_+) \neq \emptyset$. *Then* $\sum_{\mathfrak{p} \in \Lambda} \Delta_\mathfrak{p} \leq 1$ *and what follows is true: Let* C *denote a given ech-HF submanifold from* $\mathcal{M}_1(\hat{\Theta}_-, \hat{\Theta}_+)$. *There exists a unique* $S \in \mathcal{K}$ *such that* $(\hat{\Theta}_-, \hat{\Theta}_+) \in \hat{\mathcal{Z}}^S$; *and there exists a unique* $(1, \mathcal{C})$ *in the corresponding version of* $\mathcal{M}^*$ *whose respective components are C's intersections with the* $f^{-1}(1+z_*, 2-z_*)$ *portion of* $\mathbb{R} \times Y$ *and the various* $\mathfrak{p} \in \Lambda$ *versions of* $\mathbb{R} \times \mathcal{H}^+_{\mathfrak{p}*}$.

*Proof of Proposition 8.2*: Delete the assertion $\sum_{\mathfrak{p} \in \Lambda} \Delta_\mathfrak{p} \leq 1$ from Proposition 8.2, and suppose that the resulting weaker proposition is false. If this is the case, then there is a sequence $\{D_n, C_n\}_{n=1,2,...}$ whose constituents will now be described. First, what is denoted by $D_n$ is a data set that can written as $((\hat{\Theta}_{n-}, \hat{\Theta}_{n+}), (z_{*n}, \delta_n, x_{0n}, R_n, J_n))$ where $z_{*n} < 1/n$ and $\delta_n < n^{-2} z_{*n}$. The latter with $x_{0n}$ and $R_n$ are suitable for defining the geometry of Y. Meanwhile, $J_n$ is an almost complex structure on the $(\delta_n, x_{0n}, R_n)$ version of $\mathbb{R} \times Y$ as described in Section 1c. In addition $J_n$ with $(z_{*n}, \delta_n, x_{0n}, R_n)$ are such that Propositions 7.1-7.3 can be invoked using any Lipshitz submanifold from $\mathcal{K}$. Meanwhile, $\hat{\Theta}_{n-}$ and $\hat{\Theta}_{n+}$ are elements in the index n version of $\hat{\mathcal{Z}}_{ech,M}$ that are defined in part by respective HF cycles that can be assumed to be independent of the index n. The $\hat{\Theta}_{n-}$ and also $\hat{\Theta}_{n+}$ elements from $\{\hat{\gamma}_\mathfrak{p}^+ \hat{\gamma}_\mathfrak{p}^-\}_{\mathfrak{p} \in \Lambda}$ are also independent of n. What is denoted by $C_n$ signifies an ech-HF submanifold defined by the data $D_n$ from a submanifold in $\mathcal{M}_1(\hat{\Theta}_{n-}, \hat{\Theta}_{n+})$ but $C_n$ is not from some $\pi_1^{-1}(1)$ element in the index n version of $\mathcal{M}^*$. No generality is lost by assuming that $C_n$'s version of the set $\{\Delta_\mathfrak{p}\}_{\mathfrak{p} \in \Lambda}$ is independent of n. The eight steps that follow derive nonsense with such a sequence.

The assertion that $\sum_{\mathfrak{p} \in \Lambda} \Delta_\mathfrak{p} \leq 1$ follows from Proposition 7.1 if $\mathcal{M}_1(\hat{\Theta}_-, \hat{\Theta}_+)$ has an ech-HF submanifold from the corresponding version of $\mathcal{M}^*$.

Step 1: The submanifold $C_n$ has a normal bundle which also inherits a holomorphic line bundle structure. Use $N_{C_n}$ to denote this bundle. There is an associated first order operator that maps sections of $N_{C_n}$ to sections of $N_{C_n} \otimes T^{0,1}C_n$. This operator also has the form depicted on the right hand side of (1.25). Use $\hat{\mathcal{D}}_{C_n}$ to denote this operator. This is Fredholm when mapping the $L^2_1$ space of sections $N_{C_n}$ to the $L^2$ space of sections of $N_C \otimes T^{0,1}C$. These respective domain and range spaces are denoted by $\mathbb{H}_1$ and $\mathbb{L}$ in what follows. Note that this Fredholm incarnation of $\hat{\mathcal{D}}_{C_n}$ has index 1 and trivial cokernel.

Step 2: What follows is a consequence of the assumption that $\{z_{*n}\}_{n=1,2,\ldots}$ has limit zero: Given $\varepsilon$, the conclusions of Proposition II.7.2 can be invoked for all sufficiently large n versions of $C_n$. This being the case, fix some small $z_* > 0$. For n large, $z_{*n}$ will be less than $z_*$. For such n, use $\mathcal{H}^+_{\mathfrak{p}*n}$ to denote the version of $\mathcal{H}^+_{\mathfrak{p}*}$ that is defined using $z_*$ and the data set $(\delta_n, x_{0n}, R_n)$. Use $\Psi_{\mathfrak{p}n}$ to denote the version of the map $\Psi_\mathfrak{p}$ that is defined using $z_*$ and the data set $(\delta_n, x_{0n}, R_n, J_n)$.

Proposition II.5.8 and Lemma II.4.7 can be invoked when n is large to conclude the following: Fix $\mathfrak{p} \in \Lambda$ and let $C_{\mathfrak{p}n}$ denote the intersection between $C_n$ and $\mathbb{R} \times \mathcal{H}^+_{\mathfrak{p}*n}$. This is a smooth, properly embedded submanifold with boundary. Moreover, it is the image via the map $\Psi_{\mathfrak{p}n}$ of a graph in the $(z_*, \delta_n, x_{0n}, R_n, J_n)$ version of $\mathbb{R} \times \mathcal{X}$ which is defined by a map from to $\mathbb{R}^2$ from the complement of $\Delta_\mathfrak{p}$ points in $\mathbb{R} \times I_*$ where $\hat{u} = 0$. This map has the form

$$(x, \hat{u}) \to (x, \hat{u}, \hat{\phi} = \varphi^{\mathfrak{p}n}(x, \hat{u}), h = \varsigma^{\mathfrak{p}n}(x, \hat{u})) .$$

(8.1)

Note that $C_{\mathfrak{p}n}$ obeys all of the requirements listed in Part 2 of Section 7a. Therefore, the fact that $C_n$ is not from the index n version of $\mathcal{M}^*$ is not due to properties of its intersection with $\mathbb{R} \times \mathcal{H}^+_{\mathfrak{p}*}$.

Step 3: With n large assumed large, introduce $C_{Sn}$ to denote $C_n$'s intersection with the $t \in [1+z_*, 2-z_*]$ part of $\mathbb{R} \times M_{\delta_n}$. When viewed in $\mathbb{R} \times [1+z_*, 2-z_*] \times \Sigma$, this $C_{Sn}$ is a smooth, properly embedded submanifold with boundary. A neighborhood of each component of the bundary of $C_{Sn}$ can be depicted as a graph of the sort that is described by PROPERTY 4 in Section 7a with the constant $z_S$ being n-dependent now. Use $(\varphi^{Sn}, \varsigma^{Sn})$ to denote $C_{Sn}$'s version of the functions $(\varphi^S, \varsigma^S)$ that appear in this PROPERTY 4. Note in particular that (7.1) holds by virtue of the fact that $C_{Sn}$ attaches seemlessly along its boundaries with the boundaries of $\cup_{\mathfrak{p} \in \Lambda} C_{\mathfrak{p}n}$ to give the surface $C_n$. It follows from this last remark that $C_n$ obeys the requirements from Part 3 in Section 7a for membership in the index n version of $\mathcal{M}^*$.

Granted this, and granted what is said in the final paragraph of Step 2, then $C_n$'s lack of membership in the index n version of $\mathcal{M}^*$ must be due to some property of $C_{Sn}$. To see what this might be, view $C_{Sn}$ as a submanifold in $\mathbb{R} \times [1+z_*, 2-z_*] \times \Sigma$, and suppose that there exists a Lipshitz submanifold, $S \in \mathcal{K}$ such that PROPERTY 1 in Part 1 of Section 7a is obeyed. If this is the case, then Proposition II.7.3 supplies a $\mathcal{K}$-compatible, and in particular, n-independent constant $c > 1$ with the following significance: If $z_* < c^{-1}$, then both PROPERTY 1 and PROPERTY 2 in Part 1 of Section 7a are obeyed with some perhaps different choice for S from $\mathcal{K}$. In any event, PROPERTY 1 impies PROPERTY 3 in Part 1 of Section 7a; and PROPERTY 4 in Part 1 of Section 7a follows if $z_* < c^{-1}$ as well.

These last observations lead to the following conclusion:

*There exists $\varepsilon_* > 0$ with the following significance: If $z_* < \varepsilon_*$, then no sufficiently large* n *version* $C_{Sn}$ *lies entirely in the radius $\varepsilon_*$ tubular neighborhood of any submanifold from $\mathcal{K}$.*

(8.2)

This assertion leads to the desired nonsense as it is proved false in the upcoming Step 6.

Step 4: The restriction of $N_{Cn}$ to $C_{Sn}$ is denoted by $N_{Sn}$ and the restriction to $C_{\mathfrak{p}n}$ is denoted by $N_{\mathfrak{p}n}$. Let $\eta$ denote a given, smooth section of the bundle $N_{Cn}$ over $C_n$. The restriction of the section $\eta$ to $C_{Sn}$ and to each $\mathfrak{p} \in \Lambda$ version of $C_{\mathfrak{p}n}$ defines the G + 1 tuple, $(\eta_S, \{\eta_\mathfrak{p}\}_{\mathfrak{p} \in \Lambda})$ with $\eta_S$ denoting a section of $N_{Sn}$ and with each $\mathfrak{p} \in \Lambda$ version of $C_{\mathfrak{p}n}$ denoting a section of $N_{\mathfrak{p}n}$. Taking this view of $C^\infty(C_n; N_{Cn})$ leads to the equivalent definition of the space $\mathbb{H}$ of $L^2_1$ space of sections of $N_{Cn}$ given in the next paragraph.

The space $\mathbb{H}$ is the completion of a subspace of $C^\infty(C_{Sn}, N_{Sn}) \oplus (\oplus_{\mathfrak{p} \in \Lambda} C^\infty(C_{\mathfrak{p}n}, N_{\mathfrak{p}n}))$. The subspace consists of elements with compact support and with boundary values as follows: Let $(\eta_S, \{\eta_\mathfrak{p}\}_{\mathfrak{p} \in \Lambda})$ denote an element in the subspace. Given $\mathfrak{p} \in \Lambda$, view a neighborhood of $C_{Sn}$ near a given critical point from $\mathfrak{p}$ as in Step 3. With this view understood, write $\eta_S$ on the cooresponding $z = z_*$ boundary as a pair of functions of x, these denoted by $(\varphi^{S'}, \varsigma^{S'})$. Meanwhile, write $\eta_\mathfrak{p}$ on the contiguous boundary of $C_{\mathfrak{p}n}$ as a pair of functions of x, these $(\varphi^{\mathfrak{p}'}, \varsigma^{\mathfrak{p}'})$. Then $\varsigma^{S'} = \varsigma^{\mathfrak{p}'}$ and $\varphi^{\mathfrak{p}'} = \varphi^{S'}$. The relevant completion of this subspace is defined by the respective $L^2_1$ norms on the spaces of compactly supported sections of $N_{Sn}$ and each $\mathfrak{p} \in \Lambda$ version of $C_{\mathfrak{p}n}$.

The range space $\mathbb{L}$ for $\hat{\mathcal{D}}_{Cn}$ can be viewed as the completion of the space of compactly supported sections of $C^{\infty}(C_{Sn}, N_{Sn} \otimes T^{0,1}S_n) \oplus (\oplus_{\mathfrak{p} \in \Lambda} C^{\infty}(C_{\mathfrak{p}n}, N_{\mathfrak{p}n} \otimes T^{0,1}C_{\mathfrak{p}n}))$ using the norm that is defined by the respective $L^2$ norms for each summand.

Use $D_{Sn}$ to denote the restriction of $\hat{\mathcal{D}}_{Cn}$ to $C_{Sn}$. Fix $\mathfrak{p} \in \Lambda$ and view $C_{\mathfrak{p}n}$ as the graph of $(\varphi^{\mathfrak{p}n}, \varsigma^{\mathfrak{p}n})$. The restriction of the operator $\hat{\mathcal{D}}_{Cn}$ to $C_{\mathfrak{p}n}$ is given by the $\mathfrak{h} = (\varphi^{\mathfrak{p}n}, \varsigma^{\mathfrak{p}n})$ version of (3.9) with the functions $a_1$, $a_2$ and $b$ defined by the index $n$ data set. The latter incarnation of $\hat{\mathcal{D}}_{Cn}$ is denoted by $D_{\mathfrak{p}n}$ in what follows.

<u>Step 5</u>: For each $\tau \in [0, 1]$, use $\mathbb{H}_\tau$ now to denote the Banach space that is obtained by the $L^2_1$ completion of the subspace of compactly supported sections in $C^{\infty}(C_{Sn}, N_{Sn} \otimes T^{0,1}S_n) \oplus (\oplus_{\mathfrak{p} \in \Lambda} C^{\infty}(C_{\mathfrak{p}n}, N_{\mathfrak{p}n} \otimes T^{0,1}C_{\mathfrak{p}n}))$ whose boundary values obey

$$\varsigma^{S'} = \tau \varsigma^{\mathfrak{p}'} \quad \text{and} \quad \varphi^{\mathfrak{p}'} = \tau \varphi^{S'}.$$
(8.3)

The operator $\mathcal{D}_{Cn} = (D_{Sn}, \{D_{\mathfrak{p}n}\}_{\mathfrak{p} \in \Lambda})$ acts as a bounded operator from each $\tau \in [0, 1]$ version of $\mathbb{H}_\tau$ to $\mathbb{L}$. The arguments for Lemma 7.15 can be used with almost no changes to prove that $\mathcal{D}_{Cn}$ defines a Fredholm operator from each $\tau \in [0, 1]$ version of $\mathbb{H}_\tau$ to $\mathbb{L}$ and that the index of each such Fredholm incarnation of $\mathcal{D}_{Cn}$ is equal to the Fredholm index of the $\tau = 1$ version, this being 1 since the $\tau = 1$ version is $\hat{\mathcal{D}}_{Cn}$.

Consider now the $\tau = 0$ version. The latter is a direct sum of $G + 1$ Fredholm operators. The first of these is $D_{Sn}$ acting on the $L^2_1$ completion of the subspace of compactly supported sections of $N_{Sn}$ that obey the following boundary condition: Write a given section on a given boundary component as a pair of functions, $(\varphi^{S'}, \varsigma^{S'})$. Then the boundary condition asserts only that $\varsigma^{S'} = 0$. Note in particular that this boundary condition makes no reference to any $\mathfrak{p} \in \Lambda$. The range space for this Fredholm operator is the $L^2$ completion of the space of compactly supported sections of $N_{Sn} \otimes T^{0,1}C_{Sn}$. Use index($D_{Sn}$) in what follows to denote the Fredholm index of this Fredholm incarnation of $D_{Sn}$

Meanwhile, each $\mathfrak{p} \in \Lambda$ labels an operator in the aforementioned direct sum. The latter is $D_{\mathfrak{p}n}$ acting on the $L^2_1$ completion of the space of compactly supported sections of $N_{\mathfrak{p}n}$ whose boundary values are as follows: Write a section on the a boundary component as $(\varphi^{\mathfrak{p}'}, \varsigma^{\mathfrak{p}'})$. Then $\varphi^{\mathfrak{p}'} = 0$. Note that this condition makes no reference to $S$ or to the other elements in $\Lambda$. The range space for this Fredholm incarnation of $D_{\mathfrak{p}n}$ is the $L^2$ completion of the space of compactly supported sections of $N_{\mathfrak{p}n} \otimes T^{0,1}C_{\mathfrak{p}n}$. This incarnation of $D_{\mathfrak{p}n}$ is described by Proposition 5.7; it has trivial cokernel and kernel dimension equal to $\Delta_{\mathfrak{p}}$.

What was said in the three previous paragraphs implies that

$$\text{index}(D_{S_n}) + \Sigma_{\mathfrak{p}\in\Lambda}\Delta_\mathfrak{p} = 1 \tag{8.4}$$

This understood, it follows that $\text{index}(D_{S_n}) \leq 1$ and that $\text{index}(D_{S_n}) \leq 0$ if any $\mathfrak{p} \in \Lambda$ version of $\Delta_\mathfrak{p}$ are positive.

Step 6: Invoke Proposition II.7.2 and use Corollary 7.2 in [L] to find a subsequence of $\{(D_n, C_n)\}_{n=1,2,\ldots}$ (hence renumbered consecutively from 1) and a sequence $\{\varepsilon_n\}_{n=1,2,\ldots} \in (0, 1)$ with the following three properties: The latter subsequence is decreasing and converges in [0, 1] to 0. Second, Proposition II.7.2 applies to each $C_n$ with $\varepsilon = \varepsilon_n$. Third, the various index n versions of the relevant broken, singular admissable sets are identical. Let $\Xi$ denote this set. The lemma that follows uses $\Xi$ to compute the number $\text{index}(D_{S_n})$.

The notation used by the lemma writes a given element in $\Xi$ as $((S, u), \vartheta_\Sigma)$ where $(S, u)$ denotes a Lipshitz subvariety and $\vartheta_\Sigma$ denotes a finite set of constant $(s, t)$ slices of $\mathbb{R} \times (1, 2) \times \Sigma$. A given slice can appear more than once in $\vartheta_\Sigma$. The lemma uses $n_\Sigma$ to denote the number of elements in $\vartheta_\Sigma$. As noted in the paragraph prior to Proposition 7.2, the subvariety pair $(S, u)$ has an associated Fredholm operator, this denoted by $D_S$.

**Lemma 8.3**: *There exists an n-independent constant $\kappa > 1$ such that if $z_* < \kappa^{-1}$ and if n is sufficiently large, then $\text{index}(D_{S_n}) = \Sigma_{((S,u),\vartheta_\Sigma)\in\Xi}(\text{index}(D_S)+2n_\Sigma)$.*

This lemma is proved momentarily. Accept it for now.

What with (8.4), this lemma implies that either $\Xi$ has just one component, and the latter is has $\vartheta_\Sigma = \emptyset$; or else there exists a non-$\mathbb{R}$ invariant Lipshitz submanifold, S, with $\text{index}(D_S) \leq 0$. As this is precluded by Proposition 7.2, it follows as a consequence that $\Xi$ has but a single element with $\vartheta_\Sigma = \emptyset$.

Granted this last conclusion, invoke the second bullet of Proposition 7.3 to see that all sufficiently large n versions of $C_n$ violate what is asserted in (8.2). This observation constitutes the desired nonsense.

Step 7: This step and Step 8 contain the

***Proof of Lemma 8.3***: Fix $Z \in ((S, u), \vartheta_\Sigma) \in \Xi$ and let $|Z|$ denote the union of $u(S)$ with the curves from $\vartheta_\Sigma$. Use the data given in Proposition II.7.2 to obtain a subsequence of $\{C_n\}_{n=1,2,\ldots}$, hence renumbered from 1, and a sequence $\{s_n\}_{n=1,2,\ldots}$ with the following property: Fix $n \in \{1, 2, \ldots\}$ and use $X_n$ to denote $[-4n, 4n] \times [1 + \frac{1}{n}, 2 - \frac{1}{n}] \times \Sigma$. View $X_n$ for the moment as a subset of $\mathbb{R} \times M_{\delta_n}$. Translate the surface $C_n$ by $s_n$ along the $\mathbb{R}$ factor of $\mathbb{R} \times [1,2] \times \Sigma$ so that each point in $C_n \cap X_n$ has distance at most $\frac{1}{n}$ from some point in $|Z| \cap X_n$, and vice versa. Furthermore, if $\mu$ is any 2-form on $X_n$ with $|\mu| \le 1$ and $|\nabla \mu| < n$, then the integral of $\mu$ over $C_n \cap X_n$ differs from $\int_{u(S) \cap X_n} \mu + \Sigma_{\Sigma' \in \vartheta_\Sigma} \int_{\Sigma'} \mu$ by less than $\frac{1}{n}$.

There exists by assumption, a constant $z_S > 0$ such that there are G components of the $t \in [1, 1+z_S]$ part of S and each is described by PROPERTY 5 in Part 1 of Section 1g. The $t \in [2-z_S, 2]$ part of S consists of an analogous set of G components. Meanwhile, there exists $s_S \ge 1$ such that the $|s| > s_S$ part of S is described by PROPERTY 6 in Part 1 of Section 1g. Assume that $z_* < 10^{-4} z_S$. Then the various components of the $s_n$-translate of $C_{Sn}$ where $t \in [1+z_*, 1+\frac{1}{2} z_S]$, where $t \in [2-\frac{1}{2} z_S, 2-z_*]$ and where $|s| \in [2s_S, 2n]$ obey the conclusions of PROPERTY 1 of Part 1 in Section 7a.

Granted what was just said, the arguments in Section 4 of [L] will write index($D_{Sn}$) as a sum of various contributions that can be readily identified with the terms in Lemma 8.3's sum. To do this, focus again on a given $((S, u), \vartheta_\Sigma) \in \Xi$. Truncate the $s_n$-translate of $C_{Sn}$ on the slices where $|s| = 3s_n$. The result, when n is large, has 2G constant s boundary arcs that run from the $t = 1+z_*$ boundary to the $t = 2-z_*$ boundary of $\mathbb{R} \times [1+z_*, 2-z_*] \times (T_- \cap T_+)$. Attach to each such arc a properly embedded, infinite strip that is a graph over either the $s \le -3s_n$ or $s \ge 3s_n$ part of $\mathbb{R} \times [1+z_*, 2-z_*]$ of a smooth map to $T_- \cap T_+$ that converges as $|s| \to \infty$ at an exponential rate in $|s|$ to the nearby $C_- \cap C_+$ point. Let $Z_n$ denote the resulting submanifold with boundary in $\mathbb{R} \times [1+z_*, 2-z_*] \times \Sigma$. This submanifold has a corresponding version of the operator in (1.25) which is Fredholm when viewed as a linear map between the $Z_n$ analogs of the Banach space domain and range spaces that were defined for $D_{Sn}$. The associated Fredholm index is denoted in what follows by index($D_{Zn}$).

Standard gluing theorems can be used to prove that index($D_{Sn}$) is the sum of the various $Z \in \Xi$ versions of index($D_{Zn}$) when n is large. See for example, Lemma 9.6 in [HT] for a statement in an analgous context but where the operator is defined on a manifold without boundary. The corresponding lemma for the case at hand is proved using arguments that differ only cosmetically. (These gluing theorems are geometric expressions of the excision property that is obeyed by the index of Fredholm elliptic operators on manifolds.)

The next step explains why index($D_{Zn}$) = index($D_S$) + $2n_\Sigma$ when n is large.

<u>Step 8</u>:  Because $t \in [1+z_*, 1+z_S] \cup [2-z_S, 2-z_*]$ part of $Z_n$ is a graph over the analogous part of S which is very close to S when $z_* \ll 1$, the arguments from Section 4 of [L] can be applied directly to $Z_n$ to prove the equality index($D_{Zn}$) = index($D_S$) + $2n_\Sigma$.

To elaborate, digress for the moment to review some basic facts. Suppose that X is a smooth 4-manifold with an almost complex structure, and suppose that $Z \subset X$ is a compact, pseudoholomorphic submanifold without boundary. Then Z has a version of the operator that is depicted on the right hand side of (1.25) which is Fredholm as a map from the $L^2_1$ space of sections of its normal bundle to the space of $L^2$ sections of the latter tensored with $T^{0,1}Z$. View Z's fundamental class as defining an element in $H_2(X; \mathbb{Z})$ and let $e_Z$ denote the image in $H^2(X; \mathbb{Z})$ of the Poincaré dual to this class. Let $c_{1X}$ denote the first Chern class of $T^{2,0}X$. The index of this operator is the value of $e_Z - c_1$ on the fundamental class of Z, this written as $(e_Z - c_1)[Z]$.

In the case at hand, $Z = \mathcal{Z}_n$ is a non-compact manifold with boundary so the preceding formula does not apply. However, the constructions in Section 4 of [L] write index($D_{\mathcal{Z}_n}$) as a sum of three terms: The first is a contribution from the boundary of $\mathcal{Z}_n$; the second is a contribution from the $|s| \gg 1$ part of $\mathcal{Z}_n$; and the third can be written as the evaluation on the fundamental class of $\mathcal{Z}_n$ of a version of $e_Z - c_1$. To elaborate, the fundamental class of $\mathcal{Z}_n$ is viewed as a class in a certain relative second homology and the version of $e_Z - c_1$ is a certain class in the second cohomology with compact support. Meanwhile, index($D_S$) has the analogous decomposition into a sum of three terms.

The fact that any large n version of $\mathcal{Z}_n$ is a graph over S near its boundary and also where $|s|$ is very large implies that the first two terms in the three term sum for index($D_{\mathcal{Z}_n}$) are the same as the first two terms in the sum for index($D_{S_n}$).

To relate the respective third terms in the three term sums for index($D_{\mathcal{Z}_n}$) and index($D_S$), identify $\mathbb{R} \times [1+z_*, 2-z_*] \times \Sigma$ with $\mathbb{R} \times [1,2] \times \Sigma$ by the diffeomorphism from $[1+z_*, 2-z_*]$ to $[1,2]$ that maps $t$ to $\frac{t-3z_*}{1-2z_*}$ so as to view $\mathcal{Z}_n$ as a properly embedded submanifold with boundary in $\mathbb{R} \times [1,2] \times \Sigma$. With this identification understood, let $[\mathcal{Z}_n]$ denote the relative fundamental class of $\mathcal{Z}_n$, this a relative class on $\mathbb{R} \times [1,2] \times \Sigma$. Let $[S]$ denote the corresponding relative second homology class given by the push-forward via $u$ of the relative fundamental class of the surface S; and let $[\Sigma]$ denote the fundamental class of $\Sigma$. Then $[\mathcal{Z}_n] = [S] + n_\Sigma[\Sigma]$. The analog of $c_1$ in this context is a certain compactly supported class that represents the first Chern class of the complex line bundle $T^{2,0}(\mathbb{R} \times [1,2] \times \Sigma))$. The latter class evaluates as $2 - 2G$ on $[\Sigma]$. The $\mathcal{Z}_n$ analog of $e_Z$ is denoted by $e_{Z_n}$ and the S analog is denoted by $e_S$. The fact that any large n version of $\mathcal{Z}_n$ is a graph over S near its boundary and also where $|s|$ is very large implies that $e_{Z_n} - e_S$ is in $H^2(\mathbb{R} \times [1,2] \times \Sigma)$ and that it evaluates on $[\mathcal{Z}_n]$ to give

$$e_{Z_n}[\mathcal{Z}_n] = e_S[S] + n_\Sigma e_\Sigma[\Sigma] + 2n_\Sigma[S] \cdot [\Sigma]$$

(8.5)

where the notation is as follows: What is denoted by $e_\Sigma$ is the image in the second cohomology of the Poincaré dual of $[\Sigma]$. Meanwhile, $[S]\cdot[\Sigma]$ denotes the intersection number between S and $\Sigma$. As $e_\Sigma[\Sigma] = 0$ and $[S]\cdot[\Sigma] = G$, this formula plus what was said about the first Chern class of the bundle $T^{2,0}(\mathbb{R}\times[1,2]\times\Sigma)$ implies that the third term in the respective third three term expressions for index($D_{Z_n}$) and index($D_{S_n}$) are such that

$$(e_{Z_n} - c_1)[Z_n] = (e_S - c_1)[S] + n_\Sigma(2 - 2G) + 2n_\Sigma G \ .$$

(8.6)

This identity with the aforementioned identity between the first two terms in the three term expressions for index($D_{Z_n}$) and index($D_{S_n}$) lead directly to the desired equality index($D_{Z_n}$) = index($D_{S_n}$) + $2n_\Sigma$.

**b) Quillen's construction and orientations**

As explained in [HS], [Hu1] and Section 9 of [HT], the differential for embedded contact homology is defined using certain dimension one moduli spaces of J-holomorphic subvarieties in $\mathbb{R}\times Y$. The definition involves a ±1 weight that is define by comparing two orientations that can be defined for these moduli spaces. The first is defined by the $\mathbb{R}$-action that is induced by the constant translations along the $\mathbb{R}$ factor of $\mathbb{R}\times Y$. The second is defined using notions that were introduced by Quillen [Q] about determinant line bundles for parametrized families of Fredholm operators. Section 9 uses what is said here and in Sections 8c-e to describe the weight that is used for the embedded contact homology differential differs by a purely S-dependent sign from the weight used in Proposition 8.1. This subsection describes the relevant version of Quillen's construction of orientations. The story is told in five parts.

*Part 1*: To say more about the Fredholm operator that is used to define the differential for embedded contact homology, fix $(\tau = 1, \mathcal{C}) \in \mathcal{M}^*$ and let C denote the corresponding ech-HF submanifold. The operator in question is the operator $\mathcal{D}_\mathcal{C}$ from Lemma 7.15 acting here on a slightly larger domain. The range space is the same as for the original. The domain is denoted here by $\mathbb{H}_S$. The space $\mathbb{H}_S$ is defined just as $\mathbb{H}_{\tau=1}$ in Part 1 of Section 7e but for the following: Elements in the S-labeled summand of $\mathbb{H}_{\tau=1}$ are required to be $L^2$ orthogonal on the $t \in [1+z_*, 2-z_*]$ part of S to the restriction of the kernel of $D_S$. This last condition is not imposed on the elements in the S-labeled summand of $\mathbb{H}_S$. In any case the $\tau = 1$ version of (7.32) is imposed. If $z_* \leq c^{-1}$ with $c \geq 1$ purely S-dependent (of $\mathcal{K}$-compatible), then there is a canonical isomorphism between $\mathbb{H}_S$ and kernel($D_S$) $\oplus \mathbb{H}_{\tau=1}$.

This version of $\mathcal{D}_C$ with domain $\mathbb{H}_S$ is denoted in what follows by $\hat{\mathcal{D}}_C$. The notation here is meant to indicate that the latter operator can be defined intrinsically as an operator on C. As noted in the proof of Proposition 8.2, it has the form as what is depicted on the right hand side of (1.25). This intrinsic definition identifes the domain Hilbert space with the space $L^2_1$ sections of the C's normal bundle and the range Hilbert space to the space of $L^2$ sections of the tensor product of this normal bundle with $T^{0,1}C$. This intrinsic definition does not reference (7.32).

*Part 2*: This part of the subsection summarizes Quillen's construction. To start, suppose for the moment that $\mathcal{H}_1$ and $\mathcal{H}_2$ are Hilbert spaces and $\mathcal{D}: \mathcal{H}_1 \to \mathcal{H}_2$ is a Fredholm operator with positive index and trivial cokernel. Define $\mathbb{D}\text{et}(\mathcal{D})$ to be the real line given by the top exterior power of the kernel of $\mathcal{D}$.

Suppose next that $\mathcal{Y}$ denotes a smooth, finite dimensional manifold and that $\mathcal{H}_1$ and $\mathcal{H}_2$ are Hilbert space bundles over $\mathcal{Y}$. (Most of what is said here generalizes readily to the case when $\mathcal{Y}$ is a Hilbert manifold.) Let $\mathcal{D}$ now denote a continous section of $\text{Hom}(\mathcal{H}_1, \mathcal{H}_2)$ whose restriction to each fiber is Fredholm. What follows describes a real line bundle that is defined over $\mathcal{Y}$ in a canonical fashion by $\mathcal{D}$. To this end, remark first that the manifold $\mathcal{Y}$ has a locally finite cover with the following property: Let $U \subset \mathcal{Y}$ denote a set from this cover. Then there exist a non-negative integer, n, and a bundle homomorphism L: $U \times \mathbb{R}^{2n} \to \mathcal{H}_2|_U$ such that

$$\mathcal{D} + L: \mathcal{H}_1|_U \oplus (U \times \mathbb{R}^{2n}) \to \mathcal{H}_2|_U$$

(8.7)

restricts to each fiber as a linear map with positive index and trivial cokernel. This understood, define the real line bundle $\mathbb{D}\text{et}|_U \to U$ to be the bundle whose fiber at any given $y \in Y$ is the top exterior power of the kernel of $(L + \mathcal{D})|_y$. As explained by Quillen, different choices for the integer n and, given n, for the homomorphism L subject to the condition that $\mathcal{D} + L$ have positive index and trivial cokernel give isomorphic versions of $\mathbb{D}\text{et}|_U$. It follows as a consequence that these line bundles over the open sets of the given cover patch together over the pairwise intersections to define a real line bundle over $\mathcal{Y}$. It also follows (by taking subdivisions) that two covers of $\mathcal{Y}$ with the requisite properties supply isomorphic bundles. This being the case, the construction just described defines from $\mathcal{D}$ a canonical real line bundle over $\mathcal{Y}$. This is the bundle $\mathbb{D}\text{et}$.

The particular version of $\mathbb{D}\text{et}$ that is used in the definition of the embedded contact homology differential is defined over the moduli spaces of ech-HF submanifolds. Let C denote the ech-HF submanifold that is associated to a given pair $(1, \mathcal{C}) \in \mathcal{M}^*$. The fiber of the relevant version of $\mathbb{D}\text{et}$ over C is the determinant line of the operator $\hat{\mathcal{D}}_C$.

*Part 3*: The versions of $\mathbb{D}\text{et}$ that arise in what follows are orientable. This part of the subsection sets up the conventions that are used in the subsequent parts that concern choices of orientations. If V and V´ denote an ordered pair of oriented vector spaces, then their direct sum has a canonical orientation that is obtained as follows: Let n and n´ denote the respective dimensions of V and V´. Choose respective basis $\{v_1, \ldots, v_n\}$ for V and $\{v_1´, \ldots, v_{n´}´\}$ for V´ such that $v_1 \wedge \cdots \wedge v_n$ defines the orientation for $\det(V) = \wedge^n V$ and $v_1´ \wedge \cdots \wedge v_{n´}´$ defines the orientation for $\det(V´) = \wedge^{n´} V´$. The orientation for $\det(V \oplus V´) = \wedge^{n+n´}(V \oplus V´)$ is defined by

$$v_1 \wedge \cdots \wedge v_n \wedge v_1´ \wedge \cdots \wedge v_{n´}´. \tag{8.8}$$

The oriented tensor product of the lines $\det(V)$ and $\det(V´)$ is defined to be the line $\det(V \oplus V´)$ with the orientation given by (8.8). This oriented line is denoted by $\det(V) \otimes \det(V´)$. The oriented lines $\det(V´) \otimes \det(V)$ and $\det(V) \otimes \det(V´)$ are isomorphic as oriented lines if and only if $nn´$ is an even number.

It follows as a consequence of what was just said that ordering issues are minimized in any given situation when one or both of V and V´ have even dimensions; this is why (8.7) uses only even dimensional Euclidean spaces. In particular, the restriction to even dimensions in (8.7) makes it easier to compare orientations for $\mathbb{D}\text{et}$.

What follows is meant to provide an abstract but relevant illustration. The vector space $\mathbb{R}^{2n}$ in here and in subsequent parts of this subsection always denotes the eponymous vector space with a standard orientation, chosen once and for all time. Suppose that $\mathcal{D}$ is a Fredholm operator with trivial cokernel and positive index. Choose an orientation for kernel($\mathcal{D}$) so as to orient the line $\mathbb{D}\text{et}(\mathcal{D})$. Let $\mathcal{H}_2$ denote the range space for $\mathcal{D}$, let n denote any given positive integer and let L: $\mathbb{R}^{2n} \to \mathcal{H}_2$ denote any given map. The kernel of $\mathcal{D}+L$ is canonically isomorphic to kernel($\mathcal{D}$) $\oplus \mathbb{R}^{2n}$, and so the oriented line det(kernel($\mathcal{D}$)) is canonically isomorphic as an oriented line to det(kernel($\mathcal{D} \oplus L$)). As a consequence, the orientation on $\mathbb{D}\text{et}(\mathcal{D})$ defined via its identfication with det(kernel($\mathcal{D}$)) is the same as that defined by its identification with det(kernel($\mathcal{D} \oplus L$)).

For example, suppose now that $\mathcal{D}$ has trivial cokernel and zero index. Fix non-negative integers n and n´ whose sum is at least 1, and fix linear maps L: $\mathbb{R}^{2n} \to \mathcal{H}_2$ and L´: $\mathbb{R}^{2n´} \to \mathcal{H}_2$. The kernel of $(\mathcal{D} + L) + L´: (\mathcal{H}_1 \oplus \mathbb{R}^{2n}) \oplus \mathbb{R}^{2n´} \to \mathcal{H}_2$ is canonically isomorphic to $\mathbb{R}^{2n} \oplus \mathbb{R}^{2n´}$ and the kernel of $(\mathcal{D}+L´)+L: (\mathcal{H}_1 \oplus \mathbb{R}^{2n´}) \oplus \mathbb{R}^{2n} \to \mathcal{H}_2$ is canonically isomorphic to $\mathbb{R}^{2n´} \oplus \mathbb{R}^{2n}$. The orientation on $\mathbb{D}\text{et}(\mathcal{D})$ that comes by identifying the latter with the oriented, top exterior power of $\mathbb{R}^{2n} \oplus \mathbb{R}^{2n´}$ is the same as that defined by identifying $\mathbb{D}\text{et}(\mathcal{D})$ with the top exterior power of $\mathbb{R}^{2n´} \oplus \mathbb{R}^{2n}$.

To end this illustration, suppose that $\mathcal{D}'$ is a second Fredholm operator with trivial cokernel and positive index. Orient kernel($\mathcal{D}'$) so as to orient the line $\mathbb{D}\mathrm{et}(\mathrm{kernel}(\mathcal{D}'))$ and so orient $\mathbb{D}\mathrm{et}(\mathcal{D}')$. Because the integer 2n in (8.7) is even, the orientation of the oriented line $\mathbb{D}\mathrm{et}(\mathcal{D}) \otimes \mathbb{D}\mathrm{et}(\mathcal{D}')$ is insensitive to the choice for n and L in (8.7) when defining $\mathbb{D}\mathrm{et}(\mathcal{D})$ subject to the constraint that $\mathcal{D}+L$ has positive index and trivial cokernel.

*Part 4*: Suppose $(\tau = 1, \mathcal{C}) \in \mathcal{M}^*$ and let C denote the corresponding ech-HF submanifold. As noted in Part 1, the operator that is used to define the embedded contact homology differential is the intrinsically defined operator $\hat{\mathcal{D}}_C$. Of particular concern with regards to the differential is a certain orientation for the real line $\mathbb{D}\mathrm{et}(\hat{\mathcal{D}}_C)$. The operator $\hat{\mathcal{D}}_C$ can be viewed as the operator $\mathcal{D}_C$ with domain space $\mathbb{H}_S = \mathrm{kernel}(D_S) \oplus \mathbb{H}_{\tau=1}$. It follows as a consequence of what is said in Parts 2 and 3 that

$$\mathbb{D}\mathrm{et}(\hat{\mathcal{D}}_C) = \mathbb{D}\mathrm{et}(D_S) \otimes \mathbb{D}\mathrm{et}(\mathcal{D}_C) .$$

(8.9)

This understood, orientations for the line $\mathbb{D}\mathrm{et}(D_S)$ and for the line $\mathbb{D}\mathrm{et}(\mathcal{D}_C)$ orient the line of interest, $\mathbb{D}\mathrm{et}(\hat{\mathcal{D}}_C)$. The rest of this Part 4 explains how to relate the line $\mathbb{D}\mathrm{et}(\mathcal{D}_C)$, and thus the line $\mathbb{D}\mathrm{et}(\hat{\mathcal{D}}_C)$ to Proposition 8.1.

Suppose that $(\tau, \mathcal{C})$ is any given pair in $\mathcal{M}^*$. Use $\mathcal{D}_C$ in this case to denote the operator from Lemma 7.15 acting on the Banach space $\mathbb{H}_\tau$. The range Banach space is $\mathcal{C}$'s version of $\mathbb{L}$. These respective Banach spaces are the fibers over $\mathcal{M}^*$ of a pair of smooth, Banach space bundles, $\mathcal{H}_1$ and $\mathcal{H}_2$. Meanwhile, the various $(\tau, \mathcal{C})$ versions of $\mathcal{D}_C$ define the fibers of a section, $\mathcal{D}$, of $\mathrm{Hom}(\mathcal{H}_1, \mathcal{H}_2)$ which is Fredholm on each fiber. Note in this regard that the smooth variation can be proved using what is said in Step 3 from Part 4 of Section 7d (with arguments that mimick those from Step 3 of Part 3 from Section 4c). The section $\mathcal{D}$ has its associated determinant line bundle, $\mathbb{D}\mathrm{et}(\mathcal{D}) \to \mathcal{M}^*$. The $(\tau = 1, \mathcal{C})$ version of $\mathbb{D}\mathrm{et}(\mathcal{D}_C)$ that appears in (8.9) is the fiber over $(1, \mathcal{C})$ of $\mathbb{D}\mathrm{et}(\mathcal{D})$.

Hold on to $\mathbb{D}\mathrm{et}(\mathcal{D})$ for a moment and reintroduce the map $\mathfrak{F}$ as in (7.55). The differential of $\mathfrak{F}$ at the point $(\tau, 0)$ defines a Fredholm map from $\mathbb{R} \times \mathbb{H}_\tau \to \mathbb{L}$ of the form $\mathcal{D}_C + L_C$ where $L_C: \mathbb{R} \to \mathbb{L}$ is a linear map. Given the choice for J, the operator $\mathcal{D}_C + L_C$ has trivial cokernel and kernel dimension equal to $\sum_{p \in \Lambda} \Delta_p + 1$. This kernel is canonically isomorphic to $T\mathcal{M}^*$. To summarize,

$$\mathrm{kernel}(\mathcal{D}_C + L_C) = \mathrm{kernel}(\mathfrak{F}_*) = T\mathcal{M}^*|_{(\tau,\mathcal{C})} .$$

(8.10)

The differential of $\pi_1$ at $(\tau, \mathcal{C})$ appears here as the restriction to kernel$(\mathcal{D}_\mathcal{C}+L_\mathcal{C})$ of the projection map from $\mathbb{R} \times \mathbb{H}_\tau$ to the $\mathbb{R}$ factor. Meanwhile, $\wedge^{top}$ kernel$(\mathcal{D}_\mathcal{C}+L_\mathcal{C})$ is the line $\mathrm{Det}(\mathcal{D}_\mathcal{C}+L_\mathcal{C})$. Thus, an orientation for $T\mathcal{M}^*|_{(\tau,\mathcal{C})}$ is induced from an orientation for $\mathrm{Det}(\mathcal{D}_\mathcal{C}+L_\mathcal{C})$ and vice-versa.

At a point $(\tau, \mathcal{C})$ where $\tau$ is a regular value of $\pi_1$, the determinant line $\mathrm{Det}(\mathcal{D}_\mathcal{C}+L_\mathcal{C})$ is isomorphic to $\mathbb{R} \otimes \mathrm{Det}(\mathcal{D}_\mathcal{C})$. As the differential of $\pi_1$ at $(\tau, \mathcal{C})$ is an isomorphism, the differential of $\pi_1$ at $(\tau, \mathcal{C})$ identifies the $\mathbb{R}$ factor in this tensor product as the oriented line $(-\infty, \infty)$. Reintroduce now the line $\mathrm{Det}(\mathcal{D}) \to \mathcal{M}^*$. It follows as a consequence of what was just said that an orientation for the line $\mathrm{Det}(\mathcal{D})$ defines an orientation for the $\pi_1 = \tau$ level set in $\mathcal{M}^*$; and in particular for the $\tau = 1$ boundary of $\mathcal{M}^*$. Note in this regard that this level set is a union of points when each $\mathfrak{p} \in \Lambda$ version of $\Delta_\mathfrak{p}$ are zero. Meanwhile an orienation of a point is, by definition, a choice of +1 or -1. The +1 appears if the orientation of $T\mathcal{M}^*$ given by $\mathrm{Det}(\mathcal{D})$ at the point in question is that given by the differential of $\pi_1$. By the same token, the orientation given by $\mathrm{Det}(\mathcal{D})$ to the $\tau = 0$ point in $\mathcal{M}^*$ is +1 if the orientation of $T\mathcal{M}^*$ at this point agrees with that given by the differential of $\pi_1$.

*Part 5*: Suppose that $\sum_{\mathfrak{p} \in \Lambda} \Delta_\mathfrak{p} = 0$. An orientation for $\mathrm{Det}(\mathcal{D})$ induces an orientation on $T\mathcal{M}^*|_{\tau=1}$ which agrees or not with that used in Proposition 8.1; but agreement or not is the same at all points in $\mathcal{M}^*|_{\tau=1}$. It follows from (8.10) that agreement occurs if and only if there is agreement for the sole point of $\mathcal{M}^*|_{\tau=0}$. Meanwhile, Section 8c describes a completely canonical orientation for $\mathrm{Det}(\mathcal{D})$ when $\sum_{\mathfrak{p} \in \Lambda} \Delta_\mathfrak{p} = 0$. As can be seen readily from the definition in Section 8c, the resulting orientation for $\mathcal{M}^*|_{\tau=0}$ agrees with that used in Proposition 8.1. This being the case, the $\mathrm{Det}(\mathcal{D})$ orientation for $\mathcal{M}^*|_{\tau=0}$ also agrees with Proposition 8.1's orientation when $\sum_{\mathfrak{p} \in \Lambda} \Delta_\mathfrak{p} = 0$.

In the case $\sum_{\mathfrak{p} \in \Lambda} \Delta_\mathfrak{p} > 0$, choose once and for all an ordering of $\Lambda_*$ up to even permutations. Such a choice orients $(\times_{\mathfrak{p} \in \Lambda_*} \mathbb{R})$. Suppose that $y \in (\times_{\mathfrak{p} \in \Lambda_*} \mathbb{R})$ is a regular value of the map $p$. Then the tangent space to $\mathcal{M}^*_y$ is isomorphic at any given point $(\tau, \mathcal{C})$ to the kernel of the $p$'s differential. Meanwhile, the normal bundle of $\mathcal{M}^*_y$ in $\mathcal{M}^*$ is mapped isomorphically by $p$'s differential to $(\times_{\mathfrak{p} \in \Lambda_*} \mathbb{R})$. This understood, Proposition 8.1's orientation for $\mathcal{M}^*_y$ and $p$'s orientation of the normal bundle to $\mathcal{M}^*_y$ orients the tangent space to any given smooth level set of $\pi_1$ and in particular the tangent space to $\mathcal{M}^*|_{\tau=1}$. Granted this last observation, (8.10) implies that the orientation of $\mathcal{M}^*|_{\tau=1}$ given by the differential of $\pi_1 \times p$ either agrees with or disagrees with the orientation induced by an orientation of $\mathrm{Det}(\mathcal{D})$; but agreement or not is the same at all points. Moreover, agreement occurs if and only if the corresponding two orientations of $\mathcal{M}^*|_{\tau=0}$ agree.

Meanwhile, Section 8d describes a completely canonical orientation for the line $\mathbb{D}\mathrm{et}(\mathcal{D}) \to \mathcal{M}^*$ which consistently orients the level sets of $\pi_1$ in $\mathcal{M}^*$. A look ahead at what is said in Section 8d shows that this $\mathbb{D}\mathrm{et}(\mathcal{D})$ orientation for $\mathcal{M}^*|_{\tau=0}$ agrees with the one defined by Proposition 8.1 and so the $\mathbb{D}\mathrm{et}(\mathcal{D})$ orientation for $\mathcal{M}^*|_{\tau=1}$ also agrees with Proposition 8.1's orientation of $\mathcal{M}^*|_{\tau=1}$

**c) The canonical orientation when $\{\Delta_\mathfrak{p} = 0\}_{\mathfrak{p} \in \Lambda}$**

This subsection describes the promised canonical orientation for $\mathbb{D}\mathrm{et}(\mathcal{D})$ in the case when all $\mathfrak{p} \in \Lambda$ versions of $\Delta_\mathfrak{p}$ are zero. The description has four parts.

*Part 1*: The Banach space $\mathbb{H}_\tau$ in this case can be viewed as a completion of a dense domain whose whose typical element is written as $(\eta_S{}', \{(\varphi^{\mathfrak{p}'}, \varsigma^{\mathfrak{p}'})\}_{\mathfrak{p} \in \Lambda})$ where $\eta'_S$ is a compactly supported section of $N_S$ over the $t \in [1+z_*, 2-z_*]$ part of S that is orthogonal to the restriction of the kernel of $D_S$. Meanwhile, each $\mathfrak{p} \in \Lambda$ version of $(\varphi^{\mathfrak{p}'}, \varsigma^{\mathfrak{p}'})$ is a compactly supported map from $\mathbb{R} \times I_*$ to $\mathbb{R}^2$. The parameter $\tau$ enters through the boundary constraint in (7.32). The Banach space norm is that induced by the $L^2_1$ norm on sections of the $t \in [1+z_*, 2-z_*]$ part of S and the $L^2_1$ norm on maps from $\mathbb{R} \times I_*$. Meanwhile, the range space $\mathbb{L}$ is the corresponding $L^2$ completion of a dense domain with typical element $(\eta_S{}^\#, \{(\varphi^{\mathfrak{p}\#}, \varsigma^{\mathfrak{p}\#})\}_{\mathfrak{p} \in \Lambda})$ where $\eta_S{}^\#$ is a section of an appropriate 2-plane bundle of the $t \in [1+z_*, 2-z_*]$ part of S, and where each $\mathfrak{p} \in \Lambda$ version of $(\varphi^{\mathfrak{p}\#}, \varsigma^{\mathfrak{p}\#})$ is a map from $\mathbb{R} \times I_*$ to $\mathbb{R}^2$. The Banach space norm is the norm induced by the $L^2$ norm on the $t \in [1+z_*, 2-z_*]$ part of S and on $\mathbb{R} \times I_*$.

As can be seen from preceding descriptions, the Banach space $\mathbb{H}_\tau$ depends only on $\tau$ (with S fixed) and the Banach space $\mathbb{L}$ depends only on S. What is done in Step 3 of Part 3 from Section 4c can be mimicked to see that the assignment to $\tau \in [0, 1]$ of the Banach space $\mathbb{H}_\tau$ defines a smooth Banach space bundle over the interval $[0, 1]$. The latter's $\pi_1$-pull-back over $\mathcal{M}^*$ is the Banach space bundle $\mathcal{H}_1$ that was defined in Part 4 of the previous subsection. Meanwhile, the bundle $\mathcal{H}_2$ from this same Part 4 is the product bundle $\mathcal{M}^* \times \mathbb{L}$.

*Part 2*: The Hilbert space bundles $\mathcal{H}_1$ and $\mathcal{H}_2$ extend over the product $\mathcal{M}^* \times [0, 1]$ as follows: The bundle $\mathcal{H}_2$ extends as the product bundle. The bundle $\mathcal{H}_1$ extends as the pull-back of via the projection to the square $[0, 1] \times [0, 1]$ of the bundle whose fiber at any given point $(\tau, r)$ is the Hilbert space $\mathbb{H}_{r\tau}$. The homomorphism $\mathcal{D}$ from $\mathcal{H}_1$ to $\mathcal{H}_2$ over $\mathcal{M}^*$ given by $(\tau, \mathcal{C}) \to \mathcal{D}_\mathcal{C}$ extends over $\mathcal{M}^* \times [0, 1]$ so that the restriction to each fiber is Fredholm. The section at a given $((\tau, \mathcal{C}), r)$ is the operator $\mathcal{D}_\mathcal{C}$ with the parameter $\tau$ in (7.32) replaced by $r\tau$. This extended bundle homomorphism is also denoted by $\mathcal{D}$. An orientation over $\mathcal{M}^* \times \{0\}$ for the associated real line bundle $\mathbb{D}\text{et}(\mathcal{D})$ defines an orientation for $\mathbb{D}\text{et}(\mathcal{D})$ over the whole of $\mathcal{M}^* \times [0, 1]$ and thus over $\mathcal{M}^* \times \{1\}$. Of course, the converse of this last assertion is also true.

To define an orientation for $\mathbb{D}\text{et}(\mathcal{D})$ over $\mathcal{M}^* \times \{0\}$, fix $(\tau, \mathcal{C}) \in \mathcal{M}^*$. The version of $\mathcal{D}_{\mathcal{C}}$ on $\mathbb{H}_0$ is a direct sum of operators. In particular, there is one for each $\mathfrak{p} \in \Lambda$. In each case, it is a version of (3.5) with coefficients that obey (3.6) and are given by (3.9). The elements in the dense domain of this $\mathfrak{p}$-summand operator are constrained to obey the boundary conditions in (3.7). There is the remaining S-labeled summand. The relevant operator here is that depicted in (7.33) with domain given by the orthogonal complement of the restriction kernel of $D_S$ to the part of S where $t \in [1+z_*, 2-z_*]$. Elements in the domain also obey the $\tau = 0$ version of the top line in (7.32).

As explained in Step 1 from Section 5d, each $\mathfrak{p} \in \Lambda$ version of (3.5) in the aforementioned direct sum is homotopic through a family of Fredholm operators to a canonical operator, this given by (3.22) with the homotopy given by (3.21). The operators in this family all have trivial kernel and cokernel. This understood, then a once and forever choice for the orientation of the determinant line of the operator in (3.22) with boundary conditions given by (3.7) orients the determinant line for each of $\mathfrak{p} \in \Lambda$ version of (3.5).

*Part 3*: This part elaborates on what was said at the end of Part 2. To start, introduce $\mathcal{O}$ to denote the set whose elements are 6-tuples of functions on $\mathbb{R} \times I_*$ that are of the form $(\mathfrak{a}_1, \mathfrak{a}_2, \mathfrak{b}_1, \mathfrak{b}_2)$ which are suitable for use in (3.5). In particular, they must obey the constraints in (3.6). No generality is lost for what follows by restricting to the case where the integrals in the third bullet are negative. The set $\mathcal{O}$ is given the topology that is induced by its inclusions into two topological function spaces of maps from $\mathbb{R} \times I_*$ to $\mathbb{R}^6$. The first topology is the $C^\infty$ Frêchet space topology with it understood that convergence means convergence in the various $C^k$ topologies on compact subsets. The second topology is the strong $C^1$ topology. The space $\mathcal{O}$ with this topology is contractible since the constraints in (3.6) form a convex set.

Each point in $\mathcal{O}$ defines an operator to which Proposition 3.1 applies. In particular, each such operator has trivial kernel and cokernel. Let $\mathbb{D}\text{et}$ denote the corresponding determinant line bundle. It follows from what is said in Part 2 of the previous subsection that this version of $\mathbb{D}\text{et}$ has a canonoical identification with the oriented line $\wedge^2 \mathbb{R}^2$. Use this orientation. Such a choice gives the determinant line of the operators parametrized by the elements in $\mathcal{O}$ a canonical orientation. Note that the space $\mathcal{O}$ and the line $\mathbb{D}\text{et}$ see nothing of the spaces M, Y and their geometry. They do depend on the parameters $z_*$ and R only to the extent that these define $I_*$. Even so, the allowed choices form a contractible set, so the orientation of $\mathbb{D}$et just described is truly canonical and universal.

   *Part 4*: As noted, the other operator that enters the direct sum giving $\mathcal{D}_\mathcal{C}$ on $\mathbb{H}_0$ is depicted in (7.33) but acting on the orthogonal complement of restriction of the kernel of $D_S$ to the $t \in [1+z_*, 2-z_*]$ part of S. The boundary conditions are given by the $\tau = 0$ version of the top line in (7.32). This operator with the domain as indicated is denoted here by $D_{\eta 0}$. Meanwhile, $D_S^\perp$ is used here to denote the restriction of the operator $D_S$, (given in (II.6.11)) to the orthogonal complement of its kernel. The boundary conditions are given by (6.12). The Fredholm version needed here is the one described in Part 3 of Section 6e.

   The operator $D_S^\perp$ has trivial kernel and cokernel. As noted previously, this is the case for $D_{\eta 0}$ when $z_* < c^{-1}$ where $c$ is purely S-dependent (or $\mathcal{K}$-compatible). It is a straightforward task to prove that $D_{\eta 0}$ is homotopic via an essentially canonical family of Fredholm operators to $D_S^\perp$ with each operator in the family having trivial kernel and cokernel. The details are omitted but for the following description of the family: The family is the concatenation of two 1-parameter families. The first modifies the operator on the fixed domain via the family parametrized by $[0,1]$ with the $\mu \in [0,1]$ version of the operator defind by (7.33) with $\mu\eta$ replacing $\eta$. The version with $\mu = 0$ is the restriction of the operator $D_S$ to the $t \in [1+z_*, 2-z_*]$ part of S. The second part of the homotopy keeps the operator fixed as $D_S$ but changes the domain Hilbert space by introducing a parameter $\mu \in [0,1]$ and restricting $D_S$ to the $[1+\mu z_*, 2-\mu z_*]$ portion of S. The boundary conditions that define any given $\mu < 1$ Hilbert space are the $t = 1+\mu z_*$ and $t = 2-\mu z_*$ analogs of those for the $\mu = 1$ member. Likewise, the orthogonality condition with regards to the kernel of $D_S$ is changed only to the extent that the orthogonality is defined by integration over the $t \in [1+\mu z_*, 2-\mu z_*]$ part of S.

   Granted all of this, it then follows that the line bundle $\mathbb{D}\text{et}(\mathcal{D})$ along $\{0\} \times \mathcal{M}^*$ has a completely canonical orientation given a choice of orientation for the determinant line of the operator $D_S^\perp$. It follows from what is said in Part 2 of the previous subsection that such an orientation is canonically defined by identifying $\mathbb{D}\text{et}(D_S^\perp)$ with $\wedge^2 \mathbb{R}^2$. Moreover, it follows as a consequence of what is said in Parts 2 and 3 of the previous section that this orientation for $\mathbb{D}\text{et}(D_S^\perp)$ is induced from an orientation for $\mathbb{D}\text{et}(D_S)$ that comes by writing $D_S$ as $D_S = D_S^\perp + L_S$ where $L_S$ is the map from kernel($D_S$) to the range of $D_S$ that sends all elements to zero.

   *Part 5*: Given what is said in Part 3 of Section 5b, the respective orientations for $\mathbb{D}\text{et}(D_S^\perp)$ and for each $\mathfrak{p} \in \Lambda$ version of the determinant line of the operator in (3.5) defines a completely canonical orientation for $\mathbb{D}et(\mathcal{D})$ along $\mathcal{M}^* \times \{0\}$. As noted above, latter defines a completely canonical orientation for $\mathbb{D}et(\mathcal{D})$ on $\mathcal{M}^* \times [0, 1]$, and thus to $\mathbb{D}et(\mathcal{D})$ along $\mathcal{M}^* \times \{1\}$. The latter is the canonical orientation promised in Part 5 of the previous subsection.

### d) Canonical orientations when $\sum_{\mathfrak{p} \in \Lambda} \Delta_\mathfrak{p} > 0$

This subsection describes the promised canonical orientation for $\mathbb{D}et(\mathcal{D}) \to \mathcal{M}^*$ in the general case. The description that follows has four parts.

*Part 1*: An orientation for $\mathbb{D}et(\mathcal{D}) \to \mathcal{M}^*$ is defined by first extending the family of operators to a family with parameter space $\mathcal{M}^* \times [0, 1]$. This extension is defined by making the domain for any given $((\tau, \mathcal{C}), r)$ version of the operator depend on the parameter for the extra $[0, 1]$ factor. The range space is kept constant, and the operator itself stays as $\mathcal{D}_\mathcal{C}$. To say more about the r-dependence of the domain, fix a pair $((\tau, \mathcal{C}), r) \in \mathcal{M}^* \times [0, 1]$. The domain for the corresponding Fredholm operator is identical to that when $r = 1$ but for one item: The boundary conditions in (7.32) replace $\tau$ with $r\tau$. This extended family defines a homomorphism over $\mathcal{M}^* \times [0, 1]$ between the corresponding extensions of the Banach space bundles $\mathcal{H}_1$ and $\mathcal{H}_2$. The respective extensions of these bundles over $\mathcal{M}^* \times [0, 1]$ are also denoted by $\mathcal{H}_1$ and $\mathcal{H}_2$, and the homomorphism between them by $\mathcal{D}$. The latter version of $\mathcal{D}$ has its associated determinant line, $\mathbb{D}et(\mathcal{D})$. Because the factor $[0, 1]$ is contractible, an orientation for $\mathbb{D}et(\mathcal{D})|_{r=0}$ canonically induces one for $\mathbb{D}et(\mathcal{D})|_{r=1}$ and thus for the line of interest, $\mathbb{D}et(\mathcal{D}) \to \mathcal{M}^*$.

*Part 2*: To define a canonical orientation for $\mathbb{D}et(\mathcal{D})$ over the $r = 0$ boundary of $\mathcal{M}^* \times [0, 1]$, focus for the moment on a given pair $((\tau, \mathcal{C} = \{C_S, \{C_\mathfrak{p}\}_{\mathfrak{p} \in \Lambda}\}), r = 0)$ on this boundary. The relevant version of (7.32) is such that there is no coupling between the summands that define the domain space for the operator $\mathcal{D}_\mathcal{C}$. As such, the kernel and cokernel of $\mathcal{D}_\mathcal{C}$ is the direct sum of the respective kernel and cokernel for the operator $D_{\eta 0}$ from Part 4 of the previous subsection, and the respective kernels and cokernels for each $C \in \{C_\mathfrak{p}\}_{\mathfrak{p} \in \Lambda}$ version of an operator that is described either by Proposition 3.1 or by Proposition 5.7.

Part 4 of the previous subsection asserts that $D_{\eta 0}$ is homotopic through a family of Fredholm operators with trivial kernel and trivial cokernel to the operator $D_S^\perp$. In particular, this same Part 4 finds that $\mathbb{D}et(D_{\eta 0})$ is canonically isomorphic to $\mathbb{D}et(D_S^\perp)$ and thus canonically oriented by the identification $\mathbb{D}et(D_S^\perp) = \wedge^2 \mathbb{R}^2$ from Part 2 of Section 8b.

Meanwhile, if $\mathfrak{p} \in \Lambda$ and $\Delta_\mathfrak{p} = 0$, then the corresponding $C = C_\mathfrak{p}$ contribution to $\mathcal{D}_\mathcal{C}$ is an operator that is described in Parts 2 and 3 of the previous subsection. In particular, these parts of Section 8c endow the corresponding determinant line with a completely canonical orientation by identifying it with $\wedge^2 \mathbb{R}^2$.

In the case when $\Delta_p > 0$, the contribution of $C_p$ to $\mathcal{D}_C$ is an operator of the sort that is describe by Proposition 5.7. Use D in what follows to denote this operator. What follows describes an absolutely canonical orientation for $\mathbb{D}et(D)$.

Proposition 5.7 asserts that D has $\Delta_p$-dimensional kernel and trivial cokernel. This being the case, its determinant line is oriented by an orientation of its kernel. To orient the kernel of D, suppose that $\mathcal{E} \subset \mathbb{R} \times \mathcal{H}^+_{p*}$ is an end of $C_p$ whose constant $s$ slices converge as $s \to \infty$ to $\hat{\gamma}^+_p$, this the ($\hat{u} = 0$, $\cos\theta = \frac{1}{\sqrt{3}}$) integral curve of $v$. Reintroduce the notation from Sections 5a and 5c so as to talk about the kernel of D on $\mathcal{E}$. As explained in Step 5 of Part 4 from Section 5c, if $\Delta_p = 1$ and $\mathfrak{m}_p = -1$, or if $\Delta_p = 2$, there is an element in the kernel of D that can be written at large values of $s_+$ on $\mathcal{E}$ as a map from the large $s_+$ part of $\mathbb{R} \times \mathbb{R}/(2\pi\mathbb{Z})$ to $\mathbb{R}^2$ that has the form

$$(s_+, \phi_+) \to c_{\mathcal{E}}(e^{-\lambda_1 s_+}, 0) + c_{\mathcal{E}} \mathfrak{e}$$

(8.11)

where $c_{\mathcal{E}} \in \mathbb{R}$ and where $|\mathfrak{e}| \leq c_0 e^{-(\lambda_1 + 1/c_0)s_+}$.

In the case $\Delta_p = 1$ and $\mathfrak{m}_p = -1$, the kernel of D is oriented by the unique element with $c_{\mathcal{E}} = 1$. In the case $\Delta_p = 1$ and $\mathfrak{m}_p = 1$, there is an analogous orientation for the kernel of D, this defined by the $c_{\mathcal{E}} = 1$ element with $c_E$ now defined by the analog of (8.11) for the end of $C_p$ whose constant $s$ slices converge to $\hat{\gamma}^-_p$ as $s \to \infty$.

If $\Delta_p = 2$, what is said in Step 5 of Part 4 of Section 5c implies that the kernel of D has a unique basis of the form $(\mathfrak{y}_+, \mathfrak{y}_-)$ with the following properties: The element $\mathfrak{y}_+$ is given by the $c_E = 1$ version of (8.11) on the end of $C_p$ whose constant $s$ slices converge to $\hat{\gamma}^+_p$ as $s \to \infty$. Meanwhile, $|\mathfrak{y}_+| \leq c_0 e^{-(\lambda_1 + 1/c_0)s}$ on the end of $C_p$ whose large $s$ slices converge to $\hat{\gamma}^-_p$ as $s \to \infty$. The situation is reversed for $\mathfrak{y}_-$; it is given by the $c_E = 1$ version of (8.11) on the end of $C_p$ whose constant $s$ slices converge to $\hat{\gamma}^-_p$ as $s \to \infty$, and it is bounded in absolute value by $c_0 e^{-(\lambda_1 + 1/c_0)s}$ on the end whose constant $s$ slices converge to $\hat{\gamma}^+_p$ as $s \to \infty$.

The canonical nature of these orientations is explained in Part 3.

*Part 3*: This part of the subsection explains the sense in which Part 2's orientation for the kernel of D is completely canonical.

To put things in a sufficiently general framework, use $\mathcal{O}$ now to denote the set whose elements have the form (D, Q, $\gamma_{p+}$, $\gamma_{p-}$, $\mathfrak{h} = (\varphi, \varsigma)$, D) with entries as follows: What is denoted by D is a data set $(z_*, \delta, x_0, R)$ that can be used to define $\mathcal{H}_p$ and the subspace $\mathcal{H}^+_{p*}$. Meanwhile, $Q \subset \mathbb{R}$ consists of a single point if $\Delta_p = 1$; it consists of two labeled points if $\Delta_p = 2$, one point labeled with a plus sign and the other with a minus sign. The pair $\gamma_{p+}$ and $\gamma_{p-}$ are integral curves of $v$ in $\mathcal{H}^+_{p*}$ with one boundary point on the $u < 0$ boundary of $\mathcal{H}^+_{p*}$ and the other boundary point is on the $u > 0$ boundary. What is denoted by $\mathfrak{h}$ signifies a smooth map from $(\mathbb{R} \times I_*)-Q$ to $\mathbb{R}^2$ of the sort that is described at the start of Section 5c; and D denotes an operator of the sort described by Proposition 5.7. The next two paragraphs define a topology for $\mathcal{O}$.

Let $\hat{o} = (D, Q, \gamma_{p+}, \gamma_{p-}, \mathfrak{h} = (\varphi, \varsigma), D)$ denote a given element in $\mathcal{O}$. A basis for the neighborhoods of $\hat{o}$ is indexed by $(\varepsilon, \mathcal{V}, \mathcal{W})$ where the notation is as follows: First, $\varepsilon$ is a positive number, and less than $10^{-4}$ times the distance between the points from Q if $\Delta_p = 2$. The definition of $\mathcal{V}$ and $\mathcal{W}$ requires the introduction of the set $Q_\varepsilon \subset \mathbb{R} \times I_*$, this the set of points with distance less than $\varepsilon$ from Q. What is denoted by $\mathcal{V}$ is an open neighborhood of $(\varphi, \varsigma)$ in the $C^\infty$-Frechet topology on the space of smooth maps from $(\mathbb{R} \times I_*)-Q_\varepsilon$ to $\mathbb{R}^2$. To describe $\mathcal{W}$, let $(\mathfrak{a}_1, \mathfrak{a}_2, \mathfrak{b}_1, \mathfrak{b}_2)$ denote the coefficient functions that define D via (3.5). What is donoted by $\mathcal{W}$ is an open neighborhood of $(\mathfrak{a}_1, \mathfrak{a}_2, \mathfrak{b}_1, \mathfrak{b}_2)$ in the $C^\infty$ Frechet topology on the space of maps from the domain $(\mathbb{R} \times I_*)-Q_\varepsilon$ to $\mathbb{R}^4$.

Let $\mathcal{U} \subset \mathcal{O}$ denote the neighorhood of $\hat{o}$ with the given indexing set. A point $\hat{o}' = (D', Q', (\gamma_{p+}', \gamma_{p-}') \mathfrak{h}', D')$ from $\mathcal{O}$ lies in $\mathcal{U}$ when the conditions listed next are met. Each entry of D' has distance less than $\varepsilon$ from the corresponding entry of D. Corresponding points from Q and Q' have distance less than $10^{-1/4}\varepsilon$ from each other; and each of the endpoints of $\gamma_{p+}'$ and $\gamma_{p-}'$ has distance less than $\varepsilon$ from the corresponding endpoint of the respective segments $\gamma_{p+}$ and $\gamma_p$. The map from $\mathfrak{h}'$ from $(\mathbb{R} \times I_*)-Q'$ to $\mathbb{R}^2$ lies in $\mathcal{V}$; and the map from this same domain to $\mathbb{R}^4$ given by sending any given $(x, \hat{u})$ to the D' version of $(\mathfrak{a}_1, \mathfrak{a}_2, \mathfrak{b}_1, \mathfrak{b}_2)$ lies in $\mathcal{W}$.

It follows from what is said in Proposition 5.1, Lemma 5.9, Lemma 3.2 and Steps 3-4 in Section 3c that the assignment to any given data set $\hat{o} \in \mathcal{O}$ of its operator D defines a continous section of a vector bundle over $\mathcal{O}$ of the form $\text{Hom}(\mathcal{H}_1, \mathcal{H}_2)$ where $\mathcal{H}_1$ and $\mathcal{H}_2$ are Banach space bundles over $\mathcal{O}$. Let $\mathcal{D}$ denote here this section. There is the corresponding determinant line bundle, $\mathbb{D}\text{et}(\mathcal{D}) \to \mathcal{O}$.

The lemma that follows makes the salient observations about $\mathcal{O}$ and $\mathbb{D}\text{et}(\mathcal{D})$.

**Lemma 8.4**: *Define the space $\mathcal{O}$ as above.*
- *The line bundle $\mathbb{D}\text{et}(\mathcal{D}) \to \mathcal{O}$ is orientable.*
- *$\mathcal{O}$ is path connected if $\Delta_p = 1$ and it has two path connected components if $\Delta_p = 2$. In the latter case the components are distinguished by whether the + labeled point from Q is greater than or less than the - labeled point.*

This lemma is proved momentarily. The next paragraph explains how this lemma leads to a canonical orientation for the versions of $\mathbb{D}et(D)$ that arise in Part 2 of this subsection.

Let ô ∈ $\mathcal{O}$ and let D denote ô's operator. Proposition 5.7 guarantees that D has trivial cokernel and $\Delta_p$-dimensional kernel. The paragraph subsequent to (5.11) at the end of Part 2 can be repeated to define a canonical basis for the kernel of D. Lemma 8.4 guarantees that the various ô ∈ $\mathcal{O}$ versions of this basis defines a canonical isomorphism between $\mathbb{D}et(\mathcal{D})$ and the oriented product bundle $\mathcal{O} \times \mathbb{R}$. It is in this sense that the versions of $\mathbb{D}et(D)$ from Part 2 have a completely canonical orientation.

*Proof of Lemma 8.4*: To prove the first bullet, look first at Proposition 5.7 to see that each ô ∈ $\mathcal{O}$ versions of D has $\Delta_p$-dimensional kernel. This being the case, these kernels fit together as ô varies in $\mathcal{O}$ to define a $\Delta_p$-dimensional vector bundle ker($\mathcal{D}$) → $\mathcal{O}$. The basis given in the paragraphs subsequent to (5.11) define an isomorphism from this bundle to the product bundle.

The proof of the second bullet has eight steps. The first step reviews some background material. In the $\Delta_p = 1$ case, the remaining seven steps construct a continuous path parameterized by [0,7] that starts at a given element ô ∈ $\mathcal{O}$ and ends at any given second element ô´ ∈ $\mathcal{O}$. If $\Delta_p = 2$, these same steps construct a continuous path parameterized by [0,6] that starts at a given element ô ∈ $\mathcal{O}$ and ends at any given second element ô´ ∈ $\mathcal{O}$ with the following property: Let Q and Q´ denote the respective ô and ô´ versions of the ± labeled points in $\mathbb{R}$. The labels of the largest points in Q and Q´ agree. These last seven steps construct the desired path as an end-to-end concatenation of seven [0,1]-parametrized paths. By way of notation, each step writes the starting point of its segment as (D, Q, $\gamma_{p+}$, $\gamma_{p-}$, $\mathfrak{h} = (\varphi, \varsigma)$, D). They all write the point ô´ in similar fashion using primes to distingish respective components that differ.

<u>Step 1</u>: Suppose that $\mathcal{E}$ is an end in ô's submanifold $C_\mathfrak{h}$ whose constant $s$ slices converge to $\hat{\gamma}_p^+$ as $s \to \infty$. Use $\mathcal{E}´$ to denote the corresponding end of $C_{\mathfrak{h}´}$. Parametrize both as in Section 5a and Proposition 5.1; use $\mathfrak{y}$ and $\mathfrak{y}´$ denote the respective parametrizing maps. The map $\mathfrak{y}$ is is defined in part by data ($c_{n=1}$, $\phi_{n=1}$) where $c_{n=1} \in \mathbb{R} - 0$ and $\phi_1 \in \mathbb{R}/2\pi\mathbb{Z}$. The constant $c_1$ can be assumed positive; for if not, there is an equivalent parametrization with $\phi_1$ replaced by $\phi_1 + \pi$ and $c_1$ replaced by $-c_1$. Let $c_1´ > 0$ and $\phi_1´$ denote the corresponding $\mathfrak{y}´$ versions of these parameters with $c_1´ > 0$ also. There are corresponding versions of $c_1$ and $c_1´$ when the respective constant $s$ slices of $\mathcal{E}$ and $\mathcal{E}´$ converge to $\hat{\gamma}_p^-$ as $s \to \infty$.

Step 2: The first path moves ô's subset Q to the corresponding ô´ subset Q´. Such a path is readily constructed from a compactly supported isotopy of $\mathbb{R} \times I_*$ that moves only the $\mathbb{R}$ coordinate of any given point, moves a neighborhood of each point from Q by a rigid translation, and takes Q to Q´. The $\tau = 1$ member of this path is now denoted by ô.
.

Step 3: This step constructs a $[0,1]$ parametrized path $\tau \to \hat{o}_\tau$ in $\mathcal{O}$ whose $\tau = 0$ member is ô and whose $\tau = 1$ member is given by $\hat{o}_1 = (D, Q, \gamma_{p+}, \gamma_{p-}, \mathfrak{h}_1 = (\phi_1, \varsigma_1), D_1)$ where $(\phi_1, \varsigma_1) = (\phi, \varsigma)$ except in very small radius disks about the point or points in Q. In a somewhat small radius disk $(\phi_1, \varsigma_1) = (\phi´, \varsigma´)$. Likewise, $D_1 = D$ on the first disk and $D_1 = D´$ on the smaller radius disk. To do this write $C_\mathfrak{h}$ and $C_{\mathfrak{h}´}$ on corresponding ends $E$ and $E´$ using the maps $\mathfrak{y}$ and $\mathfrak{y}´$ as in Step 1. Let $(c_1, \phi_1)$ and $\mathfrak{e}_1$ denote the data that appears in Proposition 5.1's depiction of $\mathfrak{y}$, and let $(c_1´, \phi_1´)$ and $\mathfrak{e}_1´$ denote the corresponding $\mathfrak{y}´$ data set. Write $\phi_1´ = \phi_1 + \iota$ with $\iota \in [0, 2\pi)$. Set $\phi_{1\tau} = \phi_1 + \tau\iota$ for $\tau \in [0,1]$. Meanwhile, set $c_{1\tau} = c_1 + \tau(c_1´ - c_1)$ and set $\mathfrak{e}_{1\tau} = \mathfrak{e}_1 + \tau(\mathfrak{e}_1´ - \mathfrak{e}_1)$ so as to define

$$\mathfrak{y}_{*\tau} = \alpha_y(e^{-\lambda_1 s_+} + e_1, 0) + c_{1\tau} e^{-\lambda_{11} s_+}(\cos(\phi_+ - \phi_{1\tau}), r_{11}\sin n(\phi_+ - \phi_{1\tau})) + \mathfrak{e}_{1\tau}. \quad (8.12)$$

The arguments from the Steps 1-3 of the proof of Proposition 5.1 can be used in an almost verbatim fashion to find $s_* > 1$ and a smooth, $[0,1]$-parametrized family of maps, $\tau \to \mathfrak{q}_{*\tau}: [s_*, \infty) \times \mathbb{R}/(2\pi\mathbb{Z})$, with the following three properties: First, $\mathfrak{y}_\tau = \mathfrak{y}_{*\tau} + \mathfrak{q}_{*\tau}$ obeys (5.7). Second, $|\mathfrak{y}_\tau - \mathfrak{y}_{*\tau}| \leq \frac{1}{1000}\min\{c_1, c_2\} e^{-\lambda_{11} s_+}$. Third, $\mathfrak{q}_{*0} = \mathfrak{q}_{*1} = 0$.

Given the family $\{\mathfrak{y}_\tau\}_{\tau\in[0,1]}$, use what is said in Lemma 5.3 to write each $\tau \in [0,1]$ version as a map, $(x, \hat{u}) \to (\phi_\tau, \varsigma_\tau)|_{(x, \hat{u})}$ with image $\mathbb{R}^2$ and domain the complement of the given point from Q in a small radius disk in $\mathbb{R} \times I_*$ about this point. Use $r$ to denote the radius of this disk and use $\chi_r$ to denote the function on $\mathbb{R} \times I_*$ given by the rule $(x, \hat{u}) \to \chi(\frac{4}{r}((x-y)^2 + \hat{u}^2)^{1/2} - 1)$. If Q has two points, do this for both.

Extend the corresponding family of maps defined in the radius $r$ disk or disks about the points in Q over the whole of $(\mathbb{R} \times I_*) - Q$ as $\mathfrak{h}$ on the complement of the disk or disks and as the relevant $(\phi_\tau, \varsigma_\tau)$ version of $\mathfrak{h}_\tau = (1 - \chi_r)(\phi, \varsigma) + \chi_r(\phi_\tau, \varsigma_\tau)$ in each disk. Use $\{\mathfrak{h}_\tau\}_{\tau\in[0,1]}$ to denote this family of maps from $(\mathbb{R} \times I_*) - Q$ to $\mathbb{R}^2$ and for each $\tau \in [0, 1]$. Use $D_\tau$ to denote the $\mathfrak{h}_\tau$ version of (3.9) and use $\hat{o}_\tau$ to denote $(D, Q, \gamma_{p+}, \gamma_{p-}, \mathfrak{h}_\tau, D_\tau)$. The assignment $\tau \to \hat{o}_\tau$ defines a continuous path in $\mathcal{O}$ with the desired $\tau = 0$ and $\tau = 1$ members. Use $\hat{o} \in \mathcal{O}$ henceforth to denote the $\tau = 1$ member of this family.

Step 4: Use cut-off functions in the manner of Parts 1-3 of Section 6a to construct a continuous, $[0, 1]$-parametrized path $\tau \to \hat{o}_\tau = (D, Q, \gamma_\tau, \gamma_{p-}, \mathfrak{h}_\tau, D_\tau)$ in $\mathcal{O}$ whose $\tau = 0$ member is ô and whose $\tau = 1$ member is such that $\gamma_{p_-}$ and $\gamma_1$ have the same respective $\phi$ angles at their endpoints, but are such that the respective change in $\phi$ differs by $-2\pi$. Meanwhile, each $\tau \in [0, 1]$ version of $\mathfrak{h}_\tau$ agrees with $\mathfrak{h}$ in the radius $\frac{1}{2}$ disk or disks centered on the points in Q. In addition, the $\hat{u} = -R - \frac{1}{2}\ln z_*$ boundary values of $\mathfrak{h}_1$ have constant $\hat{\phi}$ component; and the corresponding $\hat{u} = R + \frac{1}{2}\ln z_*$ component winds once around $\mathbb{R}/(2\pi\mathbb{Z})$ in the anti-clockwise direction as x varies from $-\infty$ to $\infty$. By way of a parenthetical remark, the $\Psi_\mathfrak{p}$ image of the graph in $\mathbb{R} \times \mathcal{X}$ of $\mathfrak{h}_\tau$ can be guaranteed J-holomorphic only where $|s| \gg 1$. The needed modifications to what is done in Parts 1-3 of Section 6a are straightforward and are omitted. What is denoted by $D_\tau$ is the $\mathfrak{h}_\tau$ version of the operator in (3.9). Use $ô \in \mathcal{O}$ now to denote the $\tau = 1$ member of this path

Step 5: Cut-off functions in the manner of Parts 1-3 of Section 6a are now used to construct a continuous, $[0,1]$-parametrized path $\tau \to ô_\tau = (D_\tau, Q, \gamma_\tau, \gamma_{p_-,\tau}, \mathfrak{h}_\tau, D_\tau)$ in $\mathcal{O}$ whose $\tau = 0$ member is ô and whose $\tau = 1$ member is such that $D_{\tau=1} = D´$. Meanwhile, the $\phi$ angles of the boundary points of both $\gamma_{p_-,\tau}$ and $\gamma_\tau$ are independent of $\tau$; Lemma 2.1 is invoked to arrange this. The restriction of $\mathfrak{h}_\tau$ to a very small disk or disks centered on the point or points in Q is also independent of $\tau$. As in the previous steps, $D_\tau$ denotes the $\mathfrak{h}_\tau$ version of (3.9). Use ô now to denote the $\tau = 1$ member of this path

Step 6: The constructions in Parts 1-3 of Section 6a are used yet again, this time to construct a continuous, $[0,1]$-parametrized path $\tau \to ô_\tau = (D´, Q, \gamma_\tau, \gamma_{p_-,\tau}, \mathfrak{h}_\tau, D_\tau)$ which is such that $\gamma_{p_-,\tau=1} = \gamma_{p_-}´$ and $\gamma_{\tau=1}$ is the $\gamma_{p_-}´$ analog of what is denoted by $\gamma_1$ in Step 4. This is done by moving the endpoints while invoking Lemma 2.1. The map $\mathfrak{h}_\tau$ restricts to a very small radius disk or disks centered on the point or points in Q to be independent of $\tau$. As before, $D_\tau$ denotes the $\mathfrak{h}_\tau$ version of (3.9). Use $ô \in \mathcal{O}$ henceforth to denote the $\tau = 1$ member of this fifth segment.

Step 7: The construction in Step 4 is run in reverse to construct a $[0,1]$ parametrized path $\tau \to ô_\tau$ in $\mathcal{O}$ that moves $\gamma_{p_+}$ so that $ô_1 = (D´, Q, \gamma_{p_+}´, \gamma_{p_-}´, \mathfrak{h}_1, D_1)$ with the path such that each $\tau \in [0,1]$ member of $\mathfrak{h}_\tau$ is again independent of $\tau$ on some small radius disk or disks centered on the point or points in Q. The operator $D_\tau$ is the $\mathfrak{h}_\tau$ version of (3.9). As in the previous steps, use $ô \in \mathcal{O}$ to denote the $\tau = 1$ version of this path.

Step 8: This final leg of the path, denoted by $[0, 1] \to ô_\tau$, is a family of data sets that have the form $(D´, Q, \gamma_{p_+}´, \gamma_{p_-}´, \mathfrak{h}_\tau, D_\tau)$. The $\tau = 1$ member is $ô´$. The family is defined using a suitable 1-parameter family of cut-off functions to homotope $\mathfrak{h}$ to $\mathfrak{h}´$. The details are omitted as they contain no novelties.

**e) Canonical orientations for $\mathcal{M}_{\mathfrak{p}\pm}$**

Proposition II.3.4 introduces for each $\mathfrak{p} \in \Lambda$ a pair of moduli spaces, $\mathcal{M}_{\mathfrak{p}+}$ and $\mathcal{M}_{\mathfrak{p}-}$ whose consitutents are embedded disks in the part of the $\hat{u} = 0$ locus of $\mathbb{R} \times \mathcal{H}_\mathfrak{p}$ where $1 - 3\cos^2\theta < 0$. Those in $\mathcal{M}_{\mathfrak{p}+}$ sit where $\cos\theta > \frac{1}{\sqrt{3}}$ and their constant $s$ slices converge as $s \to \infty$ in an isotopic fashion to $\hat{\gamma}_\mathfrak{p}^+$. Those in $\mathcal{M}_{\mathfrak{p}-}$ sit where $\cos\theta < -\frac{1}{\sqrt{3}}$ and their constant $s$ slices converge as $s \to -\infty$ to $\hat{\gamma}_\mathfrak{p}^-$. What follows talks solely about $\mathcal{M}_{\mathfrak{p}+}$ as the $\mathcal{M}_{\mathfrak{p}-}$ story is identical save for changing $\theta$ to $\pi - \theta$.

Let $C \in \mathcal{M}_{\mathfrak{p}+}$ denote a given curve. By way of a reminder from Section II.3c, the curve $C$ is invariant under the action of $S^1$ that rotates the angle $\phi$ on $\mathcal{H}_\mathfrak{p}$ and so the vector field $\partial_\phi$ is tangent to $C$. The vector field depicted in Equation (II.3.10) is proportional to $J\partial_\phi$ and so is also tangent to $C$ also. This being the case, $C$ is foliated by the integral curves of the latter vector field except for the single point on $C$ where it and $\partial_\phi$ have are zero. This point is the $\theta = 0$ point on $C$ and the infimum of $s$ on $C$. The latter point is the only critical point of $s$ on $\phi$. Use $y \in \mathbb{R}$ to denote this point. The association to each curve in $\mathcal{M}_{\mathfrak{p}+}$ of the minimum of $s$ on the curve defines an $\mathbb{R}$-equivariant diffeomorphism from $\mathcal{M}_{\mathfrak{p}+}$ to $\mathbb{R}$.

Let $U_+ \subset \mathcal{H}_\mathfrak{p}$ denote the tubular neighborhood of $\hat{\gamma}_\mathfrak{p}^+$ that is described in Part 2 of Section 5a. The $s \gg y$ part of the given curve $C$ lies in $\mathbb{R} \times U_+$ and so can be described using the coordinates $(s_+, \phi_+, \theta_+, u_+)$ as defined in (5.5). As such, it appears as the graph of a function $(s_+, \phi_+) \to (s_+, \phi_+, a, b)$ where $\eta = (a, b)$ is given by (5.8) with $\alpha = \alpha_y < 0$.

Since $C$ is not $\mathbb{R}$-invariant, the normal projection of the vector field $\partial_s$ along $C$ supplies a canonical element to the kernel of $C$'s version of the operator $\mathcal{D}_C$. Denote the latter by $\eta_C$. Since the vector field in (3.10) is not proportional to $\partial_s$ along $C$, this canonical element $\eta_C$ is nowhere zero.

To say more about $\eta_C$, use (5.5) with (5.8) to identify the normal bundle to $C$ along $C \cap (\mathbb{R} \times U_+)$ with the product bundle using the 1-forms $(d\theta_+, du_+)$. Granted this identification, and with $C \cap (\mathbb{R} \times U_+)$ parametrized as a graph in the manner just described, it follows from (5.8) that $\eta_C$ appears as a map from the $s_+ \gg 1$ part of $\mathbb{R} \times \mathbb{R}/(2\pi\mathbb{Z})$ to $\mathbb{R}^2$ that can be written as

$$(s_+, \phi_+) \to -\lambda_+ \alpha_y (e^{-\lambda_1 s_+} + \mathfrak{e}, 0)$$

(8.13)

where $|\mathfrak{e}| \leq |\alpha_y| e^{-(\lambda_1 + 1/c_0) s_+}$.

A parenthetical remark subsequent to Proposition II.3.4 asserts that the cokernel of $\mathcal{D}_C$ is trivial. This assertion is proved momentarily. It implies that the kernel of $\mathcal{D}_C$ is 1-dimensional, this being the span of $\eta_C$. The identification between the kernel of $\mathcal{D}_C$ and the tangent space at C to $\mathcal{M}_{\mathfrak{p}+}$ gives the canonical section $\eta_C$, and this section defines the desired canonical orientation.

To see why cokernel($\mathcal{D}_C$) = 0, first identify the cokernel with the kernel of C's formal $L^2$ adjoint. This done, use the previously described parametrizations of C and C's normal bundle on $C \cap (\mathbb{R} \times U_+)$ to view an element in the cokernel of this adjoint operator as a map from the $s_+ \gg 1$ portion of $\mathbb{R} \times \mathbb{R}/(2\pi\mathbb{Z})$ to $\mathbb{R}^2$. To see what such an element looks like, view $\mathcal{D}_C$ on this part of C as an operator of the form $\mathfrak{D}_0 + \mathfrak{d}$ with $\mathfrak{D}_0$ as in (5.1) and with $\mathfrak{d}$ a first order operator on the space of maps from this large $s_+$ portion of $\mathbb{R} \times \mathbb{R}/(2\pi\mathbb{Z})$ to $\mathbb{R}^2$ with $\phi$-invariant coefficients that are bounded in absolute value by $c_0 |\alpha_\mathfrak{y}| e^{-(\lambda_1 + 1/c_0)s_+}$. With $\mathcal{D}_C$ written as $\mathfrak{D}_0 + \mathfrak{d}$, an element in the kernel of $\mathcal{D}_C$'s formal $L^2$ adjoint can be written on C's intersection with $\mathbb{R} \times U_+$ as

$$(s_+, \phi_+) \to c_0(0, e^{-\lambda_2 s_+} + \mathfrak{e}_0) + c_n(\mathfrak{y}_{n-}|_{(-s_+, \phi_+)} + \mathfrak{e}_n)$$

(8.14)

where the notation is as follows: First, $c_0 \in \mathbb{R}$, $c_n \in \mathbb{R}$ and one of these is not zero. Meanwhile, $\mathfrak{y}_{-n}$ is defined in (5.4). What is written as $\mathfrak{e}_0$ in (8.14) is $\phi$-independent and obeys $|\mathfrak{e}_0| \le c_0 e^{-(\lambda_2 + 1/c_0)s_+}$, and what is written as $\mathfrak{e}_n$ obeys $|\mathfrak{e}_n| \le c_0 e^{-(\lambda_{2n} + 1/c_0)s_+}$.

Granted that a cokernel element is a section of $N \otimes T^{0,1}C$, and granted that C is a disk, the claim that the cokernel of $\mathcal{D}_C$ is trivial follows from (8.14) because the latter forces any non-zero element in the kernel of C's adjoint to vanish at some point on C with positive local degree. This sort of vanishing is not possible by virtue of the fact that the formal $L^2$ adjoint of $\mathcal{D}_C$ differs from that of $\bar{\partial}$ by a zero'th order endomorphism.

**f) Canonical orientations for the $I_{ech}$ = 1 moduli spaces**

Assume that the defining data for the geometry of Y and $\mathbb{R} \times Y$ are such that what is said in Sections 1-7 and 8a-e hold. This subsection defines a canonical orientation for the $I_{ech}$ = 1 moduli spaces, $\{\mathcal{M}_1(\hat{\Theta}', \hat{\Theta})\}_{\hat{\Theta}', \hat{\Theta} \in \hat{\mathcal{Z}}_{ech,M}}$. The desired orientations are defined with the help of a given orientation for the real line bundle over $\mathcal{A}_{HF1}$ whose fiber over any given surface S is the determinant line for the operator $D_S$. This line bundle is denoted by $\mathbb{D}et(\hat{D})$. A given orientation for $\mathbb{D}et(\hat{D})$ is assumed in the four parts that follow. Also needed is a chosen ordering for the set $\Lambda$. The resulting ordered set is written as $\{\mathfrak{p}_1, \ldots, \mathfrak{p}_G\}$ when the ordering is relevant.

*Part 1*: Fix a pair $(\hat{\Theta}', \hat{\Theta}) \in \hat{\mathcal{Z}}_{ech,M} \times \hat{\mathcal{Z}}_{ech,M}$ and write these two elements respectively as $((\hat{\mathfrak{v}}', k'), (\mathfrak{k}_\mathfrak{p}', O_\mathfrak{p}')_{\mathfrak{p} \in \Lambda})$ and $((\hat{\mathfrak{v}}, k), (\mathfrak{k}_\mathfrak{p}, O_\mathfrak{p})_{\mathfrak{p} \in \Lambda})$. It follows from Propositions 8.1 and 8.2 that the moduli space $\mathcal{M}_1(\hat{\Theta}', \hat{\Theta})$ is non-empty if and only if one of the two conditions listed in the upcoming (8.15) hold. By way of notation, (8.15) introduces $\mathcal{A}_{HF1}((\hat{\upsilon}',k'),(\hat{\upsilon},k))$ to denote the space of Lipshitz submanifolds with the following three properties: If $S \in \mathcal{A}_{HF1}((\hat{\upsilon}',k'),(\hat{\upsilon},k))$, then the constant $s$ slices of S converge in an isotopic fashion as $s \to \infty$ to the arcs that comprise $\hat{\upsilon}$; and they converge in an isotopic fashion as $s \to -\infty$ to the arcs that comprise $\hat{\upsilon}'$. In addition, the surface S has intersection number $k - k'$ with the arc $\gamma^{(z_0)}$. Finally, the operator $D_S$ has Fredholm index 1. Equation (8.15) also uses $\Delta_\mathfrak{p}$ to denote the number of elements in a given $\mathfrak{p} \in \Lambda$ version of $o_\mathfrak{p}$.

- $\mathcal{A}_{HF1}((\hat{\upsilon}',k'),(\hat{\upsilon},k)) \neq \emptyset$ *and each* $\mathfrak{p} \in \Lambda$ *version of* $(\mathfrak{k}_\mathfrak{p}, o_\mathfrak{p})$ *equals* $(\mathfrak{k}_\mathfrak{p}', o_\mathfrak{p}')$.
- $(\hat{\upsilon}',k') = (\hat{\upsilon},k)$ *and there exists precisely one* $\mathfrak{p} \in \Lambda$ *such that* $(\mathfrak{k}_\mathfrak{p}, o_\mathfrak{p}) \neq (\mathfrak{k}_\mathfrak{p}', o_\mathfrak{p}')$. *In this case,* $\Delta_\mathfrak{p}' = \Delta_\mathfrak{p} - 1$ *and one of the following holds:*
  a) $\mathfrak{k}_\mathfrak{p}' = \mathfrak{k}_\mathfrak{p}$.
  b) $\mathfrak{k}_\mathfrak{p}' = \mathfrak{k}_\mathfrak{p} \pm 1$.

(8.15)

The upcoming Part 2 specifes the desired orientations for the case of the first bullet in (8.15), Part 3 considers the case of Item a) of the second bullet, and Part 4 speaks to the case of Item b) of the second bullet.

*Part 2*: Suppose that $\mathcal{M}_1(\hat{\Theta}', \hat{\Theta})$ is described by the first bullet in (8.15). Any given element $\vartheta \in \mathcal{M}_1(\hat{\Theta}', \hat{\Theta})$ is a disjoint union of components with some union from a version of $\mathcal{M}^*$ that is defined by a surface $S \in \mathcal{A}_{HF1}$. The other components are $\mathbb{R}$-invariant cylinders of the form $\mathbb{R} \times \hat{\gamma}_\mathfrak{p}^+$ or $\mathbb{R} \times \hat{\gamma}_\mathfrak{p}^-$ for various $\mathfrak{p} \in \Lambda$. The contribution to $\vartheta$ from any given $\mathfrak{p}$ depends on $o_\mathfrak{p}$; either none, one or both of these $\mathfrak{p}$-labeled cylinders can be present.

Let $\mathcal{C} = (1, C)$ denote $\mathcal{M}^*$ part of $\vartheta$. This submanifold may also be a union of components, but in any event, precisely one such component is not $\mathbb{R}$-invariant. In any event, the tangent space to the curve C is canonically identified with the kernel of the corresponding version of the operator $\hat{\mathcal{D}}_C$ and it is therefore oriented by a choice of orientation for the line $\mathbb{D}\text{et}(\hat{\mathcal{D}}_C)$. The desired orientation for the latter is supplied by (8.9) using the given orientation for $\mathbb{D}\text{et}(D_S)$ and the canonical orientations for the line $\mathbb{D}\text{et}(\mathcal{D}_C)$ given in Section 8c.

*Part 3*: Suppose that $\mathcal{M}_1(\hat{\Theta}', \hat{\Theta})$ is described by Item a) of the second bullet in (8.15). This version of $\mathcal{M}_1(\hat{\Theta}', \hat{\Theta})$ contains but a single $\mathbb{R}$-orbit. Let $\vartheta$ denote a given point on this orbit. The element $\vartheta$ has some union of $\mathbb{R}$-invariant components that define an element in a version of $\mathcal{M}^*$. The latter is defined by an HF-cycle. Part 5 in Section

II.2b associates an *orientation sign*, either +1 or -1, to each integral $\mathfrak{v}$ in this cycle. Let $N_+$ denote the number of positive orientation signs.

The given point $\vartheta$ also contains a union of $\mathbb{R}$-invariant cylinders, each of the form $\mathbb{R} \times \hat{\gamma}_{\mathfrak{p}'}^+$ or $\mathbb{R} \times \hat{\gamma}_{\mathfrak{p}'}^-$ for various $\mathfrak{p}' \in \Lambda - \mathfrak{p}$. The contribution to $\vartheta$ from any given such $\mathfrak{p}'$ depends on $o_{\mathfrak{p}'}$; either none, one or both of these $\mathfrak{p}$-labeled cylinders can be present. The nature of the remaining components depends on which of the cases listed below occur.

- $o_{\mathfrak{p}'} = \{0\}$ *and* $o_{\mathfrak{p}} = \{1\}$.
- $o_{\mathfrak{p}'} = \{-1\}$ *and* $o_{\mathfrak{p}} = \{1, -1\}$.
- $o_{\mathfrak{p}'} = \{0\}$ *and* $o_{\mathfrak{p}} = \{-1\}$.
- $o_{\mathfrak{p}'} = \{1\}$ *and* $o_{\mathfrak{p}} = \{1, -1\}$.

(8.16)

In the case of the first bullet the pair $\mathfrak{p}$ contributes to $\vartheta$ an element from Proposition II.3.4's moduli space $\mathcal{M}_{\mathfrak{p}+}$; and in the case of the second bullet, it contributes such an element and also the cylinder $\mathbb{R} \times \hat{\gamma}_{\mathfrak{p}}^-$. In either case, let $C_+$ denote the element from $\mathcal{M}_{\mathfrak{p}+}$. The tangent space at $\vartheta$ of $\mathcal{M}_1(\hat{\Theta}', \hat{\Theta})$ is canonically isomorphic to that of $\mathcal{M}_{\mathfrak{p}+}$ at $C_+$. Meanwhile $T\mathcal{M}_{\mathfrak{p}+}$ has its canonical orientation from Section 8e. Use $\mathfrak{o}$ to denote the orientation for $T\mathcal{M}_1(\hat{\Theta}', \hat{\Theta})$ that comes from this canonical orientation for $T\mathcal{M}_{\mathfrak{p}+}$. This may or may not be the desired orientation. To say if it is or not, introduce $k \in \{1, \dots, G\}$ to denote label for $\mathfrak{p}$ when $\Lambda$ is written as $\{\mathfrak{p}_1, \dots, \mathfrak{p}_G\}$ and introduce N to denote $N = \sum_{1 \leq k' < k} \Delta_{k'}$. The desired orientation for $T\mathcal{M}_1(\hat{\Theta}', \hat{\Theta})$ is $(-1)^{N_+ + N} \mathfrak{o}$.

In the case of the third bullet, the pair $\mathfrak{p}$ contributes to $\vartheta$ an element from Proposition II.3.4's moduli space $\mathcal{M}_{\mathfrak{p}-}$; and in the case of the fourth bullet, it contributes such an element and also the cylinder $\mathbb{R} \times \hat{\gamma}_{\mathfrak{p}}^+$. In either case, let $C$ denote the element from $\mathcal{M}_{\mathfrak{p}-}$. The tangent space at $\vartheta$ of $\mathcal{M}_1(\hat{\Theta}', \hat{\Theta})$ is canonically isomorphic to that of $\mathcal{M}_{\mathfrak{p}-}$ at $C$. The tangent bundle to $\mathcal{M}_{\mathfrak{p}-}$ also has a canonical orientation from Section 8e. Use $\mathfrak{o}$ now to denote the orientation for $T\mathcal{M}_1(\hat{\Theta}', \hat{\Theta})$ that comes from this canonical orientation for $T\mathcal{M}_{\mathfrak{p}-}$. Reintroduce $k \in \{1, \dots, G\}$ to denote label for $\mathfrak{p}$ when $\Lambda$ is written as $\{\mathfrak{p}_1, \dots, \mathfrak{p}_G\}$ and $N = \sum_{1 \leq k' < k} \Delta_{k'}$. The desired orientation for $T\mathcal{M}_1(\hat{\Theta}', \hat{\Theta})$ is $(-1)^{N_+ + N} \mathfrak{o}$ in the case of the third bullet in (8.16) and it is $(-1)^{N_+ + N + 1} \mathfrak{o}$ in the case of the fourth bullet.

*Part 4*: This part deals with Item b) of the second bullet in (8.15). Consider first the case where $\mathfrak{k}_{\mathfrak{p}}' = \mathfrak{k}_{\mathfrak{p}} + 1$. In this case, either the first or the second bullet in (8.16) holds. In any event, a given $\vartheta$ from $\mathcal{M}_1(\hat{\Theta}', \hat{\Theta})$ consists of a union of components. Some subset of these define an element $\mathcal{C} = (1, C)$ from a version of $\mathcal{M}^*$. The pair $\mathfrak{p}$ contributes the $\mathbb{R}$-invariant cylinder $\mathbb{R} \times \hat{\gamma}_{\mathfrak{p}}^-$ if and only if the second bullet in (8.16) is relevant. The various $\mathfrak{p}' \in \Lambda - \mathfrak{p}$ contribute either none, one or both $\mathbb{R}$-invariant cylinders from the set $\{\mathbb{R} \times \hat{\gamma}_{\mathfrak{p}'}^-, \mathbb{R} \times \hat{\gamma}_{\mathfrak{p}'}^+\}$, this depending as usual on $\mathfrak{o}_{\mathfrak{p}'}$. This understood, the tangent space to $\mathcal{M}_1(\hat{\Theta}', \hat{\Theta})$ at $\vartheta$ is canonically isomorphic to the kernel of the operator $\hat{\mathcal{D}}_C$. The latter is oriented via (8.9) using the given orientation for the relevant version of S and the canonical orientations supplied by Sections 8c and 8d for $\mathbb{D}\mathrm{et}(\mathcal{D}_C)$. Note in this regard that S is $\mathbb{R}$-invariant and so $D_S$ has trivial kernel and cokernel. Use $\mathfrak{o}$ to denote the orientation for $T\mathcal{M}_1(\hat{\Theta}', \hat{\Theta})$ at $\vartheta$ that comes via the afore-mentioned identification with the kernel of $\hat{\mathcal{D}}_C$. Use k again to denote the label for $\mathfrak{p}$ when $\Lambda$ is written as $\{1, \ldots, G\}$ and use N again to denote $\sum_{1 \leq k' < k} \Delta_{k'}$. The desired orientation for $\vartheta$'s component of $T\mathcal{M}_1(\hat{\Theta}', \hat{\Theta})$ is $(-1)^{N_+ + N} \mathfrak{o}$.

Suppose next that $\mathfrak{k}_{\mathfrak{p}}' = \mathfrak{k}_{\mathfrak{p}} - 1$. The story here is almost identical to that just told but for two salient and very much related changes. First, either the third or the fourth bullets in (8.16) hold. In the case of the fourth bullet the element $\vartheta$ contains the $\mathbb{R}$-invariant cylinder $\mathbb{R} \times \hat{\gamma}_{\mathfrak{p}}^+$. In either case, let $\mathfrak{o}$ again denote the orientation for $T\mathcal{M}_1(\hat{\Theta}', \hat{\Theta})$ at $\vartheta$ that comes via the canonical identification with the kernel of the relevant version of $\hat{\mathcal{D}}_C$ with it understood that the kernel of the latter is oriented using (8.9) as before. Reintroduce the integer N. The desired orientation for $\vartheta$'s component of $T\mathcal{M}_1(\hat{\Theta}', \hat{\Theta})$ is $(-1)^{N_+ + N} \mathfrak{o}$ if the third bullet of (8.16) is relevant, and it is $(-1)^{N_+ + N + 1} \mathfrak{o}$ if the fourth bullet in (8.16) is relevant.

**g) Coherent orientations**

The definition of the Heegard Floer differential requires the specification of an orientation for the low dimensional components of $\mathcal{A}_{HF}$. There are constraints on the choice that are described in Section 6 of [L]. A choice that obeys the constraints is said to be a *coherent system of orientations*. The definition of the embedded contact homology differential likewise requires the specification of suitably constrained orientations for the low dimensional components of $\mathcal{M}_{ech}$. Orientations that obey the latter constraints are also said to constitute a coherent system. The constraints are given in Section 9.5 of [HT2]. Proposition 8.5 in the upcoming Part 2 of this subsection makes a precise the assertion that a coherent system of orientations for $\mathcal{A}_{HF}$ leads to a suitably compatible coherent system of orientations for $\mathcal{M}_{ech}$. Part 1 of the subsection set up the needed background information.

*Part 1*: A closed integral curve of $v$ is said to be hyperbolic if the associated linearized return map in $SL(2; \mathbb{R})$ has two real eigenvalues with neither equal to 1 or -1. The integral curve is said to be positive when these eigenvalues are positive.

As noted previously, Section 9.5 in [HT2] defines what is meant by a coherent orientation for the components of $\mathcal{M}_{ech}$. There are four constraints which are labeled (OR1)–(OR4) in this section of [HT2]. Those labled (OR1) and (OR4) are normalization constraints that set the orientation in specific instances. The salient constraints are those expressed by (OR2) and (OR3). The former asserts that the orientation should be compatible with end-to-end concatenation of the subvarieties in $\mathcal{M}_{ech}$. The condition expressed by (OR3) constrains the orientation over the components of $\mathcal{M}_{ech}$ whose elements consist of disjoint unions of two or more submanifolds. This (OR3) constraint requires an a priori choice of ordering for the ends of any given element whose constant $s$ slices converge as $s \to \infty$ to positive hyperbolic integral curves of $v$. Also needed is a choice of ordering for those ends that converge as $s \to -\infty$ to positive hyperbolic integral curves of $v$.

Section 8f describes orientations for $\{\mathcal{M}_1(\hat{\Theta}', \hat{\Theta})\}_{\hat{\Theta}',\hat{\Theta} \in \hat{\mathcal{Z}}_{ech,M}}$, these being the 1-dimensional components of $\mathcal{M}_{ech}$. The definition in Section 8f makes no reference to an ordering of the relevant ends of the constituent submanifolds. Even so, Proposition 8.5 refers to the Section 8f orientations when describing a coherent system of orientations for $\mathcal{M}_{ech}$. This referral implicitly invokes the ordering given momentarily for the ends of the elements in $\{\mathcal{M}_1(\hat{\Theta}', \hat{\Theta})\}_{\hat{\Theta}',\hat{\Theta} \in \hat{\mathcal{Z}}_{ech,M}}$.

It proves useful to first define an ordering for the positive hyperbolic integral curves of $v$ in any given element from $\mathcal{Z}_{ech,M}$. To this end, let $\Theta$ denote a given element from $\mathcal{Z}_{ech,M}$ and let $\Theta^+ \subset \Theta$ denote the subset of positive, hyperbolic integral curves. Introduce $\Theta_{ech,M} \subset \Theta$ to denote the subset of closed integral curves of $v$ that cross one or more of the handles $\{\mathcal{H}_\mathfrak{p}\}_{\mathfrak{p} \in \Lambda}$. Given $\gamma \in \Theta_{ech,M}$, let $n_\gamma$ denote the smallest of the labels of those $\mathfrak{p} \in \Lambda$ with $\mathcal{H}_\mathfrak{p} \cap \gamma \neq \emptyset$. Order $\Theta_{ech,M}$ so that the corresponding ordered set of integers $\{n_\gamma\}_{\gamma \in \Theta_{ech,M}}$ is increasing. Use $\Theta^+_{ech,M} \subset \Theta_{ech,M}$ to denote the corresponding ordered subset of positive hyperbolic elements. Let $\mathfrak{o}_{\Theta,\mathfrak{p}}$ denote the subset from the set $\{\hat{\gamma}^+_\mathfrak{p}, \hat{\gamma}^-_\mathfrak{p}\}$ that come from $\Theta$. If this set has two elements, order it as just written. Use this convention to order $\Theta^+$ as

$$\Theta^+ = \{\Theta^+_{ech,M}, \mathfrak{o}_{\mathfrak{p}_1}, \ldots, \mathfrak{o}_{\mathfrak{p}_G}\}.$$

(8.17)

Let $\hat{\Theta}'$ and $\hat{\Theta}$ denote elements in $\hat{\mathcal{Z}}_{ech,M}$ and let C denote a given element in $\mathcal{M}_1(\hat{\Theta}', \hat{\Theta})$. What follows directly describes the convention for the ordering of C's ends for the case described by the first bullet in (8.15). The submanifold C may have more than one component, but only one is not $\mathbb{R}$-invariant. The remaining components are $\mathbb{R}$-invariant cylinders. Some union of the latter with the component that is not $\mathbb{R}$ invariant defines an ech-HF submanifold. Write this ech-HF submanifold as $C_1 \cup C_2$ with $C_1$ denoting the non-$\mathbb{R}$ invariant component. The desired ordering places the ends of $C_1$ before those of $C_2$. Meanwhile, the ends of $C_1$ are ordered amongst themselves so as to be consistent with the ordering of the relevant $\Theta$ or $\Theta'$ version of (8.17)'s set $\Theta^+_{\text{ech},M}$, and likewise those of $C_2$. All of the remaining $\mathbb{R}$-invariant cylinders in C are from the set $\{\mathbb{R} \times \hat{\gamma}_{\mathfrak{p}}^+, \mathbb{R} \times \hat{\gamma}_{\mathfrak{p}}^-\}_{\mathfrak{p} \in \Lambda}$. Their ends are ordered after the ends from $C_1 \cup C_2$. The union of the respective $s \gg 1$ and $s \ll -1$ ends from these sorts of $\mathbb{R}$-invariant cylinders is then ordered so as to be consistent with the ordering given by the relevant version of (8.17)'s ordered set $\{\mathfrak{o}_{\mathfrak{p}_1}, \ldots, \mathfrak{o}_{\mathfrak{p}_G}\}$.

The convention in the case described by the second bullet in (8.16) orders the relevant subset of $s \gg 1$ ends of C so as to give the ordering in (8.17), and likewise for the $s \ll -1$ ends of C.

*Part 2*: This part explains how a coherent system of orientations for $\mathcal{A}_{HF}$ leads to one for $\mathcal{M}_{\text{ech}}$ and in particular for the components of $\{\mathcal{M}_1(\hat{\Theta}', \hat{\Theta})\}_{\hat{\Theta}', \hat{\Theta} \in \hat{\mathcal{Z}}_{\text{ech},M}}$. The assertion that such is the case is given by Proposition 8.5. The proposition refers to the notion from Section 2a of a weakly compact $\mathcal{K} \subset \mathcal{A}_{HF}$. As noted in [L] and implied by Lemma II.6.6, the tangent space to $\mathcal{A}_{HF}$ at any given Lipshitz submanifold S has a canonical identification with the kernel of the corresponding operator $D_S$. This being the case, an orientation for the tangent space of $\mathcal{A}_{HF}$ is neither more nor less than an orientation for the real line bundle $\mathbb{D}\text{et}(\hat{D}) \to \mathcal{A}_{HF}$.

**Proposition 8.5**: *Suppose that a coherent system of orientations has been chosen for $\mathcal{A}_{HF}$ and thus for the line bundle $\mathbb{D}\text{et}(\hat{D})$. There exists a weakly compact set $\mathcal{K} \subset \mathcal{A}_{HF}$ that contains all elements in $\mathcal{A}_{HF1}$ and has the following significance: Choose a $\mathcal{K}$-compatible data set $\mathfrak{D} = (z_*, \delta, x_0, R)$ from the collection described in Proposition 7.2 for a suitable choice of $\kappa$, and choose the almost complex structure on $\mathbb{R} \times Y$ pursuant to the constraints in Proposition 7.1-7.3 and 8.1 and 8.2. Use this data to define $\mathcal{Z}_{\text{ech},M}$ and $\mathcal{M}_{\text{ech}}$. There exists a coherent system of orientations for $\mathcal{M}_{\text{ech}}$ whose restriction to $\{\mathcal{M}_1(\hat{\Theta}', \hat{\Theta})\}_{\hat{\Theta}', \hat{\Theta} \in \hat{\mathcal{Z}}_{\text{ech},M}}$ agrees with that defined in Section 8f using the orientations given by $\mathbb{D}\text{et}(\hat{D})$.*

**Proof of Proposition 8.5**: With the goal a proof of Theorem 1.1, coherent orientation systems are needed only for the one dimensional components and certain two dimensional components of $\mathcal{A}_{HF}$ and $\mathcal{M}_{\text{ech}}$. To keep this long paper form being even longer, the coherence for the orientations of $\mathcal{M}_{\text{ech}}$ will be verified only these relevant components.

The various $\mathfrak{p} \in \Lambda$ versions of $\hat{\gamma}_\mathfrak{p}^+$ and $\hat{\gamma}_\mathfrak{p}^-$ are all positive hyperbolic closed integral curves of $v$. Proposition II.2.7 characterizes the other positive hyperbolic integral curves of $v$ that can appear in any given element from $\mathcal{Z}_{ech,M}$. With the preceding understood, it is a straightforward task to verify that the orientations for the various elements in $\{\mathcal{M}_1(\hat{\Theta}', \hat{\Theta})\}_{\hat{\Theta}', \hat{\Theta} \in \hat{\mathcal{Z}}_{ech,M}}$ given in Section 8f using Part 1's ordering of the relevent ends the constituent submanifolds obey the (OR3) constraint from Section 9.5 in [HT2]. The (OR3) constraint is the only relevant constraint on the 1-dimensional components of $\mathcal{M}_{ech}$.

Consider next the case where $C \subset \mathcal{M}_{ech}$ is a 2-dimensional component. There are two sorts of components to consider. The first is the union of cross-sectional 2-sphere in the handle $\mathcal{H}_0$ with $\mathbb{R}$-invariant cylinders. The sphere component is oriented using the convention (OR1) from Section 9.5 of [HT2]. The constraint in (OR3) of [HT2] is obeyed as long as the ordering for the set of $s \gg 1$ ends is the same as that for the set of $s \ll -1$ ends. The other constraints in Section 9.5 of [HT2] are not relevant.

The other sort of 2-dimensional component contains submanifolds that lie entirely in $\mathbb{R} \times (M_\delta \cup (\cup_{\mathfrak{p} \in \Lambda} \mathcal{H}_\mathfrak{p}))$. The only salient constraint to consider for this case is that given by (OR2). The concern with (OR2) arises when a 2-dimensional component of $\mathcal{M}_{ech}$ has the following property: Fix $\varepsilon > 0$ and the component has two or more open subsets that are described by Proposition II.7.2 using distinct versions of $\Xi$ of the form $\{\mathcal{Z}_1, \mathcal{Z}_2\}$ where $\mathcal{Z}_1 = \{(S_1, u_1), \emptyset\}$ and $\mathcal{Z}_2 = \{(S_2, u_2), \emptyset\}$ are such that $S_1$ and $S_2$ come from $\mathcal{A}_{HF1}$. What with Lemma 9.6 in [HT2], the end to end concatentation using any such version of $\Xi$ orients the relevant component. The constraint (OR2) requires that all such orientations agree.

To see about (OR2), note that end to end concatenations of the pair $S_1$ and $S_2$ from any given such $\Xi$ supply Lipshitz submanifolds in a 2-dimensional component of $\mathcal{A}_{HF}$. The construction is described in Appendix B of [L]. Moreover, what is said in this appendix implies a Heegard Floer version of Lemma 9.6 in [HT2]. This analog orients the relevant 2-dimensional component of $\mathcal{A}_{HF}$ given orientations for $\mathbb{D}et(S_1)$ and $\mathbb{D}et(S_2)$. With the preceding understood, suppose for a moment that all relevant versions of $\Xi$ define in this way the same 2-dimensional component of $\mathcal{A}_{HF}$. The corresponding set of orientations for this component will agree if the orienations for the components of $\mathcal{A}_{HF}$ constitute a coherent system. Meanwhile, an appropriate choice for the set $\mathcal{K}$ has the following property: Fix a 2-dimensional component of $\mathcal{M}_{ech}$ of the sort under consideration and there is but one 2-dimensional components of $\mathcal{A}_{HF}$ that can arise in this manner.

Granted the preceding, a straightforward modification to what is said in Step 3 in the proof of Proposition 7.1 for a suitable version of (7.43) proves the following: If the set $\mathcal{K}$ and Proposition 7.2's constant $\kappa$ are chosen appropriately, then the isomorphism given by (8.9) is compatible with respect to end-to-end concatenations of ech-HF subvarieties on the one hand and Heegard Floer subvarieties from $\mathcal{A}_{HF1}$ on the other. This fact implies that the (OR2) constraint in Section 9.5 of [HT2] is obeyed by the 2-dimensional components of $\mathcal{M}_{ech}$.

## 9. Proof of Theorem 1.1

This section uses the results from the previous two sections to prove the assertions made by the various bullets in Theorem 1.1.

### a) The grading of the ech chain complex

This section addresses the assertion made by the fourth bullet of Theorem 1.1. The proof of this fourth bullet has five parts.

*Part 1*: Suppose that $(\hat{\Theta}', \hat{\Theta})$ is an ordered pair of elements from $\hat{\mathcal{Z}}_{ech,M}$. The grading difference $gr_{ech}(\hat{\Theta}') - gr_{ech}(\hat{\Theta})$ is equal modulo the integer $p_M$ to -1 times the ech index $I(\cdot)$ of a suitable relative 2-cycle. In particular, suppose that $k \in \{0, 1, 2...\}$ and that C is an element in $\mathcal{M}_k(\hat{\Theta}', \hat{\Theta})$. Then the aforementioned grading difference is equal to -k modulo $p_M$. The formula in Definition 4.3 of [Hu2] defines I(C).

In some of the cases considered in the subsequent parts of the proof, the integer I(C) is equal to the Fredholm index that is defined in Equation (4.3) of [Hu2], this being a consequence of the fact that all integral curves of $\nu$ from elements in $\mathcal{Z}_{ech,M}$ are hyperbolic. The equivalence between the ech index and the Fredholm index is used at times in the arguments that follow.

*Part 2*: Suppose that $(\hat{\Theta}', \hat{\Theta})$ lie over the same element, $\Theta \in \mathcal{Z}_{ech,M}$. As such, they can be written respectively as $((\hat{\upsilon}, k'), (\mathfrak{k}_p, o_p)_{p \in \Lambda})$ and $((\hat{\upsilon}, k), (\mathfrak{k}_p, o_p)_{p \in \Lambda})$. If $k > k'$, then there is an ech-subvariety in $\mathcal{M}_{2(k-k')}(\hat{\Theta}', \hat{\Theta})$ which is a union of an $\mathbb{R}$-invariant part and k - k´ distinct spheres from Proposition II.3.1's moduli space $\mathcal{M}_0$. It follows as a consequence that $gr_{ech}(\hat{\Theta}') - gr_{ech}(\hat{\Theta}) = 2(k' - k)$. Granted this fact, it is enough to consider the assertion of the fourth bullet only for those cases where k = k´.

*Part 3*: Suppose that $(\hat{\Theta}', \hat{\Theta})$ are given respectively as $((\hat{\upsilon}, k), (\mathfrak{k}_p, o_p')_{p \in \Lambda})$ and $((\hat{\upsilon}, k), (\mathfrak{k}_p, o_p)_{p \in \Lambda})$. Suppose in addition that there is but one $\mathfrak{p} \in \Lambda$ where $o_p' \neq o_p$ and that the following conditions hold:

$$o_p = \{1, -1\}; \text{ or else } o_p = \{-1\} \text{ and } o_{p'} = \{0\}; \text{ or else } o_p = \{1\} \text{ and } o_{p'} = \{0\}.$$
(9.1)

Let k denote the integer $gr(o_p) - gr(o_p')  \in \{1, 2\}$. There is in this case an element in $\mathcal{M}_k(\hat{\Theta}', \hat{\Theta})$, this being a union of $\mathbb{R}$-invariant cylinders with either one element from Proposition II.3.4's moduli space $\mathcal{M}_{p+}$, or one element from the latter's $\mathcal{M}_{p-}$ in the case k = 1, or one from each in the case k = 2. It follows as a consequence that the ech grading difference $gr_{ech}(\hat{\Theta}') - gr_{ech}(\hat{\Theta})$ is equal to $gr(o_p') - gr(o_p)$.

Given what was said in Part 3, repeated applications of this last observation justify the assertion that it is sufficient to consider the assertion of fourth bullet of Theorem 1.1 only for cases with k = k´ and with all $\mathfrak{p} \in \Lambda$ versions of $o_p'$ and $o_p$ equal.

*Part 4*: Suppose that $(\hat{\Theta}', \hat{\Theta})$ are given respectively as $((\hat{\upsilon}, k), (\mathfrak{k}_p', o_p)_{p \in \Lambda})$ and $((\hat{\upsilon}, k), (\mathfrak{k}_p, o_p)_{p \in \Lambda})$. Suppose in addition that each $\mathfrak{p} \in \Lambda$ version of $o_p = \{0\}$, and that there is but one $\mathfrak{p} \in \Lambda$ with $\mathfrak{k}_p' \ne \mathfrak{k}_p$. Let $\mathfrak{q} \in \Lambda$ denote the exception and take $\mathfrak{k}_q' = \mathfrak{k}_q + 1$. Introduce $\hat{\Theta}''$ to denote $((\hat{\upsilon}, k), (\mathfrak{k}_q, o_q = \{1\}), (\mathfrak{k}_p, o_p)_{p \in \Lambda - q})$. It follows from Proposition 8.1 that $\mathcal{M}_1(\hat{\Theta}', \hat{\Theta}'') \ne \emptyset$ and so $gr_{ech}(\hat{\Theta}') = gr_{ech}(\hat{\Theta}'') - 1$. Meanwhile, it follows from what is said in Part 3 that $\mathcal{M}_1(\hat{\Theta}, \hat{\Theta}'') \ne \emptyset$ also. Thus, $gr_{ech}(\hat{\Theta})$ also equals $gr_{ech}(\hat{\Theta}'') - 1$.

Repeated applications of this last observation justify the assertion that it is enough to consider the fourth bullet only for cases with k = k´ and with all $\mathfrak{p} \in \Lambda$ versions of $(\mathfrak{k}_p, o_p)$ equal to $(\mathfrak{k}_p', o_p')$.

*Part 5*: Fix $n \in \{0, 1, 2, \dots\}$ and suppose that $(\hat{\Theta}', \hat{\Theta})$ are given respectively by $((\hat{\upsilon}, k), (\mathfrak{k}_p, o_p)_{p \in \Lambda})$ and by $((\hat{\upsilon}', k - n), (\mathfrak{k}_p, o_p)_{p \in \Lambda})$. It follows from Lemma 4.1 in [L] that there exists n and an almost complex structure $J_{HF}'$ for $\mathbb{R} \times [1, 2] \times \Sigma$ with the following properties: First, both depend only on the Heegaard-Floer data. Second, $J_{HF}'$ obeys Lipshitz' requirements and those in Section II.6.1. Third, there is a $J_{HF}'$ version of a Lipshitz submanifold, S´, whose constant s slices converge in an isotopic fashion to the arcs in $\hat{\upsilon}$ as $s \to \infty$, and to the arcs in $\hat{\upsilon}'$ as $s \to -\infty$. Fourth, this Lipshitz submanifold has intersection number n with the $f \in [1, 2]$ part of the curve $\gamma^{(z_0)}$ that is described in the fifth bullet of Part 2 in Section 1b. Fifth, the corresponding operator $D_{S'}$ has trivial cokernel. Let $\iota$ denote the Fredholm index of $D_{S'}$.

Given that n, $J_{HF}'$ and S depend only on the Heegaard Floer data, it follows that the data set $(z_*, \delta, x_0, R)$ can be taken without lost of generality so that the following is true: There is an almost complex structure J´ for $\mathbb{R} \times Y$ that is of the sort defined in Section 1c whose restriction to $\mathbb{R} \times M_\delta$ is $J_{HF}'$, and is such that Propositions 7.1-7.3 and 8.1 and 8.2 can be invoked using the Lipshitz submanifold S´. In particular, Proposition 8.1 implies that $\mathcal{M}_1(\hat{\Theta}', \hat{\Theta}) \ne \emptyset$.

This last observation implies that $\deg_{ech}(\hat{\Theta}') - \deg_{ech}(\hat{\Theta}) = gr_{HF}(\hat{\upsilon}') - gr_{HF}(\hat{\upsilon}) - 2n$. The assertion from of the fourth bullet follows from this and what is said in Parts 2-4.

### b) The ech differential

This subsection proves the assertion made by the first bullet in Theorem 1.1. By way of a reminder the endomorphism on the chain complex $\mathbb{Z}(\hat{\mathcal{Z}}_{\text{ech},M})$ that defines the differential has the form depicted in (1.19); thus saying something about the differential requires saying something about the various integers from the relevant version of the set $\{N_{\hat{\Theta}',\hat{\Theta}}\}_{\hat{\Theta}',\hat{\Theta}\in\hat{\mathcal{Z}}_{\text{ech},M}}$. The definition of any given $N_{\hat{\Theta}',\hat{\Theta}}$ is reviewed in Part 1 of what follows. Parts 2-4 say what is needed about these integers to deduce the first bullet of Theorem 1.1.

*Part 1*: Fix $(\hat{\Theta}',\hat{\Theta}) \in \hat{\mathcal{Z}}_{\text{ech},M} \times \hat{\mathcal{Z}}_{\text{ech},M}$. The corresponding integer $N_{\hat{\Theta}',\hat{\Theta}}$ that is used in (1.19) to define the embedded contact homology differential is given by a sum that is indexed by the components of $\mathcal{M}_1(\hat{\Theta}',\hat{\Theta})/\mathbb{R}$ whereby each component contributes either +1 or -1. Whether +1 or -1 is determined by comparing two orientations of the given component. To say more, keep in mind that each component of $\mathcal{M}_1(\hat{\Theta}',\hat{\Theta})$ is a 1-dimensional manifold with a free action of $\mathbb{R}$, with the action given by the constant translations along the $\mathbb{R}$ factor of $\mathbb{R} \times Y$. The generator of this action is a nowhere zero vector field on each component. This vector field defines an orientation for each component, this denoted by $\hat{o}_{\text{ech},\mathbb{R}}$. The second orientation is given by a coherent orientation for $\mathcal{M}_{\text{ech}}$. The latter orientation is denoted by $\hat{o}_{\text{ech},Q}$. Write $\hat{o}_{\text{ech},\mathbb{R}}$ as $N^C \hat{o}_{\text{ech},Q}$ with $N^C \in \{1,-1\}$. The assignment $C \to N^C$ is a locally constant, $\{1,-1\}$-valued function on $\mathcal{M}_1(\hat{\Theta}',\hat{\Theta})$. The value of this function $N^C$ on C's component of $\mathcal{M}_1(\hat{\Theta}',\hat{\Theta})$ is the component's contribution to the sum that defines $N_{\hat{\Theta}',\hat{\Theta}}$.

*Part 2*: There is an analogous definition of the endomorphism of $\mathbb{Z}(\mathcal{Z}_{\text{HF}} \times \mathbb{Z})$ that defines the differential for Heegaard Floer homology. To elaborate, any given endomorphism of $\mathbb{Z}(\mathcal{Z}_{\text{HF}} \times \mathbb{Z})$ is defined by its action on the generating set, and so by a rule of the form

$$(\hat{\upsilon}',k') \to \sum_{(\hat{\upsilon},k)\in\mathcal{Z}_{\text{HF}}\times\mathbb{Z}} M_{(\hat{\upsilon}',k'),(\hat{\upsilon},k)} (\hat{\upsilon},k)$$

(9.2)

where $M_{(\cdot)(\cdot)}$ is in all cases an integer. The paragraphs that follow define these integers when the endomorphism in question is the differential $\partial_{\text{HF}}$.

Fix ordered pairs $(\hat{\upsilon}',k')$ and $(\hat{\upsilon},k)$ of elements from $\mathcal{Z}_{\text{HF}} \times \mathbb{Z}$. Having done so, reintroduce the subspace $\mathcal{A}_{\text{HF1}}((\hat{\upsilon}',k'),(\hat{\upsilon},k)) \subset \mathcal{A}_{\text{HF}}$ from Part 1 of Section 8f. The integer $M_{(\hat{\upsilon}',k'),(\hat{\upsilon},k)}$ for the version of (9.2) that defines the differential for Heegaard Floer homology is a sum that is indexed by the components of $\mathcal{A}_{HF1}((\hat{v}',k'),(\hat{v},k))$ with each component contributing either +1 or -1.

To define these ±1 contributions, note that $\mathcal{A}_{HF1}((\hat{v}',k'),(\hat{v},k))$ has a finite set of components, with each a copy of $\mathbb{R}$. Each such copy of $\mathbb{R}$ is a free orbit of the $\mathbb{R}$ action that comes from the constant translations along the $\mathbb{R}$ factor of $\mathbb{R} \times [1,2] \times \Sigma$. The generator of this $\mathbb{R}$ action orients each component. This orientation is denoted by $\hat{o}_{HF,\mathbb{R}}$. A second orientation is that supplied by a given coherent orientation for $\mathcal{A}_{HF}$. This orientation is denoted by $\hat{o}_{HF,Q}$. Let $S \in \mathcal{A}_{HF1}((\hat{v}',k'),(\hat{v},k))$ denote a given surface and write $\hat{o}_{HF,\mathbb{R}}$ at S as $z^S \hat{o}_{HF,Q}$ where $z^S \in \{-1, 1\}$ is constant on the component of S in $\mathcal{A}_{HF1}((\hat{v}',k'),(\hat{v},k))/\mathbb{R}$. The value of $z^S$ is this component's contribution to $M_{(\hat{v}',k'),(\hat{v},k)}$.

*Part 3*: Fix a coherent system of orientations for $\mathcal{A}_{HF}$ to define the coefficients in the version of (9.2) that defines the Heegaard Floer differential. Use this same coherent system in Proposition 8.5 to define the coherent system of orientations that is used to define the embedded contact homology differential.

Fix a pair, $(\hat{\Theta}', \hat{\Theta})$, from $\hat{\mathcal{Z}}_{ech,M}$. The corresponding $N_{\hat{\Theta}',\hat{\Theta}}$ is zero unless $\mathcal{M}_1(\hat{\Theta}', \hat{\Theta})$ is non-empty and therefore $(\hat{\Theta}', \hat{\Theta})$ is described by (8.15). Consider first the case given by the first bullet in (8.15). It follows from Proposition 7.2 that the elements in $\mathcal{M}_1(\hat{\Theta}', \hat{\Theta})$ are labeled in part by the surfaces in $\mathcal{A}_{HF1}((\hat{v}',k'),(\hat{v},k))$. Let S denote a given such surface and let $\mathcal{M}_1^S(\hat{\Theta}', \hat{\Theta}) \subset \mathcal{M}_1(\hat{\Theta}', \hat{\Theta})$ denote the corresponding subset. It follows directly from Propositions 8.1, 8.2 and 8.5 plus what is said in Part 5 of Section 8b that

$$\sum_{C \in \mathcal{M}_1^S(\hat{\Theta}',\hat{\Theta})} N^C = z^S .$$

(9.3)

This implies directly that $N_{\hat{\Theta}',\hat{\Theta}} = M_{(\hat{v}',k'),(\hat{v},k)}$.

Write $\mathbb{Z}(\hat{\mathcal{Z}}_{ech,M})$ as $\mathbb{Z}(\hat{\mathcal{Z}}_{ech,M}) = \mathbb{Z}(\mathcal{Z}_{HF} \times \mathbb{Z}) \otimes (\otimes_{\mathfrak{p} \in \Lambda} \mathbb{Z}(\mathbb{Z} \times O))$ as done in Theorem 1.1. The conclusion of the preceding paragraph implies that writing $\mathbb{Z}(\hat{\mathcal{Z}}_{ech,M})$ in this way makes $\partial_{ech}$ appear as $\partial_{ech} = \partial_{HF} + L$ with L acting solely on the $(\otimes_{\mathfrak{p} \in \Lambda} \mathbb{Z}(\mathbb{Z} \times O))$ factor.

*Part 4*: The endomorphism L is defined by those $N_{\hat{\Theta}',\hat{\Theta}}$ with $(\hat{\Theta}', \hat{\Theta})$ as described in the second bullet of (8.15). It follows as a consequence that L can be written as $\sum_{\mathfrak{p} \in \Lambda} L_\mathfrak{p}$ where $L_\mathfrak{p}$ acts only on the $\mathfrak{p} \in \Lambda$ factor of $\mathbb{Z}(\mathbb{Z} \times O)$ in $\mathbb{Z}(\mathcal{Z}_{HF} \times \mathbb{Z}) \otimes (\otimes_{\mathfrak{p} \in \Lambda} \mathbb{Z}(\mathbb{Z} \times O))$. This is because there is but one pair $\mathfrak{p} \in \Lambda$ with $(\mathfrak{k}_\mathfrak{p}, O_\mathfrak{p}) \neq (\mathfrak{k}_\mathfrak{p}', O_\mathfrak{p}')$. Moreover, the relevant version of $L_\mathfrak{p}$ acts on the given generator $(\mathfrak{k}_\mathfrak{p}, O_\mathfrak{p})$ to give an integer weighted sum of generators with the weight being zero unless $\Delta_\mathfrak{p}' = \Delta_\mathfrak{p} - 1$ in which case either Item a) or Item b) in (8.15) must occur. This being the case, it follows from Propositions 8.1 and 8.2 with what is said in Part 5 of Section 8b that the corresponding integer weight is either 1 or -1.

Consider first the case when Item a) of the second bullet is obeyed. Proposition 8.5 and what is said in Part 3 of Section 8b determine these signs:

- *In the case of the first three bullets in (8.16), the sign is* $(-1)^{N_+ + N}$.
- *In the case of the fourth bullet in (8.16), the sign is* $(-1)^{N_+ + N + 1}$.

(9.4)

Suppose next that Item b) of the second bullet holds. Proposition 8.5 and what is said in Part 4 of Section 8f determine that the sign is again given by (9.4). Note in this regard that the fourth bullet in (8.16) can occur only if $\mathfrak{k}_p' = \mathfrak{k}_p - 1$.

What is written in (9.4) is consistent with what is claimed by Theorem 1.1 if and only if the integers $N_+$ are such that $(-1)^{N_+} = \varepsilon (-1)^{\deg_{HF}(\hat{\upsilon}, k)}$ with $\varepsilon \in \{-1, 1\}$ being independent of both $\hat{\upsilon}$ and k. Since the Heegaard Floer degree changes by an even integer as k varies, it is enough to verify that this is so for any given value of k. That such is the case follows from Proposition 4.8 in [L] and Equation (9) in [L]. The latter expresses an equality that was derived by J. Rasmussen in [R].

### c) The endomorphisms from the second and third bullets of Theorem 1.1

The assertion made by the second bullet about the action of the $\mathbb{U}$ map on the chain complex $\mathbb{Z}(\mathcal{Z}_{HF} \times \mathbb{Z}) \otimes (\otimes_{p \in \Lambda} \mathbb{Z}(\mathbb{Z} \times O))$ follows directly from what is said in the first paragraph of Part 2 in Section 9d. The assertions about the endomorphisms in the third bullet of Theorem 1.1 are discussed in the subsequent two parts of this subsection. The first part briefly reviews the definitions of the coefficients that appear in the corresponding versions of (1.19).

*Part 1*: This part of the subsection explains how to the endomorphisms in the third bullet of Theorem 1.1 are defined. To this end, let $\hat{\imath}$ denote one of the cycles from the set $\{\hat{\imath}_p\}_{p \in \Lambda}\}$ and let $Q^{\hat{\imath}}_{ech}$ denote the corresponding endomorphism of $\mathbb{Z}(\hat{\mathcal{Z}}_{ech,M})$. This endomorphism is described by a version of (1.19), and thus defined by the integers $\{N_{\hat{\Theta}', \hat{\Theta}}\}_{\hat{\Theta}', \hat{\Theta} \in \hat{\mathcal{Z}}_{ech,M}}$. Let $(\hat{\Theta}', \hat{\Theta})$ denote a given pair from $\hat{\mathcal{Z}}_{ech,M}$. The corresponding integer $N_{\hat{\Theta}', \hat{\Theta}}$ is non-zero only if $\mathcal{M}_1(\hat{\Theta}', \hat{\Theta})$ is non-empty. Each component in the latter set contributes an integer to a sum whose value is $N_{\hat{\Theta}', \hat{\Theta}}$. This understood, let $C \subset \mathcal{M}_1(\hat{\Theta}', \hat{\Theta})$ denote a given element, and let [C] denote the corresponding relative 2-cycle in $H_2(Y; [\Theta] - [\Theta'])$ defined by the image of C in Y via the projection from $\mathbb{R} \times Y$. The cycle $\hat{\imath}$ has been chosen so as to be disjoint from the integral curves of $v$ that appear in elements from $\hat{\mathcal{Z}}_{ech,M}$, and so there is a well defined pairing between $\hat{\imath}$ and [C] with values in $\mathbb{Z}$. Use $\langle \hat{\imath}, [C] \rangle$ to denote this pairing. Reintroduce the sign $N^C \in \{1,-1\}$ from Part 1 of Section 9c. The contribution of C's component in $\mathcal{M}_1(\hat{\Theta}', \hat{\Theta})$ to the sum for $N_{\hat{\Theta}',\hat{\Theta}}$ is $\langle \hat{\imath}, [C] \rangle N^C$.

By way of comparison, what follows summarizes from Section 8 of [L] the Heegaard Floer version of the endomorphisms that are define by the cycles from the set $\{[\gamma^{(z_0)}], \{\hat{\imath}^{(z)}\}_{z \in \yen - z_0}\}$. Let $\hat{\imath}$ now denote one of the cycles from the latter set, and let $Q^{\hat{\imath}}_{HF}$ denote the corresponding endomorphism. This endomorphism is defined by a version of (9.2). A given coefficient $M_{(\hat{\upsilon}',k'),(\hat{\upsilon},k)}$ is non-zero only if $\mathcal{A}_{HF1}((\hat{\upsilon}',k'),(\hat{\upsilon},k))$ is non-empty. If so, then each component of this space contributes an integer to a sum whose value is $M_{(\hat{\upsilon}',k'),(\hat{\upsilon},k)}$. Let S denote a given Lipshitz surface from this space. The image of S in $\Sigma \times [0,1]$ via the projection from $\mathbb{R} \times \Sigma \times [0,1]$ has a well defined intersection pairing with $\hat{\imath}$, this denoted by $\langle \hat{\imath}, S \rangle$. Reintroduce $Z^S \in \{1,-1\}$ from Part 1 of Section 9b. The component of S contributes $\langle \hat{\imath}, S \rangle Z^S$ to the sum that computes $M_{(\hat{\upsilon}',k'),(\hat{\upsilon},k)}$.

*Part 2*: Consider first the statements made by Items a) and b) of the third bullet in Theorem 1.1. To this end, fix $\hat{\Theta}'$ and $\hat{\Theta}$ from $\hat{\mathcal{Z}}_{ech,M}$. Given $S \in \mathcal{A}_{HF1}((\hat{\upsilon}',k'),(\hat{\upsilon},k))$, reintroduce from Part 3 in Section 9b the subspace $\mathcal{M}_1^S(\hat{\Theta}', \hat{\Theta})$. Given the definitions in Part 1, the assertion made by Item b) of the third bullet in Theorem 1.1 follows from (9.3) if $\langle \hat{\imath}, [C] \rangle = \langle \hat{\imath}, S \rangle$ when $C \in \mathcal{M}_1^S(\hat{\Theta}', \hat{\Theta})$. The latter equality is a consequence of Proposition 8.2.

Consider next the statement made by Item c) of the third bullet in Theorem 1.1. To this end, fix $\mathfrak{p} \in \Lambda$ so as to see about the action of the $\hat{\imath} = \hat{\imath}_\mathfrak{p}$ version of $Q_{\hat{\imath}}$. The cycle $\hat{\imath}_\mathfrak{p}$ is disjoint from $M_\delta \cup (\cup_{\mathfrak{p}' \in \Lambda - \mathfrak{p}} \mathcal{H}_{\mathfrak{p}'})$ and as a consequence, it must act as $\mathbb{I}_\mathfrak{p} + \hat{\iota}_\mathfrak{p}$. This being the case, at issue is the precise form for $\hat{\iota}_\mathfrak{p}$. The cycle $\hat{\imath}_\mathfrak{p}$ lies in the $\cos\theta > \frac{1}{\sqrt{3}}$ part of $\mathcal{H}_\mathfrak{p}$, and as a consequence any given integer $N_{\hat{\Theta}',\hat{\Theta}}$ from the $Q_{\hat{\imath}}$ version of (1.19) is zero unless $\hat{\Theta}'$ and $\hat{\Theta}$ are described by Item a) of the second bullet in (8.15) and the first two bullets in (8.16). With this point understood, Item c) of the third bullet in Theorem 1.1 follows directly from Propositions 8.1 and 8.5 plus what is said in Part 5 of Section 8b.

Cagatay Kutluhan; Department of Mathematics, Columbia University, New York, NY 10027. *Email address*: kutluhan@umich.edu.

Yi-Jen Lee; Department of Mathematics, Purdue University, West Lafayette, IN 47907
*Email address*: yjlee@math.purdue.edu.

Clifford Henry Taubes; Department of Mathematics, Harvard University, Cambridge, MA 02138. *Email address*: chtaubes@math.harvard.edu.
195